\documentclass[10pt,epsfig]{amsart}
 \usepackage{amsbsy,amssymb,amscd,amsfonts,latexsym,amstext,delarray,
 amsmath,graphicx}
 \usepackage[dvips]{epsfig}

\def\ug{{\mathbf u}}

\def\semi{\rtimes}

\newtheorem{thm}{Theorem}[section]
\newtheorem{prop}[thm]{Proposition}
\newtheorem{cor}[thm]{Corollary}
\newtheorem{lem}[thm]{Lemma}
\newtheorem{defn}[thm]{Definition}

\numberwithin{equation}{section}

\newcommand{\ie}{{\it i.e.\/}\ }

\newcommand{\cf}{{\it cf.\/}\ }

\def\sin{{{\rm sin}}}
\def\cos{{{\rm cos}}}
\def\tan{{{\rm tan}}}

\def\cal{\mathcal}

\def\id{{\mbox{Id}}}
\def\ch{{\mbox{ch}}}
\def\im{{\mbox{Im}}}
\def\Aut{{\mbox{Aut}}}

\def\hom {{\mbox{hom}}}

\def\Trace{{\mbox{Trace}}}

\def\cala{{\cal A}}
\def\calb{{\cal B}}

\def\call{{\cal L}}
\def\calm{{\cal M}}

\def\calo{{\cal O}}

\def\calc{{\cal C}}

\def\frach{{\mathfrak h}}

\def\fracS{{\mathfrak {S}}}
  
 \def\fracsu{{\mathfrak {su}}}
 \def\fracsl{{\mathfrak {sl}}}

\def\calg{{C_{\mathrm{alg}}}}
\def\bbbone{\mbox{\rm 1\hspace {-.6em} l}}

\def\ug{{\mathbf u}}
\def\mg{{\mathbf m}}
\def\vg{{\mathbf v}}

\def\Tg{{\mathbf T}}

\def\vg{{\mathbf v}}

\def\Cb{{\mathbb C}}

\def\Rb{{\mathbb R}}

\def\Zb{{\mathbb Z}}

\def\P{{\mathbb P}}

\def\qqq{\,,\quad \forall}

\begin{document}

 \thispagestyle{empty}

\title[NON COMMUTATIVE 3-SPHERES]{{\Large\bf NON COMMUTATIVE FINITE DIMENSIONAL MANIFOLDS II:}\\
 \vspace{1.5cm} MODULI SPACE AND STRUCTURE \\
 \vspace{0.5cm} OF NON COMMUTATIVE 3-SPHERES}

\vspace{2.5cm}
\bigskip
\bigskip

\author[Connes]{Alain Connes}
\author[Dubois-Violette]{Michel Dubois-Violette}

\address{A.~Connes: Coll\`ege de France \\
3, rue d'Ulm \\ Paris, F-75005 France\\ I.H.E.S. and Vanderbilt
University.} \email{alain\@@connes.org}
\address{M.~Dubois-Violette: Laboratoire de Physique Th\'eorique, UMR 8627,
Universit\'e Paris XI,
B\^atiment 210, F-91 405 Orsay Cedex, France\\
Bonn, D-53111 Germany} \email{Michel.Dubois-Violette\@@th.u-psud.fr}

\maketitle

\baselineskip=0.5cm

\begin{abstract}

This paper contains detailed proofs of our results on the moduli
space and the structure of noncommutative 3-spheres. We develop the
notion of central quadratic form for quadratic algebras, and a
general theory which creates a bridge between noncommutative
differential geometry and its purely algebraic counterpart. It
allows to construct a morphism from  an involutive quadratic
algebras to a C*-algebra constructed from the characteristic variety
and the hermitian line bundle associated to the central quadratic
form. We apply the general theory in the case of noncommutative
3-spheres and show that the above morphism corresponds to  a natural
ramified covering  by a noncommutative 3-dimensional nilmanifold. We
then compute the Jacobian of the ramified covering  and obtain the
answer as the product of a period (of an elliptic integral) by a
rational function. We describe the real and complex moduli spaces of
noncommutative 3-spheres, relate the real one to root systems and
the complex one to the orbits of a birational cubic automorphism of
three dimensional projective space. We classify  the algebras and
establish duality relations between them.

\end{abstract}

\section{Introduction}

This paper contains detailed proofs of our results on the moduli
space and the structure of noncommutative $3$-spheres announced in
\cite{ac-mdv:2003}.

Through the analysis of a specific class of noncommutative manifolds
which arose  as solutions of a simple equation of $K$-theoretic
origin we discovered rather general structures which lie at the
intersection of two fundamental aspects of noncommutative geometry,
namely

\begin{itemize}

\item Differential Geometry

\item Algebraic Geometry

\end{itemize}

This class of noncommutative manifolds, called noncommutative
$3$-spheres, has a very rich structure both at the level of the
objects themselves as well as at the level of the moduli space which
parameterizes these geometric objects. There are two aspects in the
geometry of the moduli space :

\begin{itemize}

\item The {\em real} moduli space and its scaling foliation,
its link with the alcove structure of the root system $D_3$ and the
Morse theory of the character of the signature representation.

\item The {\em complex} moduli space and its
net of elliptic curves, its link with the iteration of a cubic
transformation of $\P_3(\Cb)$.

\end{itemize}

At the level of the structure of the noncommutative $3$-spheres our
main result is to relate them to very well understood noncommutative
nilmanifolds which fall under the framework of the early theory
developped in \cite{ac:1980} and were analysed in great detail in
\cite{aba-exel:1997}
 \cite{aba-eil-exel:1998}.

The core of the paper is to construct the corresponding map of
noncommutative spaces and compute its Jacobian. The essence of the
work is to extract from very complicated computations the general
concepts that allow not only to understand what is going on but also
to extend the construction in full generality. Thus at center stage
lies the computation of the Jacobian and the gradual simplification
of the result which at first was expressed in terms of elliptic
functions and the $9$-th power of Dedekind $\eta$ function. We shall
reach at the end of the paper (Theorem \ref{vol1}) a result of
utmost simplicity while the starting point was a computation
performed even in the trigonometric case with the help of a
computer\footnote{We are grateful to Michael Trott for his help}.
While we reach a reasonable level of conceptual understanding of the
general construction of the map, we believe that a lot remains to be
discovered for the abstract construction of the calculus as well as
for the general computation of the Jacobian, let alone in the case
of higher dimensional spheres which we do not adress here.

After recalling the basic definitions and properties of the
noncommutative $3$-spheres in section \ref{basicp} we analyse the
real moduli space in section \ref{modulispace} and exhibit a
fundamental domain in terms of alcoves of the root system $D_3$.

In section \ref{triang} we define the scaling foliation and show its
compatibility with the alcove structure of the real moduli space. We
also show that the isomorphism class of the $4$-spaces ${\mathbb
R}^4(\Lambda)$ remains constant on the leaves of the foliation. In
order to prove the converse \ie that isomorphism of the $4$-spaces
${\mathbb R}^4(\Lambda)$ implies equality of the leaves we compute
in details in section \ref{geometricdata} the geometric data of the
quadratic algebra of  ${\mathbb R}^4(\Lambda)$. This allows to
finish the proof of the converse in section \ref{flowbetter}.

We exhibit in section \ref{dualities} more subtle relations between
the $4$-spaces ${\mathbb R}^4(\Lambda)$ given by {\em dualities}. At
the level of the algebras these are obtained from the general notion
of semi-cross product of quadratic algebras. At the level of the
moduli  space these dualities shrink further the fundamental domain
and that amounts essentially to the transition from the root system
$D_3$ to the larger one $C_3$.

We then use these to describe the  $4$-spaces ${\mathbb
R}^4(\Lambda)$ for degenerate values of the parameter in section
\ref{degalg}.

In section \ref{complexmod} we analyse the complex moduli space
which appears naturally as a net of elliptic curves having eight
points in common in $\P_3(\Cb)$. We show that in the generic case
these elliptic curves are the characteristic varieties of the
algebras that their points label. Moreover the canonical
correspondence $\sigma$ is simply the restriction of a globally
defined cubic map of $\P_3(\Cb)$.
 This gives in particular a very natural
choice of generators for the algebra. As a preparation for the next
section we give the natural parameterization of the  net of elliptic
curves in terms of $\vartheta$-functions.

Section \ref{pairing} is the most technical one and contains the
root of the concepts developped in full generality in section
\ref{Cstar}. In essence what we do first is starting from unitary
representations of the Sklyanin algebra constructed by Sklyanin in
his original paper we derive a one parameter family of
$\star$-homorphisms from the  algebras of ${\mathbb R}^4(\Lambda)$
in the generic case, to the  algebras of noncommutative tori. We
then use a suitable restriction to
 the $3$-spheres $S^3 (\Lambda)$. After a lot of work on this
construction we find that we can eliminate all occurences of
$\vartheta$-functions from the formulas and  obtain a purely
algebraic formulation of the construction as a morphism to a twisted
cross product $C^*$-algebra obtained from the geometric data.

Section \ref{Cstar} describes  the abstract general construction in
the framework of involutive quadratic algebras. The key notions are
those of {\em central quadratic form} and of {\em positivity} for
such forms. It is in that section that the interaction between the
two above aspects of noncommutative geometry is manifest. In fact we
construct a bridge between the purely algebraic notions such as the
geometric data of a quadratic algebra and the world of
noncommutative geometry including the topological ($C^*$-algebraic)
and differential geometric aspects (in the sense of \cite{ac:1980},
\cite{ac:1982}). At one end of the bridge one starts with the given
involutive quadratic algebra. At the other one has the $C^*$-algebra
obtained as the twisted cross product of the charateristic variety
by the canonical correspondence. The twisting is effected by an
hermitian line bundle and the construction is a special case of a
general one due to M. Pimsner. The bridge provides a $\star$-algebra
morphism. This algebra morphism has a ``trivial part" which does not
make use of the central quadratic form and lands in a ``triangular"
subalgebra of the $C^*$-algebra. This part was well-known for quite
sometime to noncommutative algebraic geometors. The non-trivaility
of our construction lies in the involved relations coming from the
cross terms mixing generators with their adjoints.

We compute in section \ref{jacobian}  the Jacobian of the above map.
We first define what we mean by the jacobian in the sense of
noncommutative geometry where Hochschild homology replaces
differential forms. The algebraic form of the result then suggests
the existence of a calculus of purely algebraic nature allowing to
express the cyclic cohomology fundamental class in terms of the
algebraic geometry of the characteristic variety and a hermitian
structure on the canonical line bundle. This is achieved in  section
\ref{calculus} and allows to finally obtain the purest form of the
computation of the jacobian in the already mentionned Theorem
\ref{vol1}. Needless to say this is a paper of highly
``computational" nature and we tried to ease the reading by
supplying in an appendix the basic factorisations of the minors and
the sixteen theta relations which are often used in the text.

\tableofcontents

\newpage

\section{The noncommutative $3$-spheres $S^3 (\Lambda) \subset {\mathbb R}^4(\Lambda)$} \label{basicp}

We shall recall in this section the basic properties of the
noncommutative spheres $S^3 (\Lambda)$ and the corresponding
$4$-spaces ${\mathbb R}^4(\Lambda)$.

\subsection{Unitary ``up to scale''}

\noindent
\medskip

We let ${\mathcal A}$ be a unital involutive algebra and first start
with a unitary ``up to scale'' in $M_q({\mathcal A})$, {\it i.e.}
\begin{equation}
\label{basic} U \in M_q ({\mathcal A}) \, , \quad UU^* = U^* U \in
{\mathcal A} \otimes \bbbone \subset {\mathcal A} \otimes M_q
({\mathbb C}) \, .
\end{equation}
\begin{lem} \label{centerC}
Let $U\in M_q(\cala)$ satisfy  (\ref{basic}) with
\[
UU^\ast=U^\ast U=C\otimes \bbbone \in \cala \otimes M_q(\mathbb C)
\]
then $C$ is in the center of the $\ast$-algebra generated by the
matrix elements of $U$.
\end{lem}
\begin{proof}
One has $(C\otimes \bbbone) U=(UU^\ast)U=U(U^\ast U)=U(C\otimes
\bbbone)$ and $(C\otimes \bbbone)U^\ast=(U^\ast
U)U^\ast=U^\ast(UU^\ast)=U^\ast(C\otimes \bbbone)$ by associativity
in $M_q(\cala)$, which implies the result.
\end{proof}
Let
\begin{equation} \label{baseq}
\tau_\mu,\, \, \mu\in \{0,\dots,q^2-1\}
\end{equation}
be an orthonormal basis of $M_q(\mathbb C)$ for the scalar product
$\langle A\vert B\rangle=\frac{1}{q}\Trace (A^\ast B)$. Then one has
\begin{equation}\label{defz}
U=\tau_\mu z^\mu,\>\>\> z^\lambda\in \cala
\end{equation}
where we used the Einstein summation convention on ``up-down
indices". The $\ast$-algebra generated by the matrix elements of $U$
is the $\ast$-subalgebra of $\cala$ generated by the $z^\lambda,
z^{\mu\ast}$ ($\lambda,\mu \in \{0,\dots, q^2-1\})$ and if $U$ is as
in the above lemma then one has the equalities
\begin{equation}\label{casiq}
C=\sum_\mu z^\mu z^{\mu\ast}=\sum_\mu z^{\mu\ast}z^\mu
\end{equation}
for the central element $C$ of the $\ast$-algebra generated by the
$z^\mu$.

\subsection{Equation ${\rm ch}_{1/2} (U) = 0$}

\noindent
\medskip

  We now turn to the relation ${\rm ch}_{1/2} (U) = 0$, in the unreduced complex \ie using the convention of summation on repeated indices,
\begin{equation}
\label{chzero} U_{i_1}^{i_0} \otimes (U^*)_{i_0}^{i_1} -
(U^*)_{i_1}^{i_0} \otimes U_{i_0}^{i_1} = 0
\end{equation}
where the  left hand side belongs to the tensor square ${\mathcal
A}^{\otimes 2} = {\mathcal A} \otimes {\mathcal A}$. Both terms in
the  left hand side give sums of the $q^2$ terms of the form
$(z^{\mu} \otimes z^{\nu\ast}) \, (\tau_{\mu})_{i_1}^{i_0} \,
(\tau_{\nu}^\ast)_{i_0}^{i_1}$ and similarly for the other.  The sum
on $i_0 , i_1$ is thus simply Trace $(\tau_{\mu} \, \tau_{\nu}^*)$,
{\it i.e.} the Hilbert Schmidt inner product  $(\tau_{\mu} ,
\tau_{\nu})=q\langle \tau_\nu\vert \tau_\mu\rangle$. This is $0$
unless $\mu = \nu$ and is $q$ if $\mu = \nu$. Thus, up to an overall
factor of $q$ the equality (\ref{chzero}) means:
\begin{equation}
\label{sym} \sum (z^{\mu} \otimes z^{\mu\ast} - z^{\mu\ast}  \otimes
z^{\mu}) = 0 \, .
\end{equation}
\begin{lem}  Equation \eqref{sym} holds iff
there exists a {\it unitary symmetric} matrix $\Lambda\in
M_{q^2}({\mathbb C}) $ such that:
\begin{equation}
\label{lambda} z^{\mu\ast} = \Lambda_{\nu}^{\mu} \, z^{\nu} \, .
\end{equation}

\end{lem}
\bigskip

\begin{proof}
Let us first assume that the $z^{\mu}$ are linearly independent
elements of ${\mathcal A}$. Let then $\varphi_{\mu}$ be linear forms
on ${\mathcal A}$ with $\varphi_{\mu} (z^{\nu}) =
\delta^{\nu}_{\mu}$. Applying $1 \otimes \varphi_{\mu}$ to
(\ref{sym}) we get
\begin{equation}
z^{\mu\ast}= \sum z^{\nu} \, \varphi_{\mu} (z^{\nu\ast}) = \sum
\Lambda_{\nu}^{\mu} \, z^{\nu} \nonumber
\end{equation}
where the matrix $\Lambda$ is uniquely prescribed by this relation.
Then since the $z^{\mu} \otimes z^{\nu}$ are linearly independent in
${\mathcal A} \otimes {\mathcal A}$ the relation (\ref{sym}) means,
looking at the coefficient of $z^{\mu} \otimes z^{\nu}$ on both
sides,
\begin{equation}
\label{eq1.11} \Lambda_{\nu}^{\mu} = \Lambda_{\mu}^{\nu}
\end{equation}
so that the matrix $\Lambda$ is symmetric.

\medskip

  Taking the adjoint of both sides in $z^{\mu\ast} = \Lambda_{\nu}^{\mu} \, z^{\nu}$ one gets
\begin{equation}
z^{\mu} = \bar\Lambda_{\nu}^{\mu} \, z^{\nu\ast} =
\bar\Lambda_{\nu}^{\mu} \, \Lambda_{\rho}^{\nu} \, z^{\rho} =
(\Lambda^* \Lambda)_{\rho}^{\mu} \, z^{\rho} \nonumber
\end{equation}
and the linear independence of the $z^{\rho}$ thus shows that:
\begin{equation}
\label{eq1.12} \Lambda^* \Lambda = 1 \, .
\end{equation}
For the general case\footnote{The simplification of the argument of
\cite{ac-mdv:2002a} given here is due to G. Skandalis}
 note that equation \eqref{sym} is invariant
by the transformation
\begin{equation}
\tilde z^{\mu} =  W^{\mu}_{\nu} \, z^{\nu} \, . \nonumber
\end{equation}
where $W\in U(q^2)$ is a unitary matrix. Moreover \eqref{lambda}
implies a similar equation for $(\tilde z^{\mu})$ with the matrix
$$
\tilde{\Lambda}=\, \bar{W}\,\Lambda\,\bar W^{t}\,,
$$
which is still symmetric and unitary. This allows to assume that the
kernel of the linear map from ${\mathbb C}^{q^2}$ to ${\mathcal A}$
determined
 by  $e^\nu\mapsto z^{\nu}$ is the linear span of a subset $I$
of the basis vectors  $e^\nu$. In other words the non-zero $z^{\nu}$
are linearly independent and the above proof ensures the existence
of a matrix fulfilling \eqref{lambda} which is symmetric and unitary
once extended by the identity on the $e^\nu\,,\,\nu \in I$.

\end{proof}

As pointed out in \cite{ac-mdv:2002a} and as will be explained in
Part III, for the study of $(2n+1)$-dimensional spherical manifolds
one should take $q=2^n$ to be coherent in particular with the
suspension functor. In the following we shall concentrate on the
3-dimensional case (and the corresponding noncommutative $\mathbb
R^4$) which is the lowest dimensional non trivial case from the
noncommutative side and for which $\ch_{1/2}(U)=0$ is the only
$K$-homological condition. Accordingly we take $q=2$ in the
following.

\subsection{Noncommutative 3-spheres and 4-planes}

\noindent
\medskip

We now turn to the case $q=2$ and we choose as orthonormal basis of
$M_2(\mathbb C)$ the basis
\begin{equation}
\label{base2} \tau_0=\bbbone \>\> \mbox{and}\>\>
\tau_k=i\sigma_k,\>\> k\in\{1,2,3\}
\end{equation}
where the $\sigma_k$ are the Pauli matrices, {\it i.e.}
\begin{equation}
\label{pauli} \sigma_1 = \begin{bmatrix} 0 &1 \\ 1 &0 \end{bmatrix}
\, , \quad \sigma_2 = \begin{bmatrix} 0 &-i \\ i &0 \end{bmatrix} \,
, \quad \sigma_3 = \begin{bmatrix} 1 &0 \\ 0 &-1 \end{bmatrix}
\end{equation}
which satisfy $\sigma_j^* = \sigma_j$, $\sigma_j^2 = 1$  and
$\sigma_k \, \sigma_{\ell} = i \, \varepsilon_{k\ell m} \, \sigma_m$
for any permutation $(k,\ell, m)$ of (1, 2, 3), where
$\varepsilon_{123} = 1$ and $\varepsilon$ is totally antisymmetric.
This allows to write
\begin{equation}
\label{eq1.3} U = z^0 + i \, \sigma_1 \, z^1 + i \, \sigma_2 \, z^2
+  i \, \sigma_3 \, z^3, \qquad z^\mu \in {\mathcal A} \,
\end{equation}
and one has $U^* = z^{0\ast} - \, i \, \sigma_1 \, z^{1\ast}- \, i
\, \sigma_2 \, z^{2\ast} - \, i \, \sigma_3 \, z^{3\ast}$ .

\medskip

  Thus
\begin{eqnarray}
UU^* &= &z^0 \, z^{0\ast} + z^1 \, z^{1\ast} + z^2 \, z^{2\ast} + z^3 \, z^{3\ast} \nonumber \\
&+ &i \, \sigma_1 \, (z^1 \, z^{0\ast} - z^0 \, z^{1\ast}+ z^2 \, z^{3\ast} - z^3 \, z^{2\ast}) \nonumber \\
&+ &i \, \sigma_2 \, (z^2 \, z^{0\ast} - z^0 \, z^{2\ast}  + z^3 \, z^{1\ast}- z^1 \, z^{3\ast}) \nonumber \\
&+ &i \, \sigma_3 \, (z^3 \, z^{0\ast} - z^0 \, z^{3\ast} + z^1 \,
z^{2\ast}  - z^2 \, z^{1\ast} ) \, . \nonumber
\end{eqnarray}
Similarly we get,
\begin{eqnarray}
U^*U &= &z^{0\ast} \, z^0 +z^{1\ast} \, z^1 + z^{2\ast}\, z^2 + z^{3\ast} \, z^3 \nonumber \\
&+ &i \, \sigma_1 \, (z^{0\ast} \, z^1 - z^{1\ast} \, z^0 + z^{2\ast} \, z^3 - z^{3\ast} \, z^2) \nonumber \\
&+ &i \, \sigma_2 \, (z^{0\ast} \, z^2 - z^{2\ast} \, z^0 + z^{3\ast} \, z^1 - z^{1\ast} \, z^3) \nonumber \\
&+ &i \, \sigma_3 \, (z^{0\ast} \, z^3 - z^{3\ast}\, z^0 + z^{1\ast}
\, z^2 - z^{2\ast} \, z^1) \, . \nonumber
\end{eqnarray}
Thus  equation \eqref{basic} is equivalent to 7 relations which are,
\begin{equation}
\label{relcenter1} \sum z^{\mu} \, z^{\mu\ast} = \sum z^{\mu\ast} \,
z^{\mu}
\end{equation}
\begin{equation}
\label{relright} z^k \, z^{0\ast} - z^0 \, z^{k\ast} + \sum
\varepsilon_{k\ell m} \, z^{\ell} \, z^{m\ast}= 0
\end{equation}
\begin{equation}
\label{relleft} z^{0\ast} \, z^k - z^{k\ast} \, z^0 + \sum
\varepsilon_{k\ell m} \, z^{\ell\ast} \, z^m = 0 \, .
\end{equation}

\bigskip

We then let ${\mathcal S}$ be the space of unitary symmetric
matrices,
\begin{equation}
\label{eq1.30} {\mathcal S} = \{ \Lambda \in M_4 ({\mathbb C}) \, ;
\ \Lambda \Lambda^* = \Lambda^* \Lambda = 1 \, , \ \Lambda^t =
\Lambda \} \, .
\end{equation}
We define for $\Lambda \in {\mathcal S}$ the algebra $C_{\rm alg}
({\mathbb R}^4 (\Lambda))$ of coordinates on ${\mathbb R}^4
(\Lambda)$  as the algebra generated by the $z^{\mu}$, $z^{\mu*}$
with the relations {\rm  (\ref{relright}), (\ref{relleft})} together
with,
\begin{equation}
\label{adjoint} z^{\mu*} = \Lambda_{\nu}^{\mu} \, z^{\nu}
\end{equation}
Note that (\ref{relcenter1}) follows automatically from
\eqref{adjoint}.
 For $S^3 (\Lambda)$ we add the inhomogeneous relation,
\begin{equation}
\label{eq1.32} \sum z^{\mu*} \, z^{\mu} = 1 \, .
\end{equation}
By Lemma \ref{centerC} the element $C = \sum z^{\mu *} \, z^{\mu}$
is in the center of the involutive algebra $C_{\rm alg} ({\mathbb
R}^4 (\Lambda))$.

\medskip

  \begin{thm} \label{univ} Let ${\mathcal A}$ be a unital involutive algebra and $U \in M_2 ({\mathcal A})$ a unitary such that ${\rm ch}_{1/2} (U) = 0$. Then there exists $\Lambda \in {\mathcal S}$ and a homomorphism $\varphi : C_{\rm alg} (S^3 (\Lambda)) \to {\mathcal A}$ such that:
\begin{equation}
U = \varphi ( \tau_\mu \,z^\mu ) \, . \nonumber
\end{equation}
Conversely for any $\Lambda$ the unitary $U =\,\tau_\mu \,z^\mu \in
M_2 (C_{\rm alg} (S^3 (\Lambda)))$ fulfills ${\rm ch}_{1/2} (U) =
0$.
  \end{thm}

\bigskip

By construction we thus obtain involutive algebras parametrized by
$\Lambda \in {\mathcal S}$. They are endowed with a canonical
element of $H_3 ({\mathcal A} , {\mathcal A})$ (in fact of $Z_3$)
given by
\begin{equation}
\label{eq1.33} [S^3 (\Lambda)] = {\rm ch}_{3/2} (U) = U_{i_1}^{i_0}
\otimes U_{i_2}^{*i_1} \otimes U_{i_3}^{i_2} \otimes U_{i_0}^{* i_3}
- U_{i_1}^{*i_0} \otimes U_{i_2}^{i_1} \otimes U_{i_3}^{*i_2}
\otimes U_{i_0}^{i_3} \, .
\end{equation}
The operations $U \to \lambda \, U$, $U \to V_1 U V_2$, $U \to U^*$,
for $\lambda \in U(1)$ and $V_j \in SU(2)$ together with the
universality described in Theorem \ref{univ} give natural
isomorphisms between the $S^3 (\Lambda)$ as follows,

\bigskip

  \begin{prop} \label{liniso} The following define isomorphisms $S^3 (\Lambda) \to S^3 (\Lambda')$ (resp. ${\mathbb R}^4(\Lambda) \to {\mathbb R}^4(\Lambda')$)
preserving $[S^3]$ (resp. $[{\mathbb R}^4]$) in the first two cases
and changing its sign in the last case:
\begin{enumerate}
\item[1)] For $\lambda \in U(1)$, $\Lambda' = \lambda^{2} \Lambda$ one lets
$$
\varphi (z_{\Lambda'}^{\mu}) = \lambda^{-1} \, z_{\Lambda}^{\mu} \,
.
$$
\item[2)] For $V \in SO(4)$, $\Lambda' = V \Lambda V^t$ one lets
$$
\varphi (z_{\Lambda'}^{\mu}) = V_{\nu}^{\mu} \, z_{\Lambda}^{\nu} \,
.
$$
\item[3)] For $\Lambda' = \varepsilon \, \Lambda^{-1} \, \varepsilon$ one lets
$$
\varphi (z_{\Lambda'}^{\mu}) = \varepsilon_{\mu} \,
z_{\Lambda}^{\mu*} \, , \ \varepsilon_0 = 1 \, , \ \varepsilon_k =
-1 \, , \ \varepsilon_{\mu\nu} = 0 \quad \mu \ne \nu \, .
$$
\end{enumerate}
\end{prop}

\bigskip

  \begin{proof} 1) Let $U = U_{\Lambda}$ be the unitary in $M_2 (C_{\rm alg} (S^3 (\Lambda))$. Then $\lambda \, U$ still fulfills ${\rm ch}_{1/2} (\lambda \, U) = 0$ and ${\rm ch}_{3/2} (\lambda \, U) = {\rm ch}_{3/2} (U)$. With $\tilde z^{\mu} = \lambda \, z^{\mu}$ one has $(\tilde z^{\mu})^* = \bar\lambda \, z^{\mu*} = \bar\lambda \, \Lambda_{\nu}^{\mu} \, z^{\nu} = \bar\lambda^2 \, \Lambda_{\nu}^{\mu} \, \tilde z^{\nu}$. This shows that $\varphi$ is a homomorphism
\begin{equation}
\varphi : C_{\rm alg} (S^3 (\Lambda')) \to C_{\rm alg} (S^3
(\Lambda)) \nonumber
\end{equation}
and that $\varphi ([S^3 (\Lambda')]) = [S^3 (\Lambda)]$.

\medskip

  2) One has ${\rm Spin}(4)= SU(2)\times SU(2)$ and the covering map $$\pi\,:\;{\rm Spin}(4)
= SU(2)\times SU(2)\to SO(4)$$ is given for any $(u,v)\in
SU(2)\times SU(2)$ and $\xi \in {\mathbb R}^4$  by $\pi(u,v) \xi
=\eta$ with
\begin{equation}
\label{spincover}
 \tau_\mu \, \eta^\mu =\,u \,(\tau_\mu \, \xi^\mu) \, v^* \, .
\end{equation}
This equality continues to hold for the natural complex linear
extension $ V_{\nu}^{\mu} \, z^{\nu}$  of $V=\,\pi(u,v)$ to
${\mathbb C}^4$ and it follows that with the notations of assertion
2) one has  $$\tau_\mu \, V_{\nu}^{\mu} \, z_{\Lambda}^{\nu} =\,
u\,U\,v^*$$ with $U = U_{\Lambda}$ as above.  The unitary
$u\,U\,v^*$ still fulfills ${\rm ch}_{1/2} (u\,U\,v^*) = 0$ and
${\rm ch}_{3/2} (u\,U\,v^*) = {\rm ch}_{3/2} (U)$ since in $M_2
({\mathbb C})$ it is the ordinary product and trace which are
involved in the formulas for ${\rm ch}_{k/2}$. Thus, with $\tilde
z^{\mu} = V_{\nu}^{\mu} \, z^{\nu}$ we just have to check that
$\tilde z^{\mu*} = \Lambda^{'\,\mu}_{\nu} \, \tilde z^{\nu}$. One
has $\Lambda' \, \tilde z = V \Lambda V^t \, V z = V \Lambda \, z =
V z^*$ and since $V$ has {\it real} coefficients this is $(V z)^* =
(\tilde z)^*$.

\medskip

  3) Let $U = U_{\Lambda}$ as above, then $U^*$ is still unitary and fulfills ${\rm ch}_{1/2} (U^*) = 0$, ${\rm ch}_{3/2} (U^*) = - \, {\rm ch}_{3/2} (U)$. It corresponds to $\tilde z^0 = z^{0\ast}$, $\tilde z^k = - z^{k\ast}$. Then $\tilde z^{\mu*} = \varepsilon_{\mu} \, z^{\mu} = \varepsilon_{\mu} (\Lambda^{-1})_{\nu}^{\mu} \, z^{\nu*} = \varepsilon_{\mu} (\Lambda^{-1})_{\nu}^{\mu} \, \varepsilon_{\nu} \, \tilde z^{\nu}$ which gives the value of $\Lambda'$. \end{proof}

  \begin{cor} \label{liniso1}
For every $\Lambda \in {\mathcal S}$ there exists $\varphi_j\in
{\mathbb R}/\,\pi{\mathbb Z}$ and isomorphisms  $S^3 (\Lambda) \to
S_{\varphi}^3$ (resp. ${\mathbb R}^4(\Lambda) \to {\mathbb
R}_{\varphi}^4$) where $S_{\varphi}^3$ corresponds to the diagonal
matrix $$\Lambda(\varphi) = \begin{bmatrix} 1 &0 \\ 0
&e^{-2i\varphi_k} \end{bmatrix}$$
\end{cor}

\begin{proof} We just need to recall why the matrix $\Lambda$ can be diagonalized
by a $W \in SO(4)$. In fact $\Lambda$ is unitary and fulfills
$\Lambda^t = \Lambda$ {\it i.e.} $\Lambda^* = J \Lambda \, J^{-1}$
where $J$ gives the real structure of ${\mathbb C}^4$. An eigenspace
$E_{\lambda} = \{ \xi \in {\mathbb C}^4 ; \Lambda \, \xi = \lambda
\, \xi \}$ is stable by $J$ since $\Lambda \, J^{-1} \, \xi = J^{-1}
\, \Lambda^* \, \xi = J^{-1} \, \bar\lambda \, \xi = \lambda \,
J^{-1} \, \xi$ (and $J^{-1} = J$). Thus we can find an orthonormal
basis of {\it real} vectors: $J \, \xi_i = \xi_i$ which are
eigenvectors for $\Lambda$. With $e_i$ the standard basis of
${\mathbb C}^4$ the map $e_i \to \xi_i$ gives an element of $O(4)$
and we can take it in $SO(4)$. Thus $\Lambda = W D \, W^t$ with $W
\in SO(4)$.
\end{proof}

The presentation of the algebra of ${\mathbb R}_{\varphi}^4$ is
given by (\ref{relright}) (\ref{relleft}) and   the relation
(\ref{lambda}) with $\Lambda = \Lambda(\varphi)$. This gives $z^{0*}
= z^0$ and $z^{k*} = e^{-2i\varphi_k} \, z^k$ so that with
$x^0=\,z^0$, $x^k = e^{-i \varphi_k} \, z^k$ we get
\begin{equation}
\label{pres0} x^{\mu\ast} = x^\mu \qqq \mu\in \{0,1,2,3\}\, .
\end{equation}
and the six other relations give
\begin{equation}
\label{firstpres} e^{i\varphi_k} \, x^k \, x^0 - e^{-i \varphi_k} \,
x^0 \, x^k + \sum \varepsilon_{k\ell m} \, e^{i(\varphi_{\ell} -
\varphi_m)} \, x^{\ell} \, x^m = 0
\end{equation}
\begin{equation}
\label{firstpres1} e^{i\varphi_k} \, x^0 \, x^k - e^{-i \varphi_k}
\, x_k \, x^0 + \sum \varepsilon_{k\ell m} \, e^{-i(\varphi_{\ell} -
\varphi_m)} \, x^{\ell} \, x^m = 0 \, .
\end{equation}
This gives, by combining (\ref{firstpres}) and (\ref{firstpres1})
the relations

\medskip
\begin{equation}
\label{pres1} \sin (\varphi_k) \, [x^0 , x^k]_+ = i\; \cos
(\varphi_{\ell} - \varphi_m) \, [x^{\ell} , x^m]
\end{equation}
\begin{equation}
\label{pres2} \cos (\varphi_k) \, [x^0 , x^k] = i \;\sin
(\varphi_{\ell} - \varphi_m) \, [x^{\ell} , x^m]_+ \, ,
\end{equation}

\medskip

where we let $ [a , b]_+=\,a\,b+\,b\,a$ be the anticommutator.

\medskip
We shall now study in much greater detail the corresponding moduli
space.

 \bigskip
 \bigskip

\section{The real moduli space ${\mathcal M}$}\label{modulispace}

We let ${\mathcal M}$ be the moduli space of noncommutative
$3$-spheres. It is obtained from proposition \ref{liniso} 1) and 2)
as the quotient
\begin{equation}
{\mathcal M}=\,(U(1) \times SO(4))\backslash {\mathcal S}\,,
\label{defmod}
\end{equation}
 of ${\mathcal S}$ by the action of $U(1) \times SO(4)$.
This  action of $U(1) \times SO(4)$ on ${\mathcal S}$ is the
restriction of the following action of $U(4)$ on ${\mathcal S}$.
\begin{equation}
\label{u4acts} \Lambda \in {\mathcal S}\to W \, \Lambda \, W^t
\qqq\, W \in U(4)
\end{equation}
which  allows to identify ${\mathcal S}$ with the homogeneous space
\begin{equation}
\label{hom1} U(4) / O(4) \simeq {\mathcal S} \, .
\end{equation}
(Note that any $\Lambda \in {\mathcal S}$ can be written as $\Lambda
= VV^t$ for some $V \in U(4)$ since it can be diagonalized by an
orthogonal matrix). The presence of $U(1)$ in \eqref{defmod} allows
to reduce to $SU(4)$ and one obtains this way a first convenient
description of ${\mathcal M}$.

\medskip

\subsection{${\mathcal M}$ in terms of $A_3$}

\noindent
\medskip

The description of ${\mathcal M}$ in terms of the compact Lie group
$SU(4)$ is given by :

\begin{prop} \label{Normalizer}
\begin{enumerate}
\item[1)] Let $N$ be the normalizer of $SO(4)$ in $SU(4)$. Then one has a
canonical isomorphism
\begin{equation}
\label{hom2} {\mathcal M} \simeq  SO(4)\backslash SU(4) / N  \, .
\end{equation}

\item[2) ] Let ${\mathbb T} \subset SU(4)$ be the maximal torus of diagonal
matrices, and $W\subset {\rm Aut}({\mathbb T} ) $ the corresponding
Weyl group. Let $D=\,{\mathbb T} \cap N$. Then  the above restricts
to an isomorphism

\begin{equation}
\label{hom3} {\mathcal M} \simeq  W\backslash {\mathbb T}/ D \, .
\end{equation}

\item[3)] The map $u \to u^2$ from  ${\mathbb T}$ to ${\mathbb T}$
induces an isomorphism
\begin{equation}
\label{hom4}
 {\mathcal M} \simeq  W\backslash {\mathbb T}/ D  \simeq \hbox{ Space of Conjugacy Classes in } PSU(4)\,.
\end{equation}
\end{enumerate}
\end{prop}

\begin{proof}
 1)  Let $Z$ be the center of $SU(4)$, it is generated by $i$ which has order 4 and contains $-1 \in SO(4)$. The normalizer $N$ is generated by $SO(4)$ and the element
\begin{equation}
\label{eq1.39} v = \,w\,\begin{bmatrix}-1 &&&0 \\ &1 \\ &&1 \\ 0
&&&1 \end{bmatrix} \, , \qquad w = e^{2\pi i /8} \, ,
\end{equation}
which implements the outer automorphism of $SO(4)$ and whose square
$ v^2 = i$ generates $Z$.

Let $X=U(1)\backslash {\mathcal S}$ considered as an homogeneous
space on $SU(4)$ using the action \eqref{u4acts} \ie $\Lambda \in
X\to W \, \Lambda \, W^t$.
 Given $\Lambda \in {\mathcal S}$ we can find $\lambda \in U(1)$ so that $\Lambda = \lambda \, \Lambda_1$ with ${\rm Det} \, \Lambda_1 = 1$, thus with
${\mathcal S}_1 = \{ \Lambda \in {\mathcal S} , {\rm Det} \, \Lambda
=1\}$ one has,
\begin{equation}
\label{eq1.38} X=\,U(1)\backslash {\mathcal S} \simeq w^{\mathbb
Z}\backslash{\mathcal S}_1 \simeq SU(4) /N\,.
\end{equation}
Indeed the first equality follows since the action \eqref{u4acts} of
$w$ is multiplication by $i$ and the second follows by computing the
isotropy group $K$ of $1\in w^{\mathbb Z}\backslash{\mathcal S}_1$.
One has $SO(4)\subset K$ and $v\in K$ since $v^t v=v^2=w^2$. Thus
$N\subset K$. Conversely given $V \in K$ one has $VV^t = i^N$ for
some $N\in {\mathbb Z}$ and thus for a suitable power of $v$ one has
$v^k \, V \, V^t \, v^k = 1$ thus $v^k \, V \in SO(4)$. This shows
that $K=N$ and we get the first statement of the proposition since
${\mathcal M} \simeq  SO(4)\backslash X$ by construction.

Note the standard description of $N$ which is given as follows. One
lets $\theta \in {\rm Aut} (SU(4))$ be given by complex conjugation,
\begin{equation}
\label{eq1.41} \theta (u) = \bar u = J \, u \, J^{-1} \qquad \forall
\, u \in SU(4) \, .
\end{equation}
One has $\theta^2 = 1$ and the fixed points of $\theta$ give $SO(4)
= SU(4)^{\theta}$. The normalizer $N$ of $SO(4)$ is characterized
by,
\begin{equation}
\label{normalizer} u^{-1} \, \theta (u) \in Z \, ,
\end{equation}
where $Z$ is the center of $SU(4)$.

\smallskip

2) Let $\sigma$ be the map $\sigma (u) = u \, u^t $ from $SU(4)$ to
the space $P{\mathcal S}_1$ of classes of elements of ${\mathcal
S}_1$ modulo the action of $Z$ by multiplication. It follows from
1) that $\sigma$ is an isomorphism of  $X=SU(4)/N$ with $P{\mathcal
S}_1$. Given $u \in SU(4)$ we can find $V \in SO(4)$ such that
\begin{equation}
\label{diag0} u \, u^t = V D \, V^t
\end{equation}
where $D$ is a diagonal matrix. Then $\sigma (u) = \sigma (V
D^{1/2})$ where $D^{1/2}$ is a diagonal square root of $D$ with
determinant equal to $1$. Thus every element of the  coset space
$SO(4) \backslash X$ can be represented by a diagonal matrix, and
the  natural map given by inclusion
\begin{equation}
\label{diag1} W\backslash {\mathbb T}/ D \to SO(4) \backslash X
\end{equation}
is surjective.

\medskip
3) When restricted to ${\mathbb T}$ the map $\sigma$ is simply the
squaring $u\to u^2$. Moreover the equality $\sigma(u_1)=\sigma(u_2)$
in $P{\mathcal S}_1$ for $u_j \in {\mathbb T}$ just means that the
$u_j^2$ define the same element of the maximal torus ${\mathbb T}/Z$
of $PSU(4)$. Thus the result follows. Let us check that the group $D
/Z$  is $({\mathbb Z} / 2)^3$. The elements of $D$ are the $v^k \,
u$ with $v$ as in (\ref{eq1.39}) and $u$ in $\mathbb T \cap SO(4)$.
One has $v^2 \in Z$ and modulo $Z \cap SO(4) = \pm 1$ the elements
of $\mathbb T \cap SO(4)$ form the Klein group $H=({\mathbb Z} / 2
{\mathbb Z})^2$. Thus one gets $D /Z\simeq({\mathbb Z} / 2 {\mathbb
Z})^3$.
\end{proof}

\medskip

  We identify the Lie algebra of $SU(4)$ with the Lie algebra of antihermitian matrices with trace 0,
\begin{equation}
\label{eq1.45} {\rm Lie} \, (SU(4)) = \{ T \in M_4 ({\mathbb C}) \,
; \ T^* = -T \, , \ {\rm Trace} \, T = 0 \}
\end{equation}
where $\theta$ is still acting by complex conjugation.

  The diagonal matrices $D \in {\mathcal D}$ form a maximal abelian Lie subalgebra of the eigenspace,
\begin{equation}
\label{eq1.46} {\rm Lie} \, (SU(4))^- = \{ T \, , \ \theta (T) = -T
\} \, .
\end{equation}
The roots $\alpha \in \Delta$ are given by,
\begin{equation}
\label{rootsa} \alpha_{\mu , \nu} (\delta) = \delta_{\mu} -
\delta_{\nu} \, .
\end{equation}

\begin{prop}
For $\delta \in {\mathcal D}$ one has,

\begin{equation}\label{latticeg}
e^{\delta} \in N \,\Leftrightarrow\, e^{2\delta} \in
Z\,\Leftrightarrow\, \alpha_{\mu , \nu} (\delta) \in i \, \pi \,
{\mathbb Z} \qqq \mu , \nu\,.
\end{equation}
\end{prop}

\medskip
\begin{proof} The equivalence between the last two conditions is a general
fact for compact Lie groups (\cf \cite{bou:1982}). The equivalence
between the first two conditions follows from the third statement of
proposition \ref{Normalizer}.
 \end{proof}
\bigskip

We let $\Gamma$ be the lattice $\Gamma \subset {\mathcal D}$
determined by the equivalent conditions \eqref{latticeg}, and
\begin{equation}
\label{torusta} {\mathbb T}_A = {\mathcal D} / \Gamma
\end{equation}
be the quotient 3-dimensional torus.

\medskip

  We let the group $W$ of permutations of 4 elements act on ${\mathbb T}_A$ by permutations of the $\delta_{\mu}$. In fact we view it as the Weyl group of the pair $(SU(4),N)$, {\it i.e.} as the quotient,
\begin{equation}
\label{eq1.51} W = {\mathcal N}/C
\end{equation}
of the normalizer of ${\mathcal D}$ in $SO(4)$ by the centralizer of
${\mathcal D}$. Note that $v$ being diagonal is in the centralizer
of ${\mathcal D}$ so that $W$ does not change in replacing $SO(4)$
by $N$ since $v^k \, u$ normalizes ${\mathcal D}$ iff $u$ does.

\bigskip

\begin{cor}
 The map $\sigma (\delta) = e^{2\delta}$ defines an isomorphism of the quotient of ${\mathbb T}_A$
by the action of $ W$ with the moduli space ${\mathcal M} = (U(1)
\times SO (4))\backslash{\mathcal S} $.
\end{cor}

\begin{proof} This is just another way to write  the
second statement of proposition \ref{Normalizer}. \end{proof}

\medskip

\subsection{ Trigonometric parameters $\varphi$ of $S^3_\varphi$}

\noindent
\medskip

We shall now describe a convenient parametrization of the torus
${\mathbb T}_A$ which gives
 the corresponding algebras in the form of corollary \ref{liniso1}.
It is given by the map
\begin{equation}
\label{phipar} \varphi = (\varphi_1 , \varphi_2 , \varphi_3) \in
({\mathbb R} / \pi \, {\mathbb Z})^3
 \to d (\varphi)=(\alpha_0,\alpha_0 - i \, \varphi_1, \alpha_0- i \, \varphi_2,\alpha_0 - i \, \varphi_3)\,,
\quad \alpha_0 = \frac{i}{4} \sum \varphi_j\,.
\end{equation}
 One has $\alpha_{0,k} (d (\varphi)) = i \, \varphi_k$ and the definition of $\Gamma$ shows that $d$ is an isomorphism. Also

$$e^{2 d (\varphi)}\simeq \begin{bmatrix} 1 &0 \\ 0 &e^{-2i\varphi_k} \end{bmatrix}$$

up to multiplication by a scalar.

\medskip
In terms of the parameters $\varphi_j$ the twelve roots $\alpha \in
\Delta$ are the following linear forms,
\begin{equation}
\label{rootphi} (\varphi_1 , \varphi_2 , \varphi_3) \to i\,\{\pm
\varphi_j , \varphi_k - \varphi_l \} \, .
\end{equation}

\medskip
The action of the Weyl group $W$ gives the following linear
transformations of the $\varphi_j$. Arbitrary permutations of the
$\varphi_j$'s correspond to permutations of
 the last three $\alpha_k$'s. The transposition of $\alpha_0$ with $\alpha_1$ corresponds to :
\begin{equation}
\label{eq1.53} T_{01} (\varphi_1 , \varphi_2 , \varphi_3) =
(-\varphi_1 , \varphi_2 - \varphi_1 , \varphi_3 - \varphi_1) \, .
\end{equation}

\medskip

The 3-spheres $S_{\varphi}^3$ are parametrized by
\begin{equation}
\label{eq1} \varphi = (\varphi_1 , \varphi_2 , \varphi_3) \in
({\mathbb R} / \pi {\mathbb Z})^3
\end{equation}
modulo the action of the Weyl group $W \simeq S_4$ generated by the
permutation group $S_3$ of the $\varphi_j$'s and
\begin{equation}
\label{eq1.69} (\varphi_j) \to (- \, \varphi_1 \, , \ \varphi_3 -
\varphi_1 \, , \ \varphi_2 - \varphi_1) = (\varphi'_j) \, .
\end{equation}

\medskip
\subsection{${\mathcal M}$ in terms of $D_3$}\label{alcoves}

\noindent
\medskip

To obtain another very convenient parametrization of the torus
${\mathbb T}_A$ we use the isomorphism $A_3 \sim D_3$, {\it i.e.}
\begin{equation}
\label{eq1.54} SU(4) \simeq {\rm Spin} \, (6) \, .
\end{equation}
which  simply comes from the spin representation of ${\rm Spin} \,
(6)$. The Clifford algebra ${\rm Cliff}_{\mathbb C} ({\mathbb R}^6)$
has dimension $2^6$ and is a matrix algebra $M_n({\mathbb C})$ with
$n = 2^3 = 8$. We let $\gamma^{\mu}$ be the corresponding
$\gamma$-matrices, with

 $$\gamma_{\mu}^* = \gamma_{\mu}\,,\quad \gamma_{\mu}\gamma_{\nu}+
\gamma_{\nu}\gamma_{\mu} = \,2\,\delta_{\mu,\nu}$$ We then let
$$\sigma_{\mu\nu} = \frac{1}{2} \, (\gamma_{\mu} \, \gamma_{\nu} - \gamma_{\nu} \, \gamma_{\mu})$$
which span a real subspace of the Clifford algebra stable under
bracket and isomorphic (up to a factor of $2$) to the Lie algebra of
$SO(6)$ of real antisymmetric $6$ by $6$ matrices.

\medskip
\begin{prop}
\begin{enumerate}
\item[1)] The following gives a parametrization $\tau$ of
a maximal torus in ${\rm Spin}(6)$,
\begin{equation}
\label{sigmapar} \theta=(\theta_j) \in (\Rb/\,2\pi\Zb)^3
\overset{\tau}{\longrightarrow} {\rm Exp}{\sum \theta_j
\sigma_{2j-1,2j}} \in {\rm Spin}(6)\subset {\rm Cliff}_{\mathbb C}
({\mathbb R}^6)
\end{equation}

\item[2)] The half Spin representation  gives an isomorphism
$\pi \;:\;{\rm Spin}(6)\to SU(4)$.

 \item[3)] One has $\pi(\tau(\theta))=\, e^{\delta_{\theta}}$ with $\delta_{\theta}$ diagonal given by
\begin{equation}
\label{thetapar} \delta_{\theta} = i \, (\theta_1 + \theta_2 +
\theta_3 \, , \ \theta_1 - \theta_2 - \theta_3 \, , \ - \, \theta_1
+ \theta_2 - \theta_3 \, , \ - \, \theta_1 - \theta_2 + \theta_3) \,
.
\end{equation}

\end{enumerate}
\end{prop}

\begin{proof}
This is a straightforward check, since in the half Spin
representation $\pi$ one has in a suitable basis
\begin{equation}
\label{eq1.59} \sigma_{12} \to \begin{bmatrix} i &&&0 \\ &i \\ &&-i
\\ 0 &&&-i \end{bmatrix} , \quad \sigma_{34} \to \begin{bmatrix} i
&&&0 \\ &-i \\ &&i \\ 0 &&&-i \end{bmatrix} , \quad \sigma_{56} \to
\begin{bmatrix} i &&&0 \\ &-i \\ &&-i \\ 0 &&&i \end{bmatrix} ,
\end{equation}
thus $\pi(\tau (\theta)) = e^{\delta_{\theta}}$ with
$\delta_{\theta}$ given by \eqref{thetapar}.
\end{proof}

Thus the trigonometric parameters $\varphi_k$ are given in terms of
the $\theta$'s by,
\begin{equation}
\label{eq1.61} \varphi_1 = 2 \, (\theta_2 + \theta_3) \, , \
\varphi_2 = 2 \, (\theta_1 + \theta_3) \, , \ \varphi_3 = 2 \,
(\theta_1 + \theta_2) \, .
\end{equation}
The natural parameters for the maximal torus ${\mathbb T}_D$ of
$SO(6)$ are the $\psi_j = 2 \, \theta_j$ and  are  defined modulo $2
\, \pi \, {\mathbb Z}$, {\it i.e.} correspond to the Lie algebra
element
\begin{equation}
\label{maxso} \ell (\psi)=\,\psi_1 \, \beta_{12} + \psi_2 \,
\beta_{34} + \psi_3 \, \beta_{56}
\end{equation}
 where the $\beta_{ij}$
form the canonical basis of real antisymmetric matrices. The kernel
of the covering ${\rm Spin} \, (6) \to SO(6)$ corresponds to
$\theta_j = \pi$ so that for $SO(6)$ the torus ${\mathbb T}_D$ is
parametrized by the $\psi_j$'s defined modulo $2\pi$ by
\begin{equation}
\label{eq1.62} \psi \in ({\mathbb R} / 2 \, \pi \, {\mathbb Z})^3
\to e^{\ell (\psi)} \, .
\end{equation}

The transition from $\psi$ to $\varphi$'s is given then by,
\begin{equation}
\label{phitopsi} \varphi_j = \psi_k + \psi_{\ell} \, , \ 2 \, \psi_j
= \varphi_k + \varphi_{\ell} - \varphi_j \, .
\end{equation}
as well as
\begin{equation}
\label{eq5} \varphi_1 - \varphi_2 = \psi_2 - \psi_1 \, , \ \varphi_2
- \varphi_3 = \psi_3 - \psi_2 \, , \ \varphi_3 - \varphi_1 = \psi_1
- \psi_3 \, .
\end{equation}

We shall now spell in great details the basic Lie group datas for
$D_3$ and get a description of the moduli space $\calm$ in these
terms. Proposition \ref{Normalizer} 3) gives a natural isomorphism
 $${\mathcal M} \simeq   \hbox{ Space of Conjugacy Classes in } PSU(4)$$
of the moduli space ${\mathcal M}$ with the space of conjugacy
classes of elements of $PSU(4)\simeq PSO(6)$. The general theory of
compact Lie groups provides a natural triangulation of such a space
of conjugacy classes in terms of alcoves. The latter are obtained as
the connected components of the complement of the union of the
singular hyperplanes.

\medskip Our aim in this section is to describe
such a triangulation in our specific case and to exhibit the role of
the singular hyperplanes. We shall see in the next section their
natural compatibility with the scaling foliation.

\medskip
In terms of the parameters  $\varphi_j$ of the 3-spheres
$S_{\varphi}^3$ the relations which  specify the non generic
situations, are all of the form
\begin{equation}
\label{eq6} \left\{ \varphi , \alpha (\varphi) = n \, \frac{\pi}{2}
\right\} = G_{\alpha , n}
\end{equation}
where $n$ is an integer and $\alpha$ is one of the twelve ``roots"
\ie $\alpha\in  \Delta=\pm\{\varphi_1 , \varphi_2 , \varphi_3 ,
\varphi_1 - \varphi_2 , \varphi_2 - \varphi_3 , \varphi_3
-\varphi_1\}$.

\smallskip

Moreover the periodicity lattice of $\varphi$ is $(\pi {\mathbb
Z})^3$ which is specified by
\begin{equation}
\label{eq7} \{ \varphi , \alpha (\varphi) \in \pi {\mathbb Z} \, , \
\forall \, \alpha \in \Delta \} = \Gamma_{\varphi} \, .
\end{equation}
We now want to relate more precisely the above situation with
canonical objects (root systems, alcoves, chambers, affine Weyl
group, nodal vectors $\ldots$) associated to the following data
$(G,\Tg)$
\begin{equation}
\label{eq8} G = PSO (6) \, , \ \Tg = \hbox{Maximal torus} \;{\mathbb
T}_D/\pm 1\  \, .
\end{equation}
We use the natural parametrization of the Lie algebra ${\rm Lie} \,
(\Tg)$,
\begin{equation}
\label{eq9} \ell (\xi)=\,\xi_1 \, \beta_{12} + \xi_2 \, \beta_{34} +
\xi_3 \, \beta_{56}  \, .
\end{equation}
Since we used the ``squaring" $u\to u^2$ in the isomorphism of
proposition \ref{Normalizer} 3), the parameter $\psi$ that appears
in the transition $\varphi \to \psi$ of equation \eqref{phitopsi} is
related to $\xi$ by
\begin{equation}
\label{doubling} \xi=\,2\,\psi  \, .
\end{equation}
In other words the natural relation between the parameters $\varphi$
and $\xi$ is
\begin{equation}
\label{phitoxi} 2\,\varphi_j = \xi_k + \xi_{\ell} \, ,  \, \xi_j =
\varphi_k + \varphi_{\ell} - \varphi_j \, .
\end{equation}

This accounts for a factor of $\frac{1}{2}$ in \eqref{eq6} but does
not yet relate it to the equation of singular hyperplanes. To
understand this relation more precisely we shall now review briefly
the standard ingredients of the theory of alcoves for the specific
data $(G,\Tg)$.
\medskip

\subsection{ Roots $\Delta = R(G,\Tg)$}

\noindent
\medskip

\smallskip They are by definition the linear forms $\alpha$ on ${\rm Lie} \, (\Tg)$ given by eigenvalues $X_{\alpha} \in {\rm Lie} \, G_{\mathbb C}$ which fulfill:
\begin{equation}
\label{eq10} [\xi,X_{\alpha}] = \alpha (\xi) \, X_{\alpha} \qquad
\forall \, \xi \in {\rm Lie} \, \Tg \subset  {\rm Lie} \, G \, .
\end{equation}
They are complex valued as defined. In our case they are given by
(up to multiplication by $i$)
\begin{equation}
\label{eq11} (\pm \, e_{\mu} \, \pm' e_{\nu}) \, \xi = \pm \,
\xi_{\mu} \, \pm' \xi_{\nu} \, .
\end{equation}

\medskip

\subsection{Singular hyperplanes $H_{\alpha , n}$}

\noindent
\medskip

\smallskip They are given by a root $\alpha \in \Delta$ and $n \in {\mathbb Z}$, with
\begin{equation}
\label{eq12} H_{\alpha , n} = \{ \xi , \alpha (\xi) = 2 \pi in \} \,
.
\end{equation}
As we shall explain below while these hyperplanes suffice to obtain
a triangulation of the space of conjugacy classes of the simply
connected covering of $G$ we shall need the additional ones with $n
\in \frac{1}{2}{\mathbb Z}$ to describe the space of conjugacy
classes in $G$ itself.
\medskip

\subsection{ Kernel of the exponential map: $\Gamma (\Tg)$ (nodal group of $\Tg$)}

\noindent
\medskip

\smallskip In our case we are dealing with $G = PSO(6)$ and thus for $\xi \in {\rm Lie} \, \Tg$, $e^\xi$ is 1 in $G$ iff $Ad (e^\xi) = 1$ since the center of $G$ is $C(G) = \{ 1 \}$. This means exactly that all eigenvalues of $ad (\xi)$ belong to ${\rm Ker} \, (\exp) = 2 \pi i {\mathbb Z}$, thus
\begin{equation}
\label{eq13} \Gamma (\Tg) = \{ \xi , \alpha (\xi) \in 2 \pi i
{\mathbb Z} \quad \forall \, \alpha \in \Delta \} \, .
\end{equation}
Thus, after applying the change of variables \eqref{phitoxi}
\begin{equation}
\label{eq14} \Gamma_{\varphi} \simeq \, \Gamma (\Tg)
\end{equation}
which shows that the periodicity lattice we have is  the nodal group
of $\Tg$.

\medskip

\subsection{ Group of nodal vectors $N(G,\Tg) \subset \Gamma (\Tg)$}

\noindent
\medskip

\smallskip This group can be defined in terms of vectors $K_{\alpha}$ which are associated to the roots $\alpha \in \Delta$ but in our case it is simpler to use the definition as follows:
\begin{equation}
\label{eq15} N(G,\Tg) = \hbox{Kernel of} \ \exp : {\rm Lie} \, \Tg
\to \tilde G = \hbox{Universal cover of}\; G\,.
\end{equation}
In our case $\tilde G = {\rm Spin} \, 6$ and in terms of
$\gamma$-matrices the exponential map takes the form,
\begin{eqnarray}
\label{eq16}
\ell (\xi) \to &&\left( \cos \frac{1}{2} \, \xi_1 + \sin \frac{1}{2} \, \xi_1 \, \gamma_1 \gamma_2 \right) \left( \cos \frac{1}{2} \, \xi_2 + \sin \frac{1}{2} \, \xi_2 \, \gamma_3 \gamma_4 \right) \nonumber \\
&&\left( \cos \frac{1}{2} \, \xi_3 + \sin \frac{1}{2} \, \xi_3 \,
\gamma_5 \gamma_6 \right)
\end{eqnarray}
whose kernel is given by
\begin{equation}
\label{eq17} N(G,\Tg) = \{ \xi \, ; \ \xi_j = 2\pi n_j \, , \quad
n_j \in {\mathbb Z} \, , \quad n_j \ \hbox{\it even} \} \, .
\end{equation}
One has $N(G,\Tg) \subset \Gamma (\Tg)$ and the quotient, of order
4, is generated by $\xi = (\pi , \pi , \pi) \in \Gamma (\Tg)$.

\medskip

\subsection{Affine Weyl group $W_a$}

\noindent
\medskip

\smallskip The affine Weyl group $W_a$ is the group generated by the reflexions associated to the hyperplane $H_{\alpha , n}$
for $\alpha \in \Delta$ and $n \in {\mathbb Z}$.

\smallskip

One has \cite{bou:1968}  (Chapter VI, proposition 1, page 173)
\begin{equation}
\label{eq18} W_a = N(G,\Tg) \rtimes W
\end{equation}
where the Weyl group is generated by the reflexions associated to
singular hyperplanes $H_{\alpha , 0}$.

\smallskip

In our case $W = S_4$ and all its elements are of the form
\begin{equation}
\label{eq19} W = \varepsilon \, \sigma \, , \quad \sigma \in S_3 \,
, \quad \varepsilon = (\varepsilon_i) \, , \quad \varepsilon_i \in
\pm 1 \, , \quad \prod \, \varepsilon_i = 1
\end{equation}
where the action on $\xi$ is by permutation of the $\xi_j$ for
$\sigma$ (careful that $(\sigma \xi)(i) = \xi (\sigma^{-1} (i))$ to
get a covariant action) and by multiplication by $\varepsilon_i$ for
$\varepsilon$.

\smallskip

For $N(G,\Tg)$ we can check that the $\check\alpha = K_{\alpha}$ are
simply given by the vectors $\pm \, e_{\mu} \pm' e_{\nu}$ which
correspond to
\begin{equation}
\label{eq20} (\pm \, 2 \pi , \pm' \, 2\pi , 0) = K_{\pm e_1 \pm'
e_2} \, , \quad \langle K_{\alpha} , \alpha \rangle = 2 \, .
\end{equation}

\bigskip

\subsection{Affine Weyl group $W'_a$}

\noindent
\medskip

\smallskip It is by definition the semi direct product,
\begin{equation}
\label{eq21} W'_a = \Gamma (\Tg) \rtimes W \, .
\end{equation}
What matters is that it still acts on the set of singular
hyperplanes. For $\gamma \in \Gamma (\Tg)$ one has $\alpha (h +
\gamma) = \alpha (h) + \alpha (\gamma)$ and $\alpha (\gamma) \in
2\pi i{\mathbb Z}$ thus one is just shifting the $n$ in $H_{\alpha ,
n}$.

\smallskip

From \cite{bou:1982} proposition 2 (Chapter 9, page 45)  $W_a$ is a
normal subgroup of $W'_a$.

\bigskip

\subsection{ Alcoves and fundamental domain}

\noindent
\medskip

 The alcoves are the connected components of $(U H_{\alpha , n})^c \subset {\rm Lie} \, \Tg$. The chambers are the components of $(U H_{\alpha , 0})^c$.

\smallskip

By construction the alcoves are intersections of half spaces and are
thus convex polyhedra.

\smallskip

By \cite{bou:1982} the affine Weyl group $W_a$ acts simply
transitively on $\Sigma=$ the set of alcoves. Since $W'_a$ is still
acting on $\Sigma$ we can identify $\Sigma$ with the homogeneous
space
\begin{equation}
\label{eq22} \Sigma = W'_a / H_X
\end{equation}
where $H_X$ is the finite isotropy group of an alcove $X$.

\smallskip

In our case, we take the following alcove :
\begin{equation}
\label{eq23} X = \{ \xi , \xi_1 + \xi_2 \geqq 0 , \xi_2 - \xi_1
\geqq 0 , \xi_3 - \xi_2 \geqq 0 , \xi_2 + \xi_3 \leqq 2\pi \} \, .
\end{equation}

\smallskip

It is a tetrahedron with all 4-faces congruent but 2 long edges and
4 short edges.

\begin{lem} The isotropy subgroup $H_X \subset W'_a$ of $X$ is generated by $w_1 = \bigl( \left( \frac{1}{2} , \frac{1}{2} , \frac{1}{2} \right)$, $\varepsilon_{12} \, \sigma_{13} \bigl)$.
\end{lem}

\bigskip

\begin{proof} For  convenience we rescale the $\xi_j$ by $2\pi$. The coordinates of the vertices of $X$ are then $0$, $p=\left( \frac{1}{2} , \frac{1}{2} , \frac{1}{2} \right)$, $q=\left( -\frac{1}{2} , \frac{1}{2} , \frac{1}{2} \right)$, $p'=(0,0,1)$.

\smallskip

One has $w_1 (0) = p$, $w_1 (p) = p + \varepsilon_{12} \,
\sigma_{13} (p) = p + \left( - \frac{1}{2} , - \frac{1}{2} ,
\frac{1}{2} \right) = p'$, $w_1 (p') = p + \varepsilon_{12} \,
\sigma_{13} (p') = p + (-1,0,0) = q$, $w_1 (q) = p +
\varepsilon_{12} \, \sigma_{13} (q) = \left( \frac{1}{2} ,
\frac{1}{2} , \frac{1}{2} \right) + \left( -\frac{1}{2} ,
-\frac{1}{2} , -\frac{1}{2} \right) = 0$.
\end{proof}
\bigskip

We let $\sigma$ be the orthogonal reflexion around the face $(pqp')$
of $X$. One has $\sigma \in W_a$ and it is given explicitly by
\begin{equation}
\label{eq25} \sigma = ((0,1,1), \varepsilon_{23} \, \sigma_{23}) \,
.
\end{equation}
Indeed, since the face is given by $\xi_2 + \xi_3 = 2\pi$ the
corresponding nodal vector is $(0,1,1)$. One checks that $\sigma^2 =
1$ and that $\sigma$ fixes $p,q,p'$, $(\sigma (p') = (0,1,1) +
(0,-1,0) = (0,0,1)=\,p')$.

We let $Y = \sigma (X)$ be the reflexion of $X$ along that face
which yields the convex pentahedra $X \cup Y$.

\bigskip

\begin{figure}
  \qquad \qquad
\includegraphics[scale=0.6]{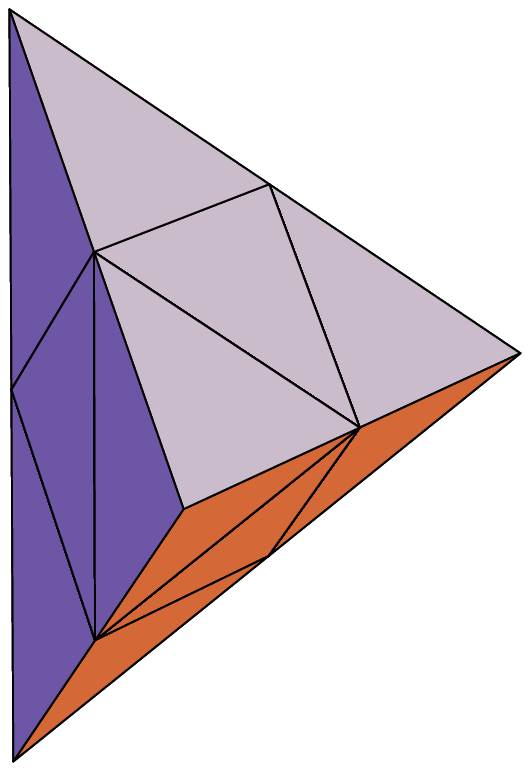}\qquad \qquad
\includegraphics[scale=0.6]{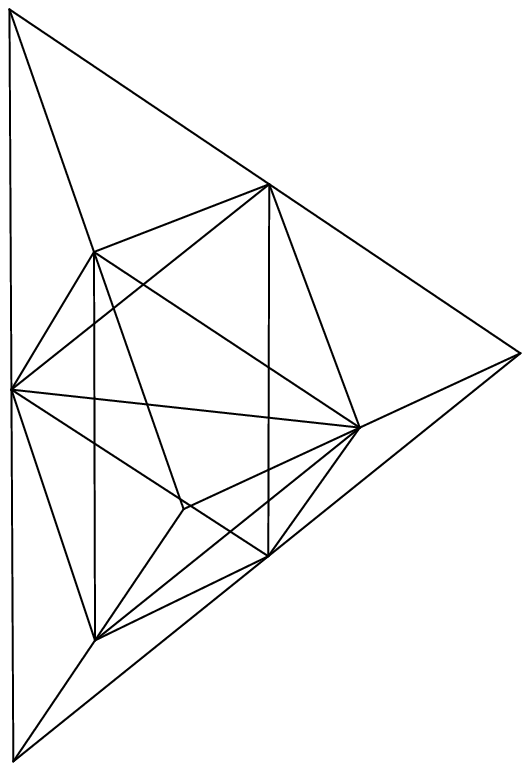}
\caption{\label{domain1} Tiling of the alcove $X$ by
$\frac{1}{2}$-alcoves}
\end{figure}

\bigskip

\begin{prop} Let $X$ be the alcove given by
\eqref{eq23} and $Y = \sigma (X)$ its reflexion along the face
$(p,q,p')$. Then
 $\frac{1}{2}(X \cup Y)$ is a fundamental domain for the action of $\Gamma (\Tg) \rtimes W$ in ${\rm Lie} \, \Tg$.
\end{prop}

\bigskip

\begin{proof} By the above lemma the action of $W'_a$ on the set of alcoves is the same as the action of $W'_a$ on $W'_a / H_X$. Let $W''_a = 2 \Gamma \rtimes W$. The left coset space $W''_a / W'_a$ is identified with the 8 elements set, $\Gamma / 2\Gamma$ where the corresponding map is given by:
\begin{equation}
\label{eq26} (\gamma , w) \in W'_a \to \hbox{Class of $w^{-1}
(\gamma)$ in $\Gamma / 2\Gamma$} \, .
\end{equation}
Indeed $(0,w^{-1}) (\gamma , w) = (w^{-1} (\gamma),1)$. We can thus
display the action (on the right) of $H_X$ as $\gamma \to
\hbox{Class of} \ \sigma_{13} \, \varepsilon_{12} (p+\gamma)$. Let
us write this transformation in terms of the $\varphi$-coordinates.
One obtains
\begin{equation}
\label{eq27} \varphi \to w \left( \varphi + \left( \frac{1}{2} ,
\frac{1}{2} , \frac{1}{2} \right)\right) \, , \quad w(\varphi) =
(-\varphi_3 , \varphi_1 - \varphi_3 , \varphi_2 - \varphi_3) \, .
\end{equation}
Thus in fact we look at the following transformation of $({\mathbb
Z} / 2)^3$,
\begin{equation}
\label{eq28} S (a_1 , a_2 , a_3) = \left( \frac{1}{2} - a_3 , a_1 -
a_3 , a_2 - a_3 \right) \, .
\end{equation}
One has $S(0) = \left( \frac{1}{2} ,0,0 \right)$, $S \left(
\frac{1}{2} , 0,0 \right) = \left( \frac{1}{2} , \frac{1}{2} , 0
\right)$, $S \left( \frac{1}{2} , \frac{1}{2} , 0 \right) = \left(
\frac{1}{2} , \frac{1}{2} , \frac{1}{2} \right)$, $ S\left(
\frac{1}{2} , \frac{1}{2} , \frac{1}{2} \right)$ $= 0$. The other
orbit is $S \left( 0,\frac{1}{2} , 0 \right) = \left( \frac{1}{2} ,
0 , \frac{1}{2} \right)$, $S \left( \frac{1}{2} , 0 , \frac{1}{2}
\right) = \left( 0,0,\frac{1}{2} \right)$, $S(0,0,\frac{1}{2}) =
\left( 0,\frac{1}{2} , \frac{1}{2} \right)$, $S \left( 0,
\frac{1}{2} , \frac{1}{2} \right) = \left( 0,\frac{1}{2} , 0
\right)$.

\smallskip

This shows that the double coset space $W''_a \backslash W'_a / H_X$
has cardinality 2 and thus that $W''_a$ just has 2 orbits in its
action on the set of alcoves. What remains is to show that $Y \notin
$ Orbit of $X$. But one has $Y = \sigma (X)$ with $\sigma$ given by
(\ref{eq25}). Thus we just need to determine the double coset of
$\sigma$, {\it i.e.} by (\ref{eq26}) the class of $\sigma_{23} \,
\varepsilon_{23} (0,1,1)$ in $\Gamma / 2\Gamma$. One just checks
that it is in the other orbit. \end{proof}

\bigskip

\smallskip

\begin{figure}
   \qquad \qquad
\includegraphics[scale=0.6]{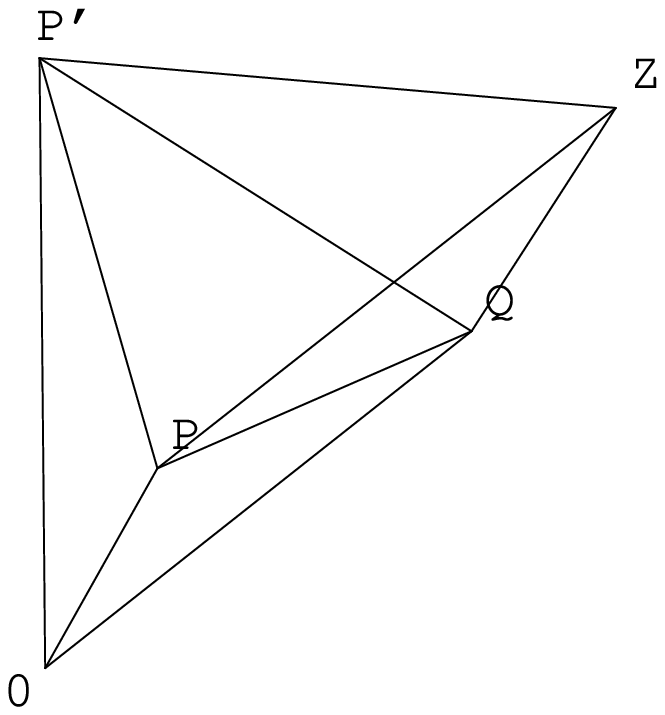}\qquad \qquad
\includegraphics[scale=0.6]{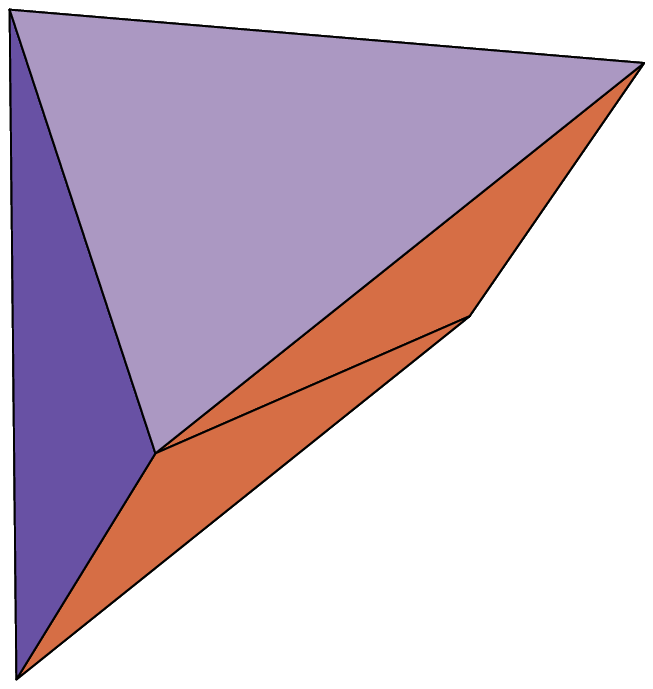}
\caption{\label{domain2} Fundamental domain }
\end{figure}

\bigskip

 We thus have a fairly simple fundamental domain for the parameter space of noncommutative $3$-spheres.
In terms of the $\varphi$'s a simple  translation of the above
result gives,

\begin{prop} \label{fundom}
1) The union $A \cup B$ of the following simplices
\begin{equation}
\label{eq29} A = \left\{ \varphi \, ; \ \frac{\pi}{2} \geqq
\varphi_1 \geqq \varphi_2 \geqq \varphi_3 \geqq 0 \right\}\,,\quad B
= \left\{ \varphi \, ; \ \varphi_3 + \frac{\pi}{2} \geqq \varphi_1
\geqq \frac{\pi}{2} \geqq \varphi_2 \geqq \varphi_3 \right\} \,,
\end{equation}
gives a fundamental domain for the action of $W$ on
$\Rb^3/\Gamma_\varphi=\,(\Rb/\pi \Zb)^3$.

2) The real moduli space ${\mathcal M}$ is obtained by glueing the
face\footnote{we use capital letters $P$, $Q$, $Q'$ for the vertices
:
 $P=(\frac{\pi}{2},\frac{\pi}{2},\frac{\pi}{2})$,
$P'=(\frac{\pi}{2},\frac{\pi}{2},0)$, $Q=(\frac{\pi}{2},0,0)$,
$Z=P+Q$. } $(P'QZ)$  to the face $(ZPP')$ by the  transformation
$\gamma \in\Gamma_\varphi\rtimes W$
$$
\gamma(\varphi_1,\varphi_2,\varphi_3)=\,(\pi-\varphi_1+\varphi_2,
\,\pi-\varphi_1+\varphi_3,\,\pi-\varphi_1) \label{Gamma}
$$
and crossing each wall, \ie each face whose supporting hyperplane
contains $0$, by the corresponding reflexion (in $W$).
\end{prop}
\smallskip

The two simplices $A$ and $B$ have a common face $(PQP')$ supported
by the hyperplane $\varphi_1=\,\frac{\pi}{2}$. In order to describe
the moduli space one needs to give the identifcations of the
boundary components. All faces of $A$ other than $(PQP')$ are walls
of chambers and thus crossing them leads to a simple reflexion on
that wall.

\bigskip

\section{The flow $F$} \label{triang}

While the space ${\mathcal M}$ is the natural moduli space for
noncommuative $3$-spheres, the moduli space of the corresponding
$4$-spaces ${\mathbb R}_{\varphi}^4$ is obtained as a space of flow
lines for a natural flow $F$ on ${\mathcal M}$.
  We define this flow $F$ as the gradient flow for the Killing metric on the Lie algebra of $SO(6)$ of the character of the virtual representation given by the
``signature". We first show that $F$ is nicely compatible with the
triangulation by alcoves and then check that the isomorphism class
of the $4$-spaces ${\mathbb R}_{\varphi}^4$ is constant along the
flow lines. We shall need for the converse to have a complete
knowledge of the geometric datas of these quadratic algebras, and
this will be obtained in section \ref{geometricdata}. Thus the
converse will be proved later on in section \ref{flowbetter}.

\subsection{Compatibility of $F$
with the triangulation by alcoves}

\noindent
\medskip

The basic compatibility of the flow $F$ with the structure of the
real moduli space ${\mathcal M}$ is given by :

\begin{prop}

\begin{enumerate}
\item[a)]
The character of the signature representation is
\begin{equation}
\label{chi} \chi( \xi)=\,{\rm Trace} \, (* \, \pi (\ell (\xi)) = -
\, 8 \prod \sin\, \xi_j \, .
\end{equation}
\item[b)] The flow $F=\nabla \;\chi(2 \psi)$ is invariant
under the action of the Weyl group $W$ and leaves each of the
singular hyperplanes $H_{\alpha,n}$, $n\in \frac{1}{2}\mathbb Z$
globally invariant.

\item[c)] In terms of the variables $\varphi_j$ one has
\begin{equation}
\label{flowphi} F=\,  \sum \,\sin (2  \varphi_j) \,\sin (\varphi_k +
\varphi_{\ell} - \varphi_j) \, \frac{\partial}{\partial \varphi_j}\,
.
\end{equation}
\end{enumerate}

\end{prop}

\medskip

\begin{proof} a) We parametrize the maximal torus
${\mathbb \Tg}_D$ of $SO(6)$ by the $\xi_j$ as in \eqref{maxso} \ie
$$\ell (\xi)=\,\xi_1 \, \beta_{12} + \xi_2 \, \beta_{34} + \xi_3 \,
\beta_{56}\,.$$ The Killing metric is simply given there (up to
scale and sign) by $\sum d \xi_j^2$.
 We identify $\wedge^3 \, {\mathbb C}^6$
with the linear span of the $e_{ijk} = \gamma_i \, \gamma_j \,
\gamma_k$, ${\rm card} \, \{ i,j,k \} = 3$ in ${\rm Cliff}_{\mathbb
C} ({\mathbb R}^6)$. We then use $\gamma_7 = \gamma_1 \, \gamma_2 \,
\gamma_3 \, \gamma_4 \, \gamma_5 \, \gamma_6$, which fulfills
$\gamma_7^2 = -1$ to define the $*$ operation by:
\begin{equation}
\label{eq1.65}
* \, e_{ijk} = \gamma_7 \, e_{ijk} \, .
\end{equation}
 One checks  that the matrix of $* \, \pi (\ell (\xi))$
 restricted to the 12 dim subspace spanned by the $e_{ijk}$ where two indices belong to one of the 3 subsets $\{ 1,2 \}$, $\{ 3,4 \}$, $\{ 5,6 \}$ is off diagonal
and hence has vanishing trace. Indeed for instance $* \, \pi (\ell
(\xi)) \, e_{123}$ is a linear combination of $e_{356}$ and
$e_{456}$ which are orthogonal to $e_{123}$.

0n the 8 dimensional subspace of $e_{ijk}$, $i \in \{ 1,2 \}$, $j
\in \{ 3,4 \}$, $k \in \{ 5,6 \}$ one gets the product of what
happens in the 2-dimensional case with a basis $e_1 , e_2$ of
$\wedge^1 \, {\mathbb C}^2$.
 One has $* \, e_1 = e_2$, $* \, e_2 = \, -e_1$ and the representation is given by $e_1 \to \cos \xi \, e_1 + \sin \xi \, e_2$, $e_2 \to \cos \xi \, e_2 - \sin \xi \, e_1$ so that
 the trace of $* \, \pi (\ell (\xi))$ is $-2 \sin \xi$. Thus we get,
\begin{equation}
\label{eq1.66} {\rm Trace} \, (* \, \pi (\ell (\xi))) = - \, 8
\;\sin \,\xi_1\; \sin \,\xi_2 \;\sin \,\xi_3 \, .
\end{equation}

b)  In terms of $SO(6)$ the Weyl group $W$ maps to the permutations
of $(\psi_1 , \psi_2 , \psi_3)$ and the kernel of this map is the
Klein subgroup which is given by the transformation,
\begin{equation}
\label{eq1.68} \psi_j \to \varepsilon_j \, \psi_j  \, \quad
\prod_1^3 \varepsilon_j = 1 \qquad (\varepsilon_j \in \{ \pm 1 \}).
\end{equation}
By construction the function $\chi(\xi)$ being a (virtual) character
is invariant by these transformations and so is the flow $X$.

c) Let us rewrite $X = \sum \frac{\partial \chi}{\partial \psi_j} \,
\frac{\partial}{\partial \, \psi_j}$ in terms of the coordinates
$\varphi_j$. One has
$$d \, \varphi_j = d \, \psi_k + d \, \psi_{\ell}\,,\quad
\sum \frac{\partial h}{\partial \varphi_j} \, d \, \varphi_j = \sum
\left( \frac{\partial h}{\partial \varphi_k} + \frac{\partial
h}{\partial \varphi_{\ell}} \right) d \, \psi_j\,,\quad
\frac{\partial h}{\partial \psi_j} = \frac{\partial h}{\partial
\varphi_k} + \frac{\partial h}{\partial \varphi_{\ell}}$$ and $X$ is
given by
\begin{eqnarray}
X = \sum \left( \frac{\partial \chi}{\partial \psi_k} + \frac{\partial \chi}{\partial \psi_{\ell}} \right) \frac{\partial}{\partial \varphi_j} &= &2 \sum (\cos 2 \, \psi_k \;\;\sin 2 \, \psi_{\ell} + \;\cos 2 \, \psi_{\ell}\; \sin 2 \, \psi_k)\; \sin 2 \, \psi_j \, \frac{\partial}{\partial \varphi_j} \nonumber \\
&= &2 \sum \sin (2 \, \varphi_j)\; \sin (\varphi_k + \varphi_{\ell}
- \varphi_j) \, \frac{\partial}{\partial \varphi_j} \, . \nonumber
\end{eqnarray}

\end{proof}

\medskip

\medskip

\subsection{Invariance of ${\mathbb R}_{\varphi}^4$ under  the flow
$F$} \label{flow}

\noindent
\medskip

We let
\begin{equation}
C_+ =\{ (\varphi , \varphi , 0)\}\,,\quad C_- =\{ ( \frac{\pi}{2} ,
\varphi , \varphi + \frac{\pi}{2} \}\,.
\end{equation}

The critical set of the flow $X$ is described as follows in terms of
the action of the Weyl group $W$,

\medskip

 \begin{lem}\label{crit}
 The critical set $C$ of $X$ is given by
\begin{equation}
C = W(C_+) \cup W(C_-) \cup W(P) \, , \qquad P = \left(
\frac{\pi}{2} , \frac{\pi}{2} , \frac{\pi}{2} \right) \, . \nonumber
\end{equation}
 \end{lem}
\medskip

One checks this directly in the $\psi$ variables.

In order to perform a change of variables we use the function,
\begin{equation}
\label{eq1.70} \delta (\varphi) = \prod \sin \varphi_j \cos
(\varphi_k - \varphi_{\ell})
\end{equation}
and we let
\begin{equation}
\label{eq1.71} D = \{ \varphi , \delta (\varphi) = 0 \} \, .
\end{equation}
By construction $D$ is invariant under the group $S_3$ of
permutations of the $\varphi_j$'s.

\bigskip

 \begin{lem} \label{crit1}
Let $\varphi \notin C$ then there exists $g \in W$ such that $g
\varphi \notin D$.
\end{lem}
\bigskip

  \begin{proof} Let us first show that the
conclusion holds if one of the $\varphi_j$ vanishes (we always work
modulo $\pi$ and all equalities below have this meaning).

Thus assume that $\varphi_3 = 0$, {\it i.e.} that $\varphi =
(\varphi_1 , \varphi_2 , 0)$.
 Then if $\varphi_1 = \varphi_2$ one is in $C$ thus
we can assume that $\varphi_1 \ne \varphi_2$. By the transformation
(\ref{eq1.69}) we get $$\varphi' = (- \, \varphi_1 , - \, \varphi_1
, \varphi_2 - \varphi_1)$$ If $\varphi_1 = 0$ we treat
$(0,\varphi_2,0) \sim (\varphi_2 ,0,0)$ by applying (\ref{eq1.69})
which gives $(- \, \varphi_2 , - \, \varphi_2 , - \, \varphi_2)$ for
which $\delta (\varphi)$ is 0 only if $\varphi_2 = 0$ in which case
we are dealing with the point $O=(0,0,0)$ which is in $C$.

Thus we can assume $\varphi_1 \ne 0$, $\varphi_2 - \varphi_1 \ne 0$
and we get  $\prod \sin \varphi'_j \ne 0$.

 One has $\varphi'_k - \varphi'_{\ell} \in\{0,\pm \varphi_2\}$ thus
the product $\prod \cos (\varphi'_k - \varphi'_{\ell})$ vanishes
only if $\varphi_2 = \frac{\pi}{2}$. We are thus dealing with
$\left( \varphi_1 , \frac{\pi}{2} , 0 \right)$ (with $\varphi_1 \ne
0$). We apply (\ref{eq1.69}) to $\left( \frac{\pi}{2} , \varphi_1, 0
\right)$ which gives
$$\varphi' = \left( - \frac{\pi}{2} , - \frac{\pi}{2} , \varphi_1 - \frac{\pi}{2} \right)$$
 If $\varphi_1 = \frac{\pi}{2}$ then one is in $C$ otherwise $\prod \sin \varphi'_j \ne 0$. One has
$\varphi'_k - \varphi'_{\ell} \in\{0,\pm \varphi_1\}$ and since
$\varphi_1 \ne \frac{\pi}{2}$ we get $\prod \cos (\varphi'_k -
\varphi'_{\ell})\ne 0$

 We have thus shown that if $\varphi_j = 0$ for some
$j$ and  $\varphi \notin C$ we can find $g \in W$ with $g \varphi
\notin D$.

Thus we can now assume that all $\varphi_j \ne 0$. Thus $\prod \sin
\varphi_j \ne 0$. If $\prod \cos (\varphi_k - \varphi_{\ell}) = 0$
we can assume $\varphi_1 - \varphi_2 = \frac{\pi}{2}$. We then apply
(\ref{eq1.69}) and get
$$ \varphi'=\,\left( - \, \varphi_1 , \varphi_3 - \varphi_1 , - \frac{\pi}{2} \right) $$
If one of the components of $\varphi'$ is 0 we are back to the
previous case thus we can assume  $\prod \sin \varphi'_j \ne 0$. The
$\varphi'_k - \varphi'_{\ell}$ give (up to sign) $\varphi_3$,
$\varphi_3 - \varphi_1 + \frac{\pi}{2}$, $\varphi_1 - \frac{\pi}{2}$
and since $\varphi_1 \ne 0$ and $\varphi'_2=\varphi_3 - \varphi_1
\ne 0$, $\prod \cos (\varphi'_k - \varphi'_{\ell}) = 0$ can occur
only if $\cos \,\varphi_3 = 0$. In that case the original $\varphi$
is $\left( \varphi_1 , \varphi_1 - \frac{\pi}{2} , \frac{\pi}{2}
\right)$ which is in $C$. \hfill  \end{proof}

\bigskip

Note then that one can always choose the $g \in W$ in the Klein
subgroup $H \subset W$. Indeed $H$ is a normal subgroup and $S_3
\subset W$ acts as the group of permutation of the $\varphi$'s and
preserves $D$. Thus for $g = \sigma \, k_1$, $\sigma \in S_3$,
$\sigma \, k_1 \, \varphi \notin D \Rightarrow k_1 \, \varphi \notin
D$.

\medskip

  Let us now use this lemma to simplify the presentation of the algebra.

  \begin{lem}\label{rescalelem} If $\delta (\varphi) \ne 0$, there exists $4$ non zero scalars $\in i^N \, {\mathbb R}^*$, such that one has
 \[
(\sin \varphi_k) \, \lambda_{\ell} \, \lambda_m + \cos
(\varphi_{\ell} - \varphi_m) \, \lambda_0 \, \lambda_k = 0
\]
and $\prod \lambda_{\mu} = - \delta (\varphi)$.
\end{lem}

\bigskip

  \begin{proof} Let us choose the square roots $$\lambda_0 = (\prod \sin \varphi_j)^{1/2}\,,\quad\lambda_k = \left( \sin \varphi_k \underset{\ell \ne k}{\prod} \cos (\varphi_k - \varphi_{\ell}) \right)^{1/2}$$ in such a way that $\prod \, \lambda_j = - \delta (\varphi)$. Note indeed that the product of the squares gives $$\prod (\sin \varphi_j)^2 \underset{k < \ell}{\prod} (\cos (\varphi_k - \varphi_{\ell}))^2 = \delta (\varphi)^2$$ and thus one can fix the $\lambda_k$ and then choose the sign of $\lambda_0$ so that it fits. Then to prove the 3 equalities one can multiply by $\lambda_0 \, \lambda_k$ which gives $(\sin \varphi_k) (- \delta (\varphi)) + \cos (\varphi_{\ell} - \varphi_m) \, \lambda_0^2 \, \lambda_k^2$. Thus one needs to check $$\sin \varphi_k \prod \cos (\varphi_{\ell_1} - \varphi_{m_1}) = \cos (\varphi_{\ell} - \varphi_m) \underset{\ell_1 \ne k}{\prod} \cos (\varphi_k - \varphi_{\ell_1}) \sin \varphi_k$$ which is an identity. \hfill \end{proof}

\bigskip

One then lets
\begin{equation}
\label{rescale1} S_{\mu} = \lambda_{\mu} \, x^{\mu} \, , (\mbox{no
summation on}\>  \mu) .
\end{equation}
Multiplying (\ref{pres1}) by $\lambda_0 \, \lambda_1 \, \lambda_2 \,
\lambda_3$ one gets
\begin{equation}
\label{sklya1} [S_{\ell} , S_m] = i \, [S_0 , S_k]_+
\end{equation}
while (\ref{pres2}) gives:
\begin{equation}
\label{eq1.76} \cos \varphi_k [S_0 , S_k] = i \, \frac{\lambda_0
\lambda_k}{\lambda_{\ell} \lambda_m} \sin (\varphi_{\ell} -
\varphi_m) \, [S_{\ell} , S_m]_+ \, .
\end{equation}
Note that $$\frac{\lambda_0 \lambda_k}{\lambda_{\ell} \lambda_m} =
\frac{\prod \lambda_{\mu}}{\lambda_{\ell}^2 \lambda_m^2} =
\frac{-\delta (\varphi)}{\lambda_{\ell}^2 \lambda_m^2}$$

\medskip

\medskip

  \begin{prop}\label{sklyaprep} If $\delta (\varphi) \ne 0$
 the quadratic algebra of ${\mathbb R}_{\varphi}^4$
admits the presentation given by {\rm (\ref{sklya1})} and
\begin{equation}
\label{sklya2} [S_0 , S_k] = i \, J_{\ell m} [S_{\ell} , S_m]_+
\end{equation}
where the $J_{\ell m}$ are given by
\begin{equation}
\label{eq1.78} J_{\ell m} = - {\rm tan} \, \varphi_k \, {\rm tan}
(\varphi_{\ell} - \varphi_m) \, .
\end{equation}
If $\varphi_k = \frac{\pi}{2}$ and $\varphi_{\ell} = \varphi_m$ the
corresponding relation gives $0=0$. If $\varphi_k = \frac{\pi}{2}$
and $\varphi_{\ell} \ne \varphi_m$ it gives $[S_{\ell} , S_m]_+ =
0$.
\end{prop}
\medskip

 \begin{proof} One needs to show that
$$\frac{- \delta (\varphi)}{\lambda_{\ell}^2 \lambda_m^2} \, \frac{\sin (\varphi_{\ell} - \varphi_m)}{\cos \varphi_k} = -{\rm tan} \, \varphi_k \, {\rm tan} (\varphi_{\ell} - \varphi_m)$$
 This amounts to the equality
$$\sin \varphi_{\ell} \;\sin \varphi_m \,\cos^2 (\varphi_{\ell} - \varphi_m)\,\cos (\varphi_{\ell} - \varphi_k) \cos (\varphi_m - \varphi_k) = \lambda_{\ell}^2 \, \lambda_m^2$$
 which follows from the definition of the $\lambda_k$'s. \hfill  \end{proof}

\bigskip

We now find a better way to write the $J_{\ell m}$. We first pass in
a faithful manner from the $\varphi_j$ to
\begin{equation}
\label{eq1.79} t_j = {\rm tan} \, \varphi_j \in {\mathbb R} \qquad
(\hbox{recall} \ \varphi_j \ne \frac{\pi}{2}) \, .
\end{equation}
Note also that $t_j \ne 0$ since we are on $\delta (\varphi) \ne 0$.

\medskip

  We then introduce the following new parameters,
\begin{equation}
\label{eq1.80} s_k (\varphi) = 1 + t_{\ell} \, t_m \, .
\end{equation}

\bigskip

 \begin{lem} \label{sigmacover}
\begin{enumerate}
\item[a)]
The map $\varphi \overset{\sigma}{\longrightarrow} s =
(s_k(\varphi))$ on $ \Omega=\,\left\{ \varphi \,|\, \cos(\varphi_k)
\ne 0 \ \forall \, k \,, \delta (\varphi) \ne 0 \right\}$ is a
double cover of the open subset of ${\mathbb R}^3$
\begin{equation}
\label{domain-s} \prod s_k \ne 0 \quad \hbox{and} \quad \prod (s_k -
1) > 0 \, .
\end{equation}
\item[b)] The map $\sigma$ is a diffeomorphism of the interior $A^\circ$ of $A$
with  $\sigma(A^\circ)=\{s \,|\,1<s_1<s_2<s_3\}$.
\item[c)] The map $\sigma$ is a diffeomorphism of $B^\circ$
with  $\sigma(B^\circ)=\{s \,|\,s_3<s_2<0, 1<s_1\}$.
\end{enumerate}
\end{lem}
\bigskip

  \begin{proof} a) The condition
$\delta (\varphi) \ne 0$ shows that $s_k \ne 0 \,,\;\forall \, k $.
Indeed ${\rm tan} \, \varphi_{\ell} \, {\rm tan} \, \varphi_m = -1$
means $\cos(\varphi_{\ell} - \varphi_m) = 0$.
 By construction $s_k \ne 1$ and $\prod (s_k - 1) = \prod t_{\ell}^2 > 0$. Knowing the $s_k$'s one gets $(\prod t_k)^2 = \prod (s_k - 1)$ and choosing the sign of the square root  gives $p=\prod t_k$ and then $t_k = p(s_k - 1)^{-1}$. Thus
 one gets a double cover and the range
is characterized by the conditions \eqref{domain-s}. The deck
transformation is simply $\varphi\to\,-\varphi$.

b) On $A^\circ=\{\varphi
\,|\,\frac{\pi}{2}>\varphi_1>\varphi_2>\varphi_3>0\}$ one has
$t_k>0$ and thus $\prod t_k>0$ so that the above map is one to one.
One checks that the range $\sigma(A^\circ)$ is given by
$\sigma(A^\circ)=\{s \,|\,1<s_1<s_2<s_3\}$.

c) On $B^\circ=\{\varphi \,|\,\frac{\pi}{2}+\varphi_3
>\varphi_1>\frac{\pi}{2}>\varphi_2>\varphi_3\}$
one has $t_1<0,t_2>0,t_3>0$ and thus $\prod t_k<0$ so that the above
map is one to one. One checks that the range $\sigma(B^\circ)$ is
given by $\sigma(B^\circ)=\{s \,|\,s_3<s_2<0, 1<s_1\}$. \hfill
\end{proof}

\smallskip

We define the transformation $\rho$ by
$$\rho (s)_k = \frac{s_{\ell} - s_m}{s_k}$$
\smallskip

 \begin{lem} \label{propor}
\begin{enumerate}
\item[a)] Let $s_k \ne 0$, then one has $\prod (1+\rho (s)_k) = \prod (1-\rho (s)_k)$.
\item[b)] If $\rho (s) = \rho (s')$ there exists $\lambda \ne 0$ with $s' = \lambda \,s$.
\item[c)] The transformation $s \to \tilde s$, $$\tilde s_k = \frac{-s_k + s_{\ell} + s_m}{s_{\ell} s_m}$$ is involutive and $\rho (\tilde s) = -\rho (s)$.
\end{enumerate}
\end{lem}

\bigskip

   \begin{proof} a) Both products give, up to the denominator $\prod s_k$, the product of $(s_k + s_{\ell} - s_m)$.

\medskip

  b) Consider the 3 linear equations (for fixed $\rho_k$) given by
\begin{equation}
\label{eq1.83} s_{\ell} - s_m - \rho_k \, s_k = 0 \, .
\end{equation}
This corresponds to the $3 \times 3$ matrix given by
\begin{equation}
\label{eq1.84} M(\rho) = \begin{bmatrix} \rho_1 &-1 &1 \\ 1 &\rho_2
&-1 \\ -1 &1 &\rho_3 \end{bmatrix} \, .
\end{equation}
The condition a) means exactly that ${\rm Det} \, (M(\rho)) = 0$.
Moreover we claim that the rank of $M(\rho)$ is 2. Indeed $\rho_2$
appears in the two minors $\begin{bmatrix} 1 &\rho_2 \\ -1 &1
\end{bmatrix}$ which gives $1+\rho_2$ and $\begin{bmatrix} -1 &1 \\
\rho_2 &-1 \end{bmatrix}$ which gives $1-\rho_2$ and one of them is
$\ne 0$.

\medskip

  This shows that the kernel is 1-dimensional and hence $\rho (s) = \rho (s')$ implies $s' = \lambda \, s$ for some $\lambda$.

\medskip

  c) The relation between $s$ and $\tilde s$ can be written as
\begin{equation}
\label{eq1.85} s_{\ell} \, \tilde s_m + s_m \, \tilde s_{\ell} = 2
\, .
\end{equation}
The determinant of the system is $2 \, s_1 \, s_2 \, s_3$ so that
(\ref{eq1.85}) determines $\tilde s_k$ uniquely. The relation is
clearly symmetric. Finally
\begin{equation}
\rho (\tilde s)_1 = \frac{\tilde s_2 - \tilde s_3}{\tilde s_1} =
\frac{(- \, s_2 + s_1 + s_3) \, s_2 - (- \, s_3 + s_1 + s_2) \,
s_3}{s_1 \, (- \, s_1 + s_2 + s_3)} = \frac{s_3 - s_2}{s_1} = - \rho
(s)_1 \, . \nonumber
\end{equation}
\hfill
 \end{proof}
\bigskip

  We now get
\begin{equation}
\label{eq1.86} \rho \circ s \, (\varphi)_k = J_{\ell m} \qquad
(\varphi_k \ne \frac{\pi}{2} \  \ \hbox{and} \ \delta (\varphi) \ne
0) \, .
\end{equation}
Indeed
 $$\frac{s_{\ell} - s_m}{s_k} = \frac{t_k t_m - t_k t_{\ell}}{1 + t_{\ell} t_m} = {\rm tan} \, \varphi_k \, {\rm tan} (\varphi_m - \varphi_{\ell}) = J_{\ell m}\,.$$

\bigskip

\begin{lem} \label{convex}
\begin{enumerate}
\item[1)] The flow $X$ fulfills
\begin{equation}
\label{eq1.87} X \, s_k (\varphi) = 4 \prod \sin \,\varphi_j \,\;
s_k (\varphi) \, .
\end{equation}
\item[2)]
Let $\varphi, \varphi' \in  A$ with $\varphi_3 >0$, $\varphi'_3 >0$.
The following conditions are equivalent:
\begin{enumerate}
\item[a)] $J_{\ell m} (\varphi') = J_{\ell m} (\varphi) \, , \quad \forall k$
\item[b)] $\varphi'$ belongs to the orbit of $\varphi$  by the flow $X$.
\end{enumerate}
\item[3)]
The same statement holds for $\varphi, \varphi' \in  B$ with
$\varphi_3+\frac{\pi}{2} >\varphi_1$, $\varphi'_3+\frac{\pi}{2}
>\varphi'_1$.
\end{enumerate}
\end{lem}

\medskip

 \begin{proof} 1) One has
\begin{equation}
X(t_k) = \sin (2 \, \varphi_k) \sin (- \, \varphi_k + \varphi_{\ell}
+ \varphi_n) \, \frac{\partial}{\partial \varphi_k} \, {\rm tan} \,
\varphi_k = 2 \, {\rm tan} \, \varphi_k \sin (- \, \varphi_k +
\varphi_{\ell} + \varphi_m) \, , \nonumber
\end{equation}
\begin{equation}
X(s_k) = X (t_{\ell}) \, t_m + t_{\ell} \, X (t_m) = 2 \, t_{\ell}
\, t_m (\sin (- \, \varphi_{\ell} + \varphi_m + \varphi_k) + \sin (-
\, \varphi_m + \varphi_{\ell} + \varphi_k)) \nonumber
\end{equation}
and using $\sin (a+b) + \sin (a-b) = 2 \sin a\, \cos b$, one gets
\begin{equation}
X(s_k) = 4 \, t_{\ell} \, t_m (\sin \varphi_k \cos (\varphi_m -
\varphi_{\ell})) = 4 \,\left( \prod \sin \varphi_j \right)\,
\frac{\cos (\varphi_m - \varphi_{\ell})}{\cos \varphi_m \cos
\varphi_{\ell}} = 4 \left( \prod \sin \varphi_j \right) (1+t_m \,
t_{\ell}) \nonumber
\end{equation}

\smallskip

2)  In fact 1) shows that the flow $X$ is up to a non-zero change of
speed the scaling flow in $\sigma(A)$. By construction $J_{\ell m}$
has homogeneity degree zero in $s_k$ thus it is preserved by $X$ and
b) $\Rightarrow$ a).

Let us show that a) $\Rightarrow$ b). We assume first that
$\frac{\pi}{2}>\varphi_1$. The same then holds for $\varphi'$ uisng
a). By Lemma \ref{propor} the equality a) implies that
$\sigma(\varphi')=\lambda \,\sigma(\varphi)$ for some non-zero
scalar $\lambda$. By Lemma \ref{sigmacover} the image
$\sigma(A^\circ)$ is convex (as well as its closure) and thus the
segment $[\sigma(\varphi),\,\sigma(\varphi')]$ is contained in
$\sigma(A)$ and its preimage under $\sigma$ is a segment in a flow
line.

On the face $Y$ determined by
 $\varphi_1 = \frac{\pi}{2}$ one has
assuming $\varphi_2 < \frac{\pi}{2} $ the equalities
 $$J_{12} = - \, {\rm tan} \, \varphi_3 \, {\rm tan} \left( \frac{\pi}{2} - \varphi_2 \right) = - \, t_3 / t_2\,,\quad
J_{31} = - \, {\rm tan} \, \varphi_2 \, {\rm tan} \left( \varphi_3 -
\frac{\pi}{2} \right) = t_2 / t_3$$ Moreover one has $X(t_2) = \cos
(\varphi_2 - \varphi_3) \, 2 \, t_2$, $X(t_3) = \cos (\varphi_2 -
\varphi_3) \, 2 \, t_3$ so that the flow $X$ restricts as the
scaling flow (up to a non-zero change of speed) in the parameters
$t_j$. Since $J_{23}=\infty$ holds iff $\varphi_1 = \frac{\pi}{2}$
for $\varphi \in  A$ we see that a) then implies
 $\varphi'_1 = \frac{\pi}{2}$ and
the proportionality $t'_j=\lambda t_j$. Since the allowed $t_j$ are
simply constrained by the inequalities $t_2\geq t_1>0$ the same
convexity argument applies. Finally let $\varphi$ with $\varphi_1 =
\varphi_2 = \frac{\pi}{2}$, $\varphi_3 \ne \frac{\pi}{2}$. Then the
only remaining parameter is $t_3$ and one has
 $X(t_3) = 2 \sin \varphi_3 \, t_3$.
Thus one is dealing with a single flow line.

The proof of 3) is similar. \hfill  \end{proof}

\newpage

\section{The geometric data of ${\mathbb R}_{\varphi}^4$}\label{geometricdata}

In this section we compute the geometric data of the quadratic
algebras of functions on ${\mathbb R}_{\varphi}^4$ for all values of
the parameter $\varphi$. These fall in eleven different classes in
each of which one gets further invariants.

\medskip

\subsection{The definition and explicit matrices}

\noindent
\medskip

We  give a list of the characteristic varieties and correspondences.
There are 11 different cases. They are described in terms of the
$\varphi$-coordinates but the result will then be translated in
invariant terms using the roots. Let us recall the definition of the
geometric data $\{E\,,\,  \sigma\,,\,{\mathcal L}\}$ for  quadratic
algebras. Let ${\mathcal A}=A(V,R)=T(V)/(R)$ be a quadratic algebra
where $V$ is a finite-dimensional complex vector space and where
$(R)$ is the two-sided ideal of the tensor algebra $T(V)$ of $V$
generated by the subspace $R$ of $V\otimes V$. Consider the subset
of $V^\ast \times V^\ast$ of  pairs $(\alpha,\beta)$ such that
\begin{equation}
\langle \omega,\alpha\otimes\beta\rangle=0,\, \, \, \alpha\not= 0,
\beta\not=0 \label{eq5.1}
\end{equation}
for any $\omega\in R$. Since $R$ is homogeneous, (\ref{eq5.1})
defines a subset
\[
\Gamma \subset P(V^\ast)\times P(V^\ast)
\]
 where $P(V^\ast)$ is the complex projective space of one-dimensional complex subspaces of $V^\ast$. Let $E_1$ and $E_2$ be the first and the second projection of $\Gamma$ in $P(V^\ast)$. It is usually assumed that they coincide i.e. that one has
\begin{equation}
E_1=E_2=E\subset P (V^\ast) \label{eq5.2}
\end{equation}
and that the correspondence $\sigma$ with graph $\Gamma$ is an
automorphism of $E$, ${\mathcal L}$ being the pull-back on $E$ of
the dual of the tautological line bundle of $P(V^\ast)$. The
algebraic variety $E$ is refered to as the characteristic variety.

The algebras we consider have 4 generators and six relations, thus
their characteristic variety is obtained as the locus of points
where a $4 \times 6$ matrix has rank less than 4. The various
matrices depend upon the choice of the parameters and are listed
below. We first give the matrix corresponding to the original
quadratic algebra of ${\mathbb R}_{\varphi}^4$,

\begin{equation}\label{matrixphi}
\begin{bmatrix}
    -\cos  ({{\varphi }_1})\, {x_1}&\cos  ({{\varphi }_1})\, {x_0}&-i \, \sin  ({{\varphi }_2}-{{\varphi }_3})\, {x_3}&-i
\, \sin  ({{\varphi }_2}-{{\varphi }_3})\, {x_2} \\
    -\cos  ({{\varphi }_2})\, {x_2}&i \, \sin  ({{\varphi }_1}-{{\varphi }_3})\, {x_3}&\cos  ({{\varphi }_2})\, {x_0}&i
\, \sin  ({{\varphi }_1}-{{\varphi }_3})\, {x_1} \\
    -\cos  ({{\varphi }_3})\, {x_3}&-i \, \sin  ({{\varphi }_1}-{{\varphi }_2})\, {x_2}&-i \, \sin  ({{\varphi
}_1}-{{\varphi }_2})\, {x_1}&\cos  ({{\varphi }_3})\, {x_0} \\
    i \, \sin  ({{\varphi }_3})\, {x_3}&-\cos  ({{\varphi }_1}-{{\varphi }_2})\, {x_2}&\cos  ({{\varphi }_1}-{{\varphi }_2})\,
{x_1}&i \, \sin  ({{\varphi }_3})\, {x_0} \\
    i \, \sin  ({{\varphi }_1})\, {x_1}&i \, \sin  ({{\varphi }_1})\, {x_0}&-\cos  ({{\varphi }_2}-{{\varphi }_3})\,
{x_3}&\cos  ({{\varphi }_2}-{{\varphi }_3})\, {x_2} \\
    i \, \sin  ({{\varphi }_2})\, {x_2}&\cos  ({{\varphi }_1}-{{\varphi }_3})\, {x_3}&i \, \sin  ({{\varphi }_2})\,
{x_0}&-\cos  ({{\varphi }_1}-{{\varphi }_3})\, {x_1}
  \end{bmatrix}
\end{equation}

When we pass to the Sklyanin algebra and eliminate the factors $i$
in replacing $ Z_0=i\,S_0 \,,\quad Z_k=\,S_k\,, $
 we get the following matrix, with $\alpha=- J_{23}$, $\beta=- J_{31}$, $\gamma=- J_{12}$,

\begin{equation}\label{matrixsklya}
\begin{bmatrix}
{z_1}&-{z_0}&\alpha \, {z_3}&\alpha \, {z_2} \\
    {z_2}&\beta \, {z_3}&-{z_0}&\beta \, {z_1} \\
    {z_3}&\gamma \, {z_2}&\gamma \, {z_1}&-{z_0} \\
    {z_3}&{z_2}&-{z_1}&{z_0} \\
    {z_1}&{z_0}&{z_3}&-{z_2} \\
    {z_2}&-{z_3}&{z_0}&{z_1}
 \end{bmatrix}
\end{equation}
The fifteen minors of the matrix \eqref{matrixphi} are listed in
factorized form in Appendix 1 and those of the matrix
\eqref{matrixsklya} in subsection \ref{minorslist} below (see
\cite{smi-sta:1992}).

\subsection{The Table}

\noindent
\medskip

We give below the list of the characteristic varieties and
correspondences in all cases. Given $z\in {\mathbb C}$
 we let $\sigma(z)$ be the conjugacy class of semi-simple
automorphisms of a curve of genus $0$ with eigenvalues
$\{z,z^{-1}\}$. By construction $\sigma(z)=\,\sigma(z^{-1})$.

The identification of the corresponding algebras $\calg(\mathbb
R^4_\varphi)$ in the nongeneric cases (cases 2 to 11 in the table
below) are described in section \ref{degalg} and may be summarized
as follows. In case 2, $\calg(\mathbb R^4_\varphi)$ is isomorphic to
a homogenized version (quadratic) $U_q(\fracsu(2))^{\hom}$  of
$U_q(\fracsu(2))$   with either $q\in \mathbb C$ with $\vert
q\vert=1$ or $q\in \mathbb R$ with $0<q<1$ or $-1<q<0$; it is
important to notice that this is a $\ast$-isomorphism, (that is the
$\fracsu(2)$ does really matter in this notation). Case 3 is
obtained by duality from case 2 as explained in Section
\ref{dualities}. In case 4, there is a missing relation so
$\calg(\mathbb R^4_\varphi)$, which corresponds formally to a
version $q=0$ of $U_q(\fracsu(2))^{\hom}$, is of exponential growth.
Case 5 is obtained by $\alpha_3$-duality (section \ref{dualities})
from case 6 which is isomorphic to a homogenized version
$U(\fracsu(2))^{\hom}$ of the universal enveloping algebra of
$\fracsu(2)$ (i.e. $q=1$ in $U_q(\fracsu(2))^{\hom}$). Case 7 is the
$\theta$-deformation studied in Part I \cite{ac-mdv:2002a} and in
\cite{ac-lan:2001} while case 8 (anti $\theta$-deformation) is
obtained by $\alpha_1$-duality (cf. section \ref{dualities}) from
case 7. Case 9 is very singular : 3 relations are missing. Case 10
is obtained by $\alpha_3$-duality (cf. section \ref{dualities}) from
case 11 which is the ordinary algebra of polynomials with 4
indeterminates $\mathbb C[x^0, x^1, x^2, x^3] $(classical case).
Note that the reality conditions above correspond to the hermiticity
of the $x^\mu$ and \underbar{not} of the Sklyanin generators
$S_\mu$.

It is important to describe the stratification in terms of the
roots. The   stratas of codimension $k$ correspond to intersections
of $k$ singular hyperplanes (up to a factor $\frac{1}{2}$) of the
form
$$
F((\alpha_1,\cdots,\alpha_k),(n_1,\cdots,n_k))=\,\cap\,\frac{1}{2}\,H_{(\alpha_j,n_j)}\,,\qquad
\alpha_j \in \Delta\,,n_j\in {\mathbb Z}\,.
$$
The two dimensional stratas are defined using a single root $\alpha
\in \Delta$ (\ie $k=1$) and since the Weyl group $W$ acts
transitively on $\Delta$ only the parity of $n$ matters which gives
the two kinds of faces $F_1$ corresponding to $n$ even and $F_2$ to
$n$ odd.

The one dimensional stratas are defined using two roots
$\alpha,\,\beta\in \Delta$ (\ie $k=2$). The roots only matter up to
sign and their relative positions is described by their angle which
(up to the   sign) can be $\frac{\pi}{2}$ in which case we write
$\alpha\perp\beta$ or  $\frac{2\pi}{3}$ in which case we write
$\alpha-\beta$. We thus have the following one dimensional stratas
 according to the parity of the $n_j$.

\begin{itemize}

\item $\alpha\perp\beta$ and $(n_1,n_2)=$(even, odd) gives the line $L$
of case 4 below.

\item $\alpha-\beta$ and $(n_1,n_2)=$(even, odd)
or $(n_1,n_2)=$(odd, odd) gives the line $L'$ of case 5  below.

\item $\alpha-\beta$ and $(n_1,n_2)=$(even, even) gives the line $L''$
of case 6  below.

\item $\alpha\perp\beta$ and $(n_1,n_2)=$(even, even) gives the line $C_+$
of case 7  below.

\item $\alpha\perp\beta$ and $(n_1,n_2)=$(odd, odd) gives the line $C_-$
of case 8  below.

\end{itemize}

To double check that the  list is complete\footnote{ we are grateful
to Marc Bellon for pointing out the subtelty of case 9) which was
incomplete in an earlier version} one can use Lemma \ref{crit} to
control the critical set $C$ and then assume by Lemma \ref{crit1}
that $\delta(\varphi)\ne 0$. Then if $\varphi \in H_{(\alpha,n)}$
and $n$ is even (resp. odd) the root $\alpha$ is one of the
differences $\varphi_k-\varphi_l$ (resp. $\varphi_k$). Thus up to
permutations of the $\varphi_k$ one obtains one of the cases 1)-6).
The complete table giving the geometric datas is the following :

\begin{center}
\begin{tabular}{|c|c|c|c|}
\hline & & & \\
\bf{Case} &\bf{Point in} \ ${\mathcal M}$ &
\bf{Characteristic variety} &\bf{Correspondence} \\
&&&\\
\hline & & & \\
\bf{1} Generic  &\ $\, \alpha(\varphi) \notin \frac{\pi}{2} \, {\mathbb Z}$ &4 points $\cup$ Elliptic curve &$({\rm id}, {\rm id}, {\rm id}, {\rm id}$, translation) \\
&&&\\ \hline & & & \\
\bf{2}$\;\;$ \bf{Even } & $\varphi_1 = \varphi_2 \,,\varphi_2
-\varphi_3
\notin \frac{\pi}{2} \, {\mathbb Z}$ &2 points, 1 line, 2 conics &$\left({\rm id},{\rm id}, {\rm id}, \sigma(\frac{i+\alpha^{1/2}}{i-\alpha^{1/2}}) , \sigma(\frac{i+\alpha^{1/2}}{i-\alpha^{1/2}}) \right)$ \\
\bf{Face} &$\varphi_j\notin \frac{\pi}{2} \, {\mathbb Z}  ,\,  \,$&&$\alpha = -J_{23}$ \\
&&&\\ \hline & & & \\
\bf{3} $\,\,$ \bf{ Odd } & $\varphi_1 = \frac{\pi}{2} , \, \varphi_2
- \varphi_3
\notin \frac{\pi}{2} \, {\mathbb Z}$ & 2 points, 1 line, 2 conics &$\Bigl({\rm id},{\rm id},  - {\rm id}, \hbox{exchange with square}$ \\
 \bf{Face} &$\varphi_k \notin \frac{\pi}{2} \, {\mathbb Z}  , \, k=2,3$ &&
$\sigma(\frac{i+\beta^{1/2}}{i-\beta^{1/2}})^2 , \, \beta = -J_{31} \Bigl)$ \\
&&&\\ \hline & & & \\
\bf{4 $\;\;\alpha\perp\beta$ } & $L=\{\left( \frac{\pi}{2} , \varphi , \varphi \right)\}$ &six lines &$({\rm id}, -{\rm id}, \hbox{cyclic permutation of 4 lines} $\\
(0,1) &&& (iso,  coarse, iso, coarse)) \\
&&&\\ \hline & & & \\
\bf{5 $\;\;\alpha-\beta$ }& $L'=\{\left( \frac{\pi}{2} , \frac{\pi}{2} , \varphi \right)\}$ &$\hbox{point} \cup P_2 ({\mathbb C})$ &(id, Symmetry of determinant $-1$) \\
 (0,1) &&&\\ &&&\\ \hline & & & \\
\bf{6 $\;\;\alpha-\beta$} & $L''=\{\left( \varphi , \varphi , \varphi \right)\}$ &$\hbox{point} \cup P_2 ({\mathbb C})$ &(id, id)\\
  (0,0) &&&\\ &&&\\ \hline & & & \\
\bf{7 $\;\;\alpha\perp\beta$ } & $C_+=\{(\varphi , \varphi , 0)\}$ & six lines & $({\rm id}, {\rm id} , \sigma(e^{\pm 2i\varphi}) ) $\\
 (0,0) &&&\\  &&&\\ \hline & & & \\
 \bf{8 $\;\;\alpha\perp\beta$ }&$C_- =\{ \left( \varphi + \frac{\pi}{2}, \frac{\pi}{2} , \varphi \right)\}$ &$\hbox{six lines}$ &$(-{\rm id} , - {\rm id}, \hbox{ exchanges with square}$ \ \\
 (1,1) &&&$\sigma(e^{\pm 4i\varphi})) $\\  &&&\\ \hline & & & \\
\bf{9} &$P =\,\left( \frac{\pi}{2} , \frac{\pi}{2} , \frac{\pi}{2} \right)$ &$P_3 ({\mathbb C})$ &$\hbox{Symmetry of determinant} \ -1$ \\
& &    &and point $\to \hbox{ line on a quadric}$ \\
 &&&\\ \hline & & & \\
\bf{10}  &$P' =\, \left( \frac{\pi}{2} , \frac{\pi}{2} , 0 \right)$ &$P_3 ({\mathbb C})$ &$\hbox{Symmetry of determinant} \ 1 $\\
&$(\hbox{in} \ C_+ \cap C_-) $& &\\
&&&\\ \hline & & & \\
\bf{11} &$0 = (0,0,0)$ &$P_3 ({\mathbb C})$ &$\hbox{id}$ \\
&$(\hbox{in} \ C_+)$& &\\
&&&\\ \hline
\end{tabular}
\end{center}

\begin{equation} \label{tablebig}
{\rm {\bf The}}\quad {\rm {\bf Geometric}}\quad {\rm {\bf Data }}
\end{equation}

\bigskip

\medskip
\subsection{Generic case}\label{minorslist}

\noindent
\medskip

 We now give the detailed description of the geometric data starting with the generic case. This case is defined by
\begin{equation}
\label{gencase} G = \left\{ \varphi \, ; \  \ \varphi_j \notin
\frac{\pi}{2}\; {\mathbb Z}\, , \ \varphi_k - \varphi_{\ell}\notin
\frac{\pi}{2}\; {\mathbb Z} \right\} \, .
\end{equation}
Then the $J_{k\ell}$ are well defined and are $\ne 0$. Let us show
that we cannot have $J_{12} = 1$, $J_{23} = -1$. Indeed ${\rm tan}
\, \varphi_3 \, {\rm tan} (\varphi_1 - \varphi_2) = -1$ means $\pi/2
- \varphi_3 = \varphi_2 - \varphi_1$ while ${\rm tan} \, \varphi_1
\, {\rm tan} (\varphi_2 - \varphi_3) = 1$ means $\frac{\pi}{2} -
\varphi_1 = \varphi_2 - \varphi_3$. This gives $\pi - \varphi_1 -
\varphi_3 = 2 \, \varphi_2 - \varphi_1 - \varphi_3$ and $2 \,
\varphi_2 = \pi$ which is not allowed by (\ref{gencase}). We can
thus apply the result of Smith-Stafford \cite{smi-sta:1992} and get:

\bigskip

  \begin{prop}\label{chargen}  For $\varphi \in G$ the geometric data of ${\mathbb R}_\varphi^4 $ is given by $4$ points and  a non degenerate elliptic curve $E$, and $\sigma$ is identity on the $4$ points and a translation of $E$ given explicitely in terms of the parameters
$\alpha_1 = \alpha$, $\alpha_2 = \beta$, $\alpha_3 = \gamma$ with
\begin{equation}
\label{minus} \alpha_k = - \, J_{\ell m} \, .
\end{equation}
by
\begin{equation}
\label{curve} E=\{z\,;\,\sum_0^3 z_j^2 = 0 \, , \ \frac{1 -
\gamma}{1+\alpha} \, z_1^2 + \frac{1+\gamma}{1-\beta} \, z_2^2 +
z_3^2 = 0 \}\, .
\end{equation}
and
\begin{equation}
\label{sigma}
\begin{bmatrix} z_0 \\ z_1 \\ z_2 \\ z_3 \end{bmatrix} \overset{\sigma}{\longrightarrow}
\begin{bmatrix}
- \, 2 \, \alpha \, \beta \, \gamma \, z_1 \, z_2 \, z_3 - z_0 (- \, z_0^2 + \beta \, \gamma \, z_1^2 + \alpha \, \gamma \, z_2^2 + \alpha \, \beta \, z_3^2) \\
2 \, \alpha \, z_0 \, z_2 \, z_3 + z_1 (z_0^2 - \beta \, \gamma \, z_1^2 + \alpha \, \gamma \, z_2^2 + \alpha \, \beta \, z_3^2) \\
2 \, \beta \, z_0 \, z_1 \, z_3 + z_2 (z_0^2 + \beta \, \gamma \, z_1^2 - \alpha \, \gamma \, z_2^2 + \alpha \, \beta \, z_3^2) \\
2 \, \gamma \, z_0 \, z_1 \, z_2 + z_3 (z_0^2 + \beta \, \gamma \,
z_1^2 + \alpha \, \gamma \, z_2^2 - \alpha \, \beta \, z_3^2)
\end{bmatrix} \, .
\end{equation}
\end{prop}
\bigskip

  \begin{proof} We rely on Smith-Stafford \cite{smi-sta:1992}. The first point is to rewrite the algebra of proposition (\ref{sklyaprep}) in the form,
\begin{equation}
\label{eq1.90} [T_0 , T_k]_+ = [T_{\ell} , T_m] \, , \ [T_0 , T_k] =
\alpha_k [T_{\ell} , T_m]_+
\end{equation}
One simply lets $T_0 = i \, S_0$, $T_k = S_k$. Then (\ref{sklya1})
means $[T_{\ell} , T_m] = [T_0 , T_k]_+$ and (\ref{sklya2}) means
$[T_0 , T_k] = - \, J_{\ell m} [T_1 , T_m]_+$.

\medskip

  Note the crucial $-$ sign in (\ref{minus}).
  By hypothesis on $\varphi \in G$ one has $\alpha \, \beta \, \gamma \ne 0$. Moreover we have seen above that for $\varphi \in G$ one cannot have $J_{12} = 1$, $J_{23} = -1$, {\it i.e.} $\alpha = -1$, $\beta = 1$ (or any cyclic transformed of that). Now in case one of the $\alpha$'s, say $\alpha$ belongs to $\pm 1$ the equality $\alpha + \beta + \gamma + \alpha \, \beta \, \gamma = 0$, {\it i.e.} $\prod (1+\alpha) = \prod (1-\alpha)$ shows that both $+1$ and $-1$ occur and since $(-1,1,x)$ is excluded the only remaining case is $(1,-1,\gamma)$ with $\gamma \notin \{ -1,0,1 \}$.

\medskip

  Thus the hypothesis of Smith-Stafford are fulfilled
and one gets from \cite{smi-sta:1992} that besides the 4 points
$(1,0,0,0) \ldots$ the characteristic variety is the curve in $P_3
({\mathbb C})$ with equations \eqref{curve}. The translation
$\sigma$ being given explicity by \eqref{sigma}.
\end{proof}

It will be useful for the computations in the degenerate case to
display the list of the 15 minors in the case of the Sklyanin
algebra, in the above parameters $(\alpha,\beta,\gamma)$. Their list
is given below (\cite{smi-sta:1992}).

\begin{equation}
\begin{matrix}
  \{-2   \, (\gamma    \, {z_1}   \, {z_2}-{z_0}   \, {z_3})   \,  (z_{0}^{2}+z_{1}^{2}+z_{2}^{2}+z_{3}^{2} ),  \\
\noalign{\vspace{0.666667ex}} \hspace{1.em} 2   \, ({z_0}   \,
{z_2}-\beta    \, {z_1}   \, {z_3})   \,
(z_{0}^{2}+z_{1}^{2}+z_{2}^{2}+z_{3}^{2} ),
  \\
\noalign{\vspace{0.666667ex}} \hspace{1.em} -2   \, ({z_0}   \,
{z_2}+{z_1}   \, {z_3})   \,  (z_{0}^{2}-\gamma    \,  (\beta    \,
z_{1}^{2}+z_{2}^{2} )+\beta
   \, z_{3}^{2} ),  \\
\noalign{\vspace{0.666667ex}} \hspace{1.em} 2   \,  (z_{0}^{2}   \,
(-(1+\gamma )   \, z_{2}^{2}+(-1+\beta )   \, z_{3}^{2} )+z_{1}^{2}
\,  ((-1+\beta
)   \, \gamma    \, z_{2}^{2}+\beta    \, (1+\gamma )   \, z_{3}^{2} ) ),  \\
\noalign{\vspace{0.666667ex}} \hspace{1.em} 2   \, (-{z_1}   \,
{z_2}+{z_0}   \, {z_3})   \,  (z_{0}^{2}-\gamma    \,  (\beta    \,
z_{1}^{2}+z_{2}^{2} )+\beta
   \, z_{3}^{2} ),     \\
\noalign{\vspace{0.666667ex}} \hspace{1.em} 2   \, ({z_0}   \,
{z_1}-\alpha    \, {z_2}   \, {z_3})   \,
(z_{0}^{2}+z_{1}^{2}+z_{2}^{2}+z_{3}^{2} ),
  \\
\noalign{\vspace{0.666667ex}} \hspace{1.em} -2   \, ({z_0}   \,
{z_1}-{z_2}   \, {z_3})   \,  (z_{0}^{2}+\gamma    \,
z_{1}^{2}-\alpha    \,  (\gamma    \,
z_{2}^{2}+z_{3}^{2} ) ),     \\
\noalign{\vspace{0.666667ex}} \hspace{1.em} -2   \, ({z_1}   \,
{z_2}+{z_0}   \, {z_3})   \,  (z_{0}^{2}+\gamma    \,
z_{1}^{2}-\alpha    \,  (\gamma    \,
z_{2}^{2}+z_{3}^{2} ) ),  \\
\noalign{\vspace{0.666667ex}} \hspace{1.em} 2   \,  (z_{0}^{2}   \,
((-1+\gamma )   \, z_{1}^{2}-(1+\alpha )   \, z_{3}^{2} )+z_{2}^{2}
\,  ((1+\alpha
)   \, \gamma    \, z_{1}^{2}+\alpha    \, (-1+\gamma )   \, z_{3}^{2} ) ),  \\
\noalign{\vspace{0.666667ex}} \hspace{1.em} 2   \,  (z_{0}^{2}   \,
((1+\beta )   \, z_{1}^{2}-(-1+\alpha )   \, z_{2}^{2} )-
((-1+\alpha )   \, \beta    \,
z_{1}^{2}+\alpha    \, (1+\beta )   \, z_{2}^{2} )   \, z_{3}^{2} ),  \\
\noalign{\vspace{0.666667ex}} \hspace{1.em} -2   \, ({z_0}   \,
{z_2}-{z_1}   \, {z_3})   \,  (z_{0}^{2}+\alpha    \,
z_{2}^{2}-\beta    \,  (z_{1}^{2}+\alpha
   \, z_{3}^{2} ) ),     \\
\noalign{\vspace{0.666667ex}} \hspace{1.em} -2   \, ({z_0}   \,
{z_1}+{z_2}   \, {z_3})   \,  (z_{0}^{2}+\alpha    \,
z_{2}^{2}-\beta    \,  (z_{1}^{2}+\alpha
   \, z_{3}^{2} ) ),     \\
\noalign{\vspace{0.666667ex}} \hspace{1.em} 2   \, ({z_0}   \,
{z_2}+\beta    \, {z_1}   \, {z_3})   \,  (z_{0}^{2}+\gamma    \,
z_{1}^{2}-\alpha    \,  (\gamma
   \, z_{2}^{2}+z_{3}^{2} ) ),     \\
\noalign{\vspace{0.666667ex}} \hspace{1.em} 2   \, ({z_0}   \,
{z_1}+\alpha    \, {z_2}   \, {z_3})   \,  (z_{0}^{2}-\gamma    \,
(\beta    \, z_{1}^{2}+z_{2}^{2} )+\beta
   \, z_{3}^{2} ),     \\
\noalign{\vspace{0.666667ex}} \hspace{1.em} 2   \, (\gamma    \,
{z_1}   \, {z_2}+{z_0}   \, {z_3})   \,  (z_{0}^{2}+\alpha    \,
z_{2}^{2}-\beta    \,  (z_{1}^{2}+\alpha
   \, z_{3}^{2} ) ) \}\\
\end{matrix}
\end{equation}

\subsection{Face $\alpha=n$ and $n$ even.
$\,F_1 =\{ \left( \varphi_1   , \varphi_1  , \varphi_3 \right)\}$}

\noindent
\medskip

In that case we have a Sklyanin algebra and with the above notations
the parameters are $\gamma=0$ while $\beta=-\alpha= \tan \varphi_1\,
\tan(\varphi_3-\varphi_1)$. The list of minors then simplifies as
follows,

\begin{equation}
\begin{matrix}
  \{{z_0} \, {z_3} \,  (z_{0}^{2}+z_{1}^{2}+z_{2}^{2}+z_{3}^{2} ),({z_0} \, {z_2}+\alpha  \, {z_1} \, {z_3}) \,  (z_{0}^{2}+z_{1}^{2}+z_{2}^{2}+z_{3}^{2} ),
 \\
\noalign{\vspace{0.666667ex}} \hspace{1.em} -({z_0} \, {z_2}+{z_1}
\, {z_3}) \,  (z_{0}^{2}-\alpha  \, z_{3}^{2} ),-\alpha  \,
z_{1}^{2} \, z_{3}^{2}-z_{0}^{2} \,
 (z_{2}^{2}+(1+\alpha ) \, z_{3}^{2} ),      \\
\noalign{\vspace{0.666667ex}} \hspace{1.em} (-{z_1} \, {z_2}+{z_0}
\, {z_3}) \,  (z_{0}^{2}-\alpha  \, z_{3}^{2} ),({z_0} \,
{z_1}-\alpha  \, {z_2} \,
{z_3}) \,  (z_{0}^{2}+z_{1}^{2}+z_{2}^{2}+z_{3}^{2} ),      \\
\noalign{\vspace{0.666667ex}} \hspace{1.em} -({z_0} \, {z_1}-{z_2}
\, {z_3}) \,  (z_{0}^{2}-\alpha  \, z_{3}^{2} ),-({z_1} \,
{z_2}+{z_0} \, {z_3}) \,
 (z_{0}^{2}-\alpha  \, z_{3}^{2} ),      \\
\noalign{\vspace{0.666667ex}} \hspace{1.em} -\alpha  \, z_{2}^{2} \,
z_{3}^{2}-z_{0}^{2} \,  (z_{1}^{2}+(1+\alpha ) \, z_{3}^{2}
),(-1+\alpha ) \,  (z_{1}^{2}+z_{2}^{2} ) \,
 (-z_{0}^{2}+\alpha  \, z_{3}^{2} ),      \\
\noalign{\vspace{0.666667ex}} \hspace{1.em} -({z_0} \, {z_2}-{z_1}
\, {z_3}) \,  (z_{0}^{2}+\alpha  \,  (z_{1}^{2}+z_{2}^{2}+\alpha  \,
z_{3}^{2} ) ),-({z_0} \,
{z_1}+{z_2} \, {z_3}) \,  (z_{0}^{2}+\alpha  \,  (z_{1}^{2}+z_{2}^{2}+\alpha  \, z_{3}^{2} ) ),  \\
\noalign{\vspace{0.666667ex}} \hspace{1.em} ({z_0} \, {z_2}-\alpha
\, {z_1} \, {z_3}) \,  (z_{0}^{2}-\alpha  \, z_{3}^{2} ),({z_0} \,
{z_1}+\alpha
 \, {z_2} \, {z_3}) \,  (z_{0}^{2}-\alpha  \, z_{3}^{2} ),      \\
\noalign{\vspace{0.666667ex}}
\hspace{1.em} {z_0} \, {z_3} \,  (z_{0}^{2}+\alpha  \,  (z_{1}^{2}+z_{2}^{2}+\alpha  \, z_{3}^{2} ) ) \}\\
\end{matrix}
\end{equation}

The detailed analysis shows that the characteristic variety is the
union of the two points with coordinates
$(z_0,z_1,z_2,z_3)=(1,0,0,0)$, $(z_0,z_1,z_2,z_3)=(0,0,0,1)$, of the
line $\{(0,z_1,z_2,0)\}$ and of the two conics obtained by
intersecting the quadric $(z_{0}^{2}+z_{1}^{2}+z_{2}^{2}+z_{3}^{2}
)=0$ with the two hyperplanes $(z_{0}^{2}-\alpha    \, z_{3}^{2}
)=0$.

The correspondence $\sigma$ is the identity on the two points, and
on the line. It restricts to the two conics and is a rational
automorphism of each. It admits two fixed points on each of them and
its derivative at the fixed points is given by the  following
complex numbers
$$
\frac{i-\sqrt \alpha}{i+\sqrt \alpha}\,,\quad \frac{i+\sqrt
\alpha}{i-\sqrt \alpha}\,.
$$

\subsection{Face $\alpha=n$ and $n$ odd.
$\,F_2 =\{ \left( \frac{\pi}{2}   , \varphi_2  , \varphi_3
\right)\}$}

\noindent
\medskip

In that case we  have a limiting case of the Sklyanin algebra where
the parameter $\alpha =\infty$ while $\gamma= -\frac{1}{\beta}$. The
list of minors can be computed directly and gives the following, up
to non-zero scalar factors,

\begin{equation}
\begin{matrix}
 \{({z_1}   \, {z_2}+\beta    \, {z_0}   \, {z_3})   \,  (z_{0}^{2}+z_{1}^{2}+z_{2}^{2}+z_{3}^{2} ),({z_0}   \, {z_2}-\beta    \,
{z_1}   \, {z_3})   \,  (z_{0}^{2}+z_{1}^{2}+z_{2}^{2}+z_{3}^{2} ),  \\
\noalign{\vspace{0.666667ex}} \hspace{1.em} (-{z_0}   \, {z_2}-{z_1}
\, {z_3})   \,  (\beta    \,  (z_{0}^{2}+z_{1}^{2}
)+z_{2}^{2}+{{\beta }^2}   \, z_{3}^{2} ),(-1+\beta
)   \,  (z_{0}^{2}+z_{1}^{2} )   \,  (-z_{2}^{2}+\beta    \, z_{3}^{2} ),  \\
\noalign{\vspace{0.666667ex}} \hspace{1.em} (-{z_1}   \, {z_2}+{z_0}
\, {z_3})   \,  (\beta    \,  (z_{0}^{2}+z_{1}^{2}
)+z_{2}^{2}+{{\beta }^2}   \, z_{3}^{2} ),-{z_2}   \, {z_3}   \,
(z_{0}^{2}+z_{1}^{2}+z_{2}^{2}+z_{3}^{2} ),({z_0}   \, {z_1}-{z_2}
\, {z_3})   \,  (-z_{2}^{2}+\beta    \, z_{3}^{2} ),
 \\
\noalign{\vspace{0.666667ex}} \hspace{1.em} ({z_1}   \, {z_2}+{z_0}
\, {z_3})   \,  (-z_{2}^{2}+\beta    \, z_{3}^{2} ),-\beta    \,
z_{0}^{2}   \, z_{3}^{2}-z_{2}^{2}   \,
 (z_{1}^{2}+(1+\beta )   \, z_{3}^{2} ),-\beta    \, z_{1}^{2}   \, z_{3}^{2}-z_{2}^{2}   \,  (z_{0}^{2}+(1+\beta )   \, z_{3}^{2} ),
 \\
\noalign{\vspace{0.666667ex}} \hspace{1.em} (-{z_0}   \, {z_2}+{z_1}
\, {z_3})   \,  (z_{2}^{2}-\beta    \, z_{3}^{2} ),(-{z_0}   \,
{z_1}-{z_2}   \, {z_3})   \,
 (z_{2}^{2}-\beta    \, z_{3}^{2} ),({z_0}   \, {z_2}+\beta    \, {z_1}   \, {z_3})   \,  (z_{2}^{2}-\beta    \, z_{3}^{2} ),
 \\
\noalign{\vspace{0.666667ex}} \hspace{1.em} {z_2}   \, {z_3}   \,
(\beta    \,  (z_{0}^{2}+z_{1}^{2} )+z_{2}^{2}+{{\beta }^2}   \,
z_{3}^{2} ),(-{z_1}   \,
{z_2}+\beta    \, {z_0}   \, {z_3})   \,  (z_{2}^{2}-\beta    \, z_{3}^{2} ) \}\\
\end{matrix}
\end{equation}

The detailed analysis shows that the characteristic variety is the
union of the two points with coordinates
$(z_0,z_1,z_2,z_3)=(0,0,1,0)$, $(z_0,z_1,z_2,z_3)=(0,0,0,1)$, of the
line $\{(z_0,z_1,0,0)\}$ and of the two conics obtained by
intersecting the quadric $(z_{0}^{2}+z_{1}^{2}+z_{2}^{2}+z_{3}^{2}
)=0$ with the two hyperplanes $(z_{2}^{2}-\beta    \, z_{3}^{2}
)=0$.

The correspondence $\sigma$ is the identity on the two points, and
is given on the line by $$\sigma(z_0,z_1,0,0)=\,(z_0,-z_1,0,0)\,. $$
It exchanges the two conics and is a rational isomorphism of one
with the other, moreover the square of $\sigma$ admits two fixed
points on each of them and its derivative at the fixed points is
given by the squares of the following complex numbers
$$
\frac{i+\sqrt \beta}{i-\sqrt \beta}\,,\quad \frac{i-\sqrt
\beta}{i+\sqrt \beta}\,.
$$

\subsection{Edge $\alpha\perp\beta$ and $(n_1,n_2)=$(even, odd).
$\,L =\{ \left( \frac{\pi}{2}   , \varphi  , \varphi \right)\}$}

\noindent
\medskip

From now on we no longer use the change of variables to the Sklyanin
algebras but we rely directly on the explicit form of the minors of
the original matrix as computed in the Appendix \ref{appendix}. Note
in particular that the parameters $x_j$ are no longer the same as
the above $z_j$ but this is irrelevant since we compute intrinsic
invariants of the quadratic algebra. In the case at hand the list of
minors simplifies (with non zero scale factors removed) to the
following,

\begin{equation}
\begin{matrix}
 \{({x_1}\, {x_2}+i \, {x_0}\, {x_3})\,  (-x_{0}^{2}+x_{1}^{2}+x_{2}^{2}+x_{3}^{2} ),  \\
\noalign{\vspace{0.666667ex}}
\hspace{1.em} ({x_0}\, {x_2}+i \, {x_1}\, {x_3})\,  (x_{0}^{2}-x_{1}^{2}-x_{2}^{2}-x_{3}^{2} ),  \\
\noalign{\vspace{0.666667ex}} \hspace{1.em} ({x_0}\, {x_2}+i \,
{x_1}\, {x_3})\,  (x_{0}^{2}-x_{1}^{2}+x_{2}^{2}+x_{3}^{2} ),\,
(x_{0}^{2}-x_{1}^{2} )\,
 (x_{2}^{2}+x_{3}^{2} ),\,   \\
\noalign{\vspace{0.666667ex}} \hspace{1.em} (-i \, {x_1}\,
{x_2}+{x_0}\, {x_3})\,  (-x_{0}^{2}+x_{1}^{2}-x_{2}^{2}-x_{3}^{2}
),0,0,0,0,0,0,0,
 \\
\noalign{\vspace{0.666667ex}}
\hspace{1.em} 0,0,0 \}\\
\end{matrix}
\end{equation}

Thus the characteristic variety consists in the six lines $\ell_j$
given in terms of free parameters $x_j$ by
\begin{equation}
\begin{matrix}
\{(0,0,x_2,x_3)\}\,,\;\{(x_0,x_1,0,0)\}\,,\;\{(x_0,x_0,x_2,-i
x_2)\}\,,\;
\{(x_0,x_0,x_2,i x_2)\}\,,\;\\
\{(x_0,-x_0,x_2,i x_2)\}\,,\;\{(x_0,-x_0,x_2,-i x_2)\}\,.\;
\end{matrix}
\end{equation}
The correspondence $\sigma$ is the identity on the first line, $-1$
on the second and permutes cyclically the four others $\ell_j$.
Passing from $\ell_4$ to $\ell_5$ or from $\ell_6$ to $\ell_3$ one
gets the coarse correspondence, while the other maps are rational
isomorphisms. The following three lines meet and their point of
intersection is mapped by the coarse correspondence to the indicated
line
$$
\ell_1\cap \ell_4 \cap \ell_5 \to \ell_5\,,\quad \ell_1\cap \ell_3
\cap \ell_6 \to \ell_3\,,\quad \ell_2\cap \ell_3 \cap \ell_4 \to
\ell_5\,,\quad \ell_2\cap \ell_5 \cap \ell_6 \to \ell_3\,.\quad
$$
This is coherent as the restriction of the relevant coarse
correspondence.

\subsection{Edge $\alpha-\beta$ and $(n_1,n_2)=$ (odd, odd) or (even, odd).
$\,L' =\{ \left( \frac{\pi}{2} , \frac{\pi}{2} , \varphi \right)\}$
}

\noindent
\medskip

In that case the list of minors simplifies (with non zero scale
factors removed) to the following,

\begin{equation}
\begin{matrix}
  \{  \, {x_0}  \, {x_3}  \,   (-x_{0}^{2}+x_{1}^{2}+x_{2}^{2}+{{\sin (\varphi )}^2}  \, x_{3}^{2}  ),  \,   \\
\noalign{\vspace{0.666667ex}} \hspace{1.em} {x_1}  \, {x_3}  \,
(-x_{0}^{2}+x_{1}^{2}+x_{2}^{2}+{{\sin (\varphi )}^2}  \, x_{3}^{2}
),  \, x_{3}^{2}  \,
({x_0}  \, {x_2}+i   \, \sin (\varphi )  \, {x_1}  \, {x_3}),  \\
\noalign{\vspace{0.666667ex}} \hspace{1.em}   (x_{0}^{2}-x_{1}^{2}
)  \, x_{3}^{2},  \, x_{3}^{2}  \, (-i   \, {x_1}  \, {x_2}+\sin
(\varphi )  \,
{x_0}  \, {x_3}),  \,   \\
\noalign{\vspace{0.666667ex}} \hspace{1.em} {x_2}  \, {x_3}  \,
(-x_{0}^{2}+x_{1}^{2}+x_{2}^{2}+{{\sin (\varphi )}^2}  \, x_{3}^{2}
),  \, x_{3}^{2}  \,
({x_0}  \, {x_1}-i   \, \sin (\varphi )  \, {x_2}  \, {x_3}),  \,   \\
\noalign{\vspace{0.666667ex}} \hspace{1.em} x_{3}^{2}  \, (i   \,
{x_1}  \, {x_2}+\sin (\varphi )  \, {x_0}  \, {x_3}),  \,
(x_{0}^{2}-x_{2}^{2}  )  \,
x_{3}^{2},  (x_{1}^{2}+x_{2}^{2}  )  \, x_{3}^{2},  \\
\noalign{\vspace{0.666667ex}} \hspace{1.em} x_{3}^{2}  \, (-{x_0}
\, {x_2}+i   \, \sin (\varphi )  \, {x_1}  \, {x_3}),x_{3}^{2}  \,
({x_0}  \,
{x_1}+i   \, \sin (\varphi )  \, {x_2}  \, {x_3}),  \,   \\
\noalign{\vspace{0.666667ex}}
\hspace{1.em} {x_1}  \, x_{3}^{3},  \, {x_2}  \, x_{3}^{3},  \, {x_0}  \, x_{3}^{3}  \}\\
\end{matrix}
\end{equation}
Thus the characteristic variety contains the hyperplane $x_{3}=0$.
For $x_{3}\neq 0$ the last three minors show that all other
coordinates vanish and this gives an additional point, not in the
above hyperplane. The correspondence $\sigma$ is the symmetry
$$
\sigma(x_0,x_1,x_2,0)=\,(-x_0,x_1,x_2,0)\,.
$$

\subsection{Edge $\alpha-\beta$ and $(n_1,n_2)=$(even, even).
$\,L'' =\{ \left( \varphi  , \varphi  , \varphi \right)\}$}

\noindent
\medskip

In that case the list of minors simplifies (with non zero scale
factors removed) to the following,

\begin{equation}
\begin{matrix}
 \{\, {x_0}\, {x_3}\,  ({{\sin (\varphi )}^2}\, x_{0}^{2}-x_{1}^{2}-x_{2}^{2}-x_{3}^{2} ),\, {x_0}\, {x_2}\,
 ({{\sin (\varphi )}^2}\, x_{0}^{2}-x_{1}^{2}-x_{2}^{2}-x_{3}^{2} ),\,   \\
\noalign{\vspace{0.666667ex}} \hspace{1.em} x_{0}^{2}\, (\sin
(\varphi )\, {x_0}\, {x_2}+ i \, {x_1}\, {x_3}),\, x_{0}^{2}\,
(x_{2}^{2}+x_{3}^{2} ),\,
  \\
\noalign{\vspace{0.666667ex}} \hspace{1.em} x_{0}^{2}\, (- i \,
{x_1}\, {x_2}+\sin (\varphi )\, {x_0}\, {x_3}),\, {x_0}\, {x_1}\,
 ({{\sin (\varphi )}^2}\, x_{0}^{2}-x_{1}^{2}-x_{2}^{2}-x_{3}^{2} ),\,   \\
\noalign{\vspace{0.666667ex}} \hspace{1.em} x_{0}^{2}\, (\sin
(\varphi )\, {x_0}\, {x_1}- i \, {x_2}\, {x_3}),\, x_{0}^{2}\, ( i
\, {x_1}\, {x_2}+\sin (\varphi )\, {x_0}\, {x_3}),\,   \\
\noalign{\vspace{0.666667ex}} \hspace{1.em} x_{0}^{2}\,
(x_{1}^{2}+x_{3}^{2} ),\, x_{0}^{2}\,  (x_{1}^{2}+x_{2}^{2} ),\,
x_{0}^{2}\, (\sin
(\varphi )\, {x_0}\, {x_2}- i \, {x_1}\, {x_3}),\,   \\
\noalign{\vspace{0.666667ex}} \hspace{1.em} x_{0}^{2}\, (\sin
(\varphi )\, {x_0}\, {x_1}+ i \, {x_2}\, {x_3}),\, x_{0}^{3}\,
{x_2},\,
x_{0}^{3}\, {x_1},\, x_{0}^{3}\, {x_3} \}\\
\end{matrix}
\end{equation}
Thus the characteristic variety contains the hyperplane $x_{0}=0$.
For $x_{0}\neq 0$ the last three minors show that all other
coordinates vanish and this gives an additional point, not in the
above hyperplane. The correspondence $\sigma$ is the identity.

\subsection{Edge $\alpha\perp\beta$ and $(n_1,n_2)=$(even, even).
$\,C_+ =\{ \left(  \varphi  , \varphi,0 \right)\}$}\label{chartheta}

\noindent
\medskip

In that case the list of minors simplifies (with non zero scale
factors removed) to the following,

\begin{equation}
\begin{matrix}
  \{  \, {x_0}  \,   (x_{1}^{2}+x_{2}^{2}  )  \, {x_3},  (x_{1}^{2}+x_{2}^{2}  )  \, (-i  \, \cos (\varphi )  \,
{x_0}  \, {x_2}+\sin (\varphi )  \, {x_1}  \, {x_3}),  \,   \\
\noalign{\vspace{0.666667ex}} \hspace{1.em} (\sin (\varphi )  \,
{x_0}  \, {x_2}+i  \, \cos (\varphi )  \, {x_1}  \, {x_3})  \,
(x_{0}^{2}+x_{3}^{2}  ),  (x_{0}^{2}  \,
x_{2}^{2}-x_{1}^{2}  \, x_{3}^{2}  ),  \,   \\
\noalign{\vspace{0.666667ex}} \hspace{1.em} {x_1}  \, {x_2}  \,
(x_{0}^{2}+x_{3}^{2}  ),  \,   (x_{1}^{2}+x_{2}^{2}  )  \, (\cos
(\varphi )  \, {x_0}  \,
{x_1}-i  \, \sin (\varphi )  \, {x_2}  \, {x_3}),  \,   \\
\noalign{\vspace{0.666667ex}} \hspace{1.em} (\sin (\varphi )  \,
{x_0}  \, {x_1}-i  \, \cos (\varphi )  \, {x_2}  \, {x_3})  \,
(x_{0}^{2}+x_{3}^{2}  ),  \,
{x_1}  \, {x_2}  \,   (x_{0}^{2}+x_{3}^{2}  ),  \,   \\
\noalign{\vspace{0.666667ex}} \hspace{1.em}   (x_{0}^{2}  \,
x_{1}^{2}-x_{2}^{2}  \, x_{3}^{2}  ),  \, \sin (4  \, \varphi )  \,
(x_{1}^{2}+x_{2}^{2}  )  \,
  (x_{0}^{2}+x_{3}^{2}  ),  \,   \\
\noalign{\vspace{0.666667ex}} \hspace{1.em}   (x_{1}^{2}+x_{2}^{2}
)  \, (i  \, \sin (\varphi )  \, {x_0}  \, {x_2}+\cos (\varphi )  \,
{x_1}  \,
{x_3}),  \,   \\
\noalign{\vspace{0.666667ex}} \hspace{1.em}   (x_{1}^{2}+x_{2}^{2}
)  \, (\sin (\varphi )  \, {x_0}  \, {x_1}+i  \, \cos (\varphi )  \,
{x_2}  \,
{x_3}),  \,   \\
\noalign{\vspace{0.666667ex}} \hspace{1.em} (\cos (\varphi )  \,
{x_0}  \, {x_2}-i  \, \sin (\varphi )  \, {x_1}  \, {x_3})  \,
(x_{0}^{2}+x_{3}^{2}  ),  \,
  \\
\noalign{\vspace{0.666667ex}} \hspace{1.em} (-i  \, \cos (\varphi )
\, {x_0}  \, {x_1}+\sin (\varphi )  \, {x_2}  \, {x_3})  \,
(x_{0}^{2}+x_{3}^{2}  ),{x_0}  \,
  (x_{1}^{2}+x_{2}^{2}  )  \, {x_3}  \}\\
\end{matrix}
\end{equation}

Thus the characteristic variety consists in the six lines $\ell_j$
given in terms of free parameters $x_j$ by
\begin{equation}
\begin{matrix}
\{(0,x_1,x_2,0)\}\,,\;\{(x_0,0,0,x_3)\}\,,\;\{(x_0,x_1,i x_1,i
x_0)\}\,,\;
\{(x_0,x_1,-i x_1,i x_0)\}\,,\;\\
\{(x_0,x_1,i x_1,-i x_0)\}\,,\;\{(x_0,x_1,-i x_1,-i x_0)\}\,.\;
\end{matrix}
\end{equation}

The correspondence $\sigma$ is the identity on the first two lines.
It preserves globally the other $\ell_j$ and induces on each of them
the rational automorphism given by multiplication by $e^{\pm 2i
\varphi}$.

\subsection{Edge $\alpha\perp\beta$ and $(n_1,n_2)=$(odd, odd).
$\,C_- =\{ \left( \frac{\pi}{2}+\varphi, \frac{\pi}{2}  ,
\varphi\right)\}$}

\noindent
\medskip

In that case the list of minors simplifies (with non zero scale
factors removed) to the following,

\begin{equation}
\begin{matrix}
 \{(\sin(\varphi)x_1x_2+i\cos(\varphi) x_0x_3)(x^2_0-x^2_2),\ x_1x_3(x^2_0-x^2_2),\\
\noalign{\vspace{0.666667ex}}
 x_0x_2(x^2_1-x^2_3),\ x^2_0x^2_3-x^2_1x^2_2,\\
 \noalign{\vspace{0.666667ex}}
 (i\cos(\varphi) x_1x_2 -\sin(\varphi)x_0x_3)(x^2_1-x^2_3),\\
\noalign{\vspace{0.666667ex}}
 (-i\sin(\varphi)x_0x_1+\cos(\varphi)x_2x_3)(x^2_0-x^2_2),\\
\noalign{\vspace{0.666667ex}}
 (\cos(\varphi)x_0x_1-i\sin(\varphi)x_2x_3)(x^2_0-x^2_2),\\
\noalign{\vspace{0.666667ex}}
 (\cos(\varphi)x_1x_2-i\sin(\varphi) x_0 x_3) (x^2_0-x^2_2),\\
\noalign{\vspace{0.666667ex}}
 \sin(4\varphi)(x^2_0-x^2_2)(x^2_1-x^2_3),\\
\noalign{\vspace{0.666667ex}}
 (x^2_0x^2_1-x^2_2x^2_3),\ x_0x_2(x^2_1-x^2_3),\\
\noalign{\vspace{0.666667ex}}
 (\cos(\varphi)x_0x_1+i\sin (\varphi)x_2x_3)(x^2_1-x^2_3),\ x_1x_3(x^2_0-x^2_2),\\
 \noalign{\vspace{0.666667ex}}
 (i \sin(\varphi)x_0x_1+\cos(\varphi)x_2x_3)(x^2_1-x^2_3),\\
\noalign{\vspace{0.666667ex}}
 (\sin(\varphi)x_1x_2-i\cos(\varphi) x_0x_3)(x^2_1-x^2_3)\}
  \end{matrix}
\end{equation}

Thus the characteristic variety consists in the six lines $\ell_j$
given in terms of free parameters $x_j$ by
\begin{equation}
\begin{matrix}
\{(0,x_1,0,x_3)\}\,,\;\{(x_0,0,x_2,0)\}\,,\;\{(x_0,x_1,x_0,x_1)\}\,,\;
\{(x_0,x_1,-x_0,-x_1)\}\,,\;\\
\{(x_0,x_1,x_0,-x_1)\}\,,\;\{(x_0,x_1,-x_0,x_1)\}\,.\;
\end{matrix}
\end{equation}

The correspondence $\sigma$ is $-1$ on the first two lines. It
permutes $\ell_3$ with $\ell_4$ and its square is the rational
automorphism multiplying the ratio $x_1/x_0$ by $e^{4i \varphi}$. It
permutes $\ell_5$ with $\ell_6$ and its square is the rational
automorphism multiplying the ratio $x_1/x_0$ by $e^{-4i \varphi}$.

\subsection{Vertex $P=(\frac{\pi}{2},\frac{\pi}{2},\frac{\pi}{2})$}

\noindent
\medskip

In that case all minors vanish identically. Thus the characteristic
variety is all projective space. The correspondence $\sigma$ is
given by
$$
\sigma((x_0,x_1,x_2,x_3))=\,(-x_0,x_1,x_2,x_3)\,,
$$
but it degenerates on the quadric
$$
{\mathcal Q}= \,\{x\,|\, x_0^2- \sum x_k^2=\,0\}\,,
$$
to a correspondence which assigns to any point $p\in {\mathcal Q}$ a
line $\ell(p)\subset {\mathcal Q}$ containing the point $\sigma(p)$
and belonging to one of the two families of lines that rule the
surface ${\mathcal Q}$.

\subsection{Vertex $P'=(\frac{\pi}{2},\frac{\pi}{2},0)$}

\noindent
\medskip

In that case all minors vanish identically. Thus the characteristic
variety is all projective space. The correspondence $\sigma$ is
given by
$$
\sigma((x_0,x_1,x_2,x_3))=\,(-x_0,x_1,-x_2,x_3)\,.
$$

\subsection{Vertex $O=(0,0,0)$}

\noindent
\medskip

In that case all minors vanish identically and the correspondence
$\sigma$ is the identity.

\section{Isomorphism classes of  ${\mathbb R}_{\varphi}^4$
and orbits of the flow $F$} \label{flowbetter}

We let as above ${\mathcal M}$ be the moduli space of oriented
non-commutative $3$-spheres and $P{\mathcal M}$ its quotient by the
symmetry given by proposition \ref{liniso} 3). For $\varphi \in
{\mathcal M}$ we view the algebra ${\mathcal A}=C_{\rm alg}({\mathbb
R}_{\varphi}^4)$ as a graded algebra \ie we endow it with the one
parameter group of automorphisms which rescale the generators
$x_\nu$,
\begin{equation}\label{grading}
\theta_\lambda \in {\rm Aut}{\mathcal A}\,,\quad \theta_\lambda
(x_\nu)=\,\lambda \,x_\nu\qqq \lambda \in {\mathbb R}^*\,.
\end{equation}

We let $P{\mathcal M}$ be the quotient of the real moduli space
$\mathcal M$ by the symmetry of proposition \ref{liniso} 3).

This section will be devoted to prove the following result:

  \begin{thm}\label{transflow}  Let $\varphi_j \in {\mathcal M}$
 the following conditions are equivalent:
\begin{enumerate}
\item[a)] The graded algebras $C_{\rm alg}({\mathbb R}_{\varphi_j}^4)$ are isomorphic.
\item[b)] $\varphi_2 \in \hbox{Flow line of} \ \varphi_1$
in  $P{\mathcal M}$.
\end{enumerate}
 \end{thm}

The proof of $b)\Rightarrow a)$ was given above in section
\ref{flow}. The converse is based on the information given by the
geometric data which is by construction an invariant of the graded
algebra. The proof will be broken up in the non-generic and generic
cases.

\subsection{Proof in the non-generic case}

\noindent
\medskip

To analyze the information given by the geometric data we can
restrict the parameters $\varphi$ to the fundamental domain $A \cup
B$ of proposition \ref{fundom}. The symmetry given by proposition
\ref{liniso} 3) is given explicitely by the transformation

\begin{equation}\label{symrho}
\rho(\varphi_1,\varphi_2,\varphi_3)=\,(\varphi_1,\varphi_1-\varphi_3,
\varphi_1-\varphi_2)
\end{equation}

The identification $\gamma$ of the face $(P'QZ)$ with the face
$(ZPP')$ (proposition \ref{fundom}) followed by the symmetry $\rho$
\eqref{symrho} gives the following symmetry $\sigma=\rho\circ\gamma$
of the face $(P'QZ)$,
\begin{equation}\label{symsigma}
\sigma(\varphi_1,\varphi_2,\varphi_3)=
\,(\pi-\varphi_1+\varphi_2,\varphi_2,\varphi_2-\varphi_3)
\end{equation}

For vertices we have three elements $0$, $P$, $P'$ modulo the action
of $\Gamma \rtimes W$. The geometric data allows to separate them.

One has a priori nine edges. They fall in five different classes 4)
5) 6)  7) 8).

 Using $\gamma$ and $\rho$ one sees
that the following four edges are equivalent:
$$[ZP]\sim [PP']\sim [QP']\sim [ZQ]$$
and are all in case 5).

Similarly, using $\rho$ one sees that the following two edges are
equivalent:
$$[OP]\sim [OQ]$$
and are both in case 6).

 The table of geometric datas shows that
the geometric data of ${\mathbb R}_{\varphi}^4$ determines in which
of the five cases 4)-8) one is. Thus when the flow is transitive in
the corresponding edge there is nothing to prove. This covers the
cases 4) 5) 6). The two cases  7) 8)  correspond to fixed points of
the flow. For $C_+$ one gets the edge in $[OP']\subset\bar A$ given
by $\{\varphi,\varphi,0\}$ with $0<\varphi<\frac{\pi}{2}$. The
geometric data gives back the set $e^{\pm 2i\varphi}$ and this
allows to recover $\varphi$. Thus distinct $\varphi$ give
non-isomorphic quadratic algebras.

For $C_-$ one gets the edge $[P'Z]\subset \bar B$ \ie
$\{(\frac{\pi}{2}+\varphi,\frac{\pi}{2},\varphi)\}$ whose interior
corresponds to $0<\varphi<\frac{\pi}{2}$. The geometric data gives
back the set $e^{\pm 4i\varphi}$. Thus there is an ambiguity
$\varphi\to \frac{\pi}{2}-\varphi$ in $\varphi$ knowing the
geometric data. To understand it let us note that in fact one checks
that $\sigma$ restricts to the edge $[P'Z]$ as $\varphi\to
\frac{\pi}{2}-\varphi$. This then accounts for the above ambiguity.

\medskip

We now have to deal with the faces. We start with those which are
odd (\ie $H_{\alpha,n}$ with $n$ odd). The identification $\gamma$
of the face $(P'QZ)$ with the face $(ZPP')$ (proposition
\ref{fundom}) shows that we just need to deal with $(ZPP')$ and with
the face $(QPP')$ which is common to $A$ and $B$.

For the face $(ZPP')$ the equation of the supporting hyperplane is
$\varphi_2=\,\frac{\pi}{2}$ and one is in case 3) with generic
elements of the form $(\varphi_1,\,\frac{\pi}{2},\varphi_3)$ where
$\varphi_3+\frac{\pi}{2}>\varphi_1>\frac{\pi}{2}>\varphi_3$.

The geometric data determines the square
$$
\sigma(\frac{i+\beta^{1/2}}{i-\beta^{1/2}})^2
$$
where $\beta=-\tan\varphi_3/\tan\varphi_1>0$. Then
$\frac{i+\beta^{1/2}}{i-\beta^{1/2}}$ is of modulus one and the
geometric data determines $\beta$ up to the ambiguity given by
$\beta \to 1/\beta$. But the face $(ZPP')$ admits the symmetry given
by
\begin{equation}\label{gammarho}
\gamma\circ \rho(\varphi_1,\varphi_2,\varphi_3)=
\,(\pi-\varphi_3,\,\pi-\varphi_2,\,\pi-\varphi_1)
\end{equation}
whose effect is precisely the transformation $\beta \to 1/\beta$.
Note that the segment joining $P$ to the middle of $[P'Z]$ is
globally invariant under the flow $X$.

For the face $(QPP')$ the equation of the supporting hyperplane is
$\varphi_1=\,\frac{\pi}{2}$ and one is in case 3) with generic
elements of the form $(\frac{\pi}{2},\,\varphi_2,\,\varphi_3)$ where
$\frac{\pi}{2}>\varphi_2>\varphi_3>0$. The geometric data determines
the square
$$
\sigma(\frac{i+\beta^{1/2}}{i-\beta^{1/2}})^2
$$
where $\beta=-\tan\varphi_2/\tan\varphi_3<0$. Then
$\frac{i+\beta^{1/2}}{i-\beta^{1/2}}$ is real and the geometric data
determines $\beta$ up to the ambiguity given by $\beta \to 1/\beta$.
But the inequality $\tan\varphi_2>\tan\varphi_3>0$ shows that in
fact $\beta \in ]-\infty, -1[$ so that the geometric data determines
$\beta$ uniquely.

To summarize we have up to symmetry only two odd faces, the
geometric data allows to decide (by $|q|=1$ or $q\in {\mathbb R}$)
in which case one is, and gives back the flow line up to the
remaining symmetries.

\medskip

Let us now consider the even faces (\ie $H_{\alpha,n}$ with $n$
even). Using $\rho$ we get the equivalence $(OPP')\sim (OQP')$. To
be able to use lemma \ref{convex} we concentrate on $(OPP')$ on
which $\varphi_3>0$. The equation of the supporting hyperplane is
$\varphi_1=\,\varphi_2$ and one is in case 2) with generic elements
of the form $(\varphi_1,\,\varphi_1,\,\varphi_3)$ where
$\frac{\pi}{2}>\varphi_1>\varphi_3>0$. The geometric data determines
$$
\sigma(\frac{i+\alpha^{1/2}}{i-\alpha^{1/2}})
$$
where $\alpha=\,\tan\,\varphi_1\,\tan(\varphi_1-\varphi_3)>0$. Then
$q=\frac{i+\alpha^{1/2}}{i-\alpha^{1/2}}$ is of modulus one. Thus
the geometric data determines $\alpha$ (since it determines
$\alpha^{1/2}$ up to sign).

There are however two other even faces namely $(OPQ)$ and $(ZPQ)$.
The equation of the supporting hyperplane is the same in both cases
and is $\varphi_2=\,\varphi_3$ and one is in case 2) with generic
elements of the form $(\varphi_1,\,\varphi_2,\,\varphi_2)$ where for
$(OPQ)$ one has $\frac{\pi}{2}>\varphi_1>\varphi_2>0$ while for
$(ZPQ)$ one gets
$\frac{\pi}{2}+\varphi_2>\varphi_1>\frac{\pi}{2}>\varphi_2$. The
geometric data determines
$$
\sigma(\frac{i+\alpha^{1/2}}{i-\alpha^{1/2}})
$$
where $\alpha=\,\tan\,\varphi_2\,\tan(\varphi_2-\varphi_1)<0$. Then
$q=\frac{i+\alpha^{1/2}}{i-\alpha^{1/2}}$ is real. This first allows
to distinguish these faces from the other even faces treated above.
Moreover one checks that
$$
\alpha\in ]-1,0[\qqq \varphi \in (OPQ) \,,\quad\alpha\in
]-\infty,-1[\qqq \varphi \in (ZPQ)\,.$$ Thus $q>0$ on $(OPQ)$ and
$q<0$ on $(ZPQ)$ which allows to distinguish these two faces from
each other. Finally on each of these faces the geometric data
determines $\alpha$ and hence the flow line of $\varphi$ using lemma
\ref{convex}.

 Thus we get the proof in all cases
except the generic case which we shall now treat in details
separately.

\medskip

\subsection{Basic notations for elliptic curves}\label{basicell}

\noindent
\medskip

 We
recall that given an elliptic curve $E$ viewed as a $1$-dimensional
complex manifold and chosing a base point $p_0 \in E$ one gets an
isomorphism of the universal cover $\tilde E$ of $E$ with base point
$p_0$
 $$\tilde E \overset{I}{\simeq} {\mathbb C}$$
given by the integral
$$
I(p)=\,\int_{p_0}^p \omega
$$
where $\omega $ is a holomorphic $(1,0)$-form. This isomorphism is
unique up to multiplication by $\lambda \in {\mathbb C}^*$.
 Let $L \subset {\mathbb C}$ be the lattice of periods
then $I^{-1}(L)$ is  the kernel of the covering map $\pi : \tilde E
\to E$ and one has an isomorphism $E \sim {\mathbb C}/L$. To
eliminate the choice of the base point we let $T(E)$ be the group of
translations of $E$ and note that the universal cover $\tilde T(E)$
identifies with the additive group ${\mathbb C}$,
$$\tilde T(E)\simeq{\mathbb C}\,,\quad  T(E)\simeq{\mathbb C}/L\,.$$
One can moreover take $L$ of the form $L=\,{\mathbb Z}+{\mathbb
Z}\tau$ where $\tau \in {\mathbb H}/\Gamma[1]$ and ${\mathbb H}$ is
the upper half plane ${\mathbb H}=\{z\in {\mathbb C}\,|\,{\rm
Im}z>0\}$ while in general $\Gamma[n]$ is the congruence subgroup of
level $n$ in $SL(2,{\mathbb Z})$. With ${\mathbb H}^*$ obtained from
${\mathbb H}$ by adjoining the rational points of the boundary one
has a canonical isomorphism ${\mathbb H}^*/\Gamma[1]\to {\mathbb
P}^1({\mathbb C})$ given by Jacobi's $j$ function. In terms of the
elliptic curve $E$ defined by the equation
$$
y^2=\,4x^3-g_2\,x-g_3\,,
$$
in ${\mathbb P}^2({\mathbb C})$ one has
$$
j(E)= 1728 \,\frac{g_2^3}{\Delta}
$$
where the discriminant is $\Delta=\,g_2^3-27\,g_3^2$.

One obtains a finer invariant $\lambda(E)$ if one has the additional
structure given by an isomorphism of abelian groups
$$
\phi: \quad ({\mathbb Z}/2{\mathbb Z})^2\to
\frac{1}{2}L/L=T_{2}(E)\,.
$$
where $T_{2}(E)$ is the group of two torsion elements of $T(E)$.
This allows to choose the module $\tau$ in a finer manner as $\tau
\in{\mathbb H}/\Gamma[2]$ and one has a canonical isomorphism
${\mathbb H}^*/\Gamma[2]\to {\mathbb P}^1({\mathbb C})$ given by
Jacobi's $\lambda$ function. In terms of the elliptic curve $E$
defined by the equation
$$
y^2=\,\prod\,(x-e_j)
$$
in ${\mathbb P}^2({\mathbb C})$ one has
$$
\lambda(E)= {\rm Cross \,Ratio}(e_1,e_2;e_3,e_4)
$$

\medskip
In more intrinsic terms the labelling of the two torsion $T_{2}(E)$
\begin{equation}
\label{eq1.95} \omega_j \in T(E) \, ,  \ \omega_1=\phi(1,0)  \, , \
\omega_3=\phi(0,1)\, , \ \omega_2=-\omega_1-\omega_3  \, .
\end{equation}
allows to define the following function on the group $T(E)$ of
translations of $E$:
\begin{equation}
\label{invpfun} F_\phi(u)=\,\frac{\wp_3 (u)}{\wp_3 (\omega_1)} \,,
\end{equation}
where $\wp_3$ is defined using a fixed isomorphism $\tilde T(E)
\simeq {\mathbb C}$, as the sum
\begin{equation}
\label{eq1.96} \wp_3 (u) = \left( \sum_{\pi (y) = u} y^{-2} \right)
- \left(\sum_{\pi (y) = \omega_3} y^{-2} \right)
\end{equation}
where one defines the sums by restricting $y$ to $\vert y \vert < R$
on both sides and then taking the limit. In standard notation with
the Weierstrass $\wp$-function given by
\begin{equation}
\label{eq1.97} \wp (z) = \frac{1}{z^2} + \sum_{\ell \in L^*} \left(
\frac{1}{(z+\ell)^2} - \frac{1}{\ell^2} \right)
\end{equation}
one gets
\begin{equation}
\label{eq1.98} \wp_3 (u) = \wp (z) - \wp (w_3) \, , \ \pi (z) = u \,
.
\end{equation}
 Note that the ratios involved in  (\ref{invpfun}) eliminate the scale factor $\lambda$ in the isomorphism $\tilde T(E) \simeq {\mathbb C}$.
With these notations one has
\begin{equation}
\label{invp} \lambda (E,\phi) = \, F(w_2)\,=
\frac{\wp(\omega_2)-\wp(\omega_3)} {\wp(\omega_1)-\wp(\omega_3)}\,.
\end{equation}

\smallskip

Finally the covering ${\mathbb H}^*/\Gamma[2]\to {\mathbb
H}^*/\Gamma[1]$ is simply given by the algebraic map
$$
\lambda \quad \to \quad 256 \,
\frac{(1-\lambda+\lambda^2)^3}{\lambda^2(1-\lambda)^2}
$$
while the group $\Sigma$ of deck transformations is the dihedral
group generated by the two symmetries
$$
u(z)= 1/z \quad v(z)=1-z\,.
$$
One has a unique anti-isomorphism $w:\,PSL(2,{\mathbb Z}/2{\mathbb
Z}) \to \Sigma$ such that
\begin{equation}
\label{invp1} \lambda (E,\phi\circ \alpha) = \, w(\alpha)(\lambda
(E,\phi))\,.
\end{equation}
Moreover $w(t)=v$ where $t$ is the transposition of $(a,b)\to
(b,a)$.
 Finally one gets,
\begin{equation}
\label{tF} F_{\phi\circ t}(u) = \, 1-F_{\phi}(u)\,.
\end{equation}

\bigskip

\subsection{The generic case}

\noindent
\medskip

We now deal with the generic case.

Let $s_j\in {\mathbb R}$ be three real numbers and  $(\alpha , \beta
, \gamma)$ be given by
\begin{equation}
\label{stoalpha} \alpha = \frac{s_3 - s_2}{s_1} \, , \ \beta =
\frac{s_1 - s_3}{s_2} \, , \ \gamma =  \frac{s_2 - s_1}{s_3} \, .
\end{equation}
 Let $ (E , \sigma)$ be the
pair of an elliptic curve and a translation associated to $(\alpha ,
\beta , \gamma)$ by proposition \ref{chargen}.

  \begin{lem} \label{trigtoell}
There exists an isomorphism $\phi: \;({\mathbb Z}/2{\mathbb Z})^2\to
T_2(E)$  such that with $e_j=\,s_j^{-1}$ one has
\begin{equation}
\label{elllambda} \lambda (E,\phi) = \hbox{Cross Ratio} \ (e_2 ,e_1
; e_3 , \infty) \, ,
\end{equation}
and with $F_\phi$ defined by \eqref{invpfun}
\begin{equation}
\label{ellF} F_\phi(\sigma) = \frac{s_1}{s_1 - s_3}
\end{equation}
\end{lem}
\bigskip

\begin{proof} The proof is straightforward using $\theta$-functions
to  parametrize the elliptic curve (\ref{chargen}) but we prefer  to
give an elementary direct proof. In order to prove \eqref{elllambda}
and \eqref{ellF} we start from proposition \ref{chargen} and replace
$(\alpha , \beta , \gamma)$ by their value. The equations of $E$
simplifies to
\begin{equation}
\label{simplerE} \sum_0^3\,x_\mu^2=0\,,\quad
\sum_1^3\,s_k\,x_k^2=0\,.
\end{equation}
We then rescale $x_1=a\,X_1$, $x_2=b\,X_2$, $x_3=\,X_3$, where
\begin{equation}
\label{rescale} a^2 =\, -\frac{s_3}{s_1}\,,\quad b^2 =\,
-\frac{s_3}{s_2}\,.
\end{equation}
After this rescaling the second equation of (\ref{simplerE}) gives
$ X_3^2 = \, X_1^2 + X_2^2 $ and one uses the standard rational
parametrization of the conic to
 parametrize the solutions of this equation as
\begin{equation}
\label{ratpar} x_1 = 2 \, a \, t \, , \ x_2 = b (1-t^2) \, , \ x_3 =
(1+t^2)\,.
\end{equation}
  One then writes the first equation in (\ref{chargen}) as
\begin{equation}
\label{eq1.105} x_0^2 + 4 \, a^2 \, t^2 + b^2 (1-t^2)^2 + (1+t^2)^2
= 0
\end{equation}
and since $\displaystyle 1+b^2 =\frac{(s_2-s_3)}{s_2}\ne 0$ this
reduces to the elliptic curve defined by the equation
\begin{equation}
\label{ellcurve1} y^2 = t^4 - 2 \, r \, t^2 + 1 \, ,
\end{equation}
where
$$
 r = - \frac{2 \, a^2 - b^2 + 1}{b^2 + 1} \, .$$
Using \eqref{rescale} one gets
\begin{equation}
\label{biq} r = \frac{s_1 (s_2 + s_3) - 2 \, s_2 \, s_3}{s_1 (s_3 -
s_2)} \, .
\end{equation}
One checks that $r\neq \pm 1$ since the $s_j$ are pairwise distinct.

Thus the roots of $\,x^4 - 2 \, r \, x^2 + 1 = 0$ are all distinct
and we write them as $\pm \, u$, $\pm \, v$ where
\begin{equation}
\label{roots} uv = 1\,,\quad u^2 + v^2 = 2r\,.
\end{equation}

\medskip

  The cross ratio of $(-v , u ; v , -u)$ is independent of the
above choices and is given by
\begin{equation}
\frac{-v-v}{-v+u} : \frac{u-v}{u+u} = - \frac{4 \, uv}{(u-v)^2} =
\frac{-4}{u^2 + v^2 - 2} = \, - \frac{2}{r - 1} =
\frac{1+b^2}{1+a^2} =  \frac{s_1}{s_2} \, \frac{s_2 - s_3}{s_1 -
s_3} \, .
\end{equation}
In other words with $e_j=\,s_j^{-1}$ one has
\begin{equation}
{\rm Cross\,Ratio}(-v , u ; v , -u)=\,{\rm Cross\,Ratio}(e_2 ,e_1 ;
e_3 , \infty)\, .
\end{equation}
Let then $\gamma \in {\rm SL}(2,{\mathbb C})$, $\gamma(X)=\frac{aX +
b}{cX + d}$ which transforms $(e_2 ,e_1 , e_3 , \infty)$ to $(-v , u
, v , -u)$. Then for a suitable choice of $\lambda \neq 0$ the
transformation,
$$t = \gamma(X)\,,\quad y = \frac{\lambda Y}{(cX + d)^2}\,.$$
gives a birational isomorphism $\tilde \gamma$ of the elliptic curve
defined by the equation
\begin{equation}
\label{ellcurve2} Y^2 = \prod_1^3 (X - e_i)
\end{equation}
with the curve (\ref{ellcurve1}). Choosing  the origin of
\eqref{ellcurve2} as the point at infinity one gets the point $P$
with coordinates $t=-u$, $y=0$ as the origin of (\ref{ellcurve1}).
One has $P = (p_0 , p_1 , p_2 , p_3)$ and $p_0=0$ since $y=0$. Thus
the other coordinates fulfill (\ref{simplerE}) in the simplified
form
\begin{equation}
\label{origin}
 p_1^2 +\, p_2^2+\,p_3^2 =\,0
\,,\quad s_1 \,p_1^2 + \,s_2\, p_2^2 +\,s_3 \,p_3^2= \,0 \,,
\end{equation}
 and in homogeneous expressions one can replace $p_k^2$
by $s_\ell-s_m$.

\medskip

  Let us now determine the translation $\sigma$.
We just compute the $t$ parameter of $\sigma (P)= (x_0 , x_1 , x_2 ,
x_3)$. The parameter $t$ of (\ref{ratpar}) is recovered in
homogeneous coordinates as,
\begin{equation}
\label{ratpar1} t = \frac{b \, x_1}{a \, x_2 + a \, b \, x_3} \, .
\end{equation}
  One first starts by simplifying (\ref{sigma}) when applied to $P$. One has $p_0 = 0$ and up to an overall scale ($-\alpha\,\beta \, \gamma$)
one gets using  \eqref{origin}, \ie the replacement $p_k^2\to
s_\ell-s_m$,

\begin{equation}
- \beta \, \gamma \, p_1^2 + \alpha \, \gamma \, p_2^2 + \alpha \,
\beta \, p_3^2 \to - \, s_1 + s_2 + s_3 \nonumber
\end{equation}
\begin{equation}
\beta \, \gamma \, p_1^2 - \alpha \, \gamma \, p_2^2 + \alpha \,
\beta \, p_3^2 \to s_1 - s_2 + s_3 \nonumber
\end{equation}
\begin{equation}
\beta \, \gamma \, p_1^2 + \alpha \, \gamma \, p_2^2 - \alpha \,
\beta \, p_3^2 \to s_1 + s_2 - s_3 \nonumber
\end{equation}
Thus the coordinates of $\sigma (P)$  are up to an overall scale
\begin{equation}
\label{eq1.109} x_1 = p_1 (-s_1 + s_2 + s_3) \, , \ x_2 = p_2 (s_1 -
s_2 + s_3) \, , \ x_3 = p_3 (s_1 + s_2 - s_3) \, ,
\end{equation}
where the $p_j$ are the coordinates of $P$. These are given by
(\ref{ratpar}) taking $t = -u$, with $u,v$ as above, thus,
\begin{equation}
\label{eq1.110} p_1 = -2 \, a \, u \, , \ p_2 = b (1-u^2) \, , \ p_3
= 1+u^2 \, .
\end{equation}
We then need to compute:
\begin{equation}
\label{eq1.111} t(\sigma) = \frac{b \, x_1}{a \, x_2 + a \, b \,
x_3} = \frac{b (-2 \, a \, u) (-s_1 + s_2 + s_3)}{a \, b (1-u^2)
(s_1 - s_2 + s_3) + a \, b (1+u^2) (s_1 + s_2 - s_3)} \, .
\end{equation}
We see that $a$ and $ b$ drop out and we remain with
\begin{equation}
\label{tsigma} t(\sigma) =  \frac{- \, u (-s_1 + s_2 + s_3)}{s_1 +
(s_2 - s_3) \, u^2} \, .
\end{equation}
We just need to compute ${\rm Cross\,Ratio}(-v , u, v , t(\sigma))$
which is the same as ${\rm Cross\,Ratio}(e_2 , e_1 ; e_3 , \tau)$
where $\tau$ is
 $\wp(\sigma)$, up to an affine transformation
which transforms the $e_j$ to $\wp(\omega_j)$. One has $$ {\rm
Cross\,Ratio}(-v , u, v , t(\sigma))=\,\frac{-v-v}{-v-t(\sigma)} :
\frac{u-v}{u-t(\sigma)} = \frac{u-t(\sigma)}{u-v} \,
\frac{2v}{v+t(\sigma)} = \frac{u^2 - u\,t(\sigma)}{u^2 - 1} \,
\frac{2}{1+u\,t(\sigma)}$$ using $u \, v=1$. The $u \, t(\sigma)$
only involves $u^2$ and one can simplify by its denominator to get
\begin{equation}\label{crfine}
{\rm Cross\,Ratio}(-v , u , v , t(\sigma)) = 2 \, \frac{(s_2 - s_3)
\, u^4 + (s_2 + s_3) \, u^2}{(s_1 - 2 \, s_3) \, u^4 + 2 \, s_3 \,
u^2 - s_1} \, .
\end{equation}
We then claim that one has:
\begin{equation}
\label{eq1.115} {\rm Cross\,Ratio}(-v , u, v, t(\sigma)) = \frac{s_2
- s_3}{s_1 - s_3} = {\rm Cross\,Ratio}(e_2 , e_1 ; e_3 , 0) \, .
\end{equation}
\ie that $\tau =0$. To see this we replace $u^4$ by $2 \, r \, u^2 -
1$ in the expression \eqref{crfine} for the cross ratio, which gives
\begin{equation}
{\rm Cross\,Ratio}(-v , u , v , t(\sigma))=  \, \frac{(2 (s_2 - s_3)
\, r + s_2 + s_3) \, u^2 - (s_2 - s_3)}{((s_1 - 2 \, s_3) \, r +  \,
s_3) \, u^2 -  (s_1 - s_3)} \, . \nonumber
\end{equation}
The computation using (\ref{biq}) shows that this fraction is
independent of $u^2$ and equal to $ \displaystyle \frac{s_2 -
s_3}{s_1 - s_3}$. This ends the proof of the lemma since we have
shown that, up to an affine tranformation,
$$
\wp(\omega_j)= e_j\,,\quad \wp(\sigma)=\,0\,.
$$
so that
$$
F_\phi(\sigma)=\,\frac{0- e_3}{e_1 - e_3}=\,\frac{s_1 }{s_1 -
s_3}\,.
$$
\end{proof}

We shall now give the proof of theorem \ref{transflow} in the
generic case. We start with the alcove $A$.

Given an elliptic curve $E$ and a translation $\sigma$ of $E$ we
claim that if there is a labelling $\phi$ of $T_2(E)$ such that
\begin{equation}\label{uniquelabel}
\lambda(E,\phi)\in ]0,1[\,,\quad F_\phi(\sigma)<0
\end{equation}
then this labelling is unique. Indeed the only element of
$PSL(2,{\mathbb Z}/2{\mathbb Z})$ which preserves the first
condition is the transposition $t$ and this does not preserve the
second by \eqref{tF}.

On $\sigma(A)=\{s \,|\,1<s_1<s_2<s_3\}$ one has (lemma
\ref{sigmacover})
$$\hbox{Cross Ratio} \ (s_2 , s_1 ; s_3 , 0)\in ]0,1[\,,
\quad \frac{s_1}{s_1 - s_3}<0 \,,$$ thus by lemma \ref{trigtoell}
there exists a unique labelling $\phi$ of $T_2(E)$ such that
\eqref{uniquelabel} holds. This gives back both $\hbox{Cross Ratio}
\ (s_2 , s_1 ; s_3 , 0)$ and $\displaystyle \frac{s_1}{s_1 - s_3}$.
The latter gives the ratio $\displaystyle\frac{s_3}{s_1}=a$ and the
former then gives $\displaystyle \frac{s_2-a s_1}{s_1 -a s_1}\frac{
s_1}{s_2} =\frac{s_2-a s_1}{(1 -a) s_2}$ which gives the ratio
$\displaystyle\frac{s_2}{s_1}=b$. Thus we recover the flow line
using the convexity of  $\sigma(A)$.

Let us now look at the alcove $B$. One has $\sigma(B)=\{s
\,|\,s_3<s_2<0, 1<s_1\}$. Thus $\displaystyle
s_2^{-1}<s_3^{-1}<s_1^{-1}$.

Exactly as above, given an elliptic curve $E$ and a translation
$\sigma$ of $E$ we claim that if there is at most one labelling
$\phi$ of $T_2(E)$ such that
\begin{equation}\label{uniquelabel1}
0< \,F_\phi(\sigma)\,<\,\lambda(E,\phi)\,<1\,.
\end{equation}
The algebra associated to the $s_j$ is unchanged, up to isomorphism,
by cyclic permutations of the $s_j$, and the same holds for the
associated geometric data. Thus  lemma \ref{trigtoell}  gives a
labelling $\phi$ such that
$$\lambda(E,\phi)=\hbox{Cross Ratio} \ (s_3 , s_2 ; s_1 , 0)\,,
\quad F_\phi(\sigma)=\frac{s_2}{s_2 - s_1} \,,$$ One checks that it
then fulfills \eqref{uniquelabel1} since with $e_j=s_j^{-1}$ one has
$$
\hbox{Cross Ratio}(s_3 , s_2 ; s_1 , 0)=\,\frac{e_1-e_3}{e_1-e_2}\,,
\quad\frac{s_2}{s_2 - s_1}=\,\frac{e_1}{e_1-e_2}\,.
$$
and
$$
1>\,\frac{e_1-e_3}{e_1-e_2} > \,\frac{e_1}{e_1-e_2}>0\,.
$$
Thus as above we recover the the ratios
$\displaystyle\frac{s_1}{s_2}$ and $\displaystyle\frac{s_3}{s_2}$
and the flow line of $\varphi$ using the convexity of  $\sigma(B)$.

Finally note that the conditions \eqref{uniquelabel} and
\eqref{uniquelabel1} are exclusive and thus allow to decide using
the geometric data wether $\varphi \in \sigma(A)$ or $\varphi \in
\sigma(B)$.

\section{Dualities}\label{dualities}

We show in this section that there are unexpected dualities between
the noncommutative spaces ${\mathbb R}_{\varphi}^4$ in the following
cases of Table \ref{tablebig}
$$A\leftrightarrow B\,,\quad 2\leftrightarrow 3\,,\quad 5\leftrightarrow 6\,,\quad 7\leftrightarrow 8\,,\quad 10\leftrightarrow 11\,,$$
Modulo these dualities the fundamental domain gets reduced from an
alcove of the root system $D_3$ to the smaller alcove of the root
system $C_3$.

\subsection{Semi-cross product.}

\noindent
\medskip

Let $\cala$ be a graded algebra and let $\alpha\in \Aut (\cala)$ be
a symmetry commuting with the grading (i.e. homogeneous of degree
$0$).
\begin{defn}
The semi-cross product $\cala(\alpha)$ of $\cala$ by $\alpha$ is the
graded vector space $\cala$ equipped with the bilinear product
$\,{.}_\alpha\,$ defined by
\end{defn}
\[
a\,{.}_\alpha\, b=a\,\alpha^n(b),\>\> \forall a \in \cala_n\,,\>
b\in \cala
\]
It is easily verified that this product is associative and that $\cala_m\,{.}_\alpha\, \cala_n\subset \cala_{m+n}$ so that $\cala(\alpha)$ is a graded algebra and that if $\cala$ is unital then $\cala(\alpha)$ is also  unital with the same unit.\\
Some basic properties of the semi-cross product are given by the
following proposition.
\begin{prop}\label{geodat}
Let $\cala$ be a graded algebra and let $\alpha$ be an automorphism of degree $0$ of $\cala$.\\
(i) If $\beta$ is an automorphism of degree $0$ of $\cala$ which
commutes with $\alpha$, then $\beta$ is also canonically an
automorphism of degree $0$ of $\cala(\alpha)$ and one has
\[
\cala(\alpha)(\beta)=\cala(\alpha\circ \beta).
\]
In particular one has
\[
\cala(\alpha)(\alpha^{-1})=\cala.
\]
(ii) If $\cala$ is a quadratic algebra, then $\cala(\alpha)$ is also
a quadratic algebra and its geometric datas $(E',\sigma', \call')$
are deduced from those $(E,\sigma,\call)$ of $\cala$ as follows
\[
E'=E,\>\> \sigma'=\alpha^t\circ \sigma,\>\> \call'=\call
\]
where $\alpha^t$ is induced by the transposed of $(\alpha\restriction \cala_1)$.\\
(iii) If $\cala$ is an involutive algebra with involution $x\mapsto
x^\ast$ homogeneous of degree $0$ (in short, if $\cala$ is a graded
$\ast$-algebra) and if $\alpha$ commutes with the involution, then
one defines an antilinear antimultiplicative mapping $x\mapsto
x^{\ast_\alpha}$ of $\cala(\alpha)$ onto $\cala(\alpha^{-1})$ by
setting $a^{\ast_\alpha}=\alpha^{-n}(a^\ast)$ for $a\in \cala_n$. In
particular if $\alpha^2=1$, then $\cala(\alpha)$ equipped with the
involution $x\mapsto x^{\ast_\alpha}$ is a graded $\ast$-algebra.
\end{prop}
\begin{proof} (i) One has for $a\in \cala_n$ and $b\in \cala$
\[
\beta(a\,{.}_\alpha\,b)=\beta(a\alpha^n(b))=\beta(a)\beta(\alpha^n(b))=\beta(a)\alpha^n(\beta(b))=\beta(a)\,{.}_\alpha\,
\beta(b)
\]
which shows that $\beta$ is an automorphism of $\cala(\alpha)$. One
has also
\[
a \,{.}_\alpha\,\beta^n(b)=a\alpha^n(\beta^n(b))=a(\alpha\circ
\beta)^n(b)= a \,{.}_{\alpha \circ \beta}\; b
\]
which implies $\cala(\alpha)(\beta)=\cala(\alpha\circ \beta)$.\\
(ii) Assume that $\alpha$ is a quadratic algebra i.e.
$\cala=A(V,R)=T(V)/(R)$ where $V$ is finite-dimensional and where
$(R)$ is the two-sided ideal of the tensor algebra $T(V)$ of $V$
generated by the subspace $R$ of $V\otimes V$. Let $\mg$ denote the
product of $\cala$ and $\mg'$ denote the product of $\cala(\alpha)$.
Since $V=\cala_1=\cala(\alpha)_1$ one has $\mg'=\mg\circ (\id
\otimes \alpha)$ on $V\otimes V$ and thus $\mg(R)=0$ is equivalent
to $\mg'((\id \otimes \alpha^{-1})(R))=0$ from which it follows
easily that $\cala(\alpha)=A(V,R')=T(V)/(R')$ with
\[
R'=(\id \otimes \alpha^{-1})R
\]
so $\cala(\alpha)$ is quadratic.\\
By definition the graph of $\sigma'$ is the subset
\[
\Gamma'\subset P(V^\ast) \times P(V^\ast)
\]
obtained from the subset of $V^\ast \times V^\ast$ of pairs
$(\omega,\pi),\> \omega\not =0, \pi\not=0$ such that
\[
\langle \omega\otimes \pi, r\rangle=0,\>\> \forall r\in R'.
\]
Since $R'=(\id\otimes \alpha^{-1})R$ we thus get $\sigma'=\alpha^t\circ \sigma$.\\
(iii) One has for $a\in \cala_n$ and $b\in \cala_m$
\[
\alpha^{-(n+m)}((a\,{.}_\alpha\,b)^\ast)=\alpha^{-(n+m)}
(\alpha^n(b^\ast)a^\ast)=\alpha^{-m}(b^\ast)\alpha^{-m}
(\alpha^{-n}(a^\ast))
\]
and thus
\[
(a\,{.}_\alpha\,
b)^{\ast_\alpha}=b^{\ast_\alpha}\,{.}_{\alpha^{-1}}\,
\,a^{\ast_\alpha}
\]
which implies (iii).
\end{proof}

More generally if $\cala$ is finitely generated in degree $1$ and
finitely presented i.e. if $\cala=T(V)/(R)$ with $V$
finite-dimensional and $(R)$ is the two-sided ideal of $T(V)$
generated by the graded subspace $R=\oplus^N_{n=2} R_n$ $(R_n\subset
V^{\otimes^n})$ of $T(V)$, one has $\cala(\alpha)=T(V)/(R(\alpha))$
with $R(\alpha)=\oplus^N_{n=2}R_n(\alpha)$,
\[
R_n(\alpha)=(\id\otimes \alpha^{-1}\otimes \dots \otimes
\alpha^{-(n-1)})R_n.
\]
In particular $\cala(\alpha)$ is a $N$-homogeneous algebra whenever
$\cala$ is $N$-homogeneous. Let us recall that
an algebra $\cala$ as above is said to be $N$-homogeneous iff $R=R_N\subset V^{\otimes^N}$ \cite{ber:2001a}, \cite{ber-mdv-wam:2003}. For these algebras, which generalize the quadratic algebras $(N=2)$, one has a direct extension of the Koszul duality of quadratic algebras as well as a natural generalization of the notion of Koszulity \cite{ber:2001a}, \cite{ber-mdv-wam:2003}. The stability of the corresponding homological notions with respect to the semi-cross product construction will be studied elsewhere.\\

The terminology semi-cross product of $\cala$ by $\alpha$ for
$\cala(\alpha)$ comes from the fact that it can be identified with a
subalgebra of the crossed product $\cala\semi_\alpha \mathbb Z$,
namely the subalgebra generated by the elements
\[
xW,\>\> x\in \cala_1
\]
where $W$ denotes the new invertible generator of the crossed
product defined by $WaW^{-1}=\alpha(a)$ for $a\in \cala$. Indeed one
has for $a\in \cala_m$ $b\in\cala_p$
\[
aW^nbW^p=a\alpha^n(b)\> W^{n+p}
\]
If $\cala$ is a graded $\ast$-algebra with $\alpha$ a $\ast$-homomorphism of degree $0$, one endows the crossed product of a structure of $\ast$-algebra by setting $W^\ast=W^{-1}$. The involution of the crossed product induces then, by restriction, the antilinear antimultiplicative mapping $\ast_\alpha:\cala(\alpha)\rightarrow \cala(\alpha^{-1})$ of $(iii)$ in the above proposition.\\
For the following application, we shall have $\alpha^2=1$ so
$\cala(\alpha)$ can then be identified with the corresponding
subalgebra of the crossed product $\cala\semi_\alpha \mathbb
Z/2\mathbb Z$ and if $\cala$ is a graded $\ast$-algebra with
$\alpha$ a $\ast$-homomorphism of degree $0$, $\cala\rtimes_\alpha
\mathbb Z/2\mathbb Z$ becomes a $\ast$-algebra by setting
$W=W^\ast(=W^{-1})$ and $\cala(\alpha)$ is a $\ast$-subalgebra.

\subsection{Application to $\mathbb R^4_\varphi$ and $S^3_\varphi$.}

\noindent
\medskip

The above construction allows to give a duality between the
following cases of Table \ref{tablebig}
$$A\leftrightarrow B\,,\quad 2\leftrightarrow 3\,,\quad 5\leftrightarrow 6\,,\quad 7\leftrightarrow 8\,,\quad 10\leftrightarrow 11\,,$$
where $A$ and $B$ are the two simplices that together complete case
$1$).

The explicit transformation on the $\varphi$-parameters is
\begin{equation}
\label{dual5} f_1(\varphi)=\,(\pi -\varphi_1,
\,\frac{\pi}{2}-\varphi_1+\varphi_2,\,
\frac{\pi}{2}-\varphi_1+\varphi_3)\,.
\end{equation}
but there are two others which give relevant additional dualities,
\begin{equation}
\label{dual6} f_2(\varphi)=\,(\frac{\pi}{2} +\varphi_2-\varphi_3,
\,\varphi_2,\, \frac{\pi}{2}-\varphi_1+\varphi_2)\,.
\end{equation}
\begin{equation}
\label{dual7} f_3(\varphi)=\,(\frac{\pi}{2} -\varphi_2+\varphi_3,
\,\frac{\pi}{2}-\varphi_1+\varphi_3,\, \varphi_3)\,.
\end{equation}
They are all involutions $f_j^2=1$ and on the fundamental domain
$A\cup B$ one has
$$
f_1(A)=\,B\,,\quad f_1(B)=\,A\,,\quad f_2(B)=\,B\,,\quad
f_3(A)=\,A\,.
$$
with $f_1(0PQP'Z)=(ZPQP'0)$, $f_2(PQP'Z)=(PQZP')$ and
$f_3(OPQP')=(P'PQO)$. Thus the duality $f_2$ (resp. $f_3$) operates
only in $B$ (resp. $A$). Moreover the transformations $f_2$ and
$f_3$ are conjugate under $f_1$ \ie $f_2=\,f_1\circ f_3\circ f_1$.

\medskip

By proposition \ref{liniso} 2) any element $v$ of the centralizer
$C\subset {\rm SO}(4)$ of the diagonal matrices in ${\rm SU}(4)$
 defines an automorphism of
 ${\mathbb R}_{\varphi}^4$ acting by
$v$ on the generators. We thus let $\alpha_j \in {\rm Aut}{\mathcal
A}$ correspond to the diagonal matrices with respective diagonals
given by
\begin{equation}
\label{signs} \alpha_1\,:\,(1,1,-1,-1)\,,\quad
\alpha_2\,:\,(1,-1,1,-1)\,,\quad \alpha_3\,:\,(1,-1,-1,1)\,,
\end{equation}
\bigskip

\begin{prop}
For $j\in\{1,2,3\}$ one has canonical $\,*$-isomorphisms
 \[
\rho_j:C_{\rm alg}(\mathbb R^4_{f_j(\varphi)})\rightarrow C_{\rm
alg} (\mathbb R^4_\varphi)(\alpha_j)
\]
which preserve the central element $\sum x_\mu^2$ and induce
corresponding isomorphisms ${S}_{f_j(\varphi)}^3 \simeq
{S}_{\varphi}^3(\alpha_j)\,.$

\end{prop}
\medskip

\begin{proof} One writes explicitely the isomorphisms as follows on the canonical
(self-adjoint) generators $x^\mu_j$ of $C_{\rm alg}(\mathbb
R^4_{f_j(\varphi)})$, (which are the canonical noncommutative
coordinates of $\mathbb R^4_{f_j(\varphi)}$)

\begin{equation}
\rho_1(x_1^0)=\,x^1\,W_1\,,\quad \rho_1(x_1^1)=\,x^0\,W_1\,  ,\quad
\rho_1(x_1^2)=-i \,x^2\,W_1\,,\quad
\rho_1(x_1^3)=-i\,x^3\,W_1\,,\nonumber
\end{equation}
\begin{equation}
\rho_2(x_2^0)=\,x^2\,W_2\,, \quad \rho_2(x_2^1)=-i\,x^3\,W_2\,,
\quad \rho_2(x_2^2)= \,x^0\,W_2\,, \quad
\rho_2(x_2^3)=i\,x^1\,W_2\,, \nonumber
\end{equation}
\begin{equation}
\rho_3(x_3^0)=\,x^3\,W_3\,, \quad \rho_3(x_3^1)=-i\,x^2\,W_3\,,
\quad \rho_3(x_3^2)= i\,x^1\,W_3\,, \quad
\rho_3(x_3^3)=\,x^0\,W_3\,, \quad  \nonumber
\end{equation}

\medskip
Using the signs in \eqref{signs} one checks that one gets
$\,*$-isomorphisms  ${C_{\rm alg}(\mathbb R}_{f_j(\varphi)}^4)
\simeq  {C_{\rm alg}(\mathbb R}_{\varphi}^4)(\alpha_j)$ and that one
has $\rho_j(\sum_\mu(x_j^\mu)^2)=\sum_\mu(x^\mu)^2$\, .
\end{proof}

\bigskip

\begin{figure}
\begin{center}

\includegraphics[scale=0.9]{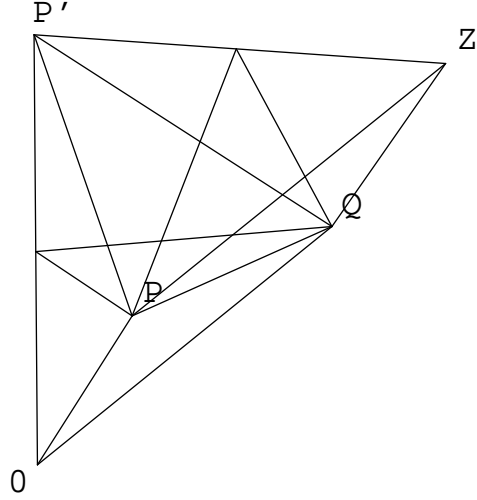}
\end{center}

\caption{\label{domain3}  The hyperplanes $f_j(\varphi)= \varphi$}
\end{figure}

\bigskip
\begin{cor}
\begin{enumerate}
\item
Let $\varphi$ be generic and $(E,\sigma,{\mathcal L})$ be the
geometric data of ${C_{\rm alg}(\mathbb R}_{\varphi}^4)$. Then the
geometric data of ${C_{\rm alg}(\mathbb R}_{f_j(\varphi)}^4)$ is
$(E,\sigma + \tau_j,{\mathcal L})$ where $\tau_j \in T_2(E)$.

\item The hyperplane $f_j(\varphi)= \varphi$
is globally invariant under the flow $Z$.

\item Let $\varphi \in A\cup B$ be generic, then it belongs
to the hyperplane $f_j(\varphi)= \varphi$ (with $j=2$ on $B$ and
$j=3$ on $A$) iff the translation $\sigma$ fulfills
$$
\sigma \in T_4(E)\,.
$$
\end{enumerate}
\end{cor}
\begin{proof} 1) Using proposition \ref{geodat} the required
equality follows if one shows that the action of $\alpha_j$ on $E$
is indeed given by a translation $\tau_j\in T_2(E)$. This can be
checked directly using $\vartheta$-functions \ie proposition
\ref{theta}.

\end{proof}

\medskip
These new symmetries suggest that one extends the group $\Gamma
(\Tg) \rtimes W$ of section \ref{alcoves} to include the $T_j$. This
amounts to adding the following new roots to the root system
$\Delta$ of section \ref{alcoves}. To the $\pm \psi_\mu\,\pm
\psi_\nu$ we adjoin the $\pm 2 \psi_\mu$. This is in fact the same
as the replacement
$$
SO(6)\to Sp(3)\,,
$$
of the compact group $SO(6)$ by the symplectic group $Sp(3)$. (\ie
of $D_3$ by $C_3$). The corresponding Weyl group is now $O(3,\mathbb
Z)$.

\smallskip
The description of $ Sp(3)$ is obtained as follows using the natural
representation $\rho\,:\,\mathbb H \to M_2(\mathbb C) $ of the field
of quaternions $\mathbb H$ as two by two matrices of the form
$\begin{bmatrix} \alpha & \beta \\ -\bar \beta &\bar \alpha
\end{bmatrix}$ with $\alpha, \beta \in \mathbb C$. By definition
$Sp(3)$ is the group of three by three matrices $Q \in M_3(\mathbb
H)$ whose image $\rho(Q)\in M_6(\mathbb C)$ is unitary. Thus it
contains as a subgroup
$$
\{  Q \in M_3(\mathbb H)\,|\, \rho(Q)\in U(6)\cap M_6(\mathbb R)\}
$$
which is isomorphic to $U(3)$ and is contained in $SO(6)$ but has
only 9 parameters.

\section{The algebras $\calg(\mathbb R^4_\varphi)$ in the nongeneric cases}\label{degalg}

In this section we shall describe the algebras $\calg(\mathbb
R^4_\varphi)$ in the nongeneric cases using the above dualities.

\subsection{$U_q(\fracsl(2))$, $U_q(\fracsu(2))$ and their homogenized versions}

\noindent
\medskip

In this section for $q\in \mathbb C_\ast$, $U_q(\fracsl(2))$ is
considered as an associative algebra while for $\vert q\vert=1$ or
$q\in \mathbb R_\ast$, $U_q(\fracsu(2))$ is considered as a
$\ast$-algebra which is a real form of $U_q(\fracsl(2))$. The
coalgebra aspect plays no role in the following and we refer to
\cite{maj:1995} for a very complete discussion of these topics. It
is convenient to start by the homogenized version. For $q\in \mathbb
C\backslash \{-1,0,1\}$ we define $U_q(\fracsl(2))^\hom$ to be the
quadratic algebra generated by 4 elements $X^+,X^-,K^+,K^-$ with
relations
\begin{equation}\label{k1}
K^+K^-=K^-K^+
\end{equation}
\begin{equation}\label{k2}
K^+X^+=qX^+K^+
\end{equation}
\begin{equation}\label{k3}
K^+X^-=q^{-1}X^-K^+
\end{equation}
\begin{equation}\label{k4}
K^-X^+=q^{-1}X^+K^-
\end{equation}
\begin{equation}\label{k5}
K^-X^-=qX^-K^-
\end{equation}
\begin{equation}\label{k6}
[X^+,X^-]=\frac{(K^+)^2-(K^-)^2}{q-q^{-1}}
\end{equation}
The ``classical limit" $U(\fracsl(2))^\hom=U_1(\fracsl(2))^\hom$ for
$q=1$ is obtained by setting $q=1+\varepsilon$ and $K^\pm=X^0\pm
\frac{\varepsilon}{2}X^3$. Letting $\varepsilon \rightarrow 0$, the
relations read then
\begin{equation}\label{x1}
[X^0,X^3]=0
\end{equation}
\begin{equation}\label{x2}
[X^0,X^+]=0
\end{equation}
\begin{equation}\label{x3}
[X^0,X^-]=0
\end{equation}
\begin{equation}\label{x4}
[X^3,X^+]=2X^0X^+
\end{equation}
\begin{equation}\label{x5}
[X^3,X^-]=-2X^0X^-
\end{equation}
\begin{equation}\label{x6}
[X^+,X^-]=X^0X^3
\end{equation}
To obtain $U_q(\fracsl(2))$, one notices that relations \eqref{k1}
\dots \eqref{k5} imply that $K^+K^-$ is central so that one may add
the inhomogeneous relation
\begin{equation}\label{kpm}
K^+K^-=\bbbone
\end{equation}
which together with relations \eqref{k1} \dots \eqref{k6} defines $U_q(\fracsl(2))$, i.e. $U_q(\fracsl(2))$ is (as associative algebra) the quotient of $U_q(\fracsl(2))^\hom$ by the two-sided ideal generated by $K^+K^--\bbbone$.\\
Similarily one notices that relations \eqref{x1} \dots \eqref{x5} imply that $X^0$ is central  and the universal enveloping algebra $U(\fracsl(2))$ is obtained by adding $X^0=1$ to the relations \eqref{x1} \dots \eqref{x6}. One has of course $\lim_{q\rightarrow 1} U_q(\fracsl(2))=U(\fracsl(2))$.\\
For $q\in \mathbb C\backslash \mathbb R$ with $\vert q\vert=1$ the
real version $U_q(\fracsu(2))^\hom$ of $U_q(\fracsl(2))^\hom$ is
obtained by endowing the algebra with the unique antilinear
antimultiplicative involution such that
\begin{equation}\label{xpm}
(X^\pm)^\ast=X^\mp
\end{equation}\label{xkpm}
\begin{equation}
(K^\pm)^\ast=K^\mp
\end{equation}
while for $q\in \mathbb R\backslash \{-1,0,1\}$,
$U_q(\fracsu(2))^\hom$ corresonds to the unique antilinear
antimultiplicative involution such that
\[
(X^\pm)^\ast=X^\mp
\]
\begin{equation}
(K^\pm)^\ast=K^\pm
\end{equation}
which gives
\[
(X^\pm)^\ast=X^\mp
\]
\begin{equation}
(X^0)^\ast=X^0
\end{equation}
\begin{equation}
(X^3)^\ast=X^3
\end{equation}
for the limiting case $q=1$, i.e. for $U(\fracsu(2))^\hom$.\\
These involutions pass to the quotient to define the $\ast$-algebras
$U_q(\fracsu(2))$ for $\vert q\vert=1$ or $q\in \mathbb R_\ast$.
Notice that the involution of $U_q(\fracsu(2))$ is obtained from the
one of $U(\fracsu(2))$ by setting $K^\pm=q^{\pm\frac{1}{2}X^3}$ (the
relations in terms of $X^\pm$ and $X^3$ differ of course).

\subsection{The algebras in the nongeneric cases}

\noindent
\medskip

We now identify the algebra $\calg(\mathbb R^4_\varphi)$ in the
cases 2 to 11 of
the table \ref{tablebig}.\\

\noindent 2. {\bf Even face} : $\varphi_1=\varphi_2=\varphi,\,
\varphi-\varphi_3 \not\in \frac{\pi}{2} \mathbb Z, \,
\varphi_k\not\in \frac{\pi}{2} \mathbb Z$. By setting
\begin{equation}\label{even0}
X^\pm = x^1\pm ix^2
\end{equation}
the relations (\ref{pres1}), (\ref{pres2}) of the algebra read then
\begin{equation}\label{even1}
[x^0,x^3]=0
\end{equation}
\begin{equation}\label{even2}
\cos(\varphi)[x^0,X^+]=\sin(\varphi-\varphi_3)(x^3X^++X^+x^3)
\end{equation}
\begin{equation}\label{even3}
\cos(\varphi-\varphi_3)[x^3,X^+]=-\sin(\varphi)(x^0X^++X^+x^3)
\end{equation}
\begin{equation}\label{even4}
\cos(\varphi)[x^0,X^-]=-\sin(\varphi-\varphi_3)(x^3X^-+X^-x^3)
\end{equation}
\begin{equation}\label{even5}
\cos(\varphi-\varphi_3)[x^3,X^-]=\sin(\varphi)(x^0X^-+X^-x^0)
\end{equation}
\begin{equation}\label{even6}
[X^+,X^-]=-2\sin(\varphi_3)(x^0x^3+x^3x^0)
\end{equation}
and since (\ref{even0}) implies $(X^+)^\ast=X^-$ one sees that
relations \eqref{even2} and \eqref{even3} are the adjoints of
relations \eqref{even4} and \eqref{even5} respectively.
 We now distinguish the following 2 regions $R$ and $R'$.\\

\noindent $R:\frac{\pi}{2}>\varphi>\varphi_3\geq 0$. By setting
\begin{equation}
K^\pm=(2\tan(\varphi_3))^{1/2}((\sin(2\varphi))^{1/2}x^0\pm
i(\sin(2(\varphi-\varphi_3)))^{1/2}x^3)
\end{equation}
the relations \eqref{even1} to \eqref{even6} become \eqref{k1} to
\eqref{k6}  with
\begin{equation}
q=\frac{1-i(\tan(\varphi-\varphi_3)
\tan(\varphi))^{1/2}}{1+i(\tan(\varphi-\varphi_3)\tan(\varphi))^{1/2}}
\end{equation}
so $\calg(\mathbb R^4_\varphi)$ coincides then with $U_q(\fracsu(2))^\hom$ for $\vert q \vert=1$,  $q\not= \pm 1$ as $\ast$-algebra (one has $(K^+)^\ast=K^-$).\\

\noindent $R':\frac{\pi}{2}>\varphi>0$ and
$\frac{\pi}{2}+\varphi>\varphi_3>\varphi$. By setting
\begin{equation}
K^\pm=(2\tan(\varphi_3))^{1/2}((\sin(2\varphi))^{1/2}x^0\pm
(\sin(2(\varphi_3-\varphi)))^{1/2} x^3)
\end{equation}
the relations \eqref{even1} to \eqref{even6} become \eqref{k1} to
\eqref{k6}  with
\begin{equation}
q=\frac{1-(\tan(\varphi_3-\varphi)\tan(\varphi))^{1/2}}{1+(\tan(\varphi_3-\varphi)\tan
(\varphi))^{1/2}}
\end{equation}
so $\calg (\mathbb R^4_\varphi)$ coincides then with $U_q(\fracsu(2))^\hom$ for $q\in ]-1,0[\cup]0,1[$ (one has $(K^\pm)^\ast=K^\pm$ here).\\

In general for case 2, it is easy to see that $\calg(\mathbb R^4_\varphi)$ is isomorphic as $\ast$-algebra either to an algebra of $R$ or of $R'$. Notice that $q>1$, for instance, is the same as $q\in ]0,1[$ by echanging $q$ and $q^{-1}$ and $K^+$ and $K^-$. Thus for case 2 the $\calg(\mathbb R^4_\varphi)$ are the $U_q(\fracsu(2))^\hom$.\\

\noindent 3. {\bf Odd face} : $\varphi_1=\frac{\pi}{2}$, $\varphi_2-\varphi_3 \not\in \frac{\pi}{2}\mathbb Z$, $\varphi_2\not\in \frac{\pi}{2}\mathbb Z$, $\varphi_3\not\in \frac{\pi}{2}\mathbb Z$. Using the analysis of Section \ref{dualities}, one sees that the cases corresponding to the plane $\varphi_1=\frac{\pi}{2}$ are in $\alpha_3$-duality with the cases corresponding to the plane $\varphi_2=\varphi_3$. On the other hand, the cases corresponding to the plane $\varphi_2=\varphi_3$ are the same as the cases corresponding to the plane $\varphi_1=\varphi_2$. So finally (taking into account the forbidden values $\frac{\pi}{2}\mathbb Z$) one sees that the case 3 (odd face) is obtained by duality (in the sense of Section 9) from the case 2 (even face)  i.e. from the $U_q(\fracsu(2))^{\hom}$.\\

\noindent 4. $\alpha \perp \beta\>\>\>   (0,1):
\varphi_1=\frac{\pi}{2},\, \varphi_2=\varphi_3=\varphi\not\in
\frac{\pi}{2}\mathbb Z$. This case is singular in the sense that one
of the 6 relations is missing namely
\[
\cos(\varphi_1)[x^0,x^1]=i\sin (\varphi_2-\varphi_3)[x^2, x^3]_+
\]
which gives ``0=0". This implies exponential growth. This is a singular limiting case of $U_q(\fracsu(2))^\hom$ for $q=0$ which separates the regions $1<q<0$ and $0<q<1$ of case 2.\\

\noindent 5. $\alpha-\beta\>\>\>  (0,1):\varphi_1=\frac{\pi}{2},\,
\varphi_2=\frac{\pi}{2},\, \varphi_3\not\in \frac{\pi}{2}\mathbb Z$.
This case is obtained by $\alpha_3$-duality (section
\ref{dualities}) from $\varphi_1=\varphi_2=\varphi_3\not\in
\frac{\pi}{2}\mathbb Z$ which is case 6 below. It corresponds to
$U_{-1}(\fracsu(2))^\hom$.
\\

\noindent 6. $\alpha-\beta\>\>\>  (0,0): \varphi_1=\varphi_2=\varphi_3\not\in \frac{\pi}{2}\mathbb Z$. The relation (\ref{pres1}), (\ref{pres2}) reduce in this case to the relations \eqref{x1} to \eqref{x6} by setting $X^\pm=x^1\pm ix^2,\, X^3=2x^3$ and $X^0=-2\sin (\varphi_1)x^0$. Thus in this case the $\ast$-algebra is isomorphic to $U(\fracsu(2))^\hom=U_1(\fracsu(2))^\hom$.\\

\noindent 7. $\alpha\perp \beta\>\>\>
(0,0):\varphi_1=\varphi_2=-\frac{1}{2}\theta\not\in
\frac{\pi}{2}\mathbb Z$, $\varphi_3=0$. This is the
``$\theta$-deformation" studied in \cite{ac-lan:2001} and
\cite{ac-mdv:2002a}. By setting
\begin{equation}
z^1=x^0+ix^3,\, z^2=x^1+ix^2,\, \bar z^1=(z^1)^\ast,\, \bar
z^2=(z^2)^\ast
\end{equation}
the relations (\ref{pres1}), (\ref{pres2}) read
\begin{equation}
\left\{
\begin{array}{l}
z^1z^2=e^{i\theta}z^2z^1\\
\bar z^1z^2=e^{-i\theta}z^2\bar z^1\\
z^1\bar z^2=e^{-i\theta}\bar z^2z^1\\
\bar z^1\bar z^2=e^{i\theta}\bar z^2\bar z^1
\end{array}\right.
\end{equation}
and we refer to Part I \cite{ac-mdv:2002a} for more details and generalizations of this algebra.\\

\noindent 8. $\alpha\perp \beta\>\>\> (1,1):\varphi_1=\frac{\pi}{2}+\varphi,\, \varphi_2=\frac{\pi}{2},\, \varphi_3=\varphi\not\in \frac {\pi}{2}\mathbb Z$. This case is obtained from the preceding case 7 ($\theta$-deformation) by $\alpha_1$-duality as explained in section \ref{dualities}.\\

\noindent 9. $\varphi_1=\varphi_2=\varphi_3=\frac{\pi}{2}$. This is the most singular case, 3 relations are missing (i.e. reduce to ``0=0") among the 6 relations (\ref{pres1}), (\ref{pres2}). The algebra has exponential growth.\\

\noindent 10. $\varphi_1=\varphi_2=\frac{\pi}{2},\, \varphi_3=0$.
This case which is at the intersection of the lines carrying case 7
and case 8 is obtained by $\alpha_3$-duality
(section \ref{dualities}) from the next case 11.\\

\noindent 11. $\varphi_1=\varphi_2=\varphi_3=0$. This is the
``classical case", the relations (\ref{pres1}), (\ref{pres2}) reduce
to $x^\mu x^\nu=x^\nu x^\mu$ for $\mu,\nu\in \{0,1,2,3\}$ so the
algebra reduces to the algebra of polynomials $\mathbb C[x^0,x^1,
x^2, x^3]$.

\bigskip
\bigskip

\noindent \underbar{Remark} : One can go much further and describe
the  $C^*$-algebras corresponding to the noncommutative spheres
$S_\varphi^3$ in all these degenerate cases. It is important in that
respect to classify the discrete series besides the obvious
continuous series of representations.

\bigskip
\bigskip

\section{The Complex Moduli Space and its Net of Elliptic Curves}
\label{complexmod}

We first explain in this section a striking
 coincidence in the generic
case between the geometric data of ${\mathbb R}_{\varphi}^4$ and the
fiber of the double cover \ref{sigmacover}. This takes place in the
real moduli space and leads us to introduce the complex moduli space
in which the equality between the two elliptic curves makes full
sense. We then describe the geometric structure of the complex
moduli space as a net of elliptic curves in three dimensional
projective space and exhibit a presentation of the algebras making
the equality ``fiber $=$ characteristic" manifest.

We showed above in the proof of lemma \ref{trigtoell} that the
geometric data of ${\mathbb R}_{\varphi}^4$ can be interpreted as
the elliptic curve
\begin{equation}
\label{ellcurve2p} E_2=\{(X,Y)\,|\,Y^2 = \prod_1^3 (X - e_i)\}\,
\end{equation}
with $e_j= \,s_j^{-1}$ and endowed with a translation $\sigma$
sending the point at $\infty$ to a point of $E_2$ whose $X$
coordinate vanishes\footnote{the two choices give isomorphic pairs
$(E_2,\sigma)$}.

One can write $E_2$ in the equivalent form,
\begin{equation}
\label{ellcurve3} E_3(\varphi)=\{(X,Y)\,|\,Y^2 = \prod_1^3 (X \,s_j
- 1)\}\,
\end{equation}
In this form this equation is the same as the one involved in the
double cover \ref{sigmacover}. More precisely one gets

\begin{figure}
\qquad \quad  \includegraphics[scale=0.7]{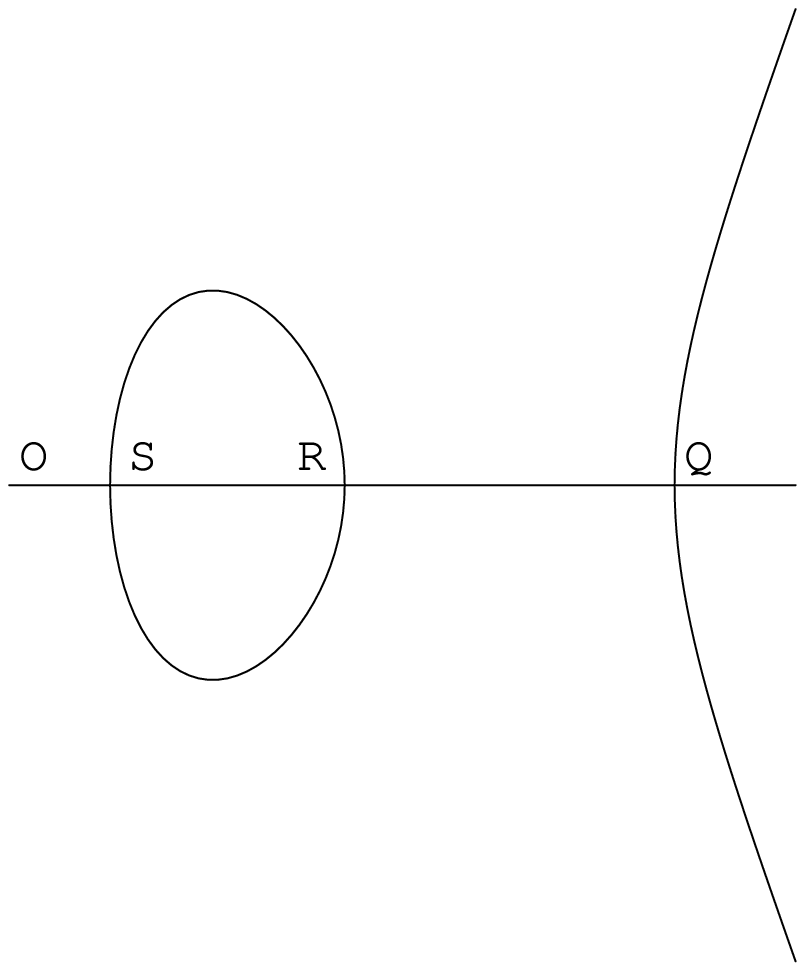}\qquad \qquad
\includegraphics[scale=0.7]{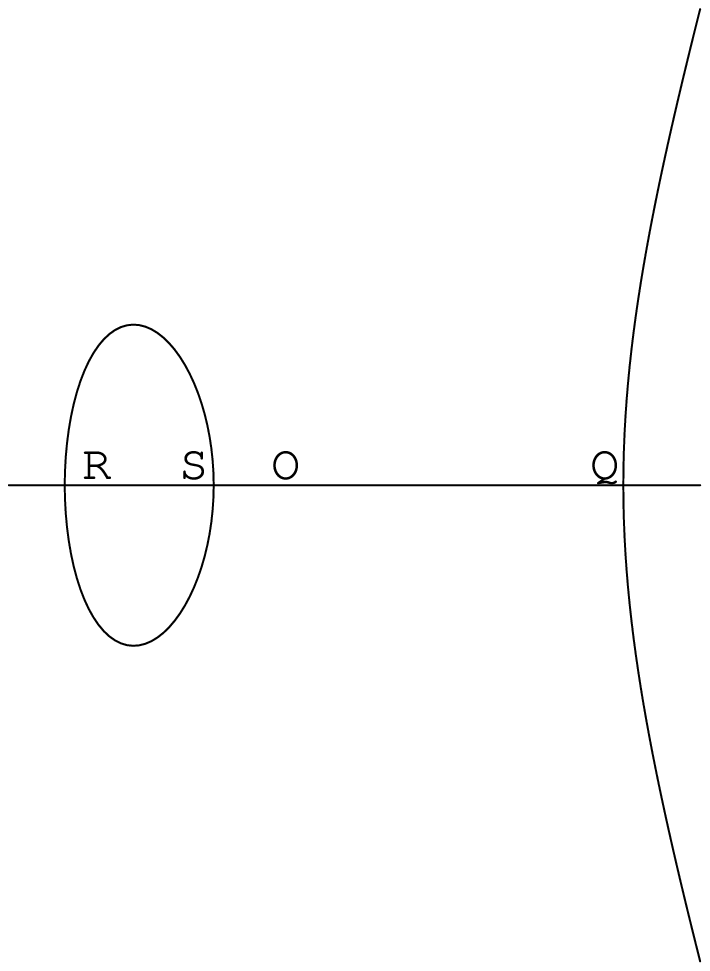}
\caption{\label{ellcurve} The elliptic curve $E_3$ and the fiber for
$\varphi \in A$ and $\varphi \in B$. }
\end{figure}

\begin{prop}\label{fiberchar1} Let $\varphi$ be generic,
and
$${\rm Fiber}\,(\varphi)=\,\{\varphi'\,{\rm generic}\,|\,
 J_{\ell m} (\varphi') = J_{\ell m} (\varphi) \, , \; \forall k\}$$
a) There is a natural isomorphism $\ell_\varphi\,:\,{\rm
Fiber}\,(\varphi)\simeq E_3(\varphi)\cap({\mathbb R}\times {\mathbb
R}^*)$ determined by the equalities
\begin{equation}
\label{ident} X(\varphi')=\,s_j(\varphi')/s_j \, , \; \forall
j\,,\quad Y(\varphi')=\,\prod \tan(\varphi'_j)\,.
\end{equation}
b) The closure  $\overline{\rm Fiber}\,(\varphi)$ is the union ${\rm
Fiber}\,(\varphi)\cup W(P)$ and $\ell_\varphi$ extends to an
isomorphism
$$  \ell_\varphi\,:\;\overline{\rm Fiber}\,(\varphi)\to
E_3(\varphi)\cap {\mathbb P}_2({\mathbb R})\,.$$
\end{prop}
\begin{proof} a) Let $\varphi'\in {\rm Fiber}\,(\varphi)$
then by lemma \ref{propor} the ratio $X=s_j(\varphi')/s_j$ is
independent of $j$. Moreover one has $$\prod_1^3 (X \,s_j -
1)=\,\left(\prod \tan(\varphi'_j)\right)^2\,,$$ thus the map is well
defined. Conversely given $(X,Y)\in E_3(\varphi)\cap({\mathbb
R}\times {\mathbb R}^*)$ one lets $s'_k=\,X\,s_k$ and  $t'_k =
Y(s'_k - 1)^{-1}$. This determines uniquely the $\varphi'_j$ such
that $\tan(\varphi'_j)=t'_j$.

b) Note that ${\rm Fiber}\,(\varphi)$ is not connected but falls in
four connected components each a flow line of the scaling flow $X$.
These correspond to the four components of
$E_3(\varphi)\cap({\mathbb R}\times {\mathbb R}^*) $ \ie of the set
of real points of $E_3(\varphi)$ which are not of two torsion. The
inverse map $\ell_\varphi^{-1}$ extends to the four points
$\{\infty,e_1,e_2,e_3\}$ of $E_3(\varphi)$ and one thus obtains a
compactification of ${\rm Fiber}\,(\varphi)$ by adding the four
points $\{P,Q,R,S\}$ in the orbit $W(P)$. The above map
$\ell_\varphi$ extends continuously to the compactification and the
four points $\{P,Q,R,S\}$ map to the four two torsion points
$\{\infty,e_1,e_2,e_3\}$ of $E_3(\varphi)$ and in that order for
$\varphi\in A$ as in Figure \ref{ellcurve}. The situation is similar
on $B$ but the order of $e_2$ and $e_3$ is reversed as well as that
of $R$ and $S$. Moreover the point $0$ is now between $Q$ and $S$
instead of being on the left of $S$ as in case $\varphi\in A$ (see
Figure \ref{ellcurve}).
\end{proof}

\medskip
In fact the action of the discrete symmetry given by the Klein group
$H\subset W$ preserves globally ${\rm Fiber}\,(\varphi)$ since it
does not alter the $J_{\ell m} (\varphi)$. In fact this discrete
symmetry admits a simple interpretation in terms of translations of
order two $\tau \in T_2(E_3)$ as follows, where we let $k_j\in H$ be
given by
$$
k_1(\varphi)=\,(-\varphi_1,\varphi_3-\varphi_1,\varphi_2-\varphi_1)
$$

\begin{prop}\label{weyltwo}
\begin{enumerate}
\item[1)] Given $\varphi'\in {\rm Fiber}\,(\varphi)$
the map $(X,Y)\to (\lambda \,X,Y)$ for $\lambda=\,s_j(\varphi')/s_j
\, , \; \forall j$ is an isomorphism  $E_3(\varphi')\simeq
E_3(\varphi)$ compatible with the isomorphisms $\ell_\varphi$ and
$\ell_{\varphi'}$ of lemma \ref{fiberchar1}.
\item[2)]
 ${\rm Fiber}\,(\varphi)$ is globally invariant
under $k_j$ and under the above isomorphism ${\rm
Fiber}\,(\varphi)\simeq {\mathbb R}^2\cap E_3(\varphi)$ the action
of $k_j$ is given by the translation by the two torsion element
$(e_j, 0)\in E_3(\varphi)$.
\end{enumerate}
\end{prop}
\begin{proof} 1) By construction the
$s_j(\varphi')$ are proportional to the $s_j$ so the conclusion
follows looking at the definition \ref{ident} of the identification
${\rm Fiber}\,(\varphi)\simeq {\mathbb R}^2\cap E_3(\varphi)$.

2) Using 1) it is enough to show that the point of $E_3(\varphi)$
associated to $k_1(\varphi)$ is obtained from $(1, t_1\,t_2\,t_3)\in
E_3(\varphi)$ (with $t_j=\tan(\varphi_j)$) by the translation of
order two associated to $(e_1,0)\in E_3(\varphi)$. One checks that
the effect of $k_1$ on the $s_j(\varphi)$ is to multiply all of them
by $X'=\,(1+ t_1^2)/(s_2\,s_3)$.
 Its
effect on $Y=\,\prod \tan(\varphi_j)$ is to replace it by
$$
Y'=-\,t_1\,(t_3-t_1)\,(t_2-t_1)/(s_2\,s_3)\,.
$$
One then checks by direct computation that the line joining the
points $(1, t_1\,t_2\,t_3)\in E_3(\varphi)$ and $(X',-Y')\in
E_3(\varphi)$ intersects $E_3(\varphi)$ in the other point
$(e_1,0)$. The result follows since colinearity of three points
$A,\, B,\, C$ on the elliptic curve $E_3$ means $A+B+C=0$ in that
abelian group, while the opposite of $A=(X,Y)$ is $-A=(X,-Y)$.
\end{proof}

\begin{figure}
\begin{center}

\includegraphics[scale=0.9]{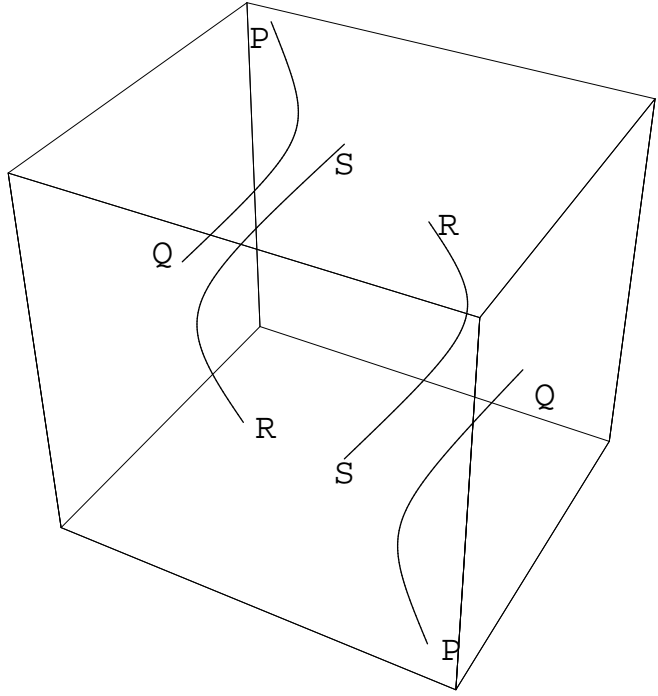}
\end{center}

\caption{\label{domain4} The flow lines }
\end{figure}

The above proposition shows that one should give a direct definition
of $E_3(\varphi)$ depending only upon the $J_{\ell m} (\varphi)$
rather than the $s_j(\varphi)$. There is also a strong reason to
define in a natural manner the complex points of $E_3(\varphi)$.
Indeed the translation $\sigma$ corresponds to points where the
coordinate $X=0$ and the corresponding equation for $Y$ is $Y^2=-1$
which does not admit real solutions.

\medskip
\subsection{Complex moduli space $\calm_{\mathbb C}$}\label{Cmoduli}

\noindent
\medskip

We shall show now that there is a natural way of
 extending the moduli space from the real to the
complex domain. The $ E_3(\varphi)$ then appear as  a net of degree
$4$ elliptic curves in $P_3(\mathbb C)$ having  $8$ points in
common. These elliptic curves will turn out to play a fundamental
role and to be  closely related to the elliptic curves of the
geometric data of the quadratic algebras which their elements label.

\smallskip

  To extend the moduli space to the
complex domain we start with the relations defining  the involutive
algebra $ C_{\mathrm{alg}}(S^3(\Lambda))$ and take for $\Lambda$ the
diagonal matrix with
\begin{equation}
\Lambda^\mu_\mu \, := \, u_\mu^{-1}
\end{equation}
where $(u_0,u_1,u_2,u_3)$ are the coordinates of
 $\ug\in  P_3(\mathbb C)$.
Using $y_{\mu}:=\Lambda^\mu_\nu z^\nu $ one obtains the homogeneous
defining relations in the form,
\begin{eqnarray}\label{uyalgebra}
u_k\, y_k\, y_0-u_0\, y_0\, y_k+u_\ell\,  y_\ell\,  y_m-u_m\, y_m\, y_\ell & = & 0   \nonumber\\
u_k\, y_0\, y_k-u_0\, y_k\, y_0+u_m\,  y_\ell\,  y_m-u_\ell\,  y_m\,
y_\ell & = & 0\label{homc}
\end{eqnarray}
for any cyclic permutation $(k,\ell,m)$ of (1,2,3). The
inhomogeneous relation becomes,
\begin{equation}\label{9.6}
\sum \, u_\mu\, y_\mu^2 =\, 1
\end{equation}
and the corresponding algebra $ C_{\mathrm{alg}}(S_{\mathbb
C}^3(\ug))$ only depends upon the class of $\ug\in  P_3(\mathbb C)$
(see the remark at the end of \ref{thetafunctions}). We let $
C_{\mathrm{alg}}(\mathbb C^4(\ug))$ be the quadratic algebra defined
by the six relations (\ref{homc}).

Taking $u_\mu\,= e^{2 i\,\varphi_\mu}$, $\varphi_0=0$,
 for all $\mu$ and $x^\mu:=e^{ i\,\varphi_\mu}y_\mu$
we obtain the defining relations of $C_{\rm alg}({\mathbb
R}_{\varphi}^4)$
 (except for ${x^\mu}^{\ast}=x^\mu$
which allows to pass from $\mathbb C^4(\ug)$ to ${\mathbb
R}_{\varphi}^4$). Thus the torus ${\mathbb T}_A$ sits naturally in
$P_3(\mathbb C)$ as
\begin{equation}\label{torusT}
{\mathbb T}_A =\, \{\ug\in  P_3(\mathbb C)\,|\;\,|u_\mu|=|u_\nu|
\qqq \mu,\,\nu\}
\end{equation}
In terms of homogeneous parameters the functions $J_{\ell
m}(\varphi)$ read as
\begin{equation}
J_{\ell m}(\varphi)=\tan (\varphi_0-\varphi_k)\tan
(\varphi_\ell-\varphi_m)
\end{equation}
for any cyclic permutation $(k,\ell, m)$ of (1,2,3), and extend to
the complex domain $\ug\in P_3(\mathbb C)$ as,
\begin{equation}\label{jlm}
J_{\ell
m}(\ug)=\frac{(u_0+u_\ell)(u_m+u_k)-(u_0+u_m)(u_k+u_\ell)}{(u_0+u_k)(u_\ell+u_m)}
\end{equation}

It follows easily from the  argument of proposition \ref{sklyaprep}
that for generic values of $\ug\in P_3(\mathbb C)$ the quadratic
algebra $ C_{\mathrm{alg}}(\mathbb C^4(\ug))$ only depends upon
$J_{k\ell}(\ug)$. We thus define the {\em complex} fiber as
\begin{equation}
F(\ug):= \{\vg\in P_3(\mathbb C)\, \vert\, J_{k\ell}(\vg)=\,
 J_{k\ell}(\ug)\}  \label{def}
\end{equation}
Let then,
\begin{equation}\label{abc}
\Phi(\ug)=( a,b,
c)=\left\{(u_0+u_1)(u_2+u_3),(u_0+u_2)(u_3+u_1),(u_0+u_3)(u_1+u_2)\right\}
\end{equation}
be the three roots of the Lagrange resolvent of the 4th degree
equation $\prod (x-u_j)=0$. We view $\Phi$ as a map
\begin{equation}
\Phi:P_3(\mathbb C)\backslash S \rightarrow P_2(\mathbb C)
\label{eq6.5}
\end{equation}
where $S$ is the following set
 of 8 points
\begin{eqnarray}
p_0=(1,0,0,0),\, \, p_1=(0,1,0,0),\, \, p_2=(0,0,1,0),\, \,
p_3=(0,0,0,1)\quad \quad
\label{eq6.3}\\
q_0=(-1,1,1,1),\, \, q_1=(1,-1,1,1),\, \, q_2=(1,1,-1,1),\, \,
q_3=(1,1,1,-1) \nonumber
\end{eqnarray}

The points $q_j$ belong to  the torus ${\mathbb T}_A$ of
\eqref{torusT}, they correspond to the orbit $W(P)=(PQRS)$.

  We extend the generic definition (\ref{def}) to arbitray
 $\ug\in P_3(\mathbb C)\backslash S$ and define  $F_\ug$
in general as the union of $S$ with the fiber $F(\ug)$ of $\Phi$
through $\ug$. It can be understood geometrically as follows.

  Let $\mathcal N$ be the net of quadrics in $P_3(\mathbb C)$ which contain $S$. Given $\ug\in P_3(\mathbb C)\backslash S$ the elements of $\mathcal N$ which contain $\ug$ form a pencil of quadrics with base locus
\begin{equation}
\cap \{Q\, \vert \, Q\in \mathcal N,\, \ug\in Q\}=Y_{\ug}
\label{eq6.6}
\end{equation}
which is an elliptic curve of degree 4 containing $S$ and $\ug$. One
has
\begin{equation}
Y_{\ug} =F_{\ug} \label{eq6.7}
\end{equation}
\begin{figure}
\begin{center}
\includegraphics[scale=0.9]{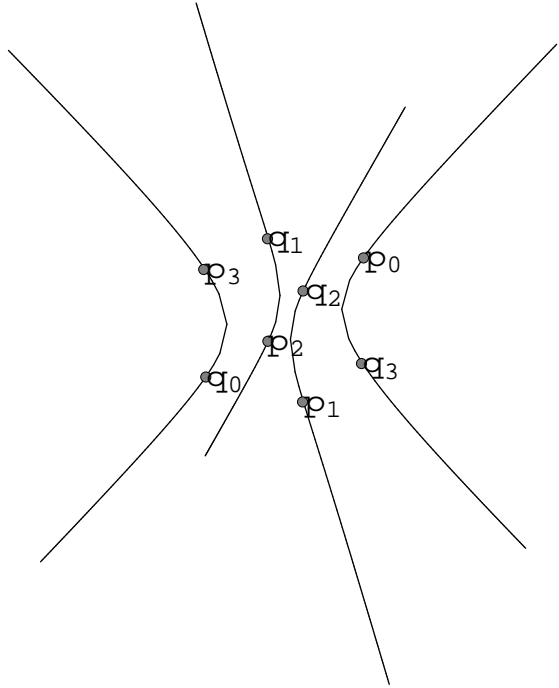}
\end{center}
\caption{\label{ellcurve10} The Elliptic Curve $F_\ug \cap
P_3(\mathbb {R})$  }
\end{figure}

\bigskip
With $\Phi(\ug)=( s_1,s_2, s_3)$ an explicit isomorphism of the
elliptic curve $F_\ug$ with the elliptic curve
$$
E_3=\,\{(X,Y)\,|\,Y^2=\,\prod(X s_j -1)\}\,
$$
of \eqref{ellcurve3} is given by
\begin{equation}\label{isoel}
(X,Y)\to \ug\,,\quad (u_k -u_0)(X\,s_k-1)-\,i\,(u_0+u_k)\,Y\,=\,0
\qqq k\,.
\end{equation}
Under this isomorphism the point at infinity in $E_3$ maps to $q_0
\in F_\ug$.

\medskip
\subsection{Notations for $\vartheta$-functions}\label{thetafunctions}

\noindent
\medskip

Let us fix our notations for elliptic $\vartheta$-functions. We fix
$\tau \in {\mathbb H}$ a complex number of stricly positive
imaginary part and let $q=\,e^{\pi i \tau}$ so that $|q|<1$. The
basic $\vartheta$-function is
\begin{equation}
\vartheta_3(z)=\,\sum_{{\mathbb Z}}\, q^{n^2}\, e^{2\pi i n z}
\label{thetafun3}
\end{equation}
It is periodic in $z$ with period $1$ and is (up to scale) the only
holomorphic section on $E=\,{\mathbb C}/L$, $L=\,{\mathbb Z}+\tau
{\mathbb Z}$ of the line bundle associated to the periodicity
conditions
$$
\xi(z+1)=\,\xi(z)\,,\quad \xi(z+\tau)=\,q^{-1}\,e^{-2\pi i
z}\,\xi(z)
$$
In particular it is equal (up to scale) to the infinite product
$$
\prod_1^\infty\,(1+q^{2n-1}\,e^{2\pi i  z})(1+q^{2n-1}\,e^{-2\pi i
z})\,,
$$
and only vanishes at $\omega_2=\frac{1}{2}+\frac{\tau}{2}$ modulo
$L$.

The three other $\vartheta$-functions are deduced from
$\vartheta_3(z)$ by the translations of order two of $E$, more
precisely one lets,
\begin{equation}
i \vartheta_1(z)=\,q^{\frac{1}{4}}\, e^{\pi i  z}\,\vartheta_3(z+
\omega_2)\,,\quad \vartheta_2(z)=\,
\vartheta_1(z+\frac{1}{2})\,,\quad \vartheta_4(z)=\,
\vartheta_3(z+\frac{1}{2})\,. \label{thetafun}
\end{equation}
They all are holomorphic sections on $E=\,{\mathbb C}/L$,
$L=\,{\mathbb Z}+\tau {\mathbb Z}$ of the line bundles associated to
the periodicity conditions of the form
$$
\xi(z+1)=\,\pm\,\xi(z)\,,\quad \xi(z+\tau)=\,\pm\, q^{-1}\,e^{-2\pi
i  z}\,\xi(z)
$$
and it follows that the linear span of their squares
$\vartheta^2_j(z)$ is a complex vector space of dimension two since
all are holomorpic sections of a line bundle of degree two. It
follows in particular that given any two theta functions the square
of any other is a linear combination of the squares of the first
two. The relevant coefficients are easy to compute from the values
$\vartheta^2_j(0)$ and if one lets
$$
k=\, \frac{\vartheta^2_2(0)}{\vartheta^2_3(0)}\,,\quad k'=\,
\frac{\vartheta^2_4(0)}{\vartheta^2_3(0)}\,,
$$
one gets
\begin{equation}
\vartheta^2_4(z)=\,k\,\vartheta^2_1(z)+\,k'\,\vartheta^2_3(z)\,,\quad
\vartheta^2_2(z)=-\,k'\,\vartheta^2_1(z)+\,k\,\vartheta^2_3(z)\,,\quad
\vartheta^2_1(z)=\,k\,\vartheta^2_4(z)-\,k'\,\vartheta^2_2(z) \,.
\label{thetarel}
\end{equation}

The $\lambda$-function of jacobi is given by $\lambda(\tau)=\,k^2$
and the $\wp$-function of Weierstrass by
\begin{equation}
\wp(z)=\,\alpha\,\frac{\vartheta^2_4(z)}{\vartheta^2_1(z)}+\,\beta\,,
\,. \label{ptheta}
\end{equation}
up to irrelevant normalization constants $\alpha$ and $\beta$. The
main identities we shall use for $\vartheta$-functions are the
sixteen theta relations recalled in  Appendix $3$.

\medskip
\subsection{ Fiber $ =$ characteristic variety, and the
birational automorphism $\sigma$ of ${\mathbb P}_3(\mathbb
C)$}\label{thetaparametrization}

\noindent
\medskip

  We shall now give, for generic values of $( a,b,c)$ a parametrization of $F_\ug$ by $\vartheta$-functions. We start with the equations for $F_\ug$
\begin{equation}
\frac{(u_0+u_1)(u_2+u_3)}{
a}=\frac{(u_0+u_2)(u_3+u_1)}{b}=\frac{(u_0+u_3)(u_1+u_2)}{c}
\label{eq6.8}
\end{equation}
and we diagonalize the above quadratic forms as follows
\begin{equation}
\begin{array}{lll}
(u_0+u_1)\,(u_2+u_3)& = & Z_0^2-Z_1^2\\
(u_0+u_2)\,(u_3+u_1) & = & Z_0^2-Z_2^2\\
 (u_0+u_3)\,(u_1+u_2)& = &Z_0^2-Z_3^2
\end{array}
\label{eq6.9}
\end{equation}
where
\begin{equation}
(Z_0,Z_1,Z_2,Z_3) = M.u \label{eq6.10}
\end{equation}
where $M$ is the involution,
\begin{equation}
M:= \frac{1}{2}\left[
\begin{array}{cccc}
1& \;1& 1&1\\
1& 1 &-1 &-1\\
1&-1& 1&-1\\
1& -1& -1& 1
\end{array}
\right]
\end{equation}
In these terms the equations for $F_\ug$ read
\begin{equation}
\frac{Z_0^2-Z_1^2}{ a}=\frac{Z_0^2-Z_2^2}{b}=\frac{Z_0^2-Z_3^2}{c}
\label{eq6.11}
\end{equation}

Let now $\tau\in \mathbb C$, $\im \,\tau >0$ and $\eta\in \mathbb C$
be such that one has, modulo projective equivalence,
\begin{equation}\label{ellparam}
( a, b, c)\sim \left(\frac{\vartheta_2(0)^2}{\vartheta_2(\eta)^2},\,
\, \frac{\vartheta_3(0)^2}{\vartheta_3(\eta)^2},\, \,
\frac{\vartheta_4(0)^2}{\vartheta_4(\eta)^2}\right)
\end{equation}
where $\vartheta_1,\vartheta_2, \vartheta_3, \vartheta_4$ are the
theta functions associated as above to $\tau$.

\smallskip
More precisely let $( a, b, c)$ be distinct non-zero complex numbers
and
\begin{equation}
\lambda=\,\frac{a}{b}\,\frac{c-b}{c-a}\,,\quad p=\,\frac{a}{a-c}
\label{ellparam1}
\end{equation}
and let $\tau\in \mathbb C$, $\im \,\tau >0$ such that
$\lambda(\tau)=\lambda$ where $\lambda$ is Jacobi's $\lambda$
function (\cf subsection \ref{basicell}). Let then $\eta\in \mathbb
C$ be such that with $\omega_1=\frac{1}{2}$,
$\omega_3=\frac{\tau}{2}$, one has
$$
\frac{\wp(\eta)-\wp(\omega_3)} {\wp(\omega_1)-\wp(\omega_3)}=\,p\,,
$$
where $\wp$ is the Weierstrass $\wp$-function  for the lattice
$L=\mathbb Z +\mathbb Z\tau\subset \mathbb C$. This last equality
does determine $\eta$ only up to sign and modulo the lattice $L$ but
this ambiguity does not affect the validity of \eqref{ellparam}
which one checks using the basic properties of $\vartheta$-functions
recalled in subsection \ref{thetafunctions}. Indeed in these terms
one has $\lambda(\tau)=\,k^2$ which gives the first equality in
\eqref{ellparam} with $( a, b, c)$ replaced by the right hand side
of \eqref{ellparam}. Using \eqref{ptheta} one gets the second
equality since for any $z$ one has
$$
\frac{\wp(z)-\wp(\omega_3)}
{\wp(\omega_1)-\wp(\omega_3)}=\,k\,\frac{\vartheta_4(z)^2}{\vartheta_1(z)^2}\,.
$$

Let then $\tau$ and $\eta$ be fixed by the above conditions, one
gets

\medskip
\begin{prop} \label{theta}
The following define isomorphisms of $\mathbb C/L$ with $F_\ug$,
\[
\varphi(z)=\left( \frac{\vartheta_1(2z)}{\vartheta_1(\eta)},\,
\frac{\vartheta_2(2z)}{\vartheta_2(\eta)},\,
\frac{\vartheta_3(2z)}{\vartheta_3(\eta)},\,
\frac{\vartheta_4(2z)}{\vartheta_4(\eta)}\right) = (Z_0,Z_1,Z_2,Z_3)
\]
and $\psi(z)=\varphi(z-\eta/2)$.
\end{prop}
 \begin{proof} Up to an affine transformation, $\varphi$ (and $\psi$ are) is the classical projective embedding of $\mathbb C/L$ in $P_3(\mathbb C)$. Thus we only need to check that the biquadratic curve $\im\,\varphi=\im\,\psi$ is given by (\ref{eq6.11}). It is thus enough to check (\ref{eq6.11}) on $\varphi(z)$. This follows from the basic relations \eqref{thetarel}
 written as
\begin{equation}\label{theta01}
\vartheta^2_3(z)\vartheta^2_2(0)=\vartheta^2_2(z)\vartheta^2_3(0)+\vartheta^2_1(z)\vartheta^2_4(0)
\end{equation}
and
\begin{equation}\label{theta02}
\vartheta^2_4(z)\vartheta^2_3(0)=\vartheta^2_1(z)\vartheta^2_2(0)+\vartheta^2_3(z)\vartheta^2_4(0)
\end{equation}
which one uses to check $\frac{Z_0^2-Z_1^2}{
a}=\frac{Z_0^2-Z_2^2}{b}$ and
$\frac{Z_0^2-Z_2^2}{b}=\frac{Z_0^2-Z_3^2}{c}$ respectively.

Let us check the first one using \eqref{theta01} to replace all
occurences of $\vartheta^2_3(2 z)$ and $\vartheta^2_3(\eta)$ by the
value given by \eqref{theta01}. One gets
$$
\frac{Z_0^2-Z_2^2}{b}=\,\frac{\vartheta^2_1(2z)}{\vartheta^2_1(\eta)}
\,\frac{\vartheta^2_3(\eta)}{\vartheta^2_3(0)}-\,\frac{\vartheta^2_3(
2 z)}{\vartheta^2_3(0)}
=\,\frac{\vartheta^2_1(2z)}{\vartheta^2_1(\eta)}\,\frac{\vartheta^2_2(\eta)}{\vartheta^2_2(0)}
+\,\frac{\vartheta^2_1(2z)\vartheta^2_4(0)}{\vartheta^2_2(0)\vartheta^2_3(0)}-(\,\frac{\vartheta^2_2(2z)}{\vartheta^2_2(0)}+\,\frac{\vartheta^2_1(2z)\vartheta^2_4(0)}{\vartheta^2_2(0)\vartheta^2_3(0)})
$$
$$
=
\,\frac{\vartheta^2_1(2z)}{\vartheta^2_1(\eta)}\,\frac{\vartheta^2_2(\eta)}{\vartheta^2_2(0)}
-\,\frac{\vartheta^2_2(2z)}{\vartheta^2_2(0)}=\,\frac{Z_0^2-Z_1^2}{a}
$$

Let us check the second one using \eqref{theta02} to replace all
occurences of $\vartheta^2_4(2 z)$ and $\vartheta^2_4(\eta)$ by the
value given by \eqref{theta02}. One gets
$$
\frac{Z_0^2-Z_3^2}{c}=\,\frac{\vartheta^2_1(2z)}{\vartheta^2_1(\eta)}
\,\frac{\vartheta^2_4(\eta)}{\vartheta^2_4(0)}-\,\frac{\vartheta^2_4(
2 z)}{\vartheta^2_4(0)}
=\,\frac{\vartheta^2_1(2z)}{\vartheta^2_1(\eta)}\,\frac{\vartheta^2_3(\eta)}{\vartheta^2_3(0)}
+\,\frac{\vartheta^2_1(2z)\vartheta^2_2(0)}{\vartheta^2_3(0)\vartheta^2_4(0)}-(\,\frac{\vartheta^2_3(2z)}{\vartheta^2_3(0)}+\,\frac{\vartheta^2_1(2z)\vartheta^2_2(0)}{\vartheta^2_3(0)\vartheta^2_4(0)})
$$
$$
\,\frac{\vartheta^2_1(2z)}{\vartheta^2_1(\eta)}\,\frac{\vartheta^2_3(\eta)}{\vartheta^2_3(0)}
-\,\frac{\vartheta^2_3(2z)}{\vartheta^2_3(0)}
=\,\frac{Z_0^2-Z_2^2}{b}
$$
 \end{proof}

  The elements of $S$ are obtained from the following values of $z$
\begin{equation}
\psi(\eta)=p_0,\, \psi(\eta+\frac{1}{2})=p_1,\,
\psi(\eta+\frac{1}{2}+\frac{\tau}{2})=p_2,\,
\psi(\eta+\frac{\tau}{2})=p_3 \label{eq6.15}
\end{equation}
and
\begin{equation}
\psi(0)=q_0,\, \psi(\frac{1}{2})=q_1,\,
\psi(\frac{1}{2}+\frac{\tau}{2})=q_2,\, \psi(\frac{\tau}{2})=q_3.
\label{eq6.16}
\end{equation}
(We used $M^{-1}.\psi$ to go back to the coordinates $u_\mu$
but note that $M^{-1}=M$ and $M(q_j)=q_j$).\\
Let $H\sim \mathbb Z_2\times \mathbb Z_2$ be the Klein subgroup of the symmetric group $\fracS_4$ acting on $P_3(\mathbb C)$ by permutation of the coordinates $(u_0,u_1,u_2,u_3)$.\\
For $\rho$ in $H$ one has $\Phi\circ \rho=\Phi$, so that $\rho$
defines for each $\ug$ an automorphism of $F_\ug$. For $\rho$ in $H$
the matrix $M\rho M^{-1}$ is diagonal with $\pm 1$ on the diagonal
and the quasiperiodicity of the $\vartheta$-functions  allows to
check that these automorphisms are translations on $F_\ug$ by the
following 2-torsion elements of $\mathbb C/L$,
\begin{equation}
\begin{array}{ll}
\rho =\left(
\begin{array}{cccc}
0 & 1 & 2 & 3\\
1 & 0 & 3 & 2
\end{array}\right) & \mbox{is translation by}\, \omega_1=\frac{1}{2},\\
\\
\rho =\left(
\begin{array}{cccc}
0 & 1 & 2 & 3\\
2 & 3 & 0 & 1
\end{array}\right) & \mbox{is translation by}\, \omega_2=\frac{1}{2}+\frac{\tau}{2}
\end{array}
\label{eq6.17}
\end{equation}
Let $\calo\subset P_3(\mathbb C)$ be the complement of the $4$
hyperplanes $\{u_\mu\,=0\}$ with $ \mu\in \{0,1,2,3\}$. Then
$(u_0,u_1,u_2,u_3)\mapsto (u_0^{-1},u_1^{-1},u_2^{-1},u_3^{-1})$
defines an involutive automorphism $I$ of $\calo$ and since one has
\begin{equation}
(u_0^{-1}+u_k^{-1})(u_\ell^{-1}+u^{-1}_m)=(u_0u_1u_2u_3)^{-1}(u_0+u_k)(u_\ell+u_m)
\label{eq6.18}
\end{equation}
it follows that $\Phi\circ I=\Phi$, so that $I$ defines for each
$\ug\in \calo\backslash \{q_0,q_1,q_2,q_3\}$ an involutive
automorphism of $F_\ug\cap\calo$ which extends canonically to
$F_\ug$. Note in fact that, as a birational map $I$ continues to
make sense on the complement of the $6$ lines
$\ell_{\mu,\nu}=\{u\,|\,u_\mu\,=0,\,u_\nu\,=0\}$ for
 $ \mu, \nu \in \{0,1,2,3\}$ using
$(u_0,u_1,u_2,u_3)\mapsto ( u_1 u_2 u_3,u_0  u_2 u_3,u_0 u_1
u_3,u_0 u_1 u_2 )$.

\medskip
\begin{prop}\label{nice}
The restriction of $I$ to $F_\ug$ is the symmetry $\psi(z)\mapsto
\psi(-z)$ around any of the points $q_\mu \in F_\ug$ in the elliptic
curve $F_\ug$.
\end{prop}
This symmetry, as well as the above translations by two torsion
elements does not refer to a choice of origin in the curve $F_\ug$.
The proof follows from identities on theta functions but it can be
seen directly using the isomorphism $E_3 \simeq F_\ug$ of
\eqref{isoel}. Indeed the symmetry around $q_0\in F_\ug$ corresponds
to the transformation $(X,Y)\to (X,-Y)$ on $E_3$ and the isomorphism
\eqref{isoel} carries this back to $(u_0,u_1,u_2,u_3)\mapsto
(u_0^{-1},u_1^{-1},u_2^{-1},u_3^{-1})$ as one checks directly
dividing each of the equations \eqref{isoel} $(u_k
-u_0)(X\,s_k-1)-\,i\,(u_0+u_k)\,Y\,=\,0$ by $u_0 u_k$ to get the
same equation but with $-Y$ instead of $Y$ for the $u_j^{-1}$.

 The torus
${\mathbb T}_A$ of \eqref{torusT} gives a covering of the real
moduli space $\calm$. For $\ug \in {\mathbb T}_A$, the point
$\Phi(\ug)$ is real with projective coordinates
\begin{equation}\label{real0}
\Phi(\ug) = (s_1,s_2,s_3)\, , \quad s_k := 1 + t_\ell\, t_m \, ,
\quad t_k:= {\rm tan}(\varphi_k- \varphi_0)
\end{equation}
The corresponding fiber $F_{\ug}$ is stable under complex
conjugation $\vg \mapsto \overline{\vg}$ and the intersection of
$F_{\ug}$ with the real moduli space is given by,
\begin{equation}\label{real}
F_{\mathbb T}(\ug)=\, F_{\ug}\cap {\mathbb T}_A=\{\vg \in F_{\ug}
\vert I(\vg)\,=\, \overline{\vg}\}
\end{equation}
The curve $F_{\ug}$ is defined over $\mathbb R$ and (\ref{real})
determines its purely \underline{imaginary} points.
 Note that  $F_{\mathbb T}(\ug) $ (\ref{real})
is  invariant under the Klein group $H$ and thus has two connected
components, we let $F_{\mathbb T}(\ug)^0 $ be the component
containing $q_0$. The real points, $ \{\vg \in F_{\ug} \vert
\vg\,=\, \overline{\vg}\}=\,F_\ug \cap P_3(\mathbb {R}) $ of
$F_{\ug}$ do play a complementary role in the characteristic variety
as we shall see below.

 Our aim now is to show that for $\ug\in P_3(\mathbb C)$ generic, there is an astute  choice of generators of the quadratic algebra $\cala_\ug= C_{\mathrm{alg}}(\mathbb C^4(\ug))$  for which  the characteristic variety $E_\ug$ actually coincides with the fiber variety $F_\ug$ and to identify the corresponding automorphism $\sigma$. Since this coincidence  no longer holds for non-generic values it is a quite remarkable fact which we first noticed by comparing the $j$-invariants of these two elliptic curves.\\

\smallskip Let  $\ug\in P_3(\mathbb C)$ be generic,
we perform the following change of generators

\begin{equation} \label{change0}
\begin{array}{llll}
y_0 & = & \sqrt{u_1-u_2} \; \sqrt{u_2-u_3}  \;\sqrt{u_3-u_1} & Y_0\\
\\
y_1 & = & \sqrt{u_0+u_2} \; \sqrt{u_2-u_3}  \;\sqrt{u_0+u_3}&  Y_1\\
\\
y_2 & = & \sqrt{u_0+u_3}  \;\sqrt{u_3-u_1}  \;\sqrt{u_0+u_1}&  Y_2\\
\\
y_3 & = & \sqrt{u_0+u_1}  \;\sqrt{u_1-u_2}  \;\sqrt{u_0+u_2}&  Y_3
\end{array}
\end{equation}

\smallskip We let $J_{\ell m}$ be as before, given by (\ref{jlm})
\begin{equation}
J_{12}=\frac{a-b}{c},\, \, \, J_{23}=\frac{b-c}{a},\, \, \,
J_{31}=\frac{c -a}{b} \label{eq5.8}
\end{equation}
with $a, b, c$ given by (\ref{abc}). Finally let $e_\nu$ be the $4$
points of $P_3(\mathbb C)$
whose homogeneous coordinates ($Z_\mu$) all vanish but one. \\

\begin{thm} \label{iden}
\smallskip 1) In terms of the $Y_\mu$, the relations of $\cala_\ug$ take the form
\begin{eqnarray}
[Y_0,Y_k]_- & = & [Y_\ell,Y_m]_+\label{eq5.6}\\
{[Y_\ell,Y_m]}_- & = & - J_{\ell m} [Y_0,Y_k]_+\label{eq5.7}
\end{eqnarray}
for any $k\in \{1,2,3\}$, $(k,\ell,m)$ being a cyclic permutation of
(1,2,3)

\smallskip 2) The characteristic variety $E_\ug$ is the union of $F_\ug$ with the $4$ points
$e_\nu$.

\smallskip 3) The automorphism $\sigma$ of the characteristic variety $E_\ug$  is given by
\begin{equation}
\psi(z) \mapsto \psi(z- \eta) \label{eq5.10}
\end{equation}
on $F_\ug$ and $\sigma=\id$ on the 4 points $e_\nu$.

\smallskip 4) The automorphism $\sigma$ is the restriction to $F_\ug$ of a birational
automorphism of $P_3(\mathbb C)$ independent of $\ug$ and defined
over $\mathbb Q$.
\end{thm}

\smallskip The similarity between the above presentation and the Sklyanin one
(cf. \eqref{sklya1}, \eqref{sklya2}) is misleading, indeed for the
latter all the characteristic varieties are contained in the same
quadric (cf. \cite{smi-sta:1992} \S 2.4)
$$
\sum x_\mu^2= 0
$$
and cant of course form a net of essentially disjoint curves.

\smallskip {\bf Proof} $\;\;$ By construction $E_\ug=\{ Z \, \vert \, {\rm Rank}\;N(Z)<4 \}$ where

\medskip
\begin{equation}
N(Z)=\left(
\begin{array}{cccc}
Z_1 & -Z_0 & Z_3 & Z_2\\
\\
Z_2 & Z_3 & -Z_0 & Z_1\\
\\
Z_3 & Z_2 & Z_1 & -Z_0\\
\\
(b-c)Z_1 & (b-c)Z_0 & -a Z_3 & a Z_2\\
\\
(c - a) Z_2 & b Z_3 & (c-a)Z_0 & -b Z_1\\
\\
(a-b) Z_3 & -c Z_2 & c Z_1 & (a-b)Z_0
\end{array}
\right) \label{eq5.9}
\end{equation}

\bigskip

\smallskip One checks that it is the union of the fiber $F_\ug$ (in the generic case)
with the above $4$ points. In fact in terms of the original
presentation \eqref{uyalgebra} \ie in terms of the $y_j$ the
characteristic variety in the generic case is the intersection of
the two quadrics
\begin{equation}
\sum y_j^2=0\,,\quad \sum u_j^2\,y_j^2=0 \,, \label{uychar}
\end{equation}

and after the change of variables \eqref{change0} it just becomes
\begin{equation}
\frac{Z_0^2-Z_1^2}{ a}=\frac{Z_0^2-Z_2^2}{b}=\frac{Z_0^2-Z_3^2}{c}
\label{charZ}
\end{equation}

as can easily be checked since only the squares $Z_j^2$ are involved
(and linearly). Note that we already knew that the fiber $F_\ug$ is
abstractly isomorphic to the elliptic curve of the characteristic
variety using lemma \ref{trigtoell}, proposition \ref{weyltwo} and
the isomorphism \eqref{isoel}. But here we have shown that their
respective embeddings in ${\mathbb P}^3({\mathbb C})$ are the same
(\ie the corresponding line bundles are the same).

The automorphism $\sigma$ of the characteristic variety $E_\ug$ is
given by definition by the equation,
\begin{equation}
N(Z)\, \sigma (Z) = 0 \label{aut}
\end{equation}
where $\sigma(Z)$ is the column vector $\sigma(Z_\mu):= M
\cdot\sigma(\ug)$ (in the variables $Z_\lambda$). One checks that
$\sigma(Z)$ is already determined by the equations in (\ref{aut})
corresponding to the first three lines in $N(Z)$ which are
independent of $a, b, c$ (see below). Thus $\sigma$ is in fact an
automorphism of $P_3(\mathbb C)$ which is the identity on the above
four points and which restricts as automorphism of $F_\ug$ for each
$\ug$ generic. One checks that $\sigma$ is the product of two
involutions which both restrict to $E_\ug$ (for $\ug$ generic)
\begin{equation}\label{twoinvolutions}
\sigma = I\circ I_0
\end{equation}
where $I$ is the involution of proposition \ref{nice} corresponding
to $u_\mu \mapsto u^{-1}_\mu$ and where $I_0$ is given by
\begin{equation}
I_0(Z_0)=-Z_0, \, \, I_0(Z_k)=Z_k \label{firstinvolution}
\end{equation}
for $k\in \{1,2,3\}$ and which restricts obviously to $E_\ug$ in
view of (\ref{eq6.9}). Both $I$ and $I_0$ are the identity on the
above four points and since $I_0$ induces the symmetry $\varphi(z)
\mapsto \varphi(-z)$ around $\varphi(0)=\psi(\eta/2)$ (proposition
\ref{theta}) one gets the result using proposition \ref{nice}.

\smallskip  The fact that $\sigma$ does not depend on the parameters $a,b,c$ plays an  important role. Explicitly we get from the first 3 equations (\ref{aut})
\begin{equation}
\sigma(Z)_\mu =\eta_{\mu\mu} (Z^3_\mu-Z_\mu\sum_{\nu\not=\mu}
Z^2_\nu-2\prod_{\lambda\not=\mu} Z_\lambda) \label{algebraicsigma}
\end{equation}
for $\mu\in \{0,1,2,3\}$, where $\eta_{00}=1$ and $\eta_{nn}=-1$ for $n\in\{1,2,3\}$. $\square$\\

\noindent \underbar{Remark} : Notice that, since the $u_\mu$ are
homogeneous coordinates on $P_3(\mathbb C)$, the generators $y_\mu$
as well as the generators $Y_\mu$ are only defined modulo a non zero
multiplicative scalar. Indeed under a change $u_\mu\mapsto \nu^2
u_\mu$ of homogeneous coordinates ($\nu\in \mathbb C\backslash
\{0\}$), the $y_\mu$ transform as $y_\mu\mapsto \nu^{-1}y_\mu$ while
the $Y_\mu$ transform as $Y_\mu\mapsto \nu^{-4}Y_\mu$ which leaves
invariant  (\ref{9.6}) and (\ref{change0}) as well as the
$x^\mu$($(x^\mu)^2=u_\mu(y_\mu)^2$). Later on, it will be convenient
to choose a normalization for the generators, e.g. in (\ref{xtoZ}).

\bigskip

\section{The map from ${\mathbb T}^2_\eta \times [0,\tau]$ to $S_\varphi^3$ and the pairing}
\label{pairing} \setcounter{equation}{0}

This section contains the main technical result of the paper \ie
both the construction of the one-parameter family of
$*$-homomorphisms from the algebra of $S_\varphi^3$ (with $\varphi$
generic) to the algebra of the non-commutative torus ${\mathbb
T}^2_\eta$ and the computation of the pairing of the image of the
Hochschild $3$-cycle with the natural ``fundamental class" for the
product ${\mathbb T}^2_\eta \times [0,\tau]$. \\

\subsection{Central elements}

\noindent
\medskip

We already saw in Section 2  (Lemma \ref{centerC}) that the algebra
$C_{\rm alg} ({\mathbb R}_{\varphi}^4)$ of ${\mathbb R}_{\varphi}^4$
contains in its center the following element,
\begin{equation}
\label{center1} Q_1=C=\,\sum (x^{\mu})^2 \, .
\end{equation}
To get another one, one first looks at the Sklyanin algebra defined
by (\ref{sklya1}) and (\ref{sklya2}) whose center contains two
natural ``Casimir" elements $C_j$,
\begin{equation}
\label{center2} C_1=\,\sum S_{\mu}^2  \,,\quad C_2=\,\sum
\,j_k\,S_{k}^2 \, ,
\end{equation}
where the $j_k$ fulfill the relations
\begin{equation}
\label{littlej} -\,\frac{j_\ell-j_m}{j_k}=\,J_{\ell m}
\end{equation}
Let us check that $C_1$ commutes with $S_0$ and $S_{\ell}$. The idea
is to only use the relation (\ref{sklya1}) so that the $J_{\ell m}$
do not interfere with the computation. This means that one treats
the commutators as follows:
\begin{equation}
[S_0 , S_k^2] = (S_0 \, S_k + S_k \, S_0) \, S_k - S_k (S_0 \, S_k +
S_k \, S_0) = i[S_k \,[S_\ell\, S_m]]\, . \nonumber
\end{equation}
Thus the sum over $k$ gives $0$ in view of the Jacobi identy.
\begin{equation}
[S_1 , S_0^2] = (S_1 \, S_0 + S_0 \, S_1) \, S_0 - S_0 (S_1 \, S_0 +
S_0 \, S_1) = - \, i (S_2 \, S_3 - S_3 \, S_2) \, S_0 + i \, S_0
(S_2 \, S_3 - S_3 \, S_2) \, . \nonumber
\end{equation}
\begin{equation}
[S_1 , S_2^2] = [S_1 , S_2] \, S_2 + S_2 \, [S_1 , S_2] = i (S_0 \,
S_3 + S_3 \, S_0) \, S_2 + i \, S_2 (S_0 \, S_3 + S_3 \, S_0)
\nonumber
\end{equation}
\begin{equation}
[S_1 , S_3^2] = [S_1 , S_3] \, S_3 + S_3 \, [S_1 , S_3] = - \, i
(S_0 \, S_2 + S_2 \, S_0) \, S_3 - i \, S_3 (S_0 \, S_2 + S_2 \,
S_0) \, . \nonumber
\end{equation}
One checks that the sum of these terms gives 0.  Using cyclic permutations the commutation with $S_k$ easily follows. \\

\noindent \underbar{Remark} : It is worth noticing here that the
relation $[C_1,S_\nu]=0$ can be written in the form
\cite{ac-mdv:2004}
\[
\sum_\mu[S_\mu,[S_\mu, S_\nu]_+]=0
\]
where it becomes apparent that it is a super Lie algebra version of
the relation defining the Yang-Mills algebra studied in
\cite{ac-mdv:2002b}. As pointed out in \cite{ac-mdv:2004} the
relation (\ref{sklya1})
 is the corresponding super analog of the self-duality relation and the fact that it implies $[C_1,S_\nu]=0$ is the content of Lemma 1 in \cite{ac-mdv:2004}.\\

\medskip

Using (\ref{rescale1}) we can then assert that in the generic case
the following element is in the center,
\begin{equation}
\label{eq2.3} (\prod \sin \,\varphi_k) \, (x^0)^2 + \sum_1^3 \cos
(\varphi_k - \varphi_{\ell}) \cos (\varphi_k - \varphi_m) \sin
\varphi_k \, (x^k)^2 \, .
\end{equation}
Substracting (\ref{center1}) multiplied by $\prod \sin \varphi_k$
and using $$\cos (\varphi_k - \varphi_{\ell})\, \cos (\varphi_k -
\varphi_m) - \sin \varphi_{\ell}\, \sin \varphi_m = \cos \varphi_k
\,\cos (- \, \varphi_k + \varphi_{\ell} + \varphi_m)$$ we then get:
\begin{equation}
\label{center22} Q_2 = \frac{1}{2}\,\sum_1^3 \sin \, 2 \varphi_k
\;\cos (- \, \varphi_k + \varphi_{\ell} + \varphi_m) \, (x^k)^2
\end{equation}
as a central element.

\medskip

  Note that to get that (\ref{center2}) is central we did not use relation (\ref{sklya2}) and thus it holds irrespective of the finiteness of the $J_{\ell m}$. Thus the case $\delta (\varphi) \ne 0$ is entirely covered to show that (\ref{center22}) is central.

\bigskip

  \begin{prop}
\begin{enumerate}
\item[1)]
Both $Q_j$ are in the center of $C_{\rm alg} ({\mathbb
R}_{\varphi}^4)$ for all values of $\varphi$.
\item[2)] Let $S_\mu=\lambda_\mu \,x^\mu$ as in (\ref{rescale1})
and $\lambda$ such that $\lambda \,j_k=(-s_k+s_\ell+s_m)/s_\ell s_m$
then
\begin{equation}
\label{relcenter} Q_1 = \frac{C_1-\lambda \,C_2}{\prod \sin
\,\varphi_k}\,,\quad Q_2=\,\lambda C_2\,.
\end{equation}

\end{enumerate}

\end{prop}

\medskip

  \begin{proof} 1) We just have to check for $Q$ and this will be done replacing it by (\ref{eq2.3}) and using instead of (\ref{sklya1}) the three relations (\ref{pres1}),
\begin{equation}
\sin \varphi_k \, [x^0 , x^k]_+ = i \cos (\varphi_{\ell} -
\varphi_m) \, [x^{\ell} , x^m] \, .
\end{equation}
We just have to repeat the same proof as for the $S_{\mu}$'s making
sure that any $[x^0 , x^k]_+$ has a $\sin \varphi_k$ as coefficient
and every $[x^{\ell} , x^m]$ a $\cos (\varphi_{\ell} - \varphi_m)$.
For the commutator with $x^0$ this follows from the term $\sin
\varphi_k \, (x^k)^2$ in (\ref{eq2.3}). For the commutator with
$x^1$ this follows from the terms $\sin \varphi_1 \, (x^0)^2$ and
$\cos (\varphi_1 - \varphi_k) \, (x^k)^2$, $k \ne 1$.

2) Note that the existence of $\lambda$ follows from lemma
\ref{propor} \ie with $\tilde s_k= (-s_k+s_\ell+s_m)/s_\ell s_m$ the
equality
$$
-\,\frac{{\tilde s}_\ell-{\tilde s}_m}{{\tilde s}_k}=\,J_{\ell m}
$$
We have already shown that $C_1=\,{\prod \sin
\,\varphi_k}\,Q_1+\,Q_2$. It remains to check that $Q_2=\,\lambda
\,C_2$. Since $\lambda \,j_k={\tilde s}_k$ this amounts to
$\lambda_k^2 \,{\tilde s}_k =\,\frac{1}{2}\, \sin \, 2 \varphi_k
\;\cos (- \, \varphi_k + \varphi_{\ell} + \varphi_m)$ \ie
$$
{\tilde s}_k=\,\frac{\cos\,\varphi_k \;\cos (- \, \varphi_k +
\varphi_{\ell} + \varphi_m)}{\cos (\varphi_k - \varphi_{\ell})\,
\cos (\varphi_k - \varphi_m)}
$$
which follows from the definition of $\tilde s_k=
(-s_k+s_\ell+s_m)/s_\ell s_m$.
\end{proof}

\medskip
In terms of the presentation    of theorem \ref{iden} one gets,

 \begin{prop} \label{centerY} Let $\cala_\ug= C_{\mathrm{alg}}(\mathbb C^4(\ug))$
 at generic $\ug$, then
 the following three linearly dependent quadratic elements
belong to the center of $\cala_\ug$,
\begin{equation} \label{qk}
Q_m\,=\,J_{k \ell}\,(Y_0^2 \, +\, Y_m^2)\,+ \, Y_k^2 \, -
\,Y_\ell^2\,.
\end{equation}
\end{prop}

\medskip
\subsection{The Hochschild cycle ${\rm ch}_{3/2} (U)$}

\noindent
\medskip

With $z^k = e^{i \varphi_k} \, x^k$ we have, with $\varphi_0 = 0$,
\begin{equation}
\label{eq3.1} U = \sum \tau_{\mu} \, z^{\mu} \, , \qquad \tau_0 = 1
\, , \quad \tau_k = i \, \sigma_k \, .
\end{equation}
We need to compute
\begin{equation}
\label{cherndef} {\rm ch}_{3/2} (U) = U_{i_0 i_1} \otimes U_{i_1
i_2}^* \otimes U_{i_2 i_3} \otimes U_{i_3 i_3}^* - U_{i_0 i_1}^*
\otimes U_{i_1 i_2} \otimes U_{i_2 i_3}^* \otimes U_{i_3 i_0} \, .
\end{equation}
The trace computation is given by:

\bigskip

  \begin{lem}  One has
\begin{equation}
\frac{1}{2} \, {\rm Trace} \, (\tau_{\alpha} \, \tau_{\beta}^* \,
\tau_{\gamma} \, \tau_{\delta}^*) = \delta_{\alpha\beta} \,
\delta_{\gamma\delta} + \delta_{\beta\gamma} \,
\delta_{\delta\alpha} - \delta_{\alpha\gamma} \,
\delta_{\beta\delta} + \varepsilon_{\alpha \beta \gamma \delta} \, .
\nonumber
\end{equation}
\end{lem}

\medskip

\begin{proof} If the set $\{ \alpha , \beta , \gamma , \delta \}$ is $\{ 0,1,2,3 \}$ we can assume $\alpha = 0$ or $\beta = 0$ by symmetry in $(\alpha , \beta , \gamma , \delta) \to (\gamma , \delta , \alpha , \beta)$. For $\alpha = 0$ we get $\frac{1}{2} \, {\rm Trace} \, (\tau_{\beta} \, \tau_{\gamma} \, \tau_{\delta})$ since the two $-$ signs from $\tau^*$ cancell. This is cyclic and antisymmetric and gives for $(1,2,3)$ using $\sigma_1 \, \sigma_2 = i \, \sigma_3$ the result $i^4 \times \frac{1}{2} \times 2 = 1$. For $\beta = 0$ we get
\begin{equation}
\frac{1}{2} \, {\rm Trace} \, (\tau_{\alpha} \, \tau_{\gamma} \,
\tau_{\delta}^* ) = -\frac{1}{2} \, {\rm Trace} \, (\tau_{\alpha} \,
\tau_{\gamma} \, \tau_{\delta}) = - \, \varepsilon_{0 \alpha \gamma
\delta} = \varepsilon_{\alpha \beta \gamma \delta} \, . \nonumber
\end{equation}
If two of the elements $\alpha , \beta , \gamma , \delta$ are equal
and the two others are different we get $0$. Thus there are 3 cases
$\alpha = \beta$, $\alpha = \gamma$, $\alpha = \delta$. One has
$\tau_{\alpha} \, \tau_{\alpha}^* = 1$, thus if $\alpha = \beta$ we
get 1. For $\alpha = \gamma$ we get $\frac{1}{2} \, {\rm Trace} \,
(\tau_{\alpha} \, \tau_{\beta}^* \, \tau_{\alpha} \,
\tau_{\beta}^*)$. The two $-$ signs cancell and give $\frac{1}{2} \,
{\rm Trace} \, (\tau_{\alpha} \, \tau_{\beta} \, \tau_{\alpha} \,
\tau_{\beta})$. We can assume $\alpha \ne \beta$. If $0 \in \{
\alpha , \beta \}$ we get $\frac{1}{2} \, {\rm Trace} \, (\tau_k^2)
= -1$. If $0 \notin \{ \alpha , \beta \}$, $\tau_{\beta} \,
\tau_{\alpha} = - \, \tau_{\alpha} \, \tau_{\beta}$ and we get again
$-1$. The case $\alpha = \delta$, $\beta = \gamma$ is as $\alpha =
\beta$, $\gamma = \delta$. Finally if all indices are equal we get
1.
\end{proof}

\bigskip

  \begin{prop} \label{chbare}
The Hochschild cycle
$${\rm ch}_{3/2} (U) \in HZ_3(C^\infty(S_\varphi^3))$$
is given (using $\varphi_0 = 0$ and up to a scalar factor) by
$$
{\rm ch}_{3/2} (U) = \sum \varepsilon_{\alpha \beta \gamma \delta}
\,\cos (\varphi_{\alpha} - \varphi_{\beta} + \varphi_{\gamma} -
\varphi_{\delta}) \, x^{\alpha} \otimes x^{\beta} \otimes x^{\gamma}
\otimes x^{\delta} $$
$$
- i \sum_{\mu , \nu} \sin 2 (\varphi_{\mu} - \varphi_{\nu}) \,
x^{\mu} \otimes x^{\nu} \otimes x^{\mu} \otimes x^{\nu} \, .
$$
\end{prop}
\bigskip

  \begin{proof} The coefficient of $x^{\alpha} \otimes x^{\beta} \otimes x^{\gamma} \otimes x^{\delta}$
in $\frac{1}{2} {\rm ch}_{3/2} (U)$
  is given by $\left( \times \, \frac{1}{2} \right)$
\begin{equation}
\label{eq3.3} \frac{1}{2} \, {\rm Trace} \, (\tau_{\alpha} \,
\tau_{\beta}^* \, \tau_{\gamma} \, \tau_{\delta}^*) \, e^{i
(\varphi_{\alpha} - \varphi_{\beta} + \varphi_{\gamma} -
\varphi_{\delta})} - \frac{1}{2} \, {\rm Trace} \, (\tau_{\alpha}^*
\, \tau_{\beta} \, \tau_{\gamma}^* \, \tau_{\delta}) \, e^{i
(-\varphi_{\alpha} + \varphi_{\beta} - \varphi_{\gamma} +
\varphi_{\delta})} \, .
\end{equation}
It is non zero only in the two cases (with cardinality denoted as
$\#$),
$$\# \, \{ \alpha , \beta , \gamma , \delta \} = 4\,,\quad
\# \, \{  \alpha , \beta , \gamma , \delta \} \leq 2$$ In the first
case we get as coefficient of $x^{\alpha} \otimes x^{\beta} \otimes
x^{\gamma} \otimes x^{\delta}$ the term
$$\varepsilon_{\alpha \beta \gamma \delta} \, e^{i (\varphi_{\alpha} - \varphi_{\beta} + \varphi_{\gamma} - \varphi_{\delta})} - \varepsilon_{\beta \gamma \delta \alpha} \, e^{-i (\varphi_{\alpha} - \varphi_{\beta} + \varphi_{\gamma} - \varphi_{\delta})} = 2
\varepsilon_{\alpha \beta \gamma \delta}  \,\cos (\varphi_{\alpha} -
\varphi_{\beta} + \varphi_{\gamma} - \varphi_{\delta})$$ since the
cyclic permutation has signature $-1$.

In the second case the terms $\delta_{\alpha\beta} \, \delta_{\gamma
\delta} \, e^{i (\varphi_{\alpha} - \varphi_{\beta} +
\varphi_{\gamma} - \varphi_{\delta})}$ and $\delta_{\beta\gamma} \,
\delta_{\delta \alpha} \, e^{i (\varphi_{\alpha} - \varphi_{\beta} +
\varphi_{\gamma} - \varphi_{\delta})}$ are just
$\delta_{\alpha\beta} \, \delta_{\gamma \delta}$ and
$\delta_{\beta\gamma} \, \delta_{\delta \alpha}$ and they cancell
with the terms coming from the second part of \ref{cherndef}
\begin{equation}
- \, \delta_{\beta\gamma} \, \delta_{\delta \alpha} \, e^{-i
(\varphi_{\alpha} - \varphi_{\beta} + \varphi_{\gamma} -
\varphi_{\delta})} = - \, \delta_{\beta\gamma} \, \delta_{\delta
\alpha} \quad \hbox{and} \quad - \delta_{\gamma\delta} \,
\delta_{\alpha\beta} \, e^{-i (\varphi_{\alpha} - \varphi_{\beta} +
\varphi_{\gamma} - \varphi_{\delta})} = - \, \delta_{\gamma\delta}
\, \delta_{\alpha\beta} \, . \nonumber
\end{equation}
Thus one remains with the following terms:
\begin{equation}
- \, \delta_{\alpha\gamma} \, \delta_{\beta\delta} \, e^{i
(\varphi_{\alpha} - \varphi_{\beta} + \varphi_{\gamma} -
\varphi_{\delta})} - (- \, \delta_{\beta\delta} \,
\delta_{\gamma\alpha} \, e^{-i (\varphi_{\alpha} - \varphi_{\beta} +
\varphi_{\gamma} - \varphi_{\delta})}) = -2 \, i \;\sin
(2(\varphi_{\alpha} - \varphi_{\beta})) \nonumber
\end{equation}
which yield the second term in the proposition.
 \end{proof}
\bigskip

We now use the rescaling (\ref{rescale1}) in the case $\delta
(\varphi) \ne 0$ and rewrite ${\rm ch}_{3/2}$ in terms of the
generators $S_{\mu}$.

\medskip

  We let
\begin{equation}
\Lambda = \overset{3}{\underset{1}{\prod}} \, ({\rm tan} \,
(\varphi_j)\; \cos (\varphi_k - \varphi_{\ell}))\label{defofL}
\end{equation}
 \begin{equation}
 s_0=0 \, ,\quad s_j=1+ {\rm tan} \, \varphi_k
{\rm tan} \, \varphi_l \qqq j\in \{1,2,3\} \, . \label{defofsj}
\end{equation}
 and
\begin{equation}
 n_0=0 \, ,\quad n_k=1 \qqq k\in \{1,2,3\} \, .
\label{defofn}
\end{equation}

\bigskip

\begin{cor}\label{propch}
One has $$
 \, \Lambda\,{\rm ch}_{3/2} = -
\sum_{\alpha , \beta , \gamma , \delta} \varepsilon_{\alpha \beta
\gamma \delta} (n_{\alpha} - n_{\beta} + n_{\gamma} - n_{\delta})
(s_{\alpha} - s_{\beta} + s_{\gamma} - s_{\delta}) \, S_{\alpha}
\otimes S_{\beta} \otimes S_{\gamma} \otimes S_{\delta}$$
$$
 + 2 \, i \,
\sum_{\mu , \nu} (-1)^{n_{\mu} - n_{\nu}} \, (s_{\mu} - s_{\nu}) \,
S_{\mu} \otimes S_{\nu} \otimes S_{\mu} \otimes S_{\nu} \, .
 \, $$
\end{cor}

\medskip
\begin{proof} We write the above formula
as $\Lambda\,{\rm ch}_{3/2} = -A+\,2iB\,.$ We let the $\lambda_\mu$
be as in lemma \ref{rescalelem} so that
\begin{equation}
\label{eq3.4} \prod_0^3 \lambda_{\mu} = - \delta (\varphi) = -
\prod_1^3 (\sin \varphi_k \cos (\varphi_{\ell} - \varphi_m))
\end{equation}
and
\begin{equation}
\label{eq3.5} \lambda_0^2 = \prod_1^3 \sin\, \varphi_j \, , \quad
\lambda_k^2 = \sin \,\varphi_k \;\cos (\varphi_k - \varphi_{\ell})\,
\cos (\varphi_k - \varphi_m) \, .
\end{equation}
One has $x^{\mu} = \frac{S_{\mu}}{ \lambda_{\mu}}$ and thus in $
\Lambda\,{\rm ch}_{3/2}$ the first terms of proposition \ref{chbare}
give
\begin{equation}
\label{eq3.6} \varepsilon_{\alpha \beta \gamma \delta}
 \, \frac{\Lambda}{\prod \lambda_{\mu}} \cos (\varphi_{\alpha} - \varphi_{\beta} + \varphi_{\gamma} - \varphi_{\delta}) S_{\alpha} \otimes S_{\beta} \otimes S_{\gamma} \otimes S_{\delta} \qquad (\varphi_0 = 0) \, .
\end{equation}
One has
\begin{equation}
\frac{\Lambda}{\prod \lambda_{\mu}} = \,- \,\frac{1}{\prod \cos
\varphi_k} \, . \nonumber
\end{equation}

Thus the presence of the term $-A$ follows from the equality

\smallskip
\begin{equation}
\frac{\cos (\varphi_{\alpha} - \varphi_{\beta} + \varphi_{\gamma} -
\varphi_{\delta})}{\prod \cos \varphi_k} = (n_{\alpha} - n_{\beta} +
n_{\gamma} - n_{\delta}) (s_{\alpha} - s_{\beta} + s_{\gamma} -
s_{\delta}) \label{ending1}
\end{equation}

which we now check. We let $t_j = {\rm tan} \, \varphi_j$. One gets
$$\frac{\cos (\varphi_1 - \varphi_2 - \varphi_3)}{\prod \cos \varphi_k} = 1 - t_2 \, t_3 + t_1 \, t_3 + t_1 \, t_2$$
and more generally for {\it any} permutation $\sigma$ of $1,2,3$ one
has
\begin{equation}
\label{eq3.7} \frac{\cos (\varphi_{\sigma (1)} - \varphi_{\sigma
(2)} - \varphi_{\sigma (3)})}{\prod \cos \varphi_k} = -
s_{\sigma(1)} + s_{\sigma(2)} + s_{\sigma(3)} \, .
\end{equation}
Let us prove \eqref{ending1}. Since $(\alpha , \beta , \gamma ,
\delta)$ is a permutation of $(0,1,2,3)$  one of the indices is $0$.

  For $\alpha = 0$ we get
$\varphi_{\alpha}=0$ and  $$\frac{\cos ( \varphi_{\alpha}-
\varphi_{\beta} + \varphi_{\gamma} - \varphi_{\delta})}{\prod \cos
\varphi_k} = - (s_{0} - s_{\beta} + s_{\gamma} - s_{\delta})$$ since
$s_0=0$. But $(n_{\alpha} - n_{\beta} + n_{\gamma} - n_{\delta})=-1$
so that \eqref{ending1} holds.

\medskip

 For $\beta = 0$ we get
$\varphi_{\beta}=0$ and  $$\frac{\cos ( \varphi_{\alpha}-
\varphi_{\beta} + \varphi_{\gamma} - \varphi_{\delta})}{\prod \cos
\varphi_k} =  (s_{\alpha} - s_{0} + s_{\gamma} - s_{\delta})$$ since
$s_0=0$. But $(n_{\alpha} - n_{\beta} + n_{\gamma} -
n_{\delta})=\,1$ so that \eqref{ending1} holds.

\medskip

  For $\gamma = 0$ we get $\varphi_{\gamma}=0$ and  $$\frac{\cos ( \varphi_{\alpha}- \varphi_{\beta} + \varphi_{\gamma} - \varphi_{\delta})}{\prod \cos \varphi_k} =
- (s_{\alpha} - s_{\beta} + s_{0} - s_{\delta})$$ since $s_0=0$. But
$(n_{\alpha} - n_{\beta} + n_{\gamma} - n_{\delta})=-1$ so that
\eqref{ending1} holds.

\medskip

Finally for $\delta = 0$, $\varphi_{\delta}=0$ and  $$\frac{\cos (
\varphi_{\alpha}- \varphi_{\beta} + \varphi_{\gamma} -
\varphi_{\delta})}{\prod \cos \varphi_k} = (s_{\alpha} - s_{\beta} +
s_{\gamma} - s_{0})$$ since $s_0=0$. But $(n_{\alpha} - n_{\beta} +
n_{\gamma} - n_{\delta})=1$ so that \eqref{ending1} holds and we
checked it in all cases.

\medskip

Let us now compute the contribution of the second
 terms of proposition \ref{chbare}.
For $i,j \in \{ 1,2,3 \}$ one has
\begin{equation}
\label{eq3.9}
 -\frac{\Lambda}{2 \lambda_i^2 \, \lambda_j^2} \,  \sin \,2 (\varphi_i - \varphi_j)
=\;s_i - s_j  \, .
\end{equation}
Indeed say with $i=2$, $j=3$, one gets
$$
\frac{\Lambda}{\lambda_2^2 \, \lambda_3^2} =\frac{\sin
\varphi_1}{\cos \varphi_1 \cos \varphi_2 \cos \varphi_3 \cos
(\varphi_2 - \varphi_3)} \, .
$$
Thus
\begin{equation}
-\frac{\Lambda}{2 \lambda_2^2 \, \lambda_3^2} \,  \sin \,2
(\varphi_2 - \varphi_3)
 =  - \frac{\sin \varphi_1 \sin (\varphi_2 - \varphi_3)}{\cos \varphi_1 \cos \varphi_2 \cos \varphi_3} = -t_1 (t_2 - t_3) = (1+t_1 \, t_3) - (1+t_1 \, t_2) = s_2 - s_3 \, .
\nonumber
\end{equation}

Next let us check that
\begin{equation}
\label{eq3.10} \frac{\Lambda}{2 \lambda_0^2 \, \lambda_k^2} \,  \sin
\,2 \varphi_k =\,s_k \, .
\end{equation}

Say with $k=1$ one has $\lambda_0^2 = \prod \sin \varphi_j$,
$\lambda_1^2 = \sin \varphi_1 \cos (\varphi_1 - \varphi_2) \cos
(\varphi_1 - \varphi_3)$ and
\begin{equation}
\frac{\Lambda}{2 \lambda_0^2 \, \lambda_1^2} \,  \sin \,2 \varphi_1
=  \frac{\cos (\varphi_2 - \varphi_3)}{\cos \varphi_2 \cos
\varphi_3} = s_1 \, . \nonumber
\end{equation}

We thus get in general,
\begin{equation}
\frac{\Lambda}{ \lambda_\mu^2 \, \lambda_\nu^2} \, \sin 2
(\varphi_{\mu} - \varphi_{\nu}) = -(-1)^{n_{\mu} - n_{\nu}} \, 2
(s_{\mu} - s_{\nu}) \, .
\end{equation}
Indeed, for $\mu , \nu \in \{ 1,2,3 \}$ this is (\ref{eq3.9}). If
both $\mu , \nu = 0$ both sides are $0$. Now both sides are
antisymmetric in $\mu , \nu$ thus one can take $\nu = 0$, $\mu \in
\{ 1,2,3 \}$. Then $n_{\mu} - n_{\nu} = 1$ and the result follows
from (\ref{eq3.10}). \end{proof}

\medskip
\subsection{Elliptic parameters}

\noindent
\medskip

Let $\varphi \in {\mathbb T}_A$ and using \eqref{ellparam} let
$\tau\in \mathbb C$, $\im \,\tau >0$ and $\eta\in \mathbb C$ and
$\lambda \in \mathbb C$ such that
\begin{equation}\label{ellparam2}
\lambda \,( s_1, s_2, s_3)=
\left(\frac{\vartheta_2(0)^2}{\vartheta_2(\eta)^2},\, \,
\frac{\vartheta_3(0)^2}{\vartheta_3(\eta)^2},\, \,
\frac{\vartheta_4(0)^2}{\vartheta_4(\eta)^2}\right)
\end{equation}
We shall call the triplet $(\tau,\eta,\lambda)$ elliptic parameters
for $\varphi$. They are not unique given  $\varphi$ but they
determine uniquely the $s_j$ and hence $\varphi$ up to an overall
sign.

\medskip

 \begin{lem}  With the above notations \eqref{ellparam1}
is equivalent to
\begin{equation}
( {\tilde s}_1, {\tilde s}_2, {\tilde s}_3)= \,\lambda
\,\left(\frac{\vartheta_2(0)
\,\vartheta_2(2\eta)}{\vartheta_2(\eta)^2},\, \,\frac{\vartheta_3(0)
\,\vartheta_3(2\eta)}{\vartheta_3(\eta)^2},\, \,
\frac{\vartheta_4(0)
\,\vartheta_4(2\eta)}{\vartheta_4(\eta)^2}\right)
 \, .
\nonumber
\end{equation}
\smallskip
\end{lem}

\medskip \begin{proof} By lemma \ref{propor} we just need to
check that with
\begin{equation}\label{sklyajk}
(j_1,j_2,j_3)=\, \left(\frac{\vartheta_2(0)
\,\vartheta_2(2\eta)}{\vartheta_2(\eta)^2},\, \,\frac{\vartheta_3(0)
\,\vartheta_3(2\eta)}{\vartheta_3(\eta)^2},\, \,
\frac{\vartheta_4(0)
\,\vartheta_4(2\eta)}{\vartheta_4(\eta)^2}\right)
\end{equation}

one has
\begin{equation}\label{sklyajk1}
({\tilde j}_1,{\tilde j}_2,{\tilde j}_3)=\,
\left(\frac{\vartheta_2(0)^2}{\vartheta_2(\eta)^2},\, \,
\frac{\vartheta_3(0)^2}{\vartheta_3(\eta)^2},\, \,
\frac{\vartheta_4(0)^2}{\vartheta_4(\eta)^2}\right)\,.
\end{equation}

\medskip
The three equalities $ j_k\,{\tilde j}_\ell +\,{\tilde j}_k
\,j_\ell=\,2 $ follow from the following identities on
$\vartheta$-functions,
$$
\vartheta_3(0)^2\,\vartheta_2(0)\,\vartheta_2(2\eta)=\,\vartheta_2(\eta)^2
\,\vartheta_3(\eta)^2 -\vartheta_1(\eta)^2\,\vartheta_4(\eta)^2\,,
$$
$$
\vartheta_2(0)^2\,\vartheta_3(0)\,\vartheta_3(2\eta)=\,\vartheta_2(\eta)^2
\,\vartheta_3(\eta)^2 +\vartheta_1(\eta)^2\,\vartheta_4(\eta)^2\,,
$$
and similarly
$$
\vartheta_4(0)^2\,\vartheta_3(0)\,\vartheta_3(2\eta)=\,\vartheta_3(\eta)^2
\,\vartheta_4(\eta)^2 -\vartheta_1(\eta)^2\,\vartheta_2(\eta)^2\,,
$$
$$
\vartheta_3(0)^2\,\vartheta_4(0)\,\vartheta_4(2\eta)=\,
\vartheta_3(\eta)^2\,\vartheta_4(\eta)^2+\,\vartheta_1(\eta)^2
\,\vartheta_2(\eta)^2 \,,
$$
and
$$
\vartheta_2(0)^2\,\vartheta_4(0)\,\vartheta_4(2\eta)=\,\vartheta_2(\eta)^2
\,\vartheta_4(\eta)^2 +\vartheta_1(\eta)^2\,\vartheta_3(\eta)^2\,,
$$
$$
\vartheta_4(0)^2\,\vartheta_2(0)\,\vartheta_2(2\eta)=\,\vartheta_2(\eta)^2
\,\vartheta_4(\eta)^2 -\vartheta_1(\eta)^2\,\vartheta_3(\eta)^2\,.
$$

\end{proof}

This lemma allows to relate the above parameters with those used by
Sklyanin and one has with the above notations (\cite{skl:1983})
\begin{equation}\label{sklyaJ}
J_{12}=\,\frac{\vartheta_1(\eta)^2\vartheta_4(\eta)^2}{\vartheta_2(\eta)^2\vartheta_3(\eta)^2}
\,,\quad
J_{23}=\,\frac{\vartheta_1(\eta)^2\vartheta_2(\eta)^2}{\vartheta_3(\eta)^2\vartheta_4(\eta)^2}
\,,\quad
J_{31}=\,-\,\frac{\vartheta_1(\eta)^2\vartheta_3(\eta)^2}{\vartheta_2(\eta)^2\vartheta_4(\eta)^2}
\,,
\end{equation}
which follows from the definition \eqref{ellparam1} of the elliptic
parameters together with \eqref{thetarel}.

\medskip
\subsection{The sphere $ S_\varphi^3$
and the noncommutative torus ${\mathbb T}^2_{\eta}$}

\noindent
\medskip

Let $\varphi \in A^\circ$ so that $ \frac{\pi}{2} > \varphi_1 >
\varphi_2 > \varphi_3 > 0$. We can then choose the elliptic
parameters $\tau$ and $\eta$ such that
\begin{equation}
\label{ellreal} \tau \in i \,\mathbb R_+\,,\quad \eta \in [0,1]\,.
\end{equation}
Then the module $q=\,e^{i\tau}\in \,]0,1[$ and the $\vartheta$
functions $\vartheta_j(z)$ are all real functions \ie fulfill
\begin{equation}
\label{selftheta} \vartheta_j(\bar z)=\,\overline{
\vartheta_j(z)}\qqq z\in \mathbb C \,.
\end{equation}
 In particular the
last elliptic parameter $\lambda$ determined by \eqref{ellparam1}
fulfills $\lambda >0$.

We shall explain in this section how to use the representations
constructed by Sklyanin \cite{skl:1983} to obtain  $*$-homomorphisms
from $C_{\rm alg} (S_\varphi^3)$ to the algebra
\begin{equation}
C^{\infty}({\mathbb T}^2_{\eta})= \,C^{\infty}({\mathbb R}/{\mathbb
Z})\rtimes_\eta {\mathbb Z}\,,
\end{equation}
obtained as the crossed product of the algebra $C^{\infty}({\mathbb
R}/{\mathbb Z})$ of smooth periodic functions by the translation
$\eta$. Recall that a generic element of $C^{\infty}({\mathbb
T}^2_{\eta})$ is of the form
$$
f=\sum_{\mathbb Z}\,f_n\;V^n
$$
while the basic algebraic rule is given by
\begin{equation}
\label{crossrule} V\,f\,V^{-1}(u)=\,f(u+\eta)\qqq u\in {\mathbb
R}/{\mathbb Z}\qqq f\in C^{\infty}({\mathbb R}/{\mathbb Z})\,.
\end{equation}

Moreover $C^{\infty}({\mathbb T}^2_{\eta})$ is an involutive algebra
with involution turning $V$ into a unitary operator.

\medskip Starting from the representations constructed in \cite{skl:1983}
and conjugating by the operator
$$
M(\xi)(u)=\, e^{-2 \pi i u\,v/\eta}\,\xi(u+\tau/4)
$$
one performs a shift in the indices of the $\vartheta$-functions
based on
$$
\vartheta_1(z+\frac{\tau}{2})=\, q^{-\frac{1}{4}}\,e^{-\pi iz}i
\vartheta_4(z)\,,\quad \vartheta_2(z+\frac{\tau}{2})=\,
q^{-\frac{1}{4}}\,e^{-\pi iz} \vartheta_3(z)\,,
$$
$$
\vartheta_3(z+\frac{\tau}{2})=\, q^{-\frac{1}{4}}\,e^{-\pi iz}
\vartheta_2(z)\,,\quad \vartheta_4(z+\frac{\tau}{2})=\,
q^{-\frac{1}{4}}\,e^{-\pi iz}i \vartheta_1(z)\,.
$$
which allows to replace the singular denominator $\vartheta_1(2u)$
by $\vartheta_4(2u)$ which no longer vanishes for $u\in {\mathbb
R}/{\mathbb Z}$.

One  obtains this way a homomorphism from the Sklyanin algebra to
$C^{\infty}({\mathbb T}^2_{\eta})$ but it is not yet unitary and to
make it so one needs to conjugate again by a multiplication operator
of the form,
$$
N(\xi)(u)=\, \bar d(u)\,\xi(u)
$$
where the function $d\in C^{\infty}({\mathbb R}/{\mathbb Z})$
fulfills the following conditions,
$$
d(u)\bar d(u)=\,\vartheta_4(2 u) \qqq u\in {\mathbb R}/{\mathbb
Z}\,.
$$
We use the identity
$$
\vartheta_3(0)^2\,\vartheta_4(0)\,\vartheta_4(2u)=\,
\vartheta_3(u)^2\,\vartheta_4(u)^2+\,\vartheta_1(u)^2
\,\vartheta_2(u)^2 \,,
$$
and thus take
\begin{equation}\label{denom}
c\,\;
d(u)=\,\vartheta_3(u)\,\vartheta_4(u)+\,i\,\vartheta_1(u)\,\vartheta_2(u)
\,,\quad c^2= \vartheta_3(0)^2\,\vartheta_4(0)\,.
\end{equation}
Note that  one has
$$
 \bar d(u)=\,d(-u)\qqq u\in {\mathbb R}/{\mathbb Z}
$$
The effect of the conjugacy $N\,.\,N^{-1}$ on simple monomials
 is  the following
$$
 \frac{f(u)}{\vartheta_4(2\,u)}\,V\to \frac{f(u)}{d(u)\,d(-u-\eta)}\,V \,,
\quad
 \frac{f(u)}{\vartheta_4(2\,u)}\,V^*\to \frac{f(u)}{d(u)\,d(-u+\eta)}\,V^*
$$

The formulas which define the images  $\rho (S_{\alpha})$ then
become, with $m\in [0,\tau]$,
\begin{equation}
\label{sktonct0} \rho (S_0) = \vartheta_1 (\eta) \,  \frac{
\vartheta_3 (2u + \eta + im)}{d(u)\,d(-u-\eta)} \, V + \vartheta_1
(\eta) \, \frac{\vartheta_3 (2u - \eta - im) }{d(u)\,d(-u+\eta)} \,
\, V^*
\end{equation}
\begin{equation}
\label{sktonct1} \rho (S_1) = -i \, \vartheta_2 (\eta) \,
\frac{\vartheta_4 (2u + \eta + im)}{d(u)\,d(-u-\eta)} \,  \, V + i
\, \vartheta_2 (\eta) \, \frac{\vartheta_4 (2u - \eta -
im)}{d(u)\,d(-u+\eta)} \,  \, V^*
\end{equation}
\begin{equation}
\label{sktonct2} \rho (S_2) = \vartheta_3 (\eta) \,
\frac{\vartheta_1 (2u + \eta + im)}{d(u)\,d(-u-\eta)} \,  \, V +
\vartheta_3 (\eta) \, \frac{\vartheta_1 (2u - \eta -
im)}{d(u)\,d(-u+\eta)} \,  \, V^*
\end{equation}
\begin{equation}
\label{sktonct3} \rho (S_3) = -\vartheta_4 (\eta) \,
\frac{\vartheta_2 (2u + \eta + im)}{d(u)\,d(-u-\eta)} \,  \, V -
\vartheta_4 (\eta) \, \frac{ \vartheta_2 (2u - \eta - im)
}{d(u)\,d(-u+\eta)} \,V^*\, .
\end{equation}

and one has

\begin{thm} The formulas \eqref{sktonct0}...\eqref{sktonct3}
define a $*$-homomorphism from the Sklyanin algebra to $C^{\infty}
({\mathbb T}_{\eta}^2)\widehat{\otimes} C^{\infty} ([0,\tau])$.
\end{thm}

\begin{proof}
Since by construction $\rho$ is conjugate to an homomorphism it is
an homomorphism and we just need to check that the images $\rho
(S_\mu)$ of the generators are self-adjoint elements of $C^{\infty}
({\mathbb T}_{\eta}^2) $ for each value of $m\in [0,\tau]$.

One has
$$
(\frac{f(u)}{d(u)\,d(-u-\eta)}\,V )^*=\,V^*\,\frac{\bar
f(u)}{d(-u)\,d(u+\eta)}
$$
since $\eta \in \mathbb R$ and $\bar d(x)=\,d(-x)$ for $x\in \mathbb
R$. Thus using \eqref{crossrule} one gets
$$
(\frac{f(u)}{d(u)\,d(-u-\eta)}\,V )^*=\,\frac{\bar
f(u-\eta)}{d(-u+\eta)\,d(u)}\,V^*
$$
 Since $m$ is real one has $$\bar \vartheta_j(x+i m)=\,\vartheta_j(x-i m)
\qqq x \in \mathbb R \qqq j$$ using \eqref{selftheta}. Thus one
checks directly the required self-adjointness of the $\rho (S_\mu)$.
\end{proof}
\medskip

To obtain a $*$-homomorphism from
 $C_{\rm alg} (S^3_\varphi)$ to $C^{\infty} ({\mathbb T}_{\eta}^2) \otimes C^{\infty} ([0,\tau])$
we need to normalize the above formulas so that the element $Q_1$ in
the center of $C_{\rm alg} ({\mathbb R}^4_\varphi)$ gets mapped to
$1$. By proposition \ref{relcenter} this amounts to introduce an
overall scaling factor given by
\begin{equation}\label{scalingfactor}
\sigma(m)=\,
(\prod\,\sin\,\varphi_j)^{1/2}\;(C_1-\lambda\,C_2)^{-1/2}\,.
\end{equation}
where the explicit values of the Casimirs $C_j$ are given from
\cite{skl:1983} by

\begin{equation}\label{casimir}
C_1=\,4\, \vartheta^2_2(im)\,,\quad C_2=\,4\, \vartheta_2(\eta
+im)\,\vartheta_2(\eta -im)\,.
\end{equation}

We can now normalize the above homomorphism $\rho$ as
\begin{equation}\label{normalrho}
\tilde \rho(S_j)=\, \sigma(m)\,\rho(S_j)\,.
\end{equation}

We then get using the change of variables
$$
S_{\mu}=\,\lambda_{\mu}\;x^{\mu}\,,
$$

\medskip
\begin{cor} \label{rhotilde} The map $\tilde \rho$
defines a $*$-homomorphism
 $$C_{\rm alg} (S^3_\varphi)\to C^{\infty} ({\mathbb T}_{\eta}^2) \widehat{\otimes} C^{\infty} ([0,\tau])\,.$$
\end{cor}
\begin{proof} Since $\varphi \in A$ one has
$\lambda_{\mu}\in \mathbb R$ and the above change of variables is a
$*$-isomorphism of $C_{\rm alg} ({\mathbb R}^4_\varphi)$ with the
Sklyanin algebra. Thus we only need to check that with the above
normalization the $*$-homomorphism $\tilde \rho$ maps the central
element $Q_1$ which determines the sphere to the element $$ 1\in
C^{\infty} ({\mathbb T}_{\eta}^2) \widehat{\otimes} C^{\infty}
([0,\tau])$$ This follows from  Equation (\ref{relcenter}).
\end{proof}

In the following we shall use the notation $C^\infty(\mathbb T
^2_\eta\times [0,\tau])$ to denote the completed tensor product
$C^\infty(\mathbb T^2_\eta) \widehat{\otimes} C^\infty([0,\tau])$.
\bigskip
\subsection{Pairing with $[{\mathbb T}_{\eta}^2]$}

\noindent
\medskip

In order to test the non-triviality of the $*$-homomorphism $\tilde
\rho$ we shall compute what will be later interpreted as its
Jacobian. In order to do this we shall pair the image

\begin{equation}
\tilde \rho_*({\rm ch}_{3/2})\in HZ_3( C^{\infty} ({\mathbb
T}_{\eta}^2\times [0,\tau]))
\end{equation}

with the natural Hochschild three cocycle obtained using the
fundamental class $[{\mathbb T}_{\eta}^2]$ introduced in
\cite{ac:1980}. Since the variable $m\in [0,\tau]$ labels the center
we shall view the above pairing as defining a function of $m$.

The basic hochschild three cocycle on $C^{\infty} ({\mathbb
T}_{\eta}^2\times[0,\tau]))$
 is given by

\begin{equation}
\tau(a_0,\cdots,a_3)=\, \sum
\,\epsilon_{ijk}\,\tau_0(a_0\,\delta_i(a_1)\,
\delta_j(a_2)\,\delta_k(a_3))\, \,.
\end{equation}
where $\tau_0$ is the trace obtained as the tensor product of the
 canonical trace $\chi$ on $C^{\infty} ({\mathbb T}_{\eta}^2)$ by
the trace on $C^{\infty} ([0,\tau])$ given by integration,

\begin{equation}
\tau_0(a ))=\,\int_0^\tau\,\chi(a (m))dm\,,\quad
\chi(f)=\int_0^1\,f(u)du\,, \quad \chi(f\,V^n)=0 \qqq n\neq 0 \,.
\end{equation}

The three basic derivations $\delta_j$ are given by
\begin{equation}
\delta_1=\,\partial/ \partial m\,,\quad \delta_2=\,\partial/
\partial u \,,\quad \delta_3=\,2\pi i\,V\,\partial/ \partial V \,,
\end{equation}
where in the last term the differentiation $V\,\partial/ \partial V
$ has the effect of multiplying by $n$ any monomial $f\,V^n$.

As the product of $\tau$ by any function $h(m)$ viewed as an element
of the center of $C^{\infty} ({\mathbb T}_{\eta}^2\times[0,\tau]))$
is still a Hochschild three cocycle we obtain a differential one
form on $[0,\tau]$ as the pairing

\begin{equation} \label{pairing0}
\omega =< {\rm ch}_{3/2},\, \tau >
\end{equation}

\smallskip
The basic lemma is then the following using the notation
\begin{equation}\label{thesigmas}
(\sigma_0,\sigma_1,\sigma_2,\sigma_3)=\,\left(0,\,\frac{\vartheta_2(0)^2}{\vartheta_2(\eta)^2},\,
\, \frac{\vartheta_3(0)^2}{\vartheta_3(\eta)^2},\, \,
\frac{\vartheta_4(0)^2}{\vartheta_4(\eta)^2}\right)
\end{equation}

\begin{lem} \label{maincomp}
With the above notations one has
$$
<\tau ,\,\sum_{\alpha , \beta , \gamma , \delta} \varepsilon_{\alpha
\beta \gamma \delta} (n_{\alpha} - n_{\beta} + n_{\gamma} -
n_{\delta}) (\sigma_{\alpha} - \sigma_{\beta} + \sigma_{\gamma} -
\sigma_{\delta}) \, \rho(S_{\alpha}) \otimes \rho(S_{\beta}) \otimes
\rho(S_{\gamma}) \otimes \rho(S_{\delta})$$
$$
- 2 \, i \, \sum_{\mu , \nu} (-1)^{n_{\mu} - n_{\nu}} \,
(\sigma_{\mu} - \sigma_{\nu}) \, \rho(S_{\mu}) \otimes \rho(S_{\nu})
\otimes \rho(S_{\mu}) \otimes \rho(S_{\nu}) >\;=
 \, $$
$$
24 \,(2\pi i)^3
\;\frac{\vartheta'_1(0)^3}{\pi^3}\;\frac{\vartheta_1(\eta)\,\vartheta_1(2im)}{
\vartheta_2(\eta)\vartheta_3(\eta)\vartheta_4(\eta)}\,\,.
$$
\end{lem}

\medskip
The precise meaning of the equality is that for any $h\in C^{\infty}
([0,\tau]))$ the evaluation of the Hochschild cocycle $\tau$ on the
product of the Hochschild cycle of the right hand side by $h$ gives
the integral
$$
\int_0^\tau\;h(m) \,g(m)\,dm
$$
where\footnote{The $q$-expansion of the fraction
$\frac{\vartheta'_1(0)}{\pi}$ has rational coefficients.}
\begin{equation} \label{gdem}
g(m)=\, 24 \,(2\pi i)^3
\;\frac{\vartheta'_1(0)^3}{\pi^3}\;\frac{\vartheta_1(\eta)\,\vartheta_1(2im)}{
\vartheta_2(\eta)\vartheta_3(\eta)\vartheta_4(\eta)}\,\,.
\end{equation}
\medskip

The proof of this is a long computation based on the ``a priori"
properties of the pairing which allow to show that the dependence in
the parameters $\eta$ and $m$ is of the expected form, while the
dependence in the module $q$ is that of a modular form. It then
follows from the explicit knowledge of enough terms in the
$q$-expansion that the above formula is valid. So far we have not
been able to eliminate completely the use of the computer to check
this validity and its understanding will only come through the
gradual simplifications below.

\medskip
We shall now show that the dependance in $m$ of the normalization
factor $\sigma(m)$ in the definition \eqref{normalrho} of the
homomorphism $\tilde \rho$ can be ignored when one computes the
pairing \eqref{pairing0}

\medskip
\begin{lem}\label{resder}
Let $\delta$, $\delta'$ be derivations of the unital algebra $\cala$
preserving a trace $\tau_0$ on $\cala$. Let $\phi_j$ be the
multilinear forms on $\cala$ given by
$$
\phi_1(a_0,a_1,a_2,a_3)=\,\tau_0(a_0\,a_1\,\delta(a_2)\,\delta'(a_3))
\,,\quad
\phi_2(a_0,a_1,a_2,a_3)=\,\tau_0(a_0\,\delta(a_1)\,a_2\,\delta'(a_3))
$$
Then for any invertible $U\in \cala$ one has
$$
\phi_j(U,U^{-1},U,U^{-1})-\,\phi_j(U^{-1},U,U^{-1},U)=\,0\,.
$$
\end{lem}

\begin{proof} One has
$$
\phi_1(U,U^{-1},U,U^{-1})=\,\tau_0(U\,U^{-1}\,\delta(U)\,\delta'(U^{-1}))
=\,-\tau_0(\delta(U)\,U^{-1}\,\delta'(U)\,U^{-1})
$$
$$
\phi_1(U^{-1},U,U^{-1},U)=\,\tau_0(U^{-1}\,U\,\delta(U^{-1})\,\delta'(U))
=\,-\tau_0(U^{-1}\,\delta(U)\,U^{-1}\,\delta'(U))
$$
thus the cyclicity of the trace proves the statement for $j=1$.
Similarly one has
$$
\phi_2(U,U^{-1},U,U^{-1})=\,\tau_0(U\,\delta(U^{-1})\,U\,\delta'(U^{-1}))
=\,\tau_0(\delta(U)\,U^{-1}\,\delta'(U)\,U^{-1})
$$
$$
\phi_2(U^{-1},U,U^{-1},U)=\,\tau_0(U^{-1}\,\delta(U)\,U^{-1}\,\delta'(U))
$$
and the cyclicity of the trace proves the statement for $j=2$.
\end{proof}

\medskip
We thus get the following result

\begin{cor}\label{rhoch}
The pairing of $\tilde \rho_*({\rm ch}_{3/2})$ with $\tau$ is given
by the differential form
$$\omega=\,-\,\frac{\sigma(m)^{4}\,g(m)}{\lambda \,\Lambda}\,dm\,,$$
with $\Lambda$ given in \eqref{defofL}, $\lambda$ by
\eqref{ellparam1}, $\sigma(m)$ by \eqref{scalingfactor} and $g(m)$
by \eqref{gdem}.
\end{cor}

\begin{proof} Using lemma \ref{maincomp}
and corollary \ref{propch} one just needs to show that the terms of
the form
$$
\delta_1(\tilde{\rho}(S_j))-\,\sigma(m)\,\delta_1(\rho(S_j))
=\,\frac{d\sigma(m)}{d m}\;\rho(S_j)
$$
do not contribute. But their total contribution is a sum of six
terms each of which is of the form
$$
\phi_j(U,U^{-1},U,U^{-1})-\,\phi_j(U^{-1},U,U^{-1},U)
$$
where $U \in \,M_2(C_{\rm alg} (S^3_\varphi))$ is the basic unitary
element while $\phi_j$ is as in lemma \ref{resder} with
$\delta,\delta'\in\,\{\delta_2,\delta_3\}$. Thus each of these terms
vanishes by lemma \ref{resder}.
\end{proof}

\bigskip
\subsection{Simplifying the
$*$-homomorphism  $\tilde \rho$} \label{simplerrho}

\noindent
\medskip

We shall make several simplifications in the formulas involved in
the construction of the $*$-homomorphism $\tilde \rho$ of corollary
\ref{rhotilde} in order to gradually eliminate all
$\vartheta$-functions and express the result in purely algebraic
terms.

The denominators involved in the construction of the
$*$-homomorphism $\tilde \rho$ are of the form
\begin{equation}\label{denom1}
d(u)\,d(-u \pm \eta)
\end{equation}

where by \eqref{denom},
$$
c\,\;
d(u)=\,\vartheta_3(u)\,\vartheta_4(u)+\,i\,\vartheta_1(u)\,\vartheta_2(u)
\,,\quad c^2= \vartheta_3(0)^2\,\vartheta_4(0)\,.
$$
Our first task will be to rewrite \eqref{denom1} as a linear form in
terms of the projective coordinates $\psi(u)$ of proposition
\ref{theta} \ie
$$
\psi(u)=\left( \frac{\vartheta_1(2 u-\eta)}{\vartheta_1(\eta)},\,
\frac{\vartheta_2(2 u-\eta)}{\vartheta_2(\eta)},\,
\frac{\vartheta_3(2 u-\eta)}{\vartheta_3(\eta)},\,
\frac{\vartheta_4(2 u-\eta)}{\vartheta_4(\eta)}\right) =
(Z_0,Z_1,Z_2,Z_3)
$$

\medskip
\begin{lem}
With the above notations one has
\begin{equation}\label{denom2}
\vartheta_3(0)\,d(u)\,d(-u +
\eta)=\,i\,\vartheta_1(\eta)\vartheta_2(\eta)\,Z_1+\,
\vartheta_3(\eta)\vartheta_4(\eta)\,Z_3
\end{equation}
\end{lem}

\begin{proof} One has
$$
c^2\,d(u)\,d(-u +
\eta)=\,(\vartheta_3(u)\,\vartheta_4(u)+\,i\,\vartheta_1(u)\,\vartheta_2(u))
(\vartheta_3(u-\eta)\,\vartheta_4(u-\eta)-\,i\,\vartheta_1(u-\eta)\,\vartheta_2(u-\eta))
$$
$$
=\,\vartheta_3(u)\,\vartheta_4(u)\,\vartheta_3(u-\eta)\,\vartheta_4(u-\eta)
+\,\vartheta_1(u)\,\vartheta_2(u)\,\vartheta_1(u-\eta)\,\vartheta_2(u-\eta)
$$
$$
+\,i\,\vartheta_1(u)\,\vartheta_2(u)\,\vartheta_3(u-\eta)\,\vartheta_4(u-\eta)
-\,i\,\vartheta_3(u)\,\vartheta_4(u)\,\vartheta_1(u-\eta)\,\vartheta_2(u-\eta)
$$
Thus using the basic addition formulas (obtained from \eqref{relt6}
and \eqref{relt15})
$$
\vartheta_3(x)\,\vartheta_4(x)\,\vartheta_3(y)\,\vartheta_4(y)
-\,\vartheta_1(x)\,\vartheta_2(x)\,\vartheta_1(y)\,\vartheta_2(y))=
\,\vartheta_3(0)\,\vartheta_4(0)\,\vartheta_3(x+y)\,\vartheta_4(x-y)
$$
$$
\vartheta_1(x)\,\vartheta_2(x)\,\vartheta_3(y)\,\vartheta_4(y)
+\,\vartheta_3(x)\,\vartheta_4(x)\,\vartheta_1(y)\,\vartheta_2(y))=
\,\vartheta_3(0)\,\vartheta_4(0)\,\vartheta_1(x+y)\,\vartheta_2(x-y)
$$
 for $x=u$, $y=\eta -u$, we get
$$
c^2\,d(u)\,d(-u + \eta)=\,\vartheta_3(0)\,\vartheta_4(0)\,(
\vartheta_3(\eta)\,\vartheta_4(2
u-\eta)+\,i\,\vartheta_1(\eta)\,\vartheta_2(2 u-\eta))\,,
$$
which gives the required equality.
\end{proof}

To simplify the numerators involved in the construction of the
$*$-homomorphism $\tilde \rho$ we pass from  generators $S_\mu$ of
the Sklyanin algebra to the generators $Y_\mu$ of Theorem \ref{iden}
by the following transformation
\begin{equation}\label{simplerho0}
S_0=\,d\,Y_2\,,\quad S_1=\,i\,Y_3\,,\quad S_2=\,d\,Y_0\,,\quad
S_3=\,-Y_1\,,
\end{equation}
where $\displaystyle
d=\,\frac{\vartheta_1(\eta)\vartheta_3(\eta)}{\vartheta_2(\eta)\vartheta_4(\eta)}\,.$

One checks that the $Y_\mu$ fulfill the presentation of Theorem
\ref{iden} using the equality \eqref{sklyaJ} $ d^2=\,-\,J_{31}\,. $
We can then reformulate the construction of the homomorphism $\rho$
in the following terms,

\medskip
\begin{lem}\label{lemrho1} With the above notations one has, up to an
overall scalar factor $\gamma$,
\begin{equation}\label{simplerho}
\rho(Y_\mu)=\,\frac{\psi_\mu(u-\,i\,m/2)}{L(u)}\,\;V^*+
\epsilon_\mu\,V\,\; \frac{\psi_\mu(u+\,i\,m/2)}{\bar L(u)}
\end{equation}
where $\epsilon=(1,1,1,-1)$ and
$$
L(u)=\, \,i\,\vartheta_1(\eta)\vartheta_2(\eta)\,\psi_1(u)+\,
\vartheta_3(\eta)\vartheta_4(\eta)\,\psi_3(u)\,,\quad \bar
L(u)=\overline{ L(u)}\,.
$$

\end{lem}

\begin{proof} One just needs to perform the  transformation
\eqref{simplerho0} on the equations
\eqref{sktonct0}....\eqref{sktonct3}. One gets an overall scalar
factor
$$
\gamma=\,\vartheta_2(\eta)\,\vartheta_4(\eta)\,\vartheta_3(0)\,.
$$
multiplying the right hand side of \eqref{simplerho} (or
equivalently dividing $L(u)$).
\end{proof}

\medskip
In order to understand \eqref{simplerho} we let
\begin{equation}\label{simplerho1}
Z=\,\psi(u-\,i\,m/2)\,,\quad Z'=\,\epsilon\;\psi(u+\,i\,m/2)\,,\quad
W= L(u)^{-1}\,V^* \,,\quad W'=\,V\,{\bar L(u)}^{-1}\,.
\end{equation}
Note that one has $Z'=\,\epsilon\,\bar Z$ and $W'=\,W^*$ but we
shall ignore that for a while and treat for instance $Z$ and $Z'$ as
independent variables. The multiplicative terms such as
$L(u)^{-1}$ do not alter the cross product rules \eqref{crossrule}
but they alter the simplification rule $V\,V^*=\,V^*\,V=\,1$. Our
next task will thus be to give a simple expression for $W\,W'$ in
terms of $(Z,Z')$.

One has by construction
\begin{equation}\label{simplerho14}
W\,W'=\,(L(u)\,\bar L(u))^{-1}\,,
\end{equation}
and we need to express the denominator in terms of $Z$ and $Z'$.
Note that we have the freedom to multiply by an arbitrary function
of $m$ since this only alters the normalization of $\rho$ which is
needed in any case to pass to $\tilde \rho$.

\medskip
\begin{lem} \label{lemrho2}
With the above notations one has,
\begin{equation}\label{simplerho2}
\nu(m)\,\,L(u)\,\bar L(u)=\,J_{23}\,(Z_0\,Z'_0+\,Z_1\,Z'_1)+\,
Z_2\,Z'_2-\,Z_3\,Z'_3\,,
\end{equation}
where $$\nu(m)=\,\frac{2\,\vartheta^2_3(i m)}{\vartheta^2_3(0)
\vartheta^2_3(\eta)\vartheta^2_4(\eta)}\,.$$
\end{lem}

\begin{proof} One has
$$
L(u)\,\bar L(u)=\,\vartheta^2_1(\eta)\vartheta^2_2(2u-\eta)\,+
\vartheta^2_3(\eta)\vartheta^2_4(2u-\eta)\,,
$$

\smallskip
thus with $$a =2u-\eta+i m \,,\quad  b =2u-\eta-i m \,,\quad
\frac{a+b}{2}=\,2u-\eta \,,\quad\frac{a-b}{2}=\,i m$$ we get
\begin{equation}\label{simplerho3}
\vartheta^2_3(i m)\,\,L(u)\,\bar L(u)=\,\vartheta^2_1(\eta)
\,\vartheta^2_3(\frac{a-b}{2}) \,\vartheta^2_2(\frac{a+b}{2})+
\vartheta^2_3(\eta)\, \,\vartheta^2_3(\frac{a-b}{2})
\,\vartheta^2_4(\frac{a+b}{2})\,.
\end{equation}
We now use the addition formulas
$$
2\,\vartheta^2_3(\frac{a-b}{2})
\,\vartheta^2_2(\frac{a+b}{2})=\,\vartheta^2_2(0)
\,\vartheta_3(a)\,\vartheta_3(b)+ \,\vartheta^2_3(0)
\,\vartheta_2(a)\,\vartheta_2(b)- \,\vartheta^2_4(0)
\,\vartheta_1(a)\,\vartheta_1(b)
$$
(adding \eqref{relt9} and \eqref{relt10}) and
$$
2\,\vartheta^2_3(\frac{a-b}{2})
\,\vartheta^2_4(\frac{a+b}{2})=\,\vartheta^2_2(0)
\,\vartheta_1(a)\,\vartheta_1(b)+ \,\vartheta^2_3(0)
\,\vartheta_4(a)\,\vartheta_4(b)+ \,\vartheta^2_4(0)
\,\vartheta_3(a)\,\vartheta_3(b)
$$
(adding \eqref{relt5} and \eqref{relt6}) which allow to write
\eqref{simplerho3} as a symmetric bilinear form in $(Z,Z')$. One
then uses \eqref{thetarel} and \eqref{sklyaJ}
$$
J_{23}=\,\frac{\vartheta^2_1(\eta)\vartheta^2_2(\eta)}{
\vartheta^2_3(\eta)\vartheta^2_4(\eta)}\,,$$ to obtain the required
equality.
\end{proof}

\medskip
\begin{prop} \label{simplerhofinal} With the above notations one has, up to an
overall scalar factor $\delta(m)$,
\begin{equation}\label{simplerho4}
\rho(Y_\mu)=\,Z_\mu\,\;W + W'\,\; Z'_\mu
\end{equation}
with algebraic rules given by
\begin{equation}\label{simplerho5}
Z_i\,W \,\; W'\,Z'_j =\,\frac{Z_i\,Z'_j}{Q(Z,Z')}\,,\quad
W\,f(Z,Z')=\,f(\sigma(Z), \sigma^{-1}(Z'))\,W\,,
\end{equation}
where $\sigma$ is the translation by $-\eta$ as in Theorem
\ref{iden} and
$$
Q(Z,Z')=\,\,J_{23}\,(Z_0\,Z'_0+\,Z_1\,Z'_1)+\,
Z_2\,Z'_2-\,Z_3\,Z'_3\,.
$$
\end{prop}

\begin{proof}
The first equality follows from lemma \ref{lemrho1} and the
definition \ref{simplerho1} of $Z,Z',W,W'$. The first algebraic rule
follows from \eqref{simplerho14} and lemma \ref{lemrho2}.

To obtain the second we need to understand the transformation
$$
\epsilon \;\psi(u+i m/2) \to \epsilon \;\psi(u-\eta + i m/2)\,,
$$
and to compare it with $\sigma^{-1}$ where $\sigma$ is the
translation by $-\eta$ as in Theorem \ref{iden}.

By construction $\sigma$ is the product \eqref{twoinvolutions} of
two involutions   $\sigma =I\circ I_0$ where $I_0$ just alters the
sign of $Z_0$ (\cf \eqref{firstinvolution}). Thus
$\sigma^{-1}=\,I_0\circ I= \,I_0\circ \sigma  \circ I_0$ and to show
that the above tranformation is $\sigma^{-1}$ it is enough to show
that $\sigma$ commutes with $I_0\circ I_3$ where
$I_3(Z)=\,\epsilon\,Z$. This follows from the commutation of
translations on the elliptic curve and can be checked directly using
\eqref{algebraicsigma}.
\end{proof}

\section{Algebraic geometry and $C^*$-algebras}\label{Cstar}

In this section we shall develop the basic relation between
noncommutative differential geometry in the sense of \cite{ac:1980}
and noncommutative algebraic geometry. This will be obtained by
abstracting  the results of  proposition \ref{simplerhofinal} of
subsection \ref{simplerrho} and giving a general construction,
independent of $\vartheta$-functions, of a homomorphism from a
quadratic algebra to a crossed product algebra constructed from the
geometric data.

\medskip

\subsection{Central Quadratic Forms and Generalised Cross-Products}

\noindent
\medskip

 Let $\cala=A(V,R)=T(V)/(R)$  be a  quadratic algebra.
Its geometric data $\{E\,,\,  \sigma\,,\,\call\}$
 is defined in
such a way that $\cala$ maps homomorphically to a cross-product
algebra obtained from sections of  powers of the line bundle $\call$
on powers of  the correspondence $\sigma$ (\cite{art-tat-vdb:1990}).

 We shall begin by a purely algebraic result
which considerably refines the above homomorphism and lands in a
richer cross-product. We use the notations of section
\ref{geometricdata} for general quadratic algebras.

\begin{defn} \label{cent}
Let $Q \in S^2(V)$ be a symmetric bilinear form on $V^\ast$ and $C$
a component of $E \times E$. We shall say that $Q$ is
\underline{central} on $C$ iff for all ($Z,\,Z'$) in $C$ and
$\omega\in R$ one has,
\begin{equation} \label{defcentral}
\omega(Z,Z')\, Q(\sigma(Z'),\sigma^{-1}(Z))+Q(Z,Z')\,
\omega(\sigma(Z'),\sigma^{-1}(Z)) =0
\end{equation}
\end{defn}

\medskip By construction the space of symmetric bilinear form on $V^\ast$
which are central on $C$ is a linear subspace of $S^2(V)$.
 Let $C$ be a component of $E \times E$
globally invariant under the map
\begin{equation}
\tilde{\sigma}(Z,Z'):=\,(\sigma(Z),\sigma^{-1}(Z')) \label{siginv}
\end{equation}
 Given a quadratic form  $Q$ central and not identically zero
on the component $C$, we define as follows an algebra $C_Q$ as a
generalised cross-product of the ring $\mathcal R$ of meromorphic
 functions on $C$
by the transformation $\tilde{\sigma}$. Let $L$, $L'\in V$ be such
that $L(Z)\,L'(Z')$ does not vanish identically on $C$. We adjoin
two generators $W_L$ and $W'_{L'}$ which besides the usual
cross-product rules,
\begin{equation}
W_L \,f = (f\circ \tilde{\sigma}) \;W_L \,, \quad W'_{L'} \,f =
(f\circ \tilde{\sigma}^{-1}) \;W'_{L'}\,, \quad \forall f\in
\mathcal R \label{cropro}
\end{equation}
fulfill the following relations,
\begin{equation}
W_L \, W'_{L'}:=\pi(Z,Z')\, , \qquad W'_{L'}\, W_L :=
\pi(\sigma^{-1}(Z),\sigma(Z')) \label{crossed}
\end{equation}
where the function $\pi(Z,Z')$ is given by the ratio,
\begin{equation}
\pi(Z,Z'):=\frac{L(Z)\,L'(Z')}{Q(Z,Z')} \label{ratio}
\end{equation}
The a priori dependence on $L$, $L'$ is eliminated by the rules,
\begin{equation}
W_{L_2}:=\frac{L_2(Z)}{L_1(Z)}\,W_{L_1}  \qquad
W'_{L'_2}:=\,W'_{L'_1} \, \frac{L'_2(Z')}{L'_1(Z')} \label{compare}
\end{equation}
which allow to adjoin all $W_L$ and $W'_{L'}$ for $L$ and $L'$ not
identically zero on the projections of $C$, without changing the
algebra and provides an intrinsic definition of $C_Q$.

\smallskip Our first result is

\begin{lem}\label{alg0}
Let $Q$ be  central and not identically zero on the component $C$.

\smallskip (i) The following equality defines a homomorphism $\rho$:
$\cala \mapsto C_Q$
\begin{equation} \label{morphism}
\sqrt{2}\;\rho(Y) :=  \frac{Y(Z)}{L(Z)}\,W_L+
W'_{L'}\,\frac{Y(Z')}{L'(Z')}\, , \qquad \forall Y \in V
\end{equation}
\smallskip (ii) If $\sigma^4 \neq \bbbone$, then $\rho(Q)=1$ where
$Q$ is viewed as an element of $T(V)/(R)$.
\end{lem}

\begin{proof} (i) Formula \eqref{morphism} is independent of $L$, $L'$ using
 (\ref{compare}) and reduces (up to $\sqrt{2}$) to $W_Y+ W'_{Y}$
when $Y$ is non-trivial on the two projections of $C$.
 It is enough to check that the
$\rho(Y)\in C_{Q}$ fulfill the quadratic relations
 $\omega\in R$. Let $\omega\in R$
$$
\omega(Z,Z')=\,\sum \omega_{ij}\,Y_i(Z)\,Y_j(Z')
$$
viewed as a bilinear form on  $V^\ast$. One has

\smallskip
$$
2\,\sum \omega_{ij}\,\rho(Y_i)\,\rho(Y_j)= \,\sum
\omega_{ij}\,(W_{Y_i}+ W'_{Y_i})(W_{Y_j}+ W'_{Y_j})= \,\sum
\omega_{ij}\,Y_i(Z)\,Y_j(\sigma(Z)\,W^2 \;+
$$
$$
\sum \omega_{ij}\,\frac{Y_i(Z)\,Y_j(Z')}{Q(Z,Z')}+\, \sum
\omega_{ij}\,\frac{Y_i(\sigma(Z'))\,Y_j(\sigma^{-1}(Z))}{Q(\sigma^{-1}(Z),
\sigma(Z'))}+\,\sum
\omega_{ij}\,W^{'\,2}\,Y_i(\sigma^{-1}(Z'))\,Y_j(Z')\,
$$

\smallskip
where $$W^2=\,\frac{1}{L(Z)L(\sigma(Z))}W_L^2\,,\quad
W^{'\,2}=W_{L'}^{'\,2}\frac{1}{L'(\sigma^{-1}(Z'))L'(Z')}\,.$$

The vanishing of the terms in $W^2$ and in $W^{'2}$ is automatic by
construction of the correspondence $\sigma$ \ie the equality
$$
\omega(Z,\sigma(Z))=\,0\qqq Z\in E\,.
$$

 The sum of the middle terms
is just
$$
\frac{\omega(Z,Z')}{Q(Z,Z')}+\,\frac{\omega(\sigma(Z'),\sigma^{-1}(Z))
}{Q(\sigma^{-1}(Z), \sigma(Z'))}=\,0\,,
$$
as follows from definition \ref{cent} and the symmetry of $Q$.

\smallskip
(ii) The above computation shows that $\rho(Q)=1$ provided one can
show that $Q(Z,\sigma(Z))=0$ and $Q(\sigma^{-1}(Z'),Z')=0$ for all
$Z, Z'$ in the projections $E$, $E'$ of $C$. We assume that
$\sigma^4(Z)$ is not identically equal to $Z$
 on each connected component of
$E$ (resp. $E'$) and use \eqref{defcentral} with $Z'=\sigma(Z)$. The
first term vanishes and we get
$$
\omega(\sigma^2(Z),\sigma^{-1}(Z))\,Q(Z,\sigma(Z))=\,0\qqq \omega
\in R\,.
$$
Thus if $Q(Z,\sigma(Z))$ does not vanish identically on  a given
connected component $E_1$ of $E$ one gets that
$$
\omega(\sigma^2(Z),\sigma^{-1}(Z))=\,0\qqq Z\in E_1\,,\quad \omega
\in R\,,
$$
  so that $\sigma^{-1}(Z)=\,\sigma^3(Z)$ for all $Z\in E_1$
which contradicts the hypothesis.
\end{proof}

\medskip Let $\cala_\ug= C_{\mathrm{alg}}(\mathbb C^4(\ug))$
 at generic $\ug$, then by proposition \ref{centerY}
the center of $\cala_\ug$ contains
 the three linearly dependent quadratic elements
\begin{equation}
Q_m\,:=\,J_{k \ell}\,(Y_0^2 \, +\, Y_m^2)\,+ \, Y_k^2 \, -
\,Y_\ell^2
\end{equation}

\begin{prop} \label{central}
Let $\cala_\ug= C_{\mathrm{alg}}(\mathbb C^4(\ug))$  at generic
$\ug$, then each $Q_m$ is central on $F_\ug \times F_\ug$ ($\subset
E_\ug \times E_\ug$).
\end{prop}

\begin{proof} One uses (\ref{twoinvolutions}) to check the algebraic identity.
\end{proof}

Together with lemma \ref{alg0} this yields non trivial homomorphisms
of $\cala_\ug$ whose unitarity will be analysed in the next section.
Note that for a general quadratic algebra $\cala=A(V,R)=T(V)/(R)$
and a quadratic form $Q \in S^2(V)$, such that $Q \in $
Center($\cala$), it does not automatically follow that $Q$ is
central on $E \times E$. For instance Proposition \ref{central} no
longer holds on $F_\ug \times \{e_\nu\}$ where $e_\nu$ is any of the
four points of $ E_\ug$ not in $F_\ug$. In fact let us describe in
some details what happens in the case of the $\theta$-deformations
\ie  $C_+=\{(\varphi,\varphi,0)\}$ (case $7$). We take the notations
of subsection \ref{chartheta} to write the characteristic variety as
the union of six lines $\ell_j$.

\begin{prop} \label{centraltheta}
Let $ C_{\mathrm{alg}}({\mathbb R}_\varphi^4)$  for $\varphi \in
C_+$, and
 $Q_k$ be defined by \eqref{center1}
and \eqref{center22}.
\begin{enumerate}
\item
Each $Q_k$ is central on $\ell_i \times \ell_j$ provided $i$ and $j$
belong to the same subsets $I=\{1,2\}$ and $J=\{3,4,5,6\}$.

 \item The bilinear form $Q_1$ does not vanish identically on $\ell_i \times \ell_j$ iff $\,i=j$
for $i,j\in I$ and iff $\,i\neq j$ for $i,j\in J$.

 \item The bilinear form $Q_2$ vanishes identically on
$\ell_2 \times \ell_j$ and $\ell_j \times \ell_2$ for all $j$.

\end{enumerate}

\end{prop}

This is proved by direct computations. Note that since $Q_1$ fails
to be central on $\ell_1 \times \ell_3$ for instance, it was crucial
to ``localize" the notion of central quadratic form to components of
the square $E\times E$ of the characteristic variety $E$. It is of
course also crucial to check the non-vanishing of $Q$ when applying
lemma \ref{alg0}, and the component $\ell_2$ does not work for $Q_2$
in that respect.

The precise table for the vanishing of the form \({Q_2}\) is the
following where $\neq $ at $(i,j)$ means that \({Q_2}\) does not
vanish identically on $\ell_i \times \ell_j$,

\bigskip
\begin{center}\begin{tabular}[c]{cccccc}
     $\neq $ &0&$\neq $ &$\neq $ &$\neq $ &$\neq $  \\
    0&0&0&0&0&0 \\
    $\neq $ &0&0&$\neq $ &$\neq $ &0 \\
    $\neq $ &0&$\neq $ &0&0&$\neq $  \\
    $\neq $ &0&$\neq $ &0&0&$\neq $  \\
    $\neq $ &0&0&$\neq $ &$\neq $ &0
 \end{tabular}\end{center}
\bigskip

\subsection{Positive Central Quadratic Forms on Quad\-ratic $\ast$-Algebras}

\noindent
\medskip

The algebra $\cala_\ug$, $\ug \in T$ is by construction a {\sl
quadratic $\ast$-algebra} i.e. a quadratic complex algebra
$\cala=A(V,R)$ which is also a $\ast$-algebra with involution
$x\mapsto x^\ast$ preserving the subspace $V$ of generators.
Equivalently one can take the generators of $\cala$ (spanning $V$)
to be hermitian elements of $\cala$. In such a case the complex
finite-dimensional vector space $V$ has a real structure given by
the antilinear involution $v\mapsto j( v)$ obtained by restriction
of $x\mapsto x^\ast$. Since one has $(xy)^\ast=y^\ast x^\ast$ for
$x,y\in \cala$, it follows that the set $R$ of relations satisfies
\begin{equation}
(j \otimes j)( R)=t(R) \label{eq5.3}
\end{equation}
 in $V\otimes V$ where $t:V\otimes V\rightarrow V\otimes V$ is the transposition $v\otimes w \mapsto t(v\otimes w)=w\otimes v$. This implies

\begin{lem}\label{conj}
 The characteristic variety is stable under the
involution $Z\mapsto j( Z)$  and one has
$$
\sigma( j(Z)) =\, j( \sigma^{-1}(Z))
$$
\end{lem}

\smallskip We let $C$ be as above an invariant component of $E \times E$
we say that $C$ is $j$-{\em real} when it is globally invariant
under the involution
\begin{equation}\label{tildej}
\tilde j(Z,\,Z'):=( j( Z'),\, j( Z))
\end{equation}
By lemma \ref{conj} this involution commutes with the automorphism
$\tilde{\sigma}$ (\ref{siginv}) and
 one has

\medskip

\begin{prop} \label{central3} Let $C$ be a $j$-real
 invariant component of $E \times E$ and $Q$ central on $C$
be such that
\begin{equation}\label{selfad}
\overline{ Q(\tilde j(Z,\,Z'))}=\,Q(Z,Z')\qqq (Z,Z')\in C\,,
\end{equation}
\begin{enumerate}
\item The following turns the
 cross-product $C_Q$ into a $\ast$-algebra,
\begin{equation}
 f^\ast(Z, \,Z'):=\overline{ f(\tilde j(Z,\,Z'))
}\,,\qquad (W_L)^\ast = W'_{j( L)}\,,\qquad (W'_{L'})^{\ast} =W_{j(
L')} \label{invol}
\end{equation}
\item The homomorphism
$\rho$ of lemma \ref{alg0} is a $\ast$-homomorphism.
\end{enumerate}
\end{prop}

\medskip
\begin{proof} We used the transpose of $j$ to define $j( L)$ in \eqref{invol}
by
\begin{equation}
 j(L)(Z)=\,\overline{L(j( Z))}\,, \qquad \forall Z \in V^\ast.
\label{inv2}
\end{equation}
The compatibility of the involution with \eqref{cropro} follows from
the commutation of $\tilde{j}$ with $\tilde{\sigma}$.

Its compatibility with \eqref{crossed} follows from
$$
\pi^\ast(Z,Z'):=\,(\frac{L(j(Z'))\,L'(j(Z))}{Q(\tilde j(Z,\,Z'))})^-
=\,\frac{j(L')(Z)\,j(L)(Z')}{Q(Z,\,Z')}\,.
$$

To check 2) one writes for $Y\in V$,
$$
\rho(Y)^*=\,(W_Y + \,W'_Y)^* =\,W'_{j(Y)}+\,W_{j(Y)}=\,\rho(j(Y)\,.
$$
\end{proof}

\medskip We have treated so far $Z$ and $Z'$ as independent variables.
We shall now restrict the above construction to the graph of $j$ \ie
to $\{(Z,Z')\in C \,\vert\, Z'= j( Z)\}$. Composing $\rho$ with the
restriction to the subset $K= \{Z\,|\,(Z,j(Z))\in C \}$ one obtains
in fact a $\ast$-homomorphism $\theta$ of $\cala=A(V,R)$ to a
twisted cross-product $C^\ast$-algebra, $C(K)
\times_{\sigma,\,\call} \mathbb {Z} $ which involves the full
geometric data $(E,\sigma, \call)$ and encodes the central quadratic
form $Q$ as a Hermitian metric on $\call$ provided $Q$ fulfills the
following {\em positivity}.

\medskip
\begin{defn}  \label{central4}
 Let $C$ be a $j$-real
 invariant component of $E \times E$ and $Q$ central on $C$.
Then $Q$ is positive on $C$ iff  it fullfills \eqref{selfad} and
$$
Q(Z,j(Z))> 0 \qqq Z\in K\,.
$$
\end{defn}

\medskip
One can then endow the line bundle $\call$ dual of the tautological
line bundle on $P(V^*)$ with the Hermitian metric defined by
\begin{equation}\label{herm}
\langle f\,L,\,g\, L'\rangle_Q(Z) =\,f(Z)\,\overline{g(Z)} \,\frac{
L(Z)\,\overline{ L'(Z)}}{Q(Z,\,j(Z))} \qquad L, L' \in V,\quad Z \in
K\, \qqq f,g \in C(K)\,.
\end{equation}

We view $f\,L$ and $g\,L'$ as sections of $\call$ and the right hand
side of the formula as a function on $K$ which expresses their inner
product $\langle f\, L, g\,L'\rangle$. This defines a Hermitian
metric on the restriction of $\call$  to $K$.

\medskip Before we proceed we need to describe the general notion
due to Pimsner \cite{pims:1997} of twisted cross product. Given a
compact space $K$, an homeomorphism $\sigma$ of $K$ and a hermitian
line bundle $\call$ on $K$ we define the  $C^\ast$-algebra $C(K)
\times_{\sigma,\,\call} \mathbb {Z} $ as the twisted cross-product
of $C(K)$ by the Hilbert $C^*$-bimodule associated to $\call$ and
$\sigma$ (\cite{aba-eil-exel:1998}, \cite{pims:1997}).

We let for each $n \geq 0$, $\call^{\sigma^n}$ be the hermitian line
bundle pullback of $\call$ by $\sigma^n$ and (cf.
\cite{art-tat-vdb:1990}, \cite{smi-sta:1992})
\begin{equation}
\call_n := \call \otimes \call^{\sigma} \otimes \cdots \otimes
\call^{\sigma^{n-1}} \label{gene2}
\end{equation}
We first define a $\ast$-algebra as the linear span of the monomials
\begin{equation}
\xi \, W^n\, , \quad W^{\ast n} \, \eta^\ast \,,\quad \xi\,,\eta \in
C(K,\call_n) \label{gene}
\end{equation}
with product given as in (\cite{art-tat-vdb:1990},
\cite{smi-sta:1992}) for $(\xi_1 \, W^{n_1})\,(\xi_2 \, W^{n_2})$ so
that
\begin{equation}
(\xi_1 \, W^{n_1})\,(\xi_2 \, W^{n_2}):= (\xi_1 \otimes
(\xi_2\circ{\sigma^{n_1}}) )\, W^{n_1+n_2} \label{gene3}
\end{equation}
We use the hermitian structure of $\call_n $ to give meaning to the
products $\eta^\ast \,\xi$ and $\xi \;\eta^\ast$ for $\xi\,,\eta \in
C(K,\call_n)$. The product then extends uniquely to an associative
product of $\ast$-algebra fulfilling the following additional rules
\begin{equation}
(W^{\ast k} \, \eta^\ast)\,( \xi \, W^k):= \, (\eta^\ast\, \xi)\circ
\sigma^{-k}\,,\qquad ( \xi \, W^k)\,(W^{\ast k} \, \eta^\ast)\,:= \,
\xi \;\eta^\ast \label{gene1}
\end{equation}

 The $C^\ast$-norm of $C(K) \times_{\sigma,\,\call} \mathbb {Z} $ is defined as for ordinary cross-products and due to the amenability of the group $\mathbb {Z} $
there is no distinction between the reduced and maximal norms. The
latter is obtained as the supremum of the norms in involutive
representations in Hilbert space. The natural positive conditional
expectation on the subalgebra $C(K)$ shows that the $C^\ast$-norm
restricts to the usual sup norm on $C(K)$.

 To lighten notations
in the next statement we  abreviate $j(Z)$ as $\bar Z$, but one
should take care that in general the expression for  $j(Z)$ can
differ from $\bar Z$ for instance with the notations of subsection
\ref{simplerrho} one gets $j(Z)=\epsilon \,\bar Z$.

\medskip

\begin{thm}\label{C*}
 Let $K \subset E$ be a compact
$\sigma$-invariant subset and $Q$ be central and strictly positive
on $\{(Z,\,\bar Z);\, Z\in K\}$. Let $\call$ be the restriction to
$K$ of the dual of the tautological line bundle on $P(V^\ast)$
endowed with the hermitian metric $\langle\;,\; \rangle_Q$.

 (i) The equality $\sqrt{2}\,\theta(Y):= Y\, W + W^\ast\,\bar Y^\ast$
yields a $\ast$-homomorphism $$\theta:\cala=A(V,R) \mapsto C(K)
\times_{\sigma,\,\call} \mathbb {Z} $$

 (ii) For any $Y \in V$
the $C^\ast$-norm of $\theta(Y)$ fulfills
$${\rm Sup}_K \|Y\|\leq \sqrt{2}\| \,\theta(Y)\|
\leq 2\,{\rm Sup}_K \|Y\| $$

 (iii) If $\sigma^4 \neq \bbbone$, then $\theta(Q)= 1$ where
$Q$ is viewed as an element of $T(V)/(R)$.
\end{thm}

\medskip
\begin{proof} (i)
The subset $\tilde{K}=\{(Z,\,j( Z));\, Z\in K\}\subset C$ is
globally invariant under $\tilde \sigma$ by lemma \ref{conj}.
Moreover $\tilde j$ defined in \eqref{tildej} is the identity on
$\tilde{K}$. Each $L\in V$ defines a section of $\call$ and hence an
element $L\,W\in C(K) \times_{\sigma,\,\call} \mathbb {Z} $. The
definition \eqref{herm} of the hermitian structure of $\call$ then
shows that the elements $L\,W$ and $W^*\, j(L')^*$ of $C(K)
\times_{\sigma,\,\call} \mathbb {Z} $ fulfill the same algebraic
rules \eqref{cropro}, \eqref{crossed} as the $W_L$ and $W'_{L'}$
while the involution of $C(K) \times_{\sigma,\,\call} \mathbb {Z} $
is the restriction of the involution of proposition \ref{selfad}.
Thus the conclusion follows from lemma \ref{alg0}.

\smallskip (ii) Since $(YW)(YW)^*=Y^*Y$ the $C^\ast$-norm
of $YW$ is  ${\rm Sup}_K \|Y\| $. It follows that $\sqrt{2}\|
\,\theta(Y)\| \leq 2\,{\rm Sup}_K \|Y\| $. For any complex number
$u$ of modulus one the map $\xi W^n\to u^n\,\xi W^n$ extends to a
$*$-automorphism of $C(K) \times_{\sigma,\,\call} \mathbb {Z} $. It
follows taking $u=i$ that $\|Y\, W - W^\ast\,\bar Y^\ast\|=\,\|Y\, W
+ W^\ast\,\bar Y^\ast\|$ and ${\rm Sup}_K \|Y\|\leq \sqrt{2}\|
\,\theta(Y)\|$.

\smallskip (iii) follows from lemma \ref{alg0}.

\end{proof}

\medskip

 We shall now apply this general result to the
algebras $C_{\mathrm{alg}}({\mathbb R}^4_\varphi)$. We take the
quadratic form
\begin{equation}
 Q(X,\,X'):=\sum X^\mu\,X^{\prime \mu}
\label{quad}
\end{equation}
in the $x$-coordinates, so that $Q$ is the canonical central element
defining the sphere $S^3_\varphi$ by the equation $Q=1$. Proposition
\ref{central} shows that in the generic case \ie for $\varphi \in A
\cup B$, the quadratic form $Q$ is central on $F_\varphi \times
F_\varphi$ with obvious notations. The positivity of $Q$ is
automatic since in the $x$-coordinates the involution $j_\varphi$
 coming from the involution of
the quadratic $\ast$-algebra $C_{\mathrm{alg}}({\mathbb
R}^4_\varphi)$ is simply complex conjugation $j_\varphi(Z)= \bar Z$,
so that $Q(X,\,j_\varphi(X))>0$ for $X \neq 0$. We thus get,

\begin{cor} \label{II}
Let $K \subset F_\varphi $  be a compact $\sigma$-invariant subset.
The homomorphism $\theta$ of Theorem \ref{C*} is a unital
 $\ast$-homomorphism from $ C_{\mathrm{alg}}(S^3_\varphi)$
to the cross-product $ C^{\infty}(K) \times_{\sigma,\,\call} \mathbb
{Z} $.
\end{cor}

\smallskip This applies in particular to $K=F_\varphi$. It follows  that one obtains a
non-trivial $C^\ast$-algebra $C^\ast(S^3_\varphi)$ as the completion
of  $ C_{\mathrm{alg}}(S^3_\varphi)$ for the semi-norm,
\begin{equation}
 \| P \|:= { \rm Sup}\| \,\pi(P) \|
\label{norm}
\end{equation}
where $\pi$ varies through all unitary representations of $
C_{\mathrm{alg}}(S^3_\varphi)$. It was clear from the start that
\eqref{norm} defines a finite $C^\ast$-semi-norm on $
C_{\mathrm{alg}}(S^3_\varphi)$ since the equation of the sphere
$\sum (x^\mu)^2=1$ together with the self-adjointness
$x^\mu=\,x^{\mu\ast}$ show that in any unitary representation one
has
$$
\|\, \pi(x^\mu) \|\leq 1\qqq \mu\,.
$$
What the above corollary gives is a lower bound for the
$C^\ast$-norm such as that given by statement (ii) of Theorem
\ref{C*} on the linear subspace $V$ of generators.

\smallskip To analyse the compact
$\sigma$-invariant subsets of $F_\varphi $
 for generic $\varphi$, we distinguish the
{\em even} case which
 corresponds to all $s_k$ having the
same sign (\cf Figure \ref{ellcurve}) (and holds for instance for
$\varphi \in A$) from the {\em odd} case when all $s_k$ dont have
the same sign. First note that in all cases the real curve
$F_\varphi \cap P_3(\mathbb {R})$ is non empty (it contains $p_0$),
and has two connected components since it is invariant under the
Klein group $H$ (\ref{eq6.17}).

\bigskip

\begin{figure}
\begin{center}
\includegraphics[scale=0.9]{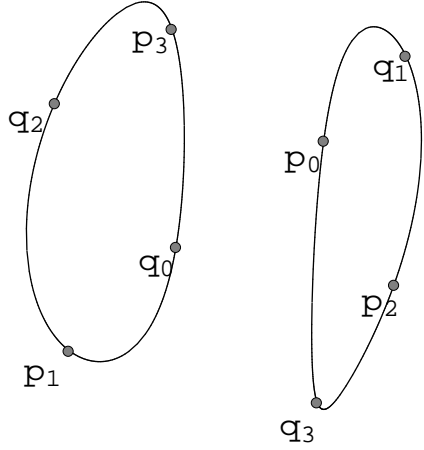}
\end{center}
\caption{\label{ellcurve8} The Elliptic Curve $F_\varphi \cap
P_3(\mathbb {R})$ (odd case )  }
\end{figure}

In the even case $\sigma$ preserves each of the two connected
components of the real curve $F_\varphi \cap P_3(\mathbb {R})$. In
the odd case it permutes them (cf. Figure \ref{ellcurve8}).

\begin{prop} \label{siginvar} Let $\varphi $ be generic and even.

\smallskip (i) Each connected component of $F_\varphi \cap P_3(\mathbb {R})$
is a minimal compact $\sigma$-invariant subset.

\smallskip (ii) Let $K \subset F_\varphi$ be a compact
$\sigma$-invariant subset, then $K$ is the sum in the elliptic curve
$F_\varphi$ with origin $p_0$ of $K_{\mathbb T}=K \cap F_{\mathbb
T}(\varphi)^0 $ (cf. \ref{real}) with the component $ C_\varphi$ of
$F_\varphi \cap P_3(\mathbb {R})$
 containing $p_0$.

\smallskip (iii) The cross-product $ C(F_\varphi) \times_{\sigma,\,\call} \mathbb {Z} $
is isomorphic to the mapping torus of the automorphism $\beta $ of
the noncommutative torus ${\mathbb T}_{\eta}^2 = C_\varphi
\times_\sigma \mathbb {Z} $ acting on the generators by the  matrix
$\left[
\begin{array}{cc}
1& 4\\
0& 1
\end{array}
\right] $.
\end{prop}

\medskip
\begin{proof} (i) This holds if we assume that $\varphi  $
is ``generic" so that the elliptic parameter $\eta$ fullfills $\eta
\notin \mathbb Q$. The diophantine approximations of $\eta$ will
play an important role later on.

\smallskip
(ii) By construction the abelian compact group $F_\varphi$ is the
product ${\mathbb T}_1 \times {\mathbb T}_2$ of the one-dimensional
tori ${\mathbb T}_1=\, C_\varphi$ and ${\mathbb T}_2=\,F_{\mathbb
T}(\varphi)^0$, \ie the component of $F_{\mathbb T}(\varphi)$
containing $q_0$ (cf. \ref{real}). The translation $\sigma$ is
$\eta\times {\rm Id}$ and the action of $\eta$ is minimal on
${\mathbb T}_1= C_\varphi$.

\smallskip (iii)
The isomorphism class of the cross-product $ C(F_\varphi)
\times_{\sigma,\,\call} \mathbb {Z} $ depends of $\call$ only
through its class as a hermitian line bundle on the two torus
$F_\varphi$. This class is entirely specified by the first Chern
class $c_1(\call)$. By construction one gets $c_1(\call)=4$ since
the space of holomorphic sections of $\call$ is the $4$-dimensional
space $V$.

Let $U_j$ be the generators of the algebra $C({\mathbb T}_{\eta}^2)$
where the presentation of the algebra is
$$
U_1\,U_2=\,e^{2 \pi i \eta}\,U_2\,U_1\,.
$$
For any integer $k$ let
 $\beta_k
$ be the automorphism of $C({\mathbb T}_{\eta}^2)$ acting on the
generators $U_j$ by
\begin{equation}
\beta_k(U_1):= U_1\,,\qquad \beta_k(U_2):= U_1^k\,U_2\,.
\end{equation}
By construction the mapping torus $T(\beta_k)$
 of
the automorphism $\beta_k $ is given by the algebra $C(T(\beta_k))$
of continuous maps $s \in \mathbb R \, \mapsto x(s)\in C({\mathbb
T}_{\eta}^2)$ such that $ x(s+1)=\beta_k( x(s))\,,\; \forall s\in
\mathbb R $. We just need to show that $C(T(\beta_4))$ is isomorphic
to $ C(F_\varphi) \times_{\sigma,\,\call} \mathbb {Z} $ and this
follows from the general isomorphism
\begin{equation}\label{maptor}
C(T(\beta_k))\simeq C({\mathbb T}_1 \times {\mathbb T}_2)
 \times_{\eta \times {\rm Id}\, ,\,\call} \mathbb {Z}\,,
\end{equation}
(with $T_j = {\mathbb R}/{\mathbb Z}$) for any hermitian line bundle
$\call$ on ${\mathbb T}_1 \times {\mathbb T}_2$ with $c_1(\call)=k$.
To check this one chooses $\call$ so that its continuous sections
$C({\mathbb T}_1 \times {\mathbb T}_2,\call)$
 are scalar functions $f(u,m)$
with $u, \,m \in {\mathbb R}$ such that
$$
f(u+1,m)=\,f(u,m)\,,\quad f(u,m+1)=\,e^{2 \pi i k\,u}\,f(u,m)\qqq
u,\,m \in {\mathbb R}\,.
$$
while its hermitian metric is given by
$$
\langle f,g\rangle(u,m)=\,f(u,m)\,\overline{g(u,m)}\qqq u,\,m \in
{\mathbb R}\,.
$$
One defines a map
$$
\alpha\;:\quad C(T(\beta_k))\to C({\mathbb T}_1 \times {\mathbb
T}_2)
 \times_{\eta \times {\rm Id},\,\call} \mathbb {Z}\,,
$$
by writing  for $x \in C(T(\beta_k))$, $x=(x(s))$, $x(s)\in
C({\mathbb T}_{\eta}^2)$ the Fourier expansion
$$
x(s)=\,\sum x(s,n)\,U_2^n\,, \quad x(s,n)\in C({\mathbb T}_1)\,.
$$
Then the $x(s,n)\in C({\mathbb T}_1)$ define sections
$$
x_n\in \,C({\mathbb T}_1 \times {\mathbb T}_2,\call_n)
$$
and one just lets
$$
 \alpha(x) =\, \sum  x_n\,W^n \in
C({\mathbb T}_1 \times {\mathbb T}_2)
 \times_{\eta \times {\rm Id} ,\,\call} \mathbb {Z}
  \qqq x\in C(T(\beta_k))\,.
$$
One then checks that this gives the required isomorphism
\eqref{maptor}.
\end{proof}

\begin{cor} \label{H3}
Let $\varphi  $ be generic and even, then $F_\varphi
\times_{\sigma,\,\call} \mathbb {Z} $ is a noncommutative
$3$-manifold with an elliptic
 action of the three dimensional
Heisenberg Lie algebra $\frach_3$ and an invariant trace $\tau$.
\end{cor}

\begin{proof} This follows \footnote{It justifies the terminology
``nilmanifold"}  from proposition \ref{siginvar} (iii). One can
construct directly the action of $\frach_3$ on
$C^{\infty}(F_\varphi) \times_{\sigma,\,\call} \mathbb {Z}$ by
choosing a constant (translation invariant) curvature connection
$\nabla$, compatible with the metric, on the hermitian line bundle
$\call$ on $F_\varphi $ (viewed in the $C^{\infty}$-category not in
the holomorphic one). The two covariant differentials $\nabla_j$
corresponding to the two vector fields $X_j$ on $F_\varphi $
generating the translations of the elliptic curve, give a natural
extension of $X_j$ as the unique derivations $\delta_j$ of
$C^{\infty}(F_\varphi) \times_{\sigma,\,\call} \mathbb {Z}$
fulfilling the rules,
\begin{eqnarray}
\delta_j(f)& = & X_j(f)\,,\quad \forall f \in C^{\infty}(F_\varphi)\nonumber\\
\delta_j(\xi \, W)& = & \nabla_j(\xi) \, W\,,\quad \forall \xi \in
C^{\infty}(F_\varphi\,,\call)
\end{eqnarray}
We let  $\delta$ be the unique derivation of $C^{\infty}(F_\varphi)
\times_{\sigma,\,\call} \mathbb {Z}$ corresponding to the grading by
powers of $W$.
 It vanishes on $C^{\infty}(F_\varphi)$ and
fulfills
\begin{equation} \label{delta}
\delta (\xi \, W^k)= i\,k\,\xi \, W^k \qquad \delta(W^{\ast k} \,
\eta^\ast)= -i\,k\,W^{\ast k} \, \eta^\ast
 \end{equation}
Let $i\,\kappa$ be the constant curvature of the connection
$\nabla$, one gets
\begin{equation}
[\delta_1,\, \delta_2]=\,\kappa\, \delta \,,\quad [\delta,\,
\delta_j]=\,0
\end{equation}
which provides the required action of the Lie algebra $\frach_3$ on
$C^{\infty}(F_\varphi) \times_{\sigma,\,\call} \mathbb {Z}$.
\end{proof}

It follows that one is exactly in the framework developped in
\cite{ac:1980}. We refer to \cite{rieff:1989} and
\cite{aba-exel:1997} where these noncommutative manifolds were
analysed in terms of crossed products by Hilbert $C^*$-bimodules.

 Integration on the translation invariant volume form $dv$
of $F_\varphi$ gives the $\frach_3$-invariant trace $\tau$,
\begin{eqnarray}
 \label{trace}
\tau(f)& = & \int f dv\,,\quad \forall f \in C^{\infty}(F_\varphi)\nonumber\\
\tau(\xi \, W^k)& = &\tau(W^{\ast k} \, \eta^\ast)\,=\,0\,,\quad
\forall k\neq 0
 \end{eqnarray}
It follows in particular that the results of \cite{ac:1980} apply to
obtain the calculus. In particular the following gives the
``fundamental class" as a $3$-cyclic cocycle,
\begin{equation}\label{3trace}
\tau_3(a_0,\,a_1,\,a_2 ,\,a_3)=\,\sum
\epsilon_{ijk}\,\tau(a_0\,\delta_i(a_1)\,\delta_j(a_2)
\,\delta_k(a_3))
\end{equation}
where the $\delta_j$ are the above derivations with
$\delta_3:=\delta$.

 We shall in fact describe the
same calculus in greater generality in the last section which will
be devoted to the computation of the Jacobian of the homomorphism
$\theta$ of corollary \ref{II}.

 Similar results hold in the odd case. Then $F_\varphi \cap P_3(\mathbb {R})$
is a minimal compact $\sigma$-invariant subset, any  compact
$\sigma$-invariant subset $K \subset F_\varphi$ is the sum in the
elliptic curve  $F_\varphi$ with origin $p_0$ of $F_\varphi \cap
P_3(\mathbb {R})$ with $K_{\mathbb T}=K \cap F_{\mathbb
T}(\varphi)^0 $ but the latter is automatically invariant under the
subgroup $H_0 \subset H$ of order $2$ of the Klein group $H$
(\ref{eq6.17})
\begin{equation}
H_0 :=\{ h \in H \vert \, h (F_{\mathbb T}(\varphi)^0 )=F_{\mathbb
T}(\varphi)^0\}
\end{equation}

\bigskip

  The group law in $F_\varphi$ is described geometrically as
follows. It involves the point $q_0$. The sum $z =x + y$ of two
points $x$ and $y$ of $F_\varphi$ is $z=I_0(w)$ where $w$ is the 4th
point of intersection of
 $F_\varphi$ with the plane determined by the three points $\{q_0, x, y\}$.
It commutes by construction with complex conjugation so that
$\overline{x+y}=\, \overline{ x}\,+\, \bar y \,,\quad \forall x,\,y
\in F_\varphi$.

  By lemma \ref{conj} the translation $\sigma$ is
imaginary for the canonical  involution $j_\varphi$. In terms of the
coordinates $Z_\mu$ this involution is  described as follows, using
(\ref{change0}) (multiplied by $e^{i (\pi/4
-\varphi_1-\varphi_2-\varphi_3)}2^{-3/2}$) to change variables.
Among the $3$ real numbers
\begin{eqnarray}
\lambda_k & = &\, \cos\varphi_\ell \,
\cos\varphi_m\,\sin(\varphi_\ell-\varphi_m)\,,\qquad k\in
\{1,2,3\}\nonumber
\end{eqnarray}
 two  have the same sign $\epsilon$ and one, $\lambda_k$,
$k\in \{1,2,3\}$, the opposite sign. Then
\begin{equation}
 j_\varphi=\,\epsilon\,I_k \,\circ \, c
\label{involE}
\end{equation}
where $c$ is  complex conjugation on the real elliptic curve
$F_\varphi$ (section 3) and  $I_\mu$  the involution
\begin{equation}
I_\mu(Z_\mu)=-Z_\mu, \, \, I_\mu(Z_\nu)=Z_\nu\,,\quad \nu \neq \mu
\label{imu}
\end{equation}
The index $k$ and the sign $\epsilon$ remain constant when $\varphi$
varies in each of the four components of the complement of the four
points $q_\mu$ in $F_{\mathbb T}(\varphi)$. The sign $\epsilon$
matters for the action of $j_\varphi$ on linear forms as in
(\ref{inv2}), but is irrelevant for the action
 on $F_\varphi$. Each involution $I_\mu$ is a symmetry
$z \mapsto p-z$ in the elliptic curve $F_\varphi$ and the products
$I_\mu \circ I_\nu$ form the Klein subgroup $H$ (\ref{eq6.17})
acting by translations of order two on  $F_\varphi$.

 The quadratic form $Q$ of (\ref{quad})
 is given in the new coordinates by,
\begin{equation} \label{q13}
Q=\;(\prod \cos^2\varphi_{\ell}) \, \sum \, t_k\, s_k \,\,Q_k \,
\end{equation}
 with $s_k := 1 + t_\ell\, t_m
\, , \: t_k:= {\rm tan}\,\varphi_k$
 and $Q_k$ defined by (\ref{qk}).

Let us assume that $0<\varphi_1<\varphi_2<\varphi_3<\pi/2$ for
instance, then the relation between the $x^\mu$ and the $Y_\mu$ is
given with the appropriate normalization of the $Y_\mu$ by
$x^\mu=\,\rho_\mu\,Y_\mu$ where
\begin{eqnarray} \label{xtoZ}
\quad\rho_0^2
&=&-\sin(\varphi_1-\varphi_2)\,\sin(\varphi_1-\varphi_3)\,\sin(\varphi_2-\varphi_3)\,,\quad
\rho_1^2 =\,
-\cos \varphi_2\, \cos \varphi_3\,\sin(\varphi_2-\varphi_3)\,,\\
\rho_2^2 &=&\cos \varphi_1\, \cos
\varphi_3\,\sin(\varphi_1-\varphi_3)\,,\quad \rho_3^2=\,-\cos
\varphi_1\, \cos \varphi_2\,\sin(\varphi_1-\varphi_2))\,.\nonumber
\end{eqnarray}

All the $\rho_\mu$ are real except for $\rho_2$ which is purely
imaginary and the  involution $j$ is $I_2\circ c$. One checks
directly that
$$
\sum \rho_\mu^2\,Y_\mu^2=\,(\prod \cos^2\varphi_{\ell}) \, \sum \,
t_k\, s_k \,\,Q_k\,.
$$

\medskip
We can now compare the $*$-homomorphism $\tilde \rho$ of section
\ref{simplerrho} with the $*$-homomorphism obtained from a positive
central quadratic form, one gets with the constants $s$ and $b$
given by,
$$
s=-s_2\,\prod\sin \varphi_j \, ,\quad b^{2}=\prod \cos \varphi_j\,
\cos(\varphi_k-\varphi_\ell)\,.$$

\medskip

\begin{prop} \label{ztoy}
Let  $0<\varphi_1<\varphi_2<\varphi_3<\pi/2$.

(i) The $*$-homomorphism $\tilde \rho$ is the $*$-homomorphism
associated to the central quadratic form $Q'$,
$$
s\,Q'=\,Q_1+\,Q_3+\,s_2\,Q_2\,,
$$
which is positive on $E$ for the involution $I_3\circ c$.

(ii) Let
$$
\beta(Y_0)=\,i\,Y_2\,,\quad \beta(Y_1)=\,\sqrt{J_{12}} \,Y_3
\,,\quad \beta(Y_2)=-\,\sqrt{J_{12}}\sqrt{J_{23}} \,Y_0 \,,\quad
\beta(Y_3)=\,i\,\sqrt{J_{23}} \,Y_1\,,
$$
the map $b\,\beta$ gives a $*$-isomorphism sending the form $Q'$
into $Q$ and the involution $I_3\circ c$ into $I_2\circ c$.
\end{prop}

\medskip

\begin{proof} (i) By proposition \ref{simplerhofinal} it is enough to
show that $Q'$ corresponds to $\sum\,(x^\mu)^2$ under the
composition of the transformations $S_\mu=\,\lambda_\mu\,x^\mu$ of
lemma \ref{rescalelem} and \eqref{simplerho0}
$$
S_0=\,d\,Y_2\,,\quad S_1=\,i\,Y_3\,,\quad S_2=\,d\,Y_0\,,\quad
S_3=\,-Y_1\,,
$$
where $d^2=-J_{31}$.

(ii) The map $b\,\beta$ is obtained as the composition of the
isomorphisms \eqref{simplerho0}, $S_\mu=\,\lambda_\mu\,x^\mu$ of
lemma \ref{rescalelem} and $x^\mu=\,\rho_\mu\,Y_\mu$ with $\rho_\mu$
given in \eqref{xtoZ}. Thus the answer follows since each of these
maps is a $*$-isomorphism and the image of $Q'$ in the $\,x^\mu$
variables is simply $\sum\,(x^\mu)^2$.

\end{proof}

 Let $\varphi  $ be generic and even
and  $v \in F_{\mathbb T}(\varphi)^0$. Let $K(v)= \,v +  C_\varphi $
be
 the minimal compact $\sigma$-invariant subset containing $v$
(Proposition \ref{siginvar} (ii)). By Corollary \ref{II} we get a
homomorphism,
\begin{equation} \label{ncu}
\theta_v \;:\; C_{\mathrm{alg}}(S^3_\varphi) \mapsto
C^\infty({\mathbb T}_{\eta}^2)
\end{equation}
whose non-triviality will be proved below in corollary \ref{nct}. We
shall first show (Theorem \ref{gal}) that it transits through the
cross-product of the field $K_q$ of meromorphic functions on the
elliptic curve by the subgroup of its Galois group Aut$_{\mathbb
C}(K_q)$ generated by the translation $\sigma$.

\medskip
For $Z = \,v + z\,$, $ z \in C_\varphi$, one has using
(\ref{involE}) and (\ref{real}),
 \begin{equation}
 j_\varphi(Z)
= \,I_\mu(Z-v)-\,I(v) \label{invol2}
\end{equation}
which is a rational function $r(v,\,Z)$. Fixing $\varphi,\,v$ and
substituting $Z$ and $Z'=r(v,\,Z)$ in the  formulas (\ref{crossed})
and (\ref{ratio}) of lemma \ref{alg0} with $L$ real such that  $0
\notin L(K(v))$, $L'=\,\epsilon\,L \circ I_\mu$ and $Q$ given by
(\ref{q13}) we obtain rational formulas for a homomorphism
$\tilde{\theta_v}$ of $ C_{\mathrm{alg}}(S^3_\varphi)$ to the
generalised cross-product of the field $K_q$ of meromorphic
functions $f(Z)$ on the elliptic curve $F_\varphi$ by $\sigma$. The
generalised  cross-product rule (\ref{crossed}) is given by $W \,
W':= \gamma(Z)$ where $\gamma$ is a rational function. Similarly $W'
\, W:=\gamma(\sigma^{-1}(Z))$. Using integration on the cycle $K(v)$
to obtain a trace, together with corollary \ref{II}, we  get,

\begin{thm} \label{gal}
 The homomorphism $\theta_v : C_{\mathrm{alg}}(S^3_\varphi) \mapsto
C^\infty({\mathbb T}_{\eta}^2)$ factorises with a homomorphism
$\tilde{\theta_v}: C_{\mathrm{alg}}(S^3_\varphi) \mapsto K_q
\times_\sigma \mathbb {Z}$
 to  the  generalised cross-product of the field $K_q$
of meromorphic functions on the elliptic curve $F_\varphi$
 by the subgroup of the Galois group Aut$_{\mathbb C}(K_q)$
generated by $\sigma$. Its image generates the hyperfinite factor of
type II$_1$ after weak closure relative to the trace given by
integration on the cycle $K(v)$.
\end{thm}

\medskip
 Elements of $K_q$
with poles on $K(v)$ are unbounded and give elements of the regular
ring of affiliated operators, but all elements of $\theta_v (
C_{\mathrm{alg}}(S^3_\varphi))$ are regular on $K(v)$. The above
generalisation of the cross-product rules (\ref{crossed}) with the
rational formula for $W \, W':= \gamma(Z)$  is similar to the
introduction of $2$-cocycles in the standard Brauer theory of
central simple algebras.

\section{The Jacobian of the Covering of $ S^3_\varphi$}\label{jacobian}

\smallskip In this section we shall analyse the morphism of $\ast$-algebras
\begin{equation}
\theta \,:\, C_{\mathrm{alg}}(S^3_\varphi)\mapsto C^\infty(F_\varphi
\times_{\sigma,\,\call} \mathbb {Z})
\end{equation}
of Corollary \ref{II}, by computing its Jacobian in the sense of
noncommutative differential geometry (\cite{ac:1982}). We postpone
the analysis at the $C^\ast$-level, in particular the role of the
discrete series, to another forthcoming publication.

\smallskip The usual Jacobian  of a smooth map
$\varphi :\, M \mapsto N$ of manifolds is obtained as the
 ratio $\displaystyle \varphi^\ast(\,\omega_N)/\omega_M$ of the
pullback of the volume form $\omega_N$ of the target manifold $N$
with the volume form $\omega_M$ of the source manifold $M$. In
noncommutative  geometry, differential forms $\omega$ of degree $d$
become Hochschild classes $\;\tilde{\omega} \in HH_d(\cala)\,,
\;\cala= C^\infty(M)$. Moreover since one works  with the dual
formulation in terms of algebras, the pullback
$\varphi^\ast(\omega_N)$ is replaced by the pushforward
$\varphi^{\,t}_\ast(\tilde{\omega}_N)$ under the corresponding
transposed morphism of algebras $\;\varphi^{\,t}(f):= f \circ
\varphi\,, \quad \forall f \in C^\infty(N)$.

\smallskip
The noncommutative sphere $S^3_\varphi$ admits a canonical ``volume
form" given by the Hochschild  $3$-cycle $\ch_{\frac{3}{2}}(U)$. Our
goal is to compute the push-forward,
\begin{equation}\label{goal}
\theta_\ast (\ch_{\frac{3}{2}}(U))\,\in HH_3(C^\infty(F_\varphi
\times_{\sigma,\,\call} \mathbb {Z}))
\end{equation}

\smallskip
 Let $\varphi$ be generic and even. The noncommutative manifold
 $F_\varphi \times_{\sigma,\,\call} \mathbb {Z}$ is, by Corollary \ref{H3},
a noncommutative $3$-dimensional nilmanifold isomorphic to the
mapping torus of an automorphism of the noncommutative $2$-torus
$T^2_\eta$. Its Hochschild homology is easily computed using the
corresponding result for the noncommutative torus (\cite{ac:1982}).
It admits in particular a canonical volume form $V\in
HH_3(C^\infty(F_\varphi \times_{\sigma,\,\call} \mathbb {Z})) $
which corresponds to the natural class in $HH_2(C^\infty(T^2_\eta)
)$ (\cite{ac:1982}). The volume form $V$ is obtained in the
cross-product
 $F_\varphi \times_{\sigma,\,\call} \mathbb {Z}$
from the translation invariant $2$-form $dv$ on $F_\varphi$.

\smallskip To compare $\,\,\theta_\ast (\ch_{\frac{3}{2}}(U))$ with $V$
we shall pair it with the  $3$-dimensional Hochschild
 cocycle $\tau_h\in HH^3(C^\infty(F_\varphi \times_{\sigma,\,\call} \mathbb {Z})) $
given, for any element $h$ of the center of $C^\infty(F_\varphi
\times_{\sigma,\,\call} \mathbb {Z})$, by
\begin{equation}\label{tauh}
\tau_h(a_0,\,a_1,\,a_2 ,\,a_3)=\tau_3(h\,a_0,\,a_1,\,a_2 ,\,a_3)
\end{equation}
where $\tau_3 \in HC^3(C^\infty(F_\varphi \times_{\sigma,\,\call}
\mathbb {Z})) $ is the fundamental class in cyclic cohomology
defined by (\ref{3trace}).

\smallskip By proposition \ref{chbare} the component $\ch_{\frac{3}{2}}(U)$ of the Chern
character  is given by,
\begin{eqnarray}
\ch_{\frac{3}{2}}(U) = &  &\sum
\epsilon_{\alpha\beta\gamma\delta}\cos(\varphi_{\alpha}-\varphi_{\beta}+\varphi_{\gamma}-\varphi_{\delta})\:\:x^\alpha\,dx^\beta\,dx^\gamma\,dx^\delta
\, - \nonumber
\\&i& \, \sum \sin 2(\varphi_\mu-\varphi_\nu)\:\:x^\mu\,dx^\nu\,dx^\mu\,dx^\nu
\end{eqnarray}
where $\varphi_0 :=0$. In terms of the $Y_\mu$ one gets,
\begin{eqnarray} \label{beauty}
\ch_{\frac{3}{2}}(U)& = & \kappa\,\sum
\delta_{\alpha\beta\gamma\delta}\,(s_{\alpha}-s_{\beta}+s_{\gamma}-s_{\delta})\:\:Y_\alpha\,dY_\beta\,dY_\gamma\,dY_\delta
\, + \nonumber
\\&& \,\kappa\, \sum \epsilon_{\alpha\beta\gamma\delta}\;(s_\alpha\,-s_\beta)\:\:\,Y_\gamma\,dY_\delta\,dY_\gamma\,dY_\delta
\end{eqnarray}
where $s_0:=0$, $s_k := 1 + t_\ell\, t_m \, , \: t_k:= {\rm
tan}\,\varphi_k$ and
\begin{equation}
\delta_{\alpha\beta\gamma\delta}
=\epsilon_{\alpha\beta\gamma\delta}\,(n_{\alpha}-n_{\beta}+n_{\gamma}-n_{\delta})
\end{equation}
with $n_0=0$ and $n_k=1$. The normalization factor is
\begin{equation}
\kappa =\,i \;\prod \cos^2(\varphi_{k}) \,
\sin(\varphi_{\ell}-\varphi_{m})
\end{equation}
Formula (\ref{beauty}) shows that, up to normalization,
$\ch_{\frac{3}{2}}(U) $ only depends on the fiber $F_\varphi$ of
$\varphi$.

\smallskip Let $\varphi $ be generic and  even, there is a similar formula
in the odd case. In our case the involutions $I$ and $I_0$ are
conjugate by a real translation $\kappa$ of the elliptic curve
$F_\varphi$ and we let $F_\varphi(0)$ be one of the two connected
components of,
\begin{equation}\label{hh}
\{ Z \in F_\varphi \,\vert \,I_0(Z)= \bar Z \}
\end{equation}
By Proposition \ref{siginvar} we can identify the center of
$C^\infty(F_\varphi \times_{\sigma,\,\call} \mathbb {Z}) $ with
$C^\infty(F_\varphi(0))$. We assume for simplicity that $\varphi_j
\in [0,\frac{\pi}{2}]$ are in cyclic order $\varphi_k < \varphi_l
<\varphi_m$ for some
 $k \in\{1,2,3\}$.

\medskip

\begin{thm} \label{vol}
Let $h \in$ Center $(C^\infty(F_\varphi \times_{\sigma,\,\call}
\mathbb {Z}) )\sim C^\infty(F_\varphi(0))$. Then
$$
\langle \ch_{\frac{3}{2}}(U), \tau_h \rangle=\,6 \,\pi\,
\Omega\,\int_{F_\varphi(0)}h(Z)\, dR(Z)
$$
where $\Omega$ is the period given by the  elliptic integral of the
first kind,
$$
\Omega=\,\int_{C_\varphi} \, \frac{Z_k dZ_0- Z_0 dZ_k}{Z_\ell Z_m}
$$
and  $R$ the rational fraction,
$$
R(Z)=\, t_k\;\frac{Z_m^2}{Z_m^2+  \,c_k\,Z_l^2}
$$
 with $c_k={\rm tg}(\varphi_l)\;{\rm cot}(\varphi_k-\varphi_\ell)$.
\end{thm}

\medskip We assume that $0<\varphi_1<\varphi_2<\varphi_3<\pi/2$
\ie that $k=1$. We start from corollary \ref{rhoch} and express the
result
\begin{equation}\label{start}
\omega=\,-\,\frac{\,\sigma(m)^{4}\,g(m)}{\lambda \,\Lambda}\,dm\,,
\end{equation}
in terms of the trigonometric parameters $\varphi_j$ and the
coordinates $Z_\mu$ of $Z$.

One has (\cf \eqref{scalingfactor})
\begin{equation}\label{scalingfactor1}
\sigma(m)^4=\,
(\prod\,\sin\,\varphi_j)^{2}\;(C_1-\lambda\,C_2)^{-2}\,,
\end{equation}
and we begin by giving a better formula for $C_1-\lambda\,C_2$.

\begin{lem} One has
\begin{equation}\label{c1-c2}
\,C_1-\lambda \,C_2=\; b_1\,\vartheta^2_1(im)+\,
 b_2\,\vartheta^2_2(im)\,,
\end{equation}
where
\begin{equation}\label{c1c2}
 b_1=\,4\;\frac{\vartheta^2_1(\eta)}{s_1\;\vartheta^2_2(\eta)}\,,\quad
 b_2=\,4\;\frac{s_1\,-1}{s_1}\,.
\end{equation}
\end{lem}

\begin{proof}  One uses
 \eqref{casimir}, $\displaystyle \lambda=\,\frac{\vartheta^2_2(0)}{s_1\;\vartheta^2_2(\eta)}$
and the identity
$$
\vartheta^2_2(0)\,\vartheta_2(\eta+i m)\,\vartheta_2(\eta-i m)
=\,\vartheta^2_2(i m)\,\vartheta^2_2(\eta)-\,\vartheta^2_1(i
m)\,\vartheta^2_1(\eta)\,.
$$

\end{proof}

The $m$-dependent terms are understood from the following lemma

\begin{lem} Let $ b_j$ be arbitrary constants, then
\begin{equation}\label{derivative}
\vartheta^2_3(0)\vartheta^2_4(0)\,
\frac{d}{du}\;\frac{\vartheta^2_1(u)}{ b_1\,\vartheta^2_1(u)+\,
 b_2\,\vartheta^2_2(u)}=\, b_2\, \frac{\vartheta'_1(0)^3}{\pi^2}\;
\frac{\vartheta_1(2u)}{( b_1\,\vartheta^2_1(u)+\,
 b_2\,\vartheta^2_2(u))^2}\,.
\end{equation}
\end{lem}

\begin{proof}
This follows from the classical identity
\begin{equation} \label{dsn}
\frac{d}{du}\;{\rm sn}(u)=\,{\rm cn}(u)\,{\rm dn}(u)\,,
\end{equation}
which implies
$$
\frac{d}{du}\;\frac{\vartheta^2_1(u)}{\vartheta^2_2(u)}=\,2\,\pi\,\vartheta^2_2(0)\,
\frac{\vartheta_1(u)\vartheta_2(u)\vartheta_3(u)\vartheta_4(u)}{\vartheta^4_2(u)}
=\,\frac{\vartheta'_1(0)^3}{\pi^2\,\vartheta^2_3(0)\vartheta^2_4(0)}
\,\frac{\vartheta_1(2u)}{\vartheta^4_2(u)}\,,
$$
using the duplication formula
$$
\vartheta_2(0)\vartheta_3(0)\vartheta_4(0)\vartheta_1(2\,u)=\,2\,
\vartheta_1(u)\vartheta_2(u)\vartheta_3(u)\vartheta_4(u)
$$
and the Jacobi derivative formula
$$
\frac{\vartheta'_1(0)}{\pi}=\,\vartheta_2(0)\vartheta_3(0)\vartheta_4(0)
$$
\end{proof}

Taking $u= im$ in \eqref{derivative} this  allows to write $\omega$
as the differential of
\begin{equation}\label{derivative1}
R(m)=\;\frac{c\;\vartheta^2_1(im)}{ b_1\,\vartheta^2_1(im)+\,
 b_2\,\vartheta^2_2(im)}
\end{equation}
using  \eqref{gdem}
\begin{equation} \label{gdem1}
g(m)=\, 24 \,(2\pi i)^3
\;\frac{\vartheta'_1(0)^3}{\pi^3}\;\frac{\vartheta_1(\eta)\,\vartheta_1(2im)}{
\vartheta_2(\eta)\vartheta_3(\eta)\vartheta_4(\eta)}\,\,.
\end{equation}
The constant $c$ is uniquely  determined and will be simplified

later, we get so far,
\begin{equation}\label{scalingfactor2}
c=\,24\,(2\pi)^3\,\frac{(\prod\,\sin\,\varphi_j)^{2}\,\vartheta_1(\eta)
\,\vartheta^2_3(0)\vartheta^2_4(0)}{\pi\, b_2\,\lambda \,\Lambda\,
\vartheta_2(\eta)\vartheta_3(\eta)\vartheta_4(\eta)}\,.
\end{equation}

\begin{lem}\label{uniqueform}

(i)  The differential form
\begin{equation}
 \chi:=\frac{Z_k dZ_0- Z_0 dZ_k}{s_k\,Z_\ell Z_m}
\end{equation}
is independent of $k$ and is, up to scale, the only holomorphic form
of type $(1,0)$ on $F_\varphi$.

(ii) One has
\begin{equation}\label{period}
2\,\pi\,\vartheta^2_4(0)\,
=\,\frac{\vartheta_1(\eta)\,\vartheta_4(\eta)}{\vartheta_2(\eta)\,\vartheta_3(\eta)}
\int_{C_\varphi} \, \frac{Z_3 dZ_0-Z_0 dZ_3 }{Z_1 Z_2}
\end{equation}
\end{lem}

\medskip
\begin{proof}  (i) Recall that the equations
defining $F_\varphi$ are
\begin{equation}\label{defeqell}
\frac{Z_0^2-Z_1^2}{
s_1}=\frac{Z_0^2-Z_2^2}{s_2}=\frac{Z_0^2-Z_3^2}{s_3}\,.
\end{equation}
One gets the required independence by differentiation.

\smallskip

(ii) Let us check  \eqref{period}. The factor $\displaystyle
\frac{\vartheta_1(\eta)\,\vartheta_4(\eta)}{\vartheta_2(\eta)\,\vartheta_3(\eta)}$
allows to replace the $Z_j$ by $\displaystyle \vartheta_{j+1}(2z)$
so that the right hand side gives using \eqref{dsn}
$$
\int_0^1 \, \frac{\vartheta_{4}(2z) d\vartheta_{1}(2z)
-\vartheta_{1}(2z) d\vartheta_{4}(2z)
}{\vartheta_{2}(2z)\vartheta_{3}(2z) }
=\,2\,\pi\,\vartheta^2_4(0)\,.
$$
\end{proof}

 In terms of $Z=(Z_j)$ one has
\begin{equation}\label{fracR}
R(m)=\,c\,s_1\,\frac{\vartheta^2_1(\eta)}{4\,\vartheta^2_2(\eta)}
\,\frac{Z_0^2}{a_1\,Z_0^2+\,a_2\,Z_1^2}
\end{equation}
where
\begin{equation}\label{coeffa}
a_1=\, \frac{\vartheta^4_1(\eta)}{\vartheta^4_2(\eta)}\,,\quad
a_2=\,s_1-1\,.
\end{equation}
One can express the coefficient $a_1$ in trigonometric terms using

\begin{lem}
\begin{equation}\label{elltotrig}
a_1=\,\frac{\vartheta^4_1(\eta)}{\vartheta^4_2(\eta)}=\,
 \frac{(s_1-s_2)(s_1-s_3) }{s_2 s_3}
\end{equation}
\end{lem}

\begin{proof} By homogeneity of the right hand side one can replace the
$s_j$ by the $\sigma_j$ of \eqref{thesigmas} One then gets
$$\frac{(s_1-s_2)(s_1-s_3) }{s_2 s_3}=\,-\,J_{12}\,J_{31}
$$
and the result follows from \eqref{sklyaJ}.
\end{proof}

\medskip
We now  simplify the product $\displaystyle
c\,s_1\,\frac{\vartheta^2_1(\eta)}{4\,\vartheta^2_2(\eta)}$
replacing $\displaystyle \vartheta^2_4(0)$ in \eqref{scalingfactor2}
by \eqref{period} and using (\cf \eqref{defofL})
\begin{equation}
\Lambda = \overset{3}{\underset{1}{\prod}} \, ({\rm tan} \,
(\varphi_j)\; \cos (\varphi_k - \varphi_{\ell}))\,,\quad
\lambda=\,\frac{\vartheta^2_3(0)}{s_2\;\vartheta^2_3(\eta)}\,.
\label{defofL1}
\end{equation}
By elementary computations and using once more \eqref{elltotrig} one
gets
\begin{equation}\label{scale1}
c\,s_1\,\frac{\vartheta^2_1(\eta)}{4\,\vartheta^2_2(\eta)}=\,
6\,\pi\,t_1\,a_1\;\,\frac{s_1}{s_3}\, \int_{C_\varphi} \, \frac{Z_3
dZ_0-Z_0 dZ_3 }{Z_1 Z_2}
\end{equation}
We thus obtain so far using the elementary equality
$$
\frac{a_2}{a_1}=\,{\rm cot}(\varphi_1-\varphi_2)\,{\rm
cot}(\varphi_1-\varphi_3)
$$
and  lemma \ref{uniqueform} which allows to eliminate the term
$\displaystyle \frac{s_1}{s_3}$, using the definition of the period
$\Omega$ the following formula for the rational fraction,
\begin{equation}\label{startr}
R(Z)=\, t_1\;\frac{Z_0^2}{Z_0^2+  \,a\,Z_1^2}
\end{equation}
 with $a={\rm cot}(\varphi_1-\varphi_2)\;{\rm cot}(\varphi_1-\varphi_3)$.

\medskip
What we have computed so far is the image of  ${\rm ch}_{3/2}$ under
the $*$-homomorphism $\tilde \rho$. By proposition \ref{ztoy} (i)
this amounts to the image of ${\rm ch}_{3/2}$ under the
$*$-homomorphism associated to the central quadratic form $Q'$. By
proposition \ref{ztoy} (ii) we get the result for $Q$ using the
isomorphism $\beta$ and this gives the following formula for $R(Z)$,
\begin{equation}\label{finalr}
R(Z)=\, t_1\;\frac{Z_3^2}{Z_3^2+  \,c_1\,Z_2^2}
\end{equation}
 with $c_1={\rm tg}(\varphi_2)\;{\rm cot}(\varphi_1-\varphi_2)$.

We shall explain below in section \ref{calculus} how to perform the
transition from \eqref{startr} to \eqref{finalr} in a conceptual
manner.

\bigskip

 In fact the conceptual  understanding of the
 simplicity of the final result of Theorem \ref{vol}
is at the origin of many of the notions developped in the present
paper and in particular of the ``rational" formulation of the
calculus which will be obtained in the last section. The geometric
meaning of Theorem \ref{vol} is the computation of the Jacobian in
the sense of noncommutative geometry of the  morphism $\theta$ as
explained above. The integral $\Omega$ is (up to a trivial
normalization factor) a standard elliptic integral, it is given by
an hypergeometric function in the variable
\begin{equation}
m:=\displaystyle\frac{s_k(s_l-s_m)}{s_l(s_k-s_m)}\,
\end{equation}
or a modular form in terms of $q$.

\begin{lem} The differential of $R$ is given on $F_\varphi$ by
$dR =\, J(Z) \; \chi$ where
\begin{equation} \label{dR}
 J(Z) =\,2\,(s_l-s_m)\,c_k\,t_k\,\frac{Z_0\,Z_1
\,Z_2 \,Z_3 }{(Z_m^2+  \,c_k\,Z_l^2)^2}
\end{equation}
\end{lem}

\begin{proof} This can easily be checked using \eqref{dsn} but it is
worthwile to give a simple direct argument. One has indeed by
definition of $\chi$
$$
d\frac{Z^2_j}{Z^2_0}=\,-2 \,s_j\,\frac{Z_0\,Z_1 \,Z_2 \,Z_3
}{Z_0^4}\,\chi\,,
$$
which gives using $$s_l\,Z^2_m-\,s_m\,Z^2_l=\,(s_l-s_m)\, Z^2_0$$
the equality
$$
d\frac{Z^2_l}{Z^2_m}=\,2 \,(s_m-s_l)\,\frac{Z_0\,Z_1 \,Z_2 \,Z_3
}{Z_m^4}\,\chi\,,
$$
and the required result.
\end{proof}

\smallskip The period $\Omega$ does not vanish and $J(Z),\,Z \in F_\varphi(0)$,
only vanishes on the $4$ ``ramification points" necessarily present
due to the symmetries.

\begin{cor} \label{nct0}
The Jacobian of the map $\theta^{\,t}$ is given by the equality
$$
 \theta_\ast(\ch_{\frac{3}{2}}(U))= 3 \,\Omega\,J\,V
$$
where $J$ is the element of the center $C^\infty(F_\varphi(0))$ of
$C^\infty(F_\varphi \times_{\sigma,\,\call} \mathbb {Z})$ given by
formula (\ref{dR}).
\end{cor}

\smallskip This statement assumes that $\varphi$
is generic in the measure theoretic sense so that $\eta$ admits good
diophantine approximation (\cite{ac:1982}). It justifies in
particular the terminology of ``ramified covering" applied to
$\theta^{\,t}$. The function $J$ has only $4$ zeros on
$F_\varphi(0)$ which correspond to the ramification.

\smallskip As shown by Theorem \ref{iden}
 the algebra $\cala_\varphi$ is defined over $\mathbb R$,
i.e. admits a  natural antilinear automorphism of period two,
$\gamma$ uniquely defined by
\begin{equation}
\gamma(Y_\mu):= Y_\mu\,,\quad \forall \mu
 \end{equation}
Theorem \ref{iden} also shows that $\sigma$ is defined over $\mathbb
R$ and hence commutes with complex conjugation $c(Z)= \bar Z$. This
gives a natural real structure $\gamma$ on the algebra $C_Q$ with
$C=F_\varphi \times F_\varphi$ and $Q$ as above,
$$
\gamma(f(Z, Z')):= \overline{f(c( Z),c( Z'))}\,,\quad
  \gamma(W_L):=W_{c(L)} \,,\quad \gamma(W'_{L'}):=W'_{c(L')}
$$
One checks that the morphism $\rho$ of lemma \ref{alg0} is ``real"
i.e. that,
\begin{equation}
\gamma \circ \rho =\, \rho \circ \gamma
\end{equation}
Since $c(Z)= \bar Z$ reverses the orientation of $F_\varphi$, while
$\gamma$ preserves the orientation of $S^3_\varphi$ it follows that
$J(\bar Z)=-J(Z)$ and $J$ necessarily vanishes on $F_\varphi(0) \cap
P_3(\mathbb {R})$.

\smallskip Note also that for general $h$ one has
$\langle\ch_{\frac{3}{2}}(U), \tau_h \rangle \neq 0$ which shows
that both $\ch_{\frac{3}{2}}(U) \in HH_3$ and $ \tau_h \in HH^3$ are
non trivial Hochschild classes. These results
 hold in the smooth algebra
$C^\infty(S^3_\varphi)$ containing the closure of
$C_{\mathrm{alg}}(S^3_\varphi)$ under holomorphic functional
calculus in the $C^\ast$ algebra $C^\ast(S^3_\varphi)$. We can also
use Theorem \ref{vol} to show the non-triviality of the morphism

$\theta_v : C_{\mathrm{alg}}(S^3_\varphi) \mapsto C^\infty({\mathbb
T}_{\eta}^2)$ of (\ref{ncu}).

\begin{cor} \label{nct}
The pullback of the fundamental class $[{\mathbb T}_{\eta}^2]$ of
the noncommutative torus by the homomorphism $\theta_v :
C_{\mathrm{alg}}(S^3_\varphi) \mapsto C^\infty({\mathbb
T}_{\eta}^2)$ of (\ref{ncu})
 is non zero, $\theta_v^\ast([{\mathbb T}_{\eta}^2]) \neq 0 \in HH^2$
provided $v$ is not a ramification point.
\end{cor}

\smallskip We have shown above the non-triviality of the Hochschild homology
and cohomology groups $HH_3(C^\infty(S^3_\varphi))$ and
$HH^3(C^\infty(S^3_\varphi))$ by exhibiting specific elements with
non-zero pairing. Combining the ramified cover
$\pi\,=\,\theta^{\,t}$  with the natural spectral geometry (spectral
triple) on the noncommutative $3$-dimensional nilmanifold $F_\varphi
\times_{\sigma,\,\call} \mathbb {Z}$ yields a natural spectral
triple on $S^3_\varphi$ in the generic case.

\section{Calculus and Cyclic Cohomology}\label{calculus}

\smallskip Theorem \ref{vol} suggests the existence of a ``rational" form
of the calculus explaining the appearance of the elliptic period
$\Omega$ and the rationality of $R$. We shall show in this last
section that this indeed the case. This will allow us to get a very
simple conceptual form of the above computation of the Jacobian in
Theorem \ref{vol1} below.

\smallskip Let us first go back to the general framework of twisted cross
products of the form
\begin{equation}
\cala=C^\infty(M)\times_{\sigma,\,\call} \mathbb {Z}
\end{equation}
where $\sigma$ is a diffeomorphism of the manifold $M$. We shall
follow \cite{ac:1986b} to construct cyclic cohomology classes from
cocycles in the bicomplex of group cohomology (with group
$\Gamma=\mathbb Z$) with coefficients in de Rham currents on $M$.
The twist by the line bundle $\call$ introduces a non-trivial
interesting nuance.

\smallskip Let $\Omega(M)$ be the
algebra  of smooth differential forms  on $M$, endowed with the
action of $\mathbb Z$
\begin{equation}
\alpha_{1,k}(\omega):=\sigma^{\ast k}\omega \,,\quad \, k\in \mathbb
Z
\end{equation}
As in \cite{ac:1994} p. 219 we let $\tilde{\Omega}(M)$ be the graded
algebra obtained as the (graded) tensor product of $\Omega(M)$ by
the exterior algebra $\wedge(\mathbb C [\mathbb Z]')$ on the
  augmentation ideal $\mathbb C [\mathbb Z]'$
in the group ring $\mathbb C [\mathbb Z]$. With $[n], n\in \mathbb
Z$ the canonical basis of $\mathbb C [\mathbb Z]$, the augmentation
$\epsilon :\mathbb C [\mathbb Z]\mapsto \mathbb C$ fulfills
$\epsilon([n])=1, \forall n\,$,  and
\begin{equation}
\delta_n := [n]- [0]\,,\quad \, n\in \mathbb Z\,,\quad n\neq 0
\end{equation}
is a linear basis of $\mathbb C [\mathbb Z]'$. The left regular
representation of $\mathbb Z$ on $\mathbb C [\mathbb Z]$ restricts
to $\mathbb C [\mathbb Z]'$ and is given on the above basis by
\begin{equation}\label{secondact}
\alpha_{2,k}(\delta_n):= \delta_{n+k}\,-\delta_k  \,,\quad \, k\in
\mathbb Z
\end{equation}
It extends to an action $\alpha_2$ of $\mathbb Z$ by automorphisms
of $\wedge\mathbb C [\mathbb Z]' $. We let $\alpha= \alpha_1 \otimes
\alpha_2$ be the tensor product action of $\mathbb Z$
 on the graded tensor product
\begin{equation}\label{graded}
\tilde{\Omega}(M)=\Omega(M)\otimes \wedge\mathbb C [\mathbb Z]'\,.
\end{equation}

\smallskip We now use the hermitian line bundle $\call$ to form the
twisted cross-product
\begin{equation}
\calc := \tilde{\Omega}(M)\times_{\alpha\,,\,\call} \mathbb {Z}
\end{equation}
We let $\call_n$ be as in (\ref{gene2}) for $n>0$ and extend its
definition for $n<0$ so that $\call_{-n}$ is the pullback by
$\sigma^n$ of the dual $\hat{\call}_n$ of $\call_n$ for all $n$. The
hermitian structure gives an antilinear isomorphism $\ast :\call_n
\mapsto \hat{\call}_n$. One has for all $n,m\in \mathbb Z$ a
canonical isomorphism
\begin{equation}\label{prodl}
\call_n\otimes \,\sigma^{*m}\call_m\simeq \call_{n+m}\,,
\end{equation}
which is by construction compatible with the hermitian structures.

The algebra $\calc$ is the linear span of monomials $\xi \, W^n$
where
\begin{equation}
\xi \in C^\infty(M,\call_n) \otimes_{C^\infty(M)} \tilde{\Omega}(M)
\end{equation}
with  the product rules (\ref{gene3}), (\ref{gene1}).

\smallskip Let $\nabla$ be a hermitian connection  on $\call$. We shall
turn $\calc$ into a differential graded algebra. By functoriality
$\nabla$ gives a hermitian connection on the $\call_k$ and hence a
graded derivation
\begin{equation}\label{nabla1}
\nabla_n  : C^\infty(M,\call_n) \otimes_{C^\infty(M)}
\Omega(M)\mapsto
 C^\infty(M,\call_n) \otimes_{C^\infty(M)} \Omega(M)
\end{equation}
whose square $\nabla_n^2$ is multiplication by the curvature
$\kappa_n\in \Omega^2(M)$ of $\call_n$,
\begin{equation}\label{curvcoc}
\kappa_{n+m}  =  \kappa_n+\sigma^{\ast n}(\kappa_m)\,,\quad \forall
n,\,m\in \mathbb Z
\end{equation}
with $\kappa_1=\kappa\in \Omega^2(M)$ the curvature of $\call$. Ones
has $d\kappa_n =0$ and we extend the differential $d$ to a graded
derivation on $ \tilde{\Omega}(M)$
 by
\begin{equation}\label{ddelta}
d \delta_n  = \,\kappa_n
\end{equation}
We can then extend $\nabla_n$ uniquely to the induced module
\begin{equation}\label{indmod}
{\mathcal E}_n=\,C^\infty(M,\call_n) \otimes_{C^\infty(M)}
\tilde{\Omega}(M)
\end{equation}
by the equality
\begin{equation}
\tilde{\nabla}_n(\xi\, \omega)= \nabla_n(\xi) \,\omega+(-1)^{{\rm
deg}(\xi)}\xi \,d\omega \,,\quad \forall \omega \in
\tilde{\Omega}(M)
\end{equation}

\medskip
\begin{prop} \label{dga}

\smallskip a) The pair $(\tilde{\Omega}(M), \,d)$ is a graded differential algebra.

\smallskip b)  Let $\alpha$
be the tensor product action
 of $\mathbb Z$ then $\alpha(k)\in {\rm Aut}(\tilde{\Omega}(M), \,d)
\qqq k\in \mathbb Z$.

\smallskip c) The following equality defines
a flat connection on  the induced module ${\mathcal E}_n$ on
$\tilde{\Omega}(M)$,
\begin{equation}\label{flatd}
\nabla'_n (\xi)= \tilde{\nabla}_n (\xi)- (-1)^{{\rm deg}(\xi)}\xi
\,\delta_n\,.
\end{equation}

\smallskip d) The graded derivation $d$ of $\tilde{\Omega}(M)$ extends uniquely
to a graded derivation of $\calc$ such that,
\begin{equation}\label{extendedd}
d(\xi \, W^n)= \nabla'_n (\xi)\, W^n
\end{equation}
which turns the pair $(\calc, \,d)$ into a graded differential
algebra.
\end{prop}

\begin{proof} a) By construction $d$ is the unique extension
of the differential $d$ of $\Omega(M)$ to a graded derivation of the
graded tensor product \eqref{graded} such that \eqref{ddelta} holds.
One just needs to check that $d^2=0$ on simple tensors
$\omega\otimes \delta_n$, one gets
$$
d(\omega\otimes \delta_n)=\,d\omega\otimes \delta_n\,+\, (-1)^{{\rm
deg}(\omega)}\,\omega\,\kappa_n \otimes 1\,,$$
$$
d^2(\omega\otimes \delta_n)=\,d^2\omega\otimes \delta_n\,+\,
(-1)^{{\rm deg}(\omega)+1}\,d\omega\,\kappa_n \otimes 1\,+
(-1)^{{\rm deg}(\omega)}\,d\omega\,\kappa_n \otimes 1\,=0\,.$$

\smallskip
b) Let us check that $\alpha(k)$ commutes with the differentiation
$d$ on simple tensors $\omega\otimes \delta_n$. One has
$$
d(\alpha(k)(\omega\otimes \delta_n))=\,d(\sigma^{*k}(\omega)\otimes
( \delta_{n+k}-\delta_k))=\,\sigma^{*k}(d\omega)\otimes (
\delta_{n+k}-\delta_k)+ (-1)^{{\rm
deg}(\omega)}\,\sigma^{*k}(\omega)\,(\kappa_{n+k}- \kappa_k)\otimes
1
$$
$$
\alpha(k)(d(\omega\otimes \delta_n))=\,\sigma^{*k}(d\omega)\otimes (
\delta_{n+k}-\delta_k)+ (-1)^{{\rm
deg}(\omega)}\,\sigma^{*k}(\omega\,\kappa_{n})\otimes 1
$$
and the equality follows from \eqref{curvcoc}.

\smallskip
c) Since $\delta_n^2=0$ one gets
$$
(\nabla'_n)^2 (\xi)=\,(\tilde{\nabla}_n)^2 (\xi)-\xi\,d\delta_n=\,0
$$
since $d \delta_n  = \,\kappa_n$ is the curvature.

\smallskip

d) Since the algebra $\calc$ is the linear span of monomials $\xi \,
W^n$ the linear map $d$ is well defined and coincides with the
differential $d$ on $ \tilde{\Omega}(M)$ since $\delta_0=0$. Let us
show that the two terms of \eqref{flatd} separately define
derivations of $\calc$. By construction the connections $\nabla_n$
are compatible with the canonical isomorphisms \eqref{prodl} and the
same holds for their extensions $\tilde{\nabla}_n$ which is enough
to show that the first term of \eqref{flatd} separately defines a
derivation of $\calc$. The proof for the second term
\begin{equation}\label{secterm}
d^{\,'}(\xi \, W^n)= \,  (-1)^{{\rm deg}(\xi)}\xi \,\delta_n\,
W^n\,,
\end{equation}
follows from \eqref{secondact} and is identical to the proof of
lemma 12 chapter III of \cite{ac:1994}. The flatness (c) of the
connections $\nabla'_n$ ensures that $d^2=0$ so that $(\calc, \,d)$
is a graded differential algebra.
\end{proof}

\smallskip To construct closed graded traces on
 this differential graded algebra we follow (\cite{ac:1986b})
and consider the double complex of group cochains (with group
$\Gamma=\mathbb Z$) with coefficients in de Rham currents on $M$.
The cochains $\gamma \in C^{n,m}$ are totally antisymmetric maps
from $\mathbb Z^{n+1}$ to the space $\Omega_{-m}(M)$ of de Rham
currents of dimension $-m$, which fulfill
\begin{equation}\label{transinv}
\gamma(k_0+k,k_1+k,k_2+k,\cdots,k_n+k)=\,\sigma_{\ast}^{- k}
\gamma(k_0,k_1,k_2,\cdots,k_n)\,,\quad \forall k,\,k_j \in \mathbb Z
\end{equation}
Besides the coboundary $d_1$ of group cohomology, given by
\begin{equation}\label{done}
(d_1\gamma)(k_0,k_1,\cdots,k_{n+1})=\, \sum_0^{n+1}\,
 (-1)^{j+m}\, \gamma(k_0,k_1,\cdots,\hat{k_j},\cdots,k_{n+1})
\end{equation}
and the coboundary $d_2$ of de Rham homology,
$$
(d_2\gamma)(k_0,k_1,\cdots,k_{n})=\,b(\gamma(k_0,k_1,\cdots,k_{n}))
$$
the curvatures $\kappa_n$ generate the further coboundary $d_3$
defined on Ker $d_1$ by,
\begin{equation}
(d_3\gamma)(k_0,\cdots,k_{n+1})=\, \sum_0^{n+1}\,
 (-1)^{j+m+1}\,\kappa_{\,k_j} \gamma(k_0,\cdots,\hat{k_j},\cdots,k_{n+1})
\end{equation}
which maps Ker $d_1 \cap C^{n,m}$ to $ C^{n+1,m+2}$. Translation
invariance follows from (\ref{curvcoc}) and $\varphi_{\ast}(\omega
C)= \varphi^{\ast-1}(\omega )\varphi_{\ast}( C)$ for $C \in
\Omega_{-m}(M)$, $\omega \in \Omega^{\ast}(M)$.

\smallskip To each $\gamma \in C^{n,m}$ one associates the functional
$\tilde{\gamma}$ on $\calc$ given by,
\begin{eqnarray}
\tilde{\gamma}(\,\xi \,W^n) & = & 0\,,\quad \forall\, n\neq 0
\,, \quad\xi \in \tilde{\Omega}(M)\nonumber\\
\tilde{\gamma}(\, \omega \otimes \delta_{k_1} \cdots \delta_{k_n})&
= & \langle \omega, \gamma(0, k_1\cdots ,k_n)\rangle \,,\quad
\forall k_j \in \mathbb Z
\end{eqnarray}

One then has

\begin{lem}\label{invariance} Let $\gamma \in C^{n,m}$, then

\smallskip
(i) One has for all $\rho \in \tilde{\Omega}(M)$,
\begin{equation}\label{inv1}
\tilde{\gamma}(\rho-\,\alpha(-k)\rho)=\,-\tilde{(d_1\gamma)}(\delta_k\,\rho)\,.
\end{equation}

\smallskip (ii) One has for all $a,b \in \calc$ with $d'$ defined in \eqref{secterm},
\begin{equation}\label{inv3}
\tilde{\gamma}(a\,b-\,(-1)^{{\rm deg}(a){\rm deg}(b)}\,b\,a)
=\,\,(-1)^{{\rm deg}(a)}\tilde{(d_1\gamma)}(a\,d'b)\,.
\end{equation}

\smallskip (iii) One has for all $\rho \in \tilde{\Omega}(M)$,
\begin{equation}\label{inv3}
\tilde{\gamma}(d\rho)=\,\tilde{(d_2\gamma)}(\rho)+
\,\tilde{(d_3\gamma)}(\rho)\,.
\end{equation}

\end{lem}

\begin{proof} (i) We can assume that $\rho$ is of the form
$$
\rho=\,\omega\otimes \delta_{k_1}\cdots \delta_{k_n}
$$
The left side of \eqref{inv1} is by construction
$$
\langle \omega, \gamma(0, k_1\cdots ,k_n)\rangle-\,
\tilde{\gamma}(\sigma^{*-k}(\omega)\otimes(\delta_{k_1-k}-\delta_{-k})
\cdots (\delta_{k_n-k}-\delta_{-k}))\,.
$$
When one expands the product one gets using $\delta_{-k}^2=0$,
$$
(\delta_{k_1-k}-\delta_{-k}) \cdots
(\delta_{k_n-k}-\delta_{-k})=\,\prod \delta_{k_j-k}+
\,\delta_{-k}\,\sum \,(-1)^i\,\prod_{j\neq i} \delta_{k_j-k}
$$
and one uses the translation invariance \eqref{transinv} to write
$$
\tilde{\gamma}(\sigma^{*-k}(\omega)\otimes\,\prod \delta_{k_j-k})
=\,\langle \omega, \gamma(k, k_1\cdots ,k_n)\rangle
$$
and
$$
 \tilde{\gamma}(\sigma^{*-k}(\omega)\otimes\,(-1)^i\,
\delta_{-k}\,\prod_{j\neq i} \delta_{k_j-k}) =\,-(-1)^i\,\langle
\omega, \gamma(0,k, k_1\cdots,\hat{k_i},\cdots,k_n)\rangle
$$
One thus obtains the same terms as in
$$
-\tilde{(d_1\gamma)}(\delta_k\,\rho)=\, -\,\langle \omega,
(d_1\gamma)(0,k, k_1\cdots ,k_n)\rangle
$$
using \eqref{done} and the graded commutation of $\delta_k$ with
$\omega$ which yields a $(-1)^m$ overall sign.

\smallskip (ii) It is enough to show that for any $k\in \mathbb Z $
and $a',b'\in \tilde{\Omega}(M)$ equation \eqref{inv2} holds for
$a=a'\,W^k$ and $b=b'\,W^{-k}$. The graded commutativity of
$\tilde{\Omega}(M)$ allows to write the graded commutator in
\eqref{inv2} as $\rho-\,\alpha(-k)\rho$ where $\rho
=a'\,\alpha(k)(b')\in \tilde{\Omega}(M)$. One has $\rho=\,a\,b$ and
(i) thus shows that
$$
\tilde{\gamma}(a\,b-\,(-1)^{{\rm deg}(a){\rm deg}(b)}\,b\,a)
=\,-\,\tilde{(d_1\gamma)}(\delta_k\,a\,b)=\, -(-1)^{{\rm
deg}(a)+{\rm deg}(b)} \,\tilde{(d_1\gamma)}(a\,b\,\delta_k)
$$
One has $$a\,d^{\,'}(b)=\,a'\,\alpha(k)((-1)^{{\rm
deg}(b)}b'\,\delta_{-k})$$ thus the result follows from the equality
$$
\alpha(k)(\delta_{-k})=\,-\,\delta_k\,.
$$

\smallskip (iii) We can assume that $\rho$ is of the form
$$
\rho=\,\omega\otimes \delta_{k_1}\cdots \delta_{k_n}
$$
One has
$$
d\rho=\,d\omega\otimes \delta_{k_1}\cdots \delta_{k_n}-\, (-1)^{{\rm
deg}(\omega)}\,\sum\,(-1)^j\,\omega\;\kappa_{k_j} \otimes
\delta_{k_1}\cdots  \hat{\delta_{k_j}}\cdots  \delta_{k_n}
$$
which gives \eqref{inv3}.
\end{proof}
\medskip

To each $\gamma \in C^{n,m}$ one associates
 the $(n-m+1)$ linear form on $\cala=C^\infty(M)\times_{\sigma,\,\call} \mathbb {Z}$
given by,

\begin{eqnarray}
&\,&\Phi(\gamma)(a_0,a_1,\cdots,a_{n-m}) =\, \nonumber \\&\,&
\,\lambda_{n,m}\sum_0^{n-m}(-1)^{j(n-m-j)}\,
\tilde{\gamma}(da_{j+1}\cdots da_{n-m} \, a_0 \,da_1 \cdots
da_{j-1}\,da_j)\qquad
\end{eqnarray}
where $\lambda_{n,m}:=\frac{n!}{(n-m+1)!}$.

\begin{lem}\label{cyclic}
(i) The Hochschild coboundary $b\Phi(\gamma)$ is equal to
$\Phi(d_1\gamma)$.

\smallskip (ii) Let $\gamma \in C^{n,m}\,\cap\,$Ker $\,d_1$. Then $\Phi(\gamma)$ is a Hochschild cocycle and

\medskip

$\qquad \qquad\qquad\qquad B\Phi(\gamma)=
\Phi(d_2\gamma)+\frac{1}{n+1}\;\Phi(d_3\gamma)$
\end{lem}

\begin{proof} (i) The proof is identical to that of Theorem 14 a) Chapter III
of \cite{ac:1994}.

\smallskip (ii) By (i) and the hypothesis $d_1(\gamma)=0$
$\Phi(\gamma)$ is a Hochschild cocycle. In fact by lemma
\ref{invariance} the functional $\tilde{\gamma}$ is a graded trace
and thus the formula for $\Phi(\gamma)$ simplifies to
\begin{equation}\label{tracesimple}
\Phi(\gamma)(a_0,a_1,\cdots,a_{n-m}) =\,(n-m+1)\; \lambda_{n,m}
\tilde{\gamma}\,( a_0 \,da_1 \cdots da_{n-m} )
\end{equation}
It follows that
$$
B_0(\Phi(\gamma))(a_0,a_1,\cdots,a_{n-m-1})=\, (n-m+1)\;
\lambda_{n,m} \tilde{\gamma}\,(d\rho)\,,\quad \rho =\,a_0 \,da_1
\cdots da_{n-m-1}\,,
$$
and $B_0(\Phi(\gamma))$ is already cyclic so that
$B(\Phi(\gamma))=(n-m)B_0(\Phi(\gamma))$.

Since the coboundary $d_2$ anticommutes with $d_1$ one has
$d_2\gamma\in$ Ker$\,d_1$ and  $\Phi(d_2\gamma)$ is also a
Hochschild cocycle and is given by \eqref{tracesimple} for
$d_2\gamma$. Let us check that $d_3\gamma\in$ Ker$\,d_1$. One has up
to an overall sign,
$$
d_1\,d_3(\gamma)=\,\sum\,(-1)^i\,d_3(\gamma)(k_0,\cdots,\hat{k_i}\cdots,k_{n+2})
=\,\sum_{i,j}\, \epsilon(i,j)\,\kappa_{\,k_j}
\gamma(k_0,\cdots,\hat{k_{i'}},\cdots,\hat{k_{j'}},\cdots,k_{n+1})
$$
where $(i',j')$ is the permutation of $(i,j)$ such that $i'<j'$ and
up to an overall sign $\epsilon(i,j)$ is the product of the
signature of this permutation by $(-1)^{i+j}$. For each $j$ the
coefficient of $\kappa_{\,k_j}$ is up to an overall sign given by
$d_1(\gamma)(k_0,\cdots,\hat{k_j}\cdots,k_{n+2}) =0$, thus
$d_1\,d_3(\gamma)=0$.

The result thus follows from lemma \ref{invariance} (iii).
\end{proof}

\medskip
 We shall now show how the above general
framework allows to reformulate the calculus involved in Theorem
\ref{vol} in rational terms. We let $M$ be the elliptic curve
$F_\varphi$ where $\varphi$ is generic and even. Let then $\nabla$
be an arbitrary hermitian connection on $\call$ and $\kappa$ its
curvature. We first display a cocycle $\gamma =\sum \gamma_{n,m}\in
\sum C^{n,m}$ which reproduces the cyclic cocycle $\tau_3$.

\begin{lem}\label{dioph}
There exists a two form $\alpha$ on  $M=F_\varphi$ and a multiple
$\lambda \,dv$
 of the translation invariant two form $dv$
such that :

\medskip

 (i)
$ \qquad \qquad \kappa_n =\,n\,\lambda \,dv+(\sigma^{\ast n}\alpha
-\alpha) \,,\quad \forall n \in \mathbb Z $

\medskip

 (ii) $\quad
d_2(\gamma_{\,j})=0\,,\quad d_1(\gamma_{3})=0\,,\quad
d_1(\gamma_{\,1})+\frac{1}{2}\,d_3(\gamma_{3}) =0\,,\quad
B\Phi(\gamma_1)=0\,,$

\medskip

 where $\gamma_{\,1}\in  C^{1,0}$ and $\gamma_{3}\in  C^{1,-2}$
 are given by

\medskip

$ \gamma_{1}(k_0,k_1):=\frac{1}{2}\,(k_1-k_0)(\sigma^{\ast
k_0}\alpha+\sigma^{\ast k_1}\alpha) \,,\quad \gamma_{3}(k_0,k_1):=
k_1-k_0 \,,\quad \forall k_j \in \mathbb Z $

\medskip

 (iii)
The class of the cyclic cocycle $\Phi(\gamma_1)+\Phi(\gamma_3)$ is
equal to $\tau_3$.
\end{lem}
\bigskip

\smallskip We use the generic hypothesis in the measure theoretic sense
to solve the ``small denominator" problem in (i). In (ii) we
identify differential forms $\omega \in \Omega^{\,d}$ of degree $d$
with the dual currents of dimension $2-d$.

\smallskip It is a general principle explained in \cite{ac:1982} that a cyclic cocycle $\tau$
generates a calculus whose differential graded algebra is obtained
as the quotient of the universal one by the radical of $\tau$. We
shall now explicitely describe the reduced calculus obtained from
the cocycle of lemma \ref{dioph} {\it(iii)}.
 We  use as above the hermitian line bundle $\call$ to form the
twisted cross-product
\begin{equation}
\calb := \Omega(M)\times_{\alpha\,,\,\call} \mathbb {Z}
\end{equation}
of the algebra $\Omega(M)$ of differential forms on $M$ by the
diffeomorphism $\sigma$. Instead of having to adjoin the infinite
number of odd elements $\delta_n$ we just adjoin two $\chi$ and $X$
as follows. We let $\delta$ be the derivation of $\calb$ such that
\begin{equation}
\delta (\xi \, W^n):= i \,n\,\xi \, W^n\,,\quad \forall \xi \in
C^\infty(M,\call_n) \otimes_{C^\infty(M)} \Omega(M)
\end{equation}
We adjoin $\chi$ to $\calb$ by tensoring $\calb$ with the exterior
algebra $\wedge\{\chi\}$ generated by an element $\chi$ of degree 1,
and extend the connection $\nabla$ (\ref{nabla1}) to the unique
graded derivation $d'$ of $\Omega'=\calb \otimes \wedge\{\chi\}$
such that,
\begin{eqnarray} \label{deltaext}
d'\,\omega & = & \nabla \omega + \chi \delta(\omega)\,,\quad \forall \omega \in \calb\nonumber\\
d'\chi & = & -\,\lambda\, dv
\end{eqnarray}
with $\lambda\, dv $ as in lemma  \ref{dioph}. By construction,
every element of $\Omega'$ is of the form
\begin{equation}
y= b_0 + b_1 \,\chi \,,\quad b_j\in \calb
\end{equation}
 One does not yet have a graded differential
algebra since  $d^{\,'2}\neq 0$. However, with $\alpha$ as in lemma
\ref{dioph} one has
\begin{equation}
d^{\,'2}(x)= [\,x, \,\alpha] \,,\quad \forall x \in \Omega'=\calb
\otimes \wedge\{\chi\}
\end{equation}
and one can apply lemma $9$ p.$229$ of \cite{ac:1994} to get a
differential graded algebra by adjoining the degree $1$ element $X:=
``d1"$ fulfilling the rules
\begin{eqnarray} \label{X}
X^2= -\alpha \,,\quad x \,X\, y=0 \,,\quad \forall x,y \in \Omega'
\end{eqnarray}
and defining the differential $d$ by,
\begin{eqnarray} \label{deltaext1}
d\,x & = & d{\,'}\,x + [\,X, \,x]\,,\quad \forall x \in \Omega'
\nonumber\\
dX & = & 0
\end{eqnarray}
where $[\,X, \,x]$ is the graded commutator. It follows from lemma
$9$ p.$229$ of \cite{ac:1994} that we obtain a differential graded
algebra $\Omega^{\ast}$, generated by $\calb$, $\xi$ and $X$. In
fact using (\ref{X}) every element of $\Omega^{\ast}$ is of the form
\begin{equation} x=x_{1,1}+x_{1,2}\,X+X\,x_{2,1}+X\,x_{2,2}X\,,\quad x_{i,j}\in \Omega'
\end{equation}
and we define the functional $\int$ on $\Omega^{\ast}$ by extending
the ordinary integral,
\begin{equation}\label{step1}
\int  \omega := \int_M \omega  \,,\quad \forall \omega \in \Omega(M)
\end{equation}
first to $\calb := \Omega(M)\times_{\alpha\,,\,\call} \mathbb {Z} $
by
\begin{equation}\label{step2}
\int  \,\xi \, W^n\,:=0  \,,\quad \forall n\neq 0
\end{equation}
then to $\Omega'$ by
\begin{equation}  \label{step3}
\int  (b_0 + b_1 \,\chi):=\int b_1  \,,\quad \forall b_j \in \calb
\end{equation}
and finally to $\Omega^{\ast}$ as in lemma $9$ p.$229$ of
\cite{ac:1994},
\begin{equation}  \label{step4}
\int (x_{1,1}+x_{1,2}\,X+X\,x_{2,1}+X\,x_{2,2}X):=\int
x_{1,1}+(-1)^{{\rm deg}(\,x_{2,2})} \int x_{2,2}\, \alpha
\end{equation}

\bigskip
\begin{thm} \label{fund} Let $M=F_\varphi$,
$\nabla$, $\alpha$ be as in lemma \ref{dioph}.

\smallskip The algebra $\Omega^{\ast}$ is a differential graded algebra
containing $C^{\infty}(M)\times_{\alpha\,,\,\call} \mathbb {Z}$.

\smallskip  The functional $\int$ is a closed graded trace on  $\Omega^{\ast}$.

\smallskip  The character of the corresponding cycle on
$C^{\infty}(M)\times_{\alpha\,,\,\call} \mathbb {Z}$
$$
\tau(a_0, \cdots, a_{3}):= \,\int \,a_0
 \,da_1 \cdots \, da_{3}   \,,\quad \forall a_j \in
C^{\infty}(M)\times_{\alpha\,,\,\call} \mathbb {Z}
$$
is cohomologous to the cyclic cocycle $\tau_3$.
\end{thm}

\smallskip It is worth noticing that the above calculus fits
with \cite{ac:1980}, \cite{mdv-mic:1994}, and \cite{mdv:2001}.

\smallskip Now in our case the line bundle $\call$ is holomorphic
and we can apply Theorem \ref{fund} to its canonical  hermitian
connection $\nabla$. We take the notations of section \ref{Cstar},
with
 $C=F_\varphi \times F_\varphi$,
and $Q$ given by (\ref{q13}). This gives a particular ``rational"
form of the calculus which explains the rationality of the answer in
Theorem \ref{vol}. We first extend as follows the construction of
 $C_Q$. We let $\Omega(C,Q)$ be the  generalised cross-product of the
algebra $\Omega(C)$ of meromorphic differential forms (in $dZ$ and
$dZ'$) on $C$ by the transformation $\tilde{\sigma}$. The
generators $W_L$ and $W'_{L'}$ fulfill the  cross-product rules,
\begin{equation}
W_L \,\,\omega =  \tilde{\sigma}^\ast(\omega) \;W_L \,, \qquad
W'_{L'} \,\omega = ( \tilde{\sigma}^{-1})^\ast(\omega) \;W'_{L'}
\label{cropro1}
\end{equation}
while (\ref{crossed}) is unchanged. The connection $\nabla$ is the
restriction to the subspace $\{Z'=\bar Z\}$ of the  unique graded
derivation $\nabla$  on $\Omega(C,Q)$ which induces the usual
differential on $\Omega(C)$ and satisfies,
\begin{eqnarray} \label{rational}
\nabla W_L & = & (\,d_Z\,\log L(Z)-d_Z\,\log Q(Z,Z'))\, W_L \nonumber\\
\nabla W'_{L'} & = & W'_{L'}\,(\,d_{Z'}\,\log L'({Z'})-d_{Z'}\,\log
Q(Z,Z'))
\end{eqnarray}
where $d_Z$ and $d_{Z'}$ are the (partial) differentials relative to
the variables $Z$ and $Z'$. Note that one needs to check that the
involved differential forms such as $d_Z\,\log L(Z)-d_Z\,\log
Q(Z,Z')$ are not only invariant under the scaling transformations $Z
\mapsto \lambda Z$ but are also \underline{basic}, i.e. have zero
restriction to the fibers of the map $\mathbb C^4 \mapsto
P_3(\mathbb {C})$, in both variables $Z$ and $Z'$. By definition the
derivation $\delta_\kappa\,=\,\nabla^2$ of $\Omega(C,Q)$ vanishes on
$\Omega(C)$ and fulfills
\begin{equation}
\delta_\kappa(W_L)\,=\,\kappa\,W_L \,,\qquad
\delta_\kappa(W'_{L'})\,=\,-\,W'_{L'}\,\kappa \label{kappa}
\end{equation}
where
\begin{equation} \label{curv1}
\kappa=d_Z\,d_{Z'}\,\log Q(Z,Z')
\end{equation}
 is a basic form which when restricted to the
subspace $\{Z'=\bar Z\}$ is  the curvature.
  We let as above $\delta$ be the derivation of  $\Omega(C,Q)$ which vanishes on $\Omega(C)$
 and is such that $\delta W_L = i\,W_L$ and $\delta W'_{L'}=-i\,W'_{L'}$.
We proceed exactly as above and get the graded algebras
$\Omega'=\Omega(C,Q)\otimes \wedge\{\chi\}$ obtained by adjoining
$\chi$ and $\Omega^{\ast}$ by adjoining $X$. We define $d'$, $d$ as
in  (\ref{deltaext}) and (\ref{deltaext1})  and the integral $\int$
by integration (\ref{step1}) on the subspace $\{Z'=\bar Z\}$
followed as above by
steps (\ref{step2}), (\ref{step3}), (\ref{step4}).\\

\begin{cor} \label{fund1} Let $\rho$:
$C_{\mathrm{alg}}(S^3_\varphi)  \mapsto C_Q$ be the morphism of
lemma \ref{alg0}.
  The equality
$$
\tau_{\mathrm{alg}}(a_0, \cdots, a_{3}):= \,\int \,\rho(a_0)
 \,d^{\,'}\rho(a_1) \cdots \, d^{\,'}\rho(a_{3})
$$
defines a $3$-dimensional Hochschild cocycle $\tau_{\mathrm{alg}}$
on $C_{\mathrm{alg}}(S^3_\varphi) $.

\smallskip Let $h \in$ Center $(C^\infty(F_\varphi \times_{\sigma,\,\call} \mathbb {Z}) )\sim
C^\infty(F_\varphi(0))$. Then
$$
\langle \ch_{\frac{3}{2}}(U), \tau_h \rangle=\langle
h\,\ch_{\frac{3}{2}}(U), \tau_{\mathrm{alg}} \rangle
$$

\end{cor}

\smallskip The computation of $d'$ only involves rational fractions
in the variables $Z$, $Z'$ (\ref{rational}), and the formula
(\ref{beauty}) for  $\ch_{\frac{3}{2}}(U)$ is polynomial in the
$W_L$, $W'_{L'}$.

 We are now ready for
a better understanding and formulation
 of the result of Theorem \ref{vol}.
Indeed what the above shows is that the denominator that appears in
the rational fraction $R(Z)$ of Theorem \ref{vol} should have to do
with the central quadratic form $Q(Z,Z')$ evaluated on the pairs
$(Z,\bar Z)$. In fact the two dimensional space of central quadratic
forms provides a natural  space of functions of the form
\begin{equation}
R(Z)=\,\frac{P(Z,\bar Z)}{Q(Z,\bar Z)}
\end{equation}
and this space is one dimensional when one mods out the constant
functions.

The following lemma shows that these functions are in fact invariant
under the correspondence $\sigma$.

\medskip
\begin{lem} \label{alg4}
Let $Q$ be  central and not identically zero on the component $C$
and $P$ be  central then the function
$$
R(Z,Z')=\,\frac{P(Z, Z')}{Q(Z,Z')}
$$
is invariant under $\tilde \sigma$.

\end{lem}

\medskip

It is thus natural now to compare the differential form $dR$ with
the form that appears in Theorem \ref{vol}.

\begin{thm} \label{vol1} Let $\varphi$ be generic and even and let
$Q$ be the central quadratic form on ${\mathbb R}^4_\varphi$
defining the three sphere $S^3_\varphi$. Let $\rho_Q$ be the
associated $*$-homomorphism
$$C^\infty(S^3_\varphi)\to C^\infty(E \times_{\sigma,\,\call} \mathbb {Z}) \,.$$
 Then for any central quadratic form $P$ not proportional to $Q$ there exists
a scalar $\mu$ such that
$$
\langle \ch_{\frac{3}{2}}(U), \tau_h \rangle=\,\mu\,\int \,h(Z)\,
dR(Z)\qqq h \in {\rm Center}\; C^\infty(E \times_{\sigma,\,\call}
\mathbb {Z}) \,,
$$
where
$$
R(Z)=\,\frac{P(Z,\bar Z)}{Q(Z,\bar Z)}\,.
$$
\end{thm}

\medskip
\begin{proof} Let us show that one can interpret \eqref{startr}
in the above terms.

Thus with $a={\rm cot}(\varphi_1-\varphi_2)\;{\rm
cot}(\varphi_1-\varphi_3)$ we need to show that
$$
R(Z)=\, t_1\;\frac{Z_0^2}{Z_0^2+  \,a\,Z_1^2}
$$
is in fact the restriction to the subset $F_\varphi(0)\subset E$ of
elements $Z$ with $c(Z)=I_0(Z)$ of  a ratio of the form
$$
\frac{P(Z,j(Z))}{Q(Z,j(Z))}\,.
$$
One has $j(Z)=\,I_3(c(Z)$ and thus
\begin{equation} \label{j3}
j(Z)= \, I_3\circ I_0(Z)\qqq Z \in F_\varphi(0)\,.
\end{equation}
We let $Q'$ be the central quadratic form of proposition \ref{ztoy}
namely
$$
s\,Q'=\,Q_1+\,Q_3+\,s_2\,Q_2\,,
$$
and we let $P'$ be given by
\begin{equation} \label{formP'}
s\,P'=\,Q_1+\,Q_3\,,
\end{equation}
A simple computation using \eqref{j3} then shows that
\begin{equation} \label{rat1}
b\;\frac{P'(Z,j(Z))}{Q'(Z,j(Z))}=\, \frac{Z_0^2}{Z_0^2+
\,a\,Z_1^2}-\,\frac{1}{s_1}\qqq Z\in F_\varphi(0)\,,
\end{equation}
where
$$b=\,\frac{\sin\varphi_2\,\sin\varphi_3}{\cos(\varphi_2-\,\varphi_3)}\,.$$
This gives the required result for the central quadratic form $P'$
and since the space of central quadratic forms is two dimensional
its quotient by multiples of $Q'$ is one dimensional so that the
result holds for all non-zero elements of this quotient.
\end{proof}

\medskip
With the above ``invariant" formulation of the formulas of Theorem
\ref{vol} we can now perform the  change of variables required in
the last part of its proof \ie explain how to pass from
\eqref{startr} to \eqref{finalr}.

We let $P$ and $Q$ be the central quadratic forms obtained from $P'$
and $Q'$ by the isomorphism $\beta$ of proposition \ref{ztoy}. The
form $Q$ is given by \eqref{q13} \ie by
\begin{eqnarray}
Q &=&\;(\prod \cos^2\varphi_{\ell}) \, \sum \, t_k\, s_k \,\,Q_k \,
=
\prod \sin(\varphi_{\ell}-\varphi_m)\,Z_0^2\\
&-&\sum\,\cos\varphi_{\ell}\,\cos\varphi_m
\,\sin(\varphi_{\ell}-\varphi_m) \,Z_k^2\nonumber
\end{eqnarray}
The form $P$ is given by
\begin{equation} \label{formP}
P\;=\;\sum_1^3\;\sin\varphi_k\,\sin(\varphi_{\ell}-\varphi_m)
\,\cos(\varphi_k-\varphi_{\ell}-\varphi_m) \, \,Z_k^2
\end{equation}
The involution $j_2$ is now given by $j_2(Z)=\,I_2(c(Z)$ and one has
\begin{equation} \label{j2}
j_2(Z)= \, I_2\circ I_0(Z)\qqq Z \in F_\varphi(0)\,.
\end{equation}
A simple computation using \eqref{j2} then shows that
\begin{equation} \label{rat2}
b\;\frac{P(Z,j_2(Z))}{Q(Z,j_2(Z))}=\, \frac{Z_3^2}{Z_3^2+
\,c_1\,Z_2^2}-\, \,\frac{1}{s_1}\qqq Z\in F_\varphi(0)\,,
\end{equation}
with $c_1={\rm tg}(\varphi_2)\;{\rm cot}(\varphi_1-\varphi_2)$ and
we thus obtain the formula required by Theorem \ref{vol}.

\section{Appendix $1$: The list of minors}\label{appendix}

 We give for convenience the list of the 15 minors of the
matrix \eqref{matrixphi}, with labels the missing lines, and in
factorized form. By setting
\begin{equation}
\left\{
\begin{array}{l}
A=x^2_0+\sum^3_{k=1}\, \cos (2\varphi_k)x^2_k\\
\\
B=\sum^3_{k=1}\, \sin(2\varphi_k)x^2_k
\end{array}
\right.
\end{equation}
one sees that these minors are combinations of the form
$M_{ij}=P_{ij}\,A+Q_{ij} \,B$.

\begin{equation}
\begin{matrix}
 M_{12}=2\, (\sin  ({{\varphi }_1}-{{\varphi }_2})\, {x_1}\, {x_2}+i\, \cos  ({{\varphi }_3})\, {x_0}\,
{x_3})\,   \\
 (-\cos  ({{\varphi }_1}-{{\varphi }_2})\,  (\cos  ({{\varphi }_1}-{{\varphi }_3})\, \sin  ({{\varphi }_1})\, x_{1}^{2}+\cos  ({{\varphi }_2}-{{\varphi }_3})\,
\sin  ({{\varphi }_2})\, x_{2}^{2} )+  \\
\noalign{\vspace{0.666667ex}} \hspace{4.em} \sin  ({{\varphi }_3})\,
(\sin  ({{\varphi }_1})\, \sin  ({{\varphi }_2})\, x_{0}^{2}-\cos
({{\varphi }_1}-{{\varphi
}_3})\, \cos  ({{\varphi }_2}-{{\varphi }_3})\, x_{3}^{2} ) )\\
=(\sin  ({{\varphi }_1}-{{\varphi }_2})\, {x_1}\, {x_2}+i\, \cos
({{\varphi }_3})\, {x_0}\,
{x_3})\,   \\
(2\sin(\varphi_1)\, \sin(\varphi_2)\, \sin(\varphi_3)\, A -
(\cos(\varphi_1-\varphi_2)\,
\cos(\varphi_3)+\sin(\varphi_1+\varphi_2)\, \sin(\varphi_3))B)
\end{matrix}
\end{equation}

\begin{equation}
\begin{matrix}
 M_{13}=2\, i\, (\cos  ({{\varphi }_2})\, {x_0}\, {x_2}+i\, \sin  ({{\varphi }_1}-{{\varphi
}_3})\, {x_1}\, {x_3})\,   \\
\noalign{\vspace{0.666667ex}} \hspace{2.em}  (-\cos  ({{\varphi
}_1}-{{\varphi }_2})\,
 (\cos  ({{\varphi }_1}-{{\varphi }_3})\, \sin  ({{\varphi }_1})\, x_{1}^{2}+\cos  ({{\varphi }_2}-{{\varphi }_3})\,
\sin  ({{\varphi }_2})\, x_{2}^{2} )+  \\
\noalign{\vspace{0.666667ex}} \hspace{4.em} \sin  ({{\varphi }_3})\,
(\sin  ({{\varphi }_1})\, \sin  ({{\varphi }_2})\, x_{0}^{2}-\cos
({{\varphi }_1}-{{\varphi
}_3})\, \cos  ({{\varphi }_2}-{{\varphi }_3})\, x_{3}^{2} ) )\\
=i\, (\cos  ({{\varphi }_2})\, {x_0}\, {x_2}+i\, \sin  ({{\varphi
}_1}-{{\varphi
}_3})\, {x_1}\, {x_3})\\
(2\sin(\varphi_1)\, \sin(\varphi_2)\, \sin(\varphi_3)\, A -
(\cos(\varphi_1-\varphi_2)\,
\cos(\varphi_3)+\sin(\varphi_1+\varphi_2)\, \sin(\varphi_3))B)
\end{matrix}
\end{equation}

\begin{equation}
\begin{matrix}
 M_{14}=2\, (\sin  ({{\varphi }_2})\, {x_0}\, {x_2}+i\, \cos  ({{\varphi }_1}-{{\varphi }_3})\, {x_1}\,
{x_3})\,   \\
\noalign{\vspace{0.666667ex}} \hspace{2.em}  (\cos  ({{\varphi
}_2})\,  (\cos  ({{\varphi }_3})\, \sin  ({{\varphi }_1})\,
x_{0}^{2}+\cos  ({{\varphi }_2}-{{\varphi
}_3})\, \sin  ({{\varphi }_1}-{{\varphi }_2})\, x_{2}^{2} )+  \\
\noalign{\vspace{0.666667ex}} \hspace{4.em} \sin  ({{\varphi
}_1}-{{\varphi }_3})\,
 (-\sin  ({{\varphi }_1})\, \sin  ({{\varphi }_1}-{{\varphi }_2})\, x_{1}^{2}+\cos  ({{\varphi }_2}-{{\varphi }_3})\,
\cos  ({{\varphi }_3})\, x_{3}^{2} ) )\\
= (\sin  ({{\varphi }_2})\, {x_0}\, {x_2}+i\, \cos  ({{\varphi }_1}-{{\varphi }_3})\, {x_1}\, {x_3})\\
(2\sin(\varphi_1)\, \cos(\varphi_2)\, \cos(\varphi_3)\,
A-(\cos(\varphi_1)\, \cos(\varphi_2-\varphi_3)-\sin(\varphi_1)\,
\sin(\varphi_2+\varphi_3))\, B)
\end{matrix}
\end{equation}

\begin{equation}
\begin{matrix}
 M_{15}=-i\, \cos  ({{\varphi }_1}-{{\varphi }_2}-{{\varphi }_3})\,   \\
\noalign{\vspace{0.666667ex}} \hspace{2.em}  (\sin  (2\, ({{\varphi
}_1}-{{\varphi }_2}))\, x_{1}^{2}\, x_{2}^{2}+\sin  (2\, ({{\varphi
}_1}-{{\varphi }_3}))\,
x_{1}^{2}\, x_{3}^{2}-   x_{0}^{2}\,  (\sin  (2\, {{\varphi }_2})\, x_{2}^{2}+\sin  (2\, {{\varphi }_3})\, x_{3}^{2} ) )\\
=  i\, \cos  ({{\varphi }_1}-{{\varphi }_2}-{{\varphi }_3})\\
((\sin(2\varphi_2)\, x^2_2+\sin(2\varphi_3)\, x^2_3)\, A -
(\cos(2\varphi_2)\, x^2_2+\cos(2\varphi_3)\, x^2_3)\, B)
\end{matrix}
\end{equation}

\begin{equation}
\begin{matrix}
 M_{16}=-2\, (-i\, \cos  ({{\varphi }_1}-{{\varphi }_2})\, {x_1}\, {x_2}+\sin  ({{\varphi }_3})\, {x_0}\,
{x_3})\,   \\
\noalign{\vspace{0.666667ex}} \hspace{2.em}  (\cos  ({{\varphi
}_2})\,  (\cos  ({{\varphi }_3})\, \sin  ({{\varphi }_1})\,
x_{0}^{2}+\cos  ({{\varphi }_2}-{{\varphi
}_3})\, \sin  ({{\varphi }_1}-{{\varphi }_2})\, x_{2}^{2} )+  \\
\noalign{\vspace{0.666667ex}} \hspace{4.em} \sin  ({{\varphi
}_1}-{{\varphi }_3})\,
 (-\sin  ({{\varphi }_1})\, \sin  ({{\varphi }_1}-{{\varphi }_2})\, x_{1}^{2}+\cos  ({{\varphi }_2}-{{\varphi }_3})\,
\cos  ({{\varphi }_3})\, x_{3}^{2} ) )\\
=-(-i\, \cos  ({{\varphi }_1}-{{\varphi }_2})\, {x_1}\, {x_2}+\sin  ({{\varphi }_3})\, {x_0}\,{x_3})\\
(2\sin(\varphi_1)\, \cos(\varphi_2)\, \cos(\varphi_3)\,
A-(\cos(\varphi_1)\, \cos(\varphi_2-\varphi_3)-\sin(\varphi_1)\,
\sin(\varphi_2+\varphi_3))\, B)
\end{matrix}
\end{equation}

\begin{equation}
\begin{matrix}
 M_{23}=2\, (i\, \cos  ({{\varphi }_1})\, {x_0}\, {x_1}+\sin  ({{\varphi }_2}-{{\varphi }_3})\, {x_2}\,
{x_3})\,   \\
\noalign{\vspace{0.666667ex}} \hspace{2.em}  (-\cos  ({{\varphi
}_1}-{{\varphi }_2})\,
 (\cos  ({{\varphi }_1}-{{\varphi }_3})\, \sin  ({{\varphi }_1})\, x_{1}^{2}+\cos  ({{\varphi }_2}-{{\varphi }_3})\,
\sin  ({{\varphi }_2})\, x_{2}^{2} )+  \\
\noalign{\vspace{0.666667ex}} \hspace{4.em} \sin  ({{\varphi }_3})\,
(\sin  ({{\varphi }_1})\, \sin  ({{\varphi }_2})\, x_{0}^{2}-\cos
({{\varphi }_1}-{{\varphi
}_3})\, \cos  ({{\varphi }_2}-{{\varphi }_3})\, x_{3}^{2} ) )\\
= (i\, \cos  ({{\varphi }_1})\, {x_0}\, {x_1}+\sin  ({{\varphi }_2}-{{\varphi }_3})\, {x_2}\, {x_3})\\
(2\sin(\varphi_1)\, \sin(\varphi_2)\, \sin(\varphi_3)\, A -
(\cos(\varphi_1-\varphi_2)\,
\cos(\varphi_3)+\sin(\varphi_1+\varphi_2)\, \sin(\varphi_3))B)
\end{matrix}
\end{equation}

\begin{equation}
\begin{matrix}
 M_{24}=2\, (\sin  ({{\varphi }_1})\, {x_0}\, {x_1}-i\, \cos  ({{\varphi }_2}-{{\varphi }_3})\, {x_2}\,
{x_3})\,   \\
\noalign{\vspace{0.666667ex}} \hspace{2.em}  (\cos  ({{\varphi
}_1})\,  (\cos  ({{\varphi }_3})\, \sin  ({{\varphi }_2})\,
x_{0}^{2}-\cos  ({{\varphi }_1}-{{\varphi
}_3})\, \sin  ({{\varphi }_1}-{{\varphi }_2})\, x_{1}^{2} )+  \\
\noalign{\vspace{0.666667ex}} \hspace{4.em} \sin  ({{\varphi
}_2}-{{\varphi }_3})\,
  (\sin  ({{\varphi }_1}-{{\varphi }_2})\, \sin  ({{\varphi }_2})\, x_{2}^{2}+\cos  ({{\varphi }_1}-{{\varphi }_3})\,
\cos  ({{\varphi }_3})\, x_{3}^{2} ) )\\
= (\sin  ({{\varphi }_1})\, {x_0}\, {x_1}-i\, \cos  ({{\varphi }_2}-{{\varphi }_3})\, {x_2}\, {x_3})\\
(2\cos(\varphi_3)\, \cos(\varphi_1)\, \sin(\varphi_2)\,
A-(\cos(\varphi_3-\varphi_1)\,
\cos(\varphi_2)-\sin(\varphi_3+\varphi_1)\, \sin(\varphi_2))\, B)
\end{matrix}
\end{equation}

\begin{equation}
\begin{matrix}
 M_{25}=-2\, i\, (\cos  ({{\varphi }_1}-{{\varphi }_2})\, {x_1}\, {x_2}-i\, \sin  ({{\varphi
}_3})\, {x_0}\, {x_3})\,   \\
\noalign{\vspace{0.666667ex}} \hspace{2.em}  (\cos  ({{\varphi
}_1})\,  (-\cos  ({{\varphi }_3})\, \sin  ({{\varphi }_2})\,
x_{0}^{2}+\cos  ({{\varphi }_1}-{{\varphi
}_3})\, \sin  ({{\varphi }_1}-{{\varphi }_2})\, x_{1}^{2} )-  \\
\noalign{\vspace{0.666667ex}} \hspace{4.em} \sin  ({{\varphi
}_2}-{{\varphi }_3})\,
 (\sin  ({{\varphi }_1}-{{\varphi }_2})\, \sin  ({{\varphi }_2})\, x_{2}^{2}+\cos  ({{\varphi }_1}-{{\varphi }_3})\,
\cos  ({{\varphi }_3})\, x_{3}^{2} ) )\\
= i\, (\cos  ({{\varphi }_1}-{{\varphi }_2})\, {x_1}\, {x_2}-i\,
\sin  ({{\varphi
}_3})\, {x_0}\, {x_3})\\
(2\cos(\varphi_3)\, \cos(\varphi_1)\, \sin(\varphi_2)\,
A-(\cos(\varphi_3-\varphi_1)\,
\cos(\varphi_2)-\sin(\varphi_3+\varphi_1)\, \sin(\varphi_2))\, B)
 \end{matrix}
\end{equation}

\begin{equation}
\begin{matrix}
 M_{26}=i\, \cos  ({{\varphi }_1}-{{\varphi }_2}+{{\varphi }_3})\,   \\
\noalign{\vspace{0.666667ex}} \hspace{2.em}  (\sin  (2\, {{\varphi
}_1})\, x_{0}^{2}\, x_{1}^{2}+\sin  (2\, ({{\varphi }_1}-{{\varphi
}_2}))\, x_{1}^{2}\,
x_{2}^{2}+   (\sin  (2\, {{\varphi }_3})\, x_{0}^{2}-\sin  (2\, ({{\varphi }_2}-{{\varphi }_3}))\, x_{2}^{2} )\, x_{3}^{2} )\\
= i\, \cos  ({{\varphi }_1}-{{\varphi }_2}+{{\varphi }_3})\\
((\sin(2 \varphi_1)\, x^2_1+ \sin(2\varphi_3)\, x^2_3)\, A -
(\cos(2\varphi_1)\, x^2_1+\cos(2\varphi_3)\, x^2_3)\, B)
\end{matrix}
\end{equation}

\begin{equation}
\begin{matrix}
 M_{34}=-i\, \cos  ({{\varphi }_1}+{{\varphi }_2}-{{\varphi }_3})\,   \\
\noalign{\vspace{0.666667ex}} \hspace{2.em}  (\sin  (2\, {{\varphi
}_1})\, x_{0}^{2}\, x_{1}^{2}+\sin  (2\, {{\varphi }_2})\,
x_{0}^{2}\, x_{2}^{2}+
  (\sin  (2\, ({{\varphi }_1}-{{\varphi }_3}))\, x_{1}^{2}+\sin  (2\, ({{\varphi }_2}-{{\varphi }_3}))\, x_{2}^{2} )\,
x_{3}^{2} )\\
= -i\, \cos  ({{\varphi }_1}+{{\varphi }_2}-{{\varphi }_3})\\
((\sin(2\varphi_1)\, x^2_1 + \sin(2\varphi_2)\, x^2_2)\,
A-(\cos(2\varphi_1)\, x^2_1+\cos(2\varphi_2)\, x^2_2)\, B)
\end{matrix}
\end{equation}

\begin{equation}
\begin{matrix}
 M_{35}=2\, i\, (i\, \sin  ({{\varphi }_2})\, {x_0}\, {x_2}+\cos  ({{\varphi }_1}-{{\varphi
}_3})\, {x_1}\, {x_3})\,   \\
\noalign{\vspace{0.666667ex}} \hspace{2.em}  (\cos  ({{\varphi
}_1})\,  (-\cos  ({{\varphi }_2})\, \sin  ({{\varphi }_3})\,
x_{0}^{2}+\cos  ({{\varphi }_1}-{{\varphi
}_2})\, \sin  ({{\varphi }_1}-{{\varphi }_3})\, x_{1}^{2} )+  \\
\noalign{\vspace{0.666667ex}} \hspace{2.em} \sin  ({{\varphi
}_2}-{{\varphi }_3})\,
 (\cos  ({{\varphi }_1}-{{\varphi }_2})\, \cos  ({{\varphi }_2})\, x_{2}^{2}+\sin  ({{\varphi }_1}-{{\varphi }_3})\,
\sin  ({{\varphi }_3})\, x_{3}^{2} ) )\\
= -i\, (i\, \sin  ({{\varphi }_2})\, {x_0}\, {x_2}+\cos  ({{\varphi
}_1}-{{\varphi
}_3})\, {x_1}\, {x_3})\\
(2\cos(\varphi_1)\, \cos(\varphi_2)\, \sin(\varphi_3)\,
A-(\cos(\varphi_1-\varphi_2)\,
\cos(\varphi_3)-\sin(\varphi_1+\varphi_2)\, \sin(\varphi_3))\, B)
\end{matrix}
\end{equation}

\begin{equation}
\begin{matrix}
 M_{36}=-2\, (\sin  ({{\varphi }_1})\, {x_0}\, {x_1}+i\, \cos  ({{\varphi }_2}-{{\varphi }_3})\, {x_2}\,
{x_3})\,   \\
\noalign{\vspace{0.666667ex}} \hspace{2.em}  (\cos  ({{\varphi
}_1})\,  (-\cos  ({{\varphi }_2})\, \sin  ({{\varphi }_3})\,
x_{0}^{2}+\cos  ({{\varphi }_1}-{{\varphi
}_2})\, \sin  ({{\varphi }_1}-{{\varphi }_3})\, x_{1}^{2} )+  \\
\noalign{\vspace{0.666667ex}} \hspace{2.em} \sin  ({{\varphi
}_2}-{{\varphi }_3})\,
 (\cos  ({{\varphi }_1}-{{\varphi }_2})\, \cos  ({{\varphi }_2})\, x_{2}^{2}+\sin  ({{\varphi }_1}-{{\varphi }_3})\,
\sin  ({{\varphi }_3})\, x_{3}^{2} ) )\\
= (\sin  ({{\varphi }_1})\, {x_0}\, {x_1}+i\, \cos  ({{\varphi }_2}-{{\varphi }_3})\, {x_2}\,{x_3})\\
(2\cos(\varphi_1)\, \cos(\varphi_2)\, \sin(\varphi_3)\,
A-(\cos(\varphi_1-\varphi_2)\,
\cos(\varphi_3)-\sin(\varphi_1+\varphi_2)\, \sin(\varphi_3))\, B)
\end{matrix}
\end{equation}

\begin{equation}
\begin{matrix}
 M_{45}=-2\, i\, (\cos  ({{\varphi }_2})\, {x_0}\, {x_2}-i\, \sin  ({{\varphi }_1}-{{\varphi
}_3})\, {x_1}\, {x_3})\,   \\
\noalign{\vspace{0.666667ex}} \hspace{2.em}  (\cos  ({{\varphi
}_1})\,  (\cos  ({{\varphi }_3})\, \sin  ({{\varphi }_2})\,
x_{0}^{2}-\cos  ({{\varphi }_1}-{{\varphi
}_3})\, \sin  ({{\varphi }_1}-{{\varphi }_2})\, x_{1}^{2} )+  \\
\noalign{\vspace{0.666667ex}} \hspace{2.em} \sin  ({{\varphi
}_2}-{{\varphi }_3})\,
 (\sin  ({{\varphi }_1}-{{\varphi }_2})\, \sin  ({{\varphi }_2})\, x_{2}^{2}+\cos  ({{\varphi }_1}-{{\varphi }_3})\,
\cos  ({{\varphi }_3})\, x_{3}^{2} ) )\\
= - i\, (\cos  ({{\varphi }_2})\, {x_0}\, {x_2}-i\, \sin  ({{\varphi
}_1}-{{\varphi
}_3})\, {x_1}\, {x_3})\\
(2\cos(\varphi_3)\, \cos(\varphi_1)\, \sin(\varphi_2)\,
A-(\cos(\varphi_3-\varphi_1)\,
\cos(\varphi_2)-\sin(\varphi_3+\varphi_1)\, \sin(\varphi_2))\, B)
\end{matrix}
\end{equation}

\begin{equation}
\begin{matrix}
 M_{46}=2\, (-i\, \cos  ({{\varphi }_1})\, {x_0}\, {x_1}+\sin  ({{\varphi }_2}-{{\varphi }_3})\, {x_2}\,
{x_3})\,   \\
\noalign{\vspace{0.666667ex}} \hspace{2.em}  (\cos  ({{\varphi
}_2})\,  (\cos  ({{\varphi }_3})\, \sin  ({{\varphi }_1})\,
x_{0}^{2}+\cos  ({{\varphi }_2}-{{\varphi
}_3})\, \sin  ({{\varphi }_1}-{{\varphi }_2})\, x_{2}^{2} )+  \\
\noalign{\vspace{0.666667ex}} \hspace{2.em} \sin  ({{\varphi
}_1}-{{\varphi }_3})\,
 (-\sin  ({{\varphi }_1})\, \sin  ({{\varphi }_1}-{{\varphi }_2})\, x_{1}^{2}+\cos  ({{\varphi }_2}-{{\varphi }_3})\,
\cos  ({{\varphi }_3})\, x_{3}^{2} ) )\\
= (-i\, \cos  ({{\varphi }_1})\, {x_0}\, {x_1}+\sin  ({{\varphi }_2}-{{\varphi }_3})\, {x_2}\, {x_3})\\
 (2\sin(\varphi_1)\, \cos(\varphi_2)\, \cos(\varphi_3)\, A-(\cos(\varphi_1)\, \cos(\varphi_2-\varphi_3)-\sin(\varphi_1)\, \sin(\varphi_2+\varphi_3))\, B)
\end{matrix}
\end{equation}

\begin{equation}
\begin{matrix}
 M_{56}=2\, (\sin  ({{\varphi }_1}-{{\varphi }_2})\, {x_1}\, {x_2}-i\, \cos  ({{\varphi }_3})\, {x_0}\,
{x_3})\,   \\
\noalign{\vspace{0.666667ex}} \hspace{2.em}  (\cos  ({{\varphi
}_1})\,  (\cos  ({{\varphi }_2})\, \sin  ({{\varphi }_3})\,
x_{0}^{2}-\cos  ({{\varphi }_1}-{{\varphi
}_2})\, \sin  ({{\varphi }_1}-{{\varphi }_3})\, x_{1}^{2} )-  \\
\noalign{\vspace{0.666667ex}} \hspace{2.em} \sin  ({{\varphi
}_2}-{{\varphi }_3})\,
  (\cos  ({{\varphi }_1}-{{\varphi }_2})\, \cos  ({{\varphi }_2})\, x_{2}^{2}+\sin  ({{\varphi }_1}-{{\varphi }_3})\,
\sin  ({{\varphi }_3})\, x_{3}^{2} ) )\\
= (\sin  ({{\varphi }_1}-{{\varphi }_2})\, {x_1}\, {x_2}-i\, \cos  ({{\varphi }_3})\, {x_0}\, {x_3})\\
(2\cos(\varphi_1)\, \cos(\varphi_2)\, \sin(\varphi_3)\,
A-(\cos(\varphi_1-\varphi_2)\,
\cos(\varphi_3)-\sin(\varphi_1+\varphi_2)\, \sin(\varphi_3))\, B)
\end{matrix}
\end{equation}

\bigskip
\bigskip

\section{Appendix $2$: The sixteeen theta relations}\label{app3}

The sixteen theta relations are the following, with
$$
(w,x,y,z)=\,M\,(a,b,c,d)\,,
$$
where
$$
M:= \frac{1}{2}\left[
\begin{array}{cccc}
1& \;1& 1&1\\
1& 1 &-1 &-1\\
1&-1& 1&-1\\
1& -1& -1& 1
\end{array}
\right]
$$

\medskip

\medskip
\begin{equation}\label{relt1}
\begin{matrix}
{{\vartheta }_2}(a ) \, {{\vartheta }_2}(b ) \, {{\vartheta }_2}(c )
\, {{\vartheta }_2}(d )+{{\vartheta }_3}(a ) \, {{\vartheta
}_3}(b ) \, {{\vartheta }_3}(c ) \, {{\vartheta }_3}(d ) =  \\
\noalign{\vspace{0.958333ex}} \hspace{1.em} {{\vartheta }_2}\;( x )
\, {{\vartheta }_2}\;( y) \, {{\vartheta }_2}\;( z ) \,
{{\vartheta }_2}\;(w )+{{\vartheta }_3}\;( x ) \,  {{\vartheta
}_3}\;(y ) \, {{\vartheta }_3}\;( z ) \, {{\vartheta
}_3}\;(w )\\
\end{matrix}
\end{equation}

\begin{equation} \begin{matrix} \label{relt2}
{{\vartheta }_3}(a ) \, {{\vartheta }_3}(b ) \, {{\vartheta }_3}(c )
\, {{\vartheta }_3}(d )-{{\vartheta }_2}(a ) \, {{\vartheta
}_2}(b ) \, {{\vartheta }_2}(c ) \, {{\vartheta }_2}(d ) =  \\
\noalign{\vspace{0.958333ex}} \hspace{1.em} {{\vartheta }_1}\;( x )
\, {{\vartheta }_1}\;( y) \, {{\vartheta }_1}\;( z ) \,
{{\vartheta }_1}\;(w )+{{\vartheta }_4}\;( x ) \,  {{\vartheta
}_4}\;( y) \, {{\vartheta }_4}\;( z ) \, {{\vartheta
}_4}\;(w )\\
\end{matrix} \end{equation}

\begin{equation} \begin{matrix} \label{relt3}
{{\vartheta }_1}(a ) \, {{\vartheta }_1}(b ) \, {{\vartheta }_1}(c )
\, {{\vartheta }_1}(d )+{{\vartheta }_4}(a ) \, {{\vartheta
}_4}(b ) \, {{\vartheta }_4}(c ) \, {{\vartheta }_4}(d ) =  \\
\noalign{\vspace{0.958333ex}} \hspace{1.em} {{\vartheta }_3}\;(w )
\, {{\vartheta }_3}\;( x ) \, {{\vartheta }_3}\;( y) \,
{{\vartheta }_3}\;( z )-{{\vartheta }_2}\;(w ) \,   {{\vartheta
}_2}\;( x ) \, {{\vartheta }_2}\;( y) \, {{\vartheta
}_2}\;( z )\\
\end{matrix} \end{equation}

\begin{equation} \begin{matrix}\label{relt4}
{{\vartheta }_4}(a ) \, {{\vartheta }_4}(b ) \, {{\vartheta }_4}(c )
\, {{\vartheta }_4}(d )-{{\vartheta }_1}(a ) \, {{\vartheta
}_1}(b ) \, {{\vartheta }_1}(c ) \, {{\vartheta }_1}(d ) =  \\
\noalign{\vspace{0.958333ex}} \hspace{1.em} {{\vartheta }_4}\;(w )
\, {{\vartheta }_4}\;( x ) \, {{\vartheta }_4}\;( y) \,  {{\vartheta
}_4}\;( z )-{{\vartheta }_1}\;(w ) \,  {{\vartheta }_1}\;( x ) \,
{{\vartheta }_1}\;( y) \, {{\vartheta
}_1}\;( z )\\
\end{matrix} \end{equation}

\begin{equation} \begin{matrix}\label{relt5}
{{\vartheta }_1}(a ) \, {{\vartheta }_1}(b ) \, {{\vartheta }_2}(c )
\, {{\vartheta }_2}(d )+{{\vartheta }_3}(c ) \, {{\vartheta
}_3}(d ) \, {{\vartheta }_4}(a ) \, {{\vartheta }_4}(b ) =  \\
\noalign{\vspace{0.958333ex}} \hspace{1.em} {{\vartheta }_1}\;( x )
\, {{\vartheta }_1}\;(w ) \, {{\vartheta }_2}\;( y) \,  {{\vartheta
}_2}\;( z )+{{\vartheta }_3}\;( y) \,  {{\vartheta }_3}\;( z ) \,
{{\vartheta }_4}\;( x ) \, {{\vartheta
}_4}\;(w )\\
\end{matrix} \end{equation}

\begin{equation} \begin{matrix}\label{relt6}
{{\vartheta }_4}(a ) \, {{\vartheta }_4}(b ) \, {{\vartheta }_3}(c )
\, {{\vartheta }_3}(d )-{{\vartheta }_1}(a ) \, {{\vartheta
}_1}(b ) \, {{\vartheta }_2}(c ) \, {{\vartheta }_2}(d ) =  \\
\noalign{\vspace{0.958333ex}} \hspace{1.em} {{\vartheta }_1}\;( y)
\, {{\vartheta }_1}\;( z ) \, {{\vartheta }_2}\;( x ) \,
{{\vartheta }_2}\;(w )+{{\vartheta }_3}\;( x ) \,   {{\vartheta
}_3}\;(w ) \, {{\vartheta }_4}\;( y) \, {{\vartheta
}_4}\;( z )\\
\end{matrix} \end{equation}

\begin{equation} \begin{matrix}\label{relt7}
{{\vartheta }_1}(a ) \, {{\vartheta }_1}(b ) \, {{\vartheta }_3}(c )
\, {{\vartheta }_3}(d )+{{\vartheta }_2}(c ) \, {{\vartheta
}_2}(d ) \, {{\vartheta }_4}(a ) \, {{\vartheta }_4}(b ) =  \\
\noalign{\vspace{0.958333ex}} \hspace{1.em} {{\vartheta }_1}\;( x )
\, {{\vartheta }_1}\;(w ) \, {{\vartheta }_3}\;( y) \,  {{\vartheta
}_3}\;( z )+{{\vartheta }_2}\;( y) \, {{\vartheta }_2}\;( z ) \,
{{\vartheta }_4}\;( x ) \, {{\vartheta
}_4}\;(w )\\
\end{matrix} \end{equation}

\begin{equation} \begin{matrix}\label{relt8}
{{\vartheta }_4}(a ) \, {{\vartheta }_4}(b ) \, {{\vartheta }_2}(c )
\, {{\vartheta }_2}(d )-{{\vartheta }_1}(a ) \, {{\vartheta
}_1}(b ) \, {{\vartheta }_3}(c ) \, {{\vartheta }_3}(d ) =  \\
\noalign{\vspace{0.958333ex}} \hspace{1.em} {{\vartheta }_1}\;( y)
\, {{\vartheta }_1}\;( z ) \, {{\vartheta }_3}\;( x ) \,
{{\vartheta }_3}\;(w )+{{\vartheta }_2}\;( x ) \,   {{\vartheta
}_2}\;(w ) \, {{\vartheta }_4}\;( y) \, {{\vartheta
}_4}\;( z )\\
\end{matrix} \end{equation}

\begin{equation} \begin{matrix}\label{relt9}
{{\vartheta }_2}(c ) \, {{\vartheta }_2}(d ) \, {{\vartheta }_3}(a )
\, {{\vartheta }_3}(b )+{{\vartheta }_2}(a ) \, {{\vartheta
}_2}(b ) \, {{\vartheta }_3}(c ) \, {{\vartheta }_3}(d ) =  \\
\noalign{\vspace{0.958333ex}} \hspace{1.em} {{\vartheta }_2}\;( x )
\, {{\vartheta }_2}\;(w ) \, {{\vartheta }_3}\;( y) \,  {{\vartheta
}_3}\;( z )+{{\vartheta }_2}\;( y) \, {{\vartheta }_2}\;( z ) \,
{{\vartheta }_3}\;( x ) \, {{\vartheta
}_3}\;(w )\\
\end{matrix} \end{equation}

\begin{equation} \begin{matrix}\label{relt10}
{{\vartheta }_3}(a ) \, {{\vartheta }_3}(b ) \, {{\vartheta }_2}(c )
\, {{\vartheta }_2}(d )-{{\vartheta }_2}(a ) \, {{\vartheta
}_2}(b ) \, {{\vartheta }_3}(c ) \, {{\vartheta }_3}(d ) =  \\
\noalign{\vspace{0.958333ex}} \hspace{1.em} {{\vartheta }_1}\;( x )
\, {{\vartheta }_1}\;(w ) \, {{\vartheta }_4}\;( y) \,  {{\vartheta
}_4}\;( z )+{{\vartheta }_1}\;( y) \,   {{\vartheta }_1}\;( z ) \,
{{\vartheta }_4}\;( x ) \, {{\vartheta
}_4}\;(w )\\
\end{matrix} \end{equation}

\begin{equation} \begin{matrix}\label{relt11}
{{\vartheta }_1}(c ) \, {{\vartheta }_1}(d ) \, {{\vartheta }_4}(a )
\, {{\vartheta }_4}(b )+{{\vartheta }_1}(a ) \, {{\vartheta
}_1}(b ) \, {{\vartheta }_4}(c ) \, {{\vartheta }_4}(d ) =  \\
\noalign{\vspace{0.958333ex}} \hspace{1.em} {{\vartheta }_3}\;(w )
\, {{\vartheta }_3}\;( x ) \, {{\vartheta }_2}\;( y) \,
{{\vartheta }_2}\;( z )-{{\vartheta }_2}\;(w ) \,  {{\vartheta
}_2}\;( x ) \, {{\vartheta }_3}\;( y) \, {{\vartheta
}_3}\;( z )\\
\end{matrix} \end{equation}

\begin{equation} \begin{matrix}\label{relt12}
{{\vartheta }_4}(a ) \, {{\vartheta }_4}(b ) \, {{\vartheta }_1}(c )
\, {{\vartheta }_1}(d )-{{\vartheta }_1}(a ) \, {{\vartheta
}_1}(b ) \, {{\vartheta }_4}(c ) \, {{\vartheta }_4}(d ) =  \\
\noalign{\vspace{0.958333ex}} \hspace{1.em} {{\vartheta }_4}\;(w )
\, {{\vartheta }_4}\;( x ) \, {{\vartheta }_1}\;( y) \,  {{\vartheta
}_1}\;( z )-{{\vartheta }_1}\;(w ) \,  {{\vartheta }_1}\;( x ) \,
{{\vartheta }_4}\;( y) \, {{\vartheta
}_4}\;( z )\\
\end{matrix} \end{equation}

\begin{equation} \begin{matrix}\label{relt13}
{{\vartheta }_2}(c ) \, {{\vartheta }_2}(d ) \, {{\vartheta }_3}(a )
\, {{\vartheta }_3}(b )+{{\vartheta }_1}(c ) \, {{\vartheta
}_1}(d ) \, {{\vartheta }_4}(a ) \, {{\vartheta }_4}(b ) =  \\
\noalign{\vspace{0.958333ex}} \hspace{1.em} {{\vartheta }_2}\;( y)
\, {{\vartheta }_2}\;( z ) \, {{\vartheta }_3}\;( x ) \,
{{\vartheta }_3}\;(w )+{{\vartheta }_1}\;( y) \,  {{\vartheta
}_1}\;( z ) \, {{\vartheta }_4}\;( x ) \, {{\vartheta
}_4}\;(w )\\
\end{matrix} \end{equation}

\begin{equation} \begin{matrix}\label{relt14}
{{\vartheta }_3}(a ) \, {{\vartheta }_3}(b ) \, {{\vartheta }_2}(c )
\, {{\vartheta }_2}(d )-{{\vartheta }_4}(a ) \, {{\vartheta
}_4}(b ) \, {{\vartheta }_1}(c ) \, {{\vartheta }_1}(d ) =  \\
\noalign{\vspace{0.958333ex}} \hspace{1.em} {{\vartheta }_2}\;( x )
\, {{\vartheta }_2}\;(w ) \, {{\vartheta }_3}\;( y) \,   {{\vartheta
}_3}\;( z )+{{\vartheta }_1}\;( x ) \,  {{\vartheta }_1}\;(w ) \,
{{\vartheta }_4}\;( y) \, {{\vartheta
}_4}\;( z )\\
\end{matrix} \end{equation}

\begin{equation} \begin{matrix}\label{relt15}
{{\vartheta }_1}(d ) \, {{\vartheta }_2}(b ) \, {{\vartheta }_3}(a )
\, {{\vartheta }_4}(c )+{{\vartheta }_1}(c ) \, {{\vartheta
}_2}(a ) \, {{\vartheta }_3}(b ) \, {{\vartheta }_4}(d ) =  \\
\noalign{\vspace{0.958333ex}} \hspace{1.em} {{\vartheta }_1}\;(w )
\, {{\vartheta }_4}\;( x ) \, {{\vartheta }_2}\;( y) \,
{{\vartheta }_3}\;( z )-{{\vartheta }_4}\;(w ) \, {{\vartheta
}_1}\;( x ) \, {{\vartheta }_3}\;( y) \, {{\vartheta
}_2}\;( z )\\
\end{matrix} \end{equation}

\begin{equation} \begin{matrix}\label{relt16}
{{\vartheta }_3}(a ) \, {{\vartheta }_2}(b ) \, {{\vartheta }_4}(c )
\, {{\vartheta }_1}(d )-{{\vartheta }_2}(a ) \, {{\vartheta
}_3}(b ) \, {{\vartheta }_1}(c ) \, {{\vartheta }_4}(d ) =  \\
\noalign{\vspace{0.958333ex}} \hspace{1.em} {{\vartheta }_3}\;(w )
\, {{\vartheta }_2}\;( x ) \, {{\vartheta }_4}\;( y) \,
{{\vartheta }_1}\;( z )-{{\vartheta }_2}\;(w ) \,  {{\vartheta
}_3}\;( x ) \, {{\vartheta }_1}\;( y) \, {{\vartheta
}_4}\;( z )\\
\end{matrix} \end{equation}

\bigskip
\bigskip

\bibliographystyle{plain}

\printindex
\end{document}